\documentclass{tac}
\pdfoutput=1

\usepackage{ifthen}
\usepackage{varioref}

\usepackage{hyperref}

\usepackage{amssymb}
\usepackage{amsmath}
\usepackage{amsfonts}
\usepackage{euscript}
\usepackage[rotate,all,color]{xy}
\usepackage[figuresright]{rotating}
\xyoption{2cell}

\usepackage{graphicx}
\usepackage{color}
\usepackage{rotating}

\usepackage{alltt}

\usepackage[utf8]{inputenc}
\usepackage[T1]{fontenc}

\DeclareMathAlphabet{\mathzc}{OT1}{pzcbig}{mb}{it}
\DeclareMathAlphabet{\mathrsfs}{OT1}{rsfs}{mb}{it}
\DeclareMathAlphabet{\mathcmssbx}{OT1}{cmss}{bx}{n}

\def\cube "#1" #2 #3 #4 {\save c="#1"+#2="#10", "#1"+#3="#11", "#1"+#4="#12",
  "#10"+#3="#101", "#11"+#4="#112", "#12"+#2="#102" +#3="#1012"\restore}

\def\tria{\NoRules\save\afterPOS{\POS **{} ?<(.1) *=0{}="zero" ; ?>(.9) : %
    "zero"+/v(0,1).7ex/ ::%
    \PATH ~={**\dir{-}} '(0,1) '(1,0) '(0,-1) '(0,1) %
    \restore};p;}

\makeatletter
\def\xyboxmatrix{%
  \def\xyboxmatrixname@{}%
  \def\xyboxmatrixsetup@{}%
  \xyFN@\xyboxmatrix@%
}
\def\xyboxmatrix@{%
  \ifx"\next\DNii@{\xyboxmatrix@prefix}%
  \else\addAT@\ifx\next\DNii@{\xyboxmatrix@setup}%
  \else\DNii@{\xyboxmatrix@ii}%
  \fi\fi\nextii@%
}
\def\xyboxmatrix@prefix"#1"{%
  \def\xyboxmatrixname@{#1}%
  \xyboxmatrix@setup%
}
\def\xyboxmatrix@setup#1#{%
  \def\xyboxmatrixsetup@{#1}%
  \xyboxmatrix@ii%
}
\def\xyboxmatrix@ii#1{%
  \def\xyboxmatrix@contents{#1}%
  \edef\xyboxmatrix@iii{\noexpand\POS*+!C\noexpand\xybox{%
      \noexpand\POS\noexpand\xymatrix\xyboxmatrixsetup@{%
        \noexpand\xyboxmatrix@contents}}="\xyboxmatrixname@"}%
  \xyboxmatrix@iii
}
\makeatother

\let\oldlabelbox=\labelbox
\def\newlabelbox#1{%
  \oldlabelbox{\vcenter{\normalbaselines%
      \let\\=\cr\ialign{\hspace{-.7ex}$\labelstyle##\hfil$\crcr#1\crcr}}
  }  
}
\let\labelbox\newlabelbox

\def\ifEqString#1#2{\def\testa{#1}\def\testb{#2}%
\ifx\testa\testb}

\renewcommand{\theta}{\vartheta}

\newcommand{\BOX}{\hbox {$\sqcap$ \kern -1em $\sqcup$}}
\newcommand{\qed}{\hskip 2em \hbox{\BOX} \vskip 2ex}

\newcommand{\To}{\Longrightarrow}
\newcommand{\Tto}{\Rrightarrow}

\def\pathspc{\overrightarrow}
\def\dblpathspc#1{\overrightarrow{\overrightarrow{#1}}}

\def\dblbarspc#1{\overline{\overline{#1}}}

\newcommand{\catfont}[1]{\mathsf{#1}}
\newcommand{\gcatfont}[1]{\mathbb{#1}}
\newcommand{\scatfont}[1]{\mathfrak{#1}}
\newcommand{\functorfont}[1]{\mathrsfs{#1}}

\def\mathcatdef#1{\expandafter\def\csname#1\endcsname{\catfont{#1}}}
\def\mathgcatdef#1{\expandafter\def\csname#1\endcsname{\gcatfont{#1}}}
\def\mathscatdef#1{\expandafter\def\csname s#1\endcsname{\scatfont{#1}}}

\mathcatdef{XMod}
\mathcatdef{Ab}
\mathcatdef{Cyl}
\mathcatdef{Set}
\mathcatdef{pSet}
\mathcatdef{pXMod}
\mathcatdef{pXSet}
\mathcatdef{Cat}
\mathcatdef{Top}
\mathcatdef{Grp}
\mathcatdef{Rng}
\mathcatdef{Mod}
\mathcatdef{CommRing}
\mathcatdef{Bimod}
\mathcatdef{Bicomod}
\mathcatdef{Comon}
\mathcatdef{Span}

\mathcatdef{C}
\mathcatdef{B}
\mathcatdef{A}
\mathcatdef{D}
\mathcatdef{Z}
\mathcatdef{E}
\mathcatdef{X}
\mathcatdef{A}
\mathcatdef{Vect}
\mathcatdef{Diff}
\mathcatdef{Gs}
\mathcatdef{T}
\mathcatdef{U}
\mathcatdef{Gray}
\mathcatdef{Graph}
\mathcatdef{RGrph}
\mathcatdef{Lax}
\mathcatdef{Fib}
\mathcatdef{Sesqui}
\def\GC{\Gray\Cat}

\newcommand{\fQ}{\functorfont{Q}}

\def\Sc{\EuScript{S}}        
\def\Cc{\mathcal C}
\def\Hc{\mathcal H}

\newcommand{\laxto}{\nrightarrow}
\renewcommand{\to}{\longrightarrow}
\renewcommand{\To}{\Longrightarrow}

\newcommand{\tto}{\genfrac{}{}{0pt}{1}{\longrightarrow}{\longrightarrow}}

\newcommand{\from}{\colon}

\newcommand{\op}{{\mathrm {op}}}
\newcommand{\iso}{\cong}
\newcommand{\id}{{\mathrm {id}}}

\newcommand{\g}{{\mathfrak g}}

\mathgcatdef{K}
\mathgcatdef{G}
\mathgcatdef{H}
\mathgcatdef{F}
\mathgcatdef{I}

\mathscatdef{S}
\mathscatdef{T}
\mathscatdef{C}
\mathscatdef{D}
\mathscatdef{E}

\renewcommand{\hat}{\widehat}

\newcommand{\ten}{\mathord{\otimes} }

\newcommand{\defterm}[2][]{{\bf #2}\ifthenelse{\equal{#1}{}}{\index{#2}}{\index{#1}}}
\newcommand{\ladj}{\dashv}

\renewcommand{\epsilon}{\varepsilon}
\newcommand{\eps}{\epsilon}
\renewcommand{\phi}{\varphi}

\newcommand{\pair}[1]{{\left<#1\right>}}
\newcommand{\lhc}{\triangleleft}
\newcommand{\rhc}{\triangleright}
\newcommand{\blhc}{\mathop{\ooalign{\hfil$\square$\hfil\cr\hfil$\lhc$\hfil\cr}}}
\newcommand{\brhc}{\mathop{\ooalign{\hfil$\square$\hfil\cr\hfil$\rhc$\hfil\cr}}}
\newcommand{\qopth}{\mathcal P}

\DeclareMathAlphabet{\mathpzc}{OT1}{pzc}{m}{it}

\newcommand{\mapar}{\ar@{|->}}
\newcommand{\ear}{\ar@{=}}
\newcommand{\lar}{\ar@{}}

\newcommand{\bten}{\boxtimes}

\newtheorem{thm}{Theorem}    
\newtheorem{cor}[thm]{Corollary}
\newtheorem{lem}[thm]{Lemma}
\newtheorem{rem}[thm]{Remark}

\newtheorem{defn}[thm]{Definition}

\def\qed {{
    \parfillskip=0pt        
    \widowpenalty=10000     
    \displaywidowpenalty=10000  
    \finalhyphendemerits=0  
    \leavevmode             
    \unskip                 
    \nobreak                
    \hfil                   
    \penalty50              
    \hskip.2em              
    \null                   
    \hfill                  
    $\square$
    \par                    
  }
}                           

\newenvironment{prf}{{\textsc {Proof}}}{\qed}

\SelectTips{cm}{}

\newboolean{draft}
\setboolean{draft}{false}

\ifthenelse{
  \boolean{draft}
}
{
  \usepackage{showkeys}

}
{
}

\newcommand{\cvshead}{}

\makeatletter
\def\setcvshead{%
  \edef\dollarcatcode{\the\catcode`\$}\catcode`\$=12
  \afterassignment\docvshead\toks@}
\def\docvshead{%
  \edef\cvshead{\the\toks@}%
  \catcode`\$=\dollarcatcode\relax}
\makeatother

\def\note#1{
  \ifthenelse{\boolean{draft}}
  {
    \marginpar{
        {\bf!!!}\tiny #1
    }
  }
  {
  }
}

\newcommand{\tqref}[1]{\text{\eqref{#1}}}

\usepackage{multicol}
\usepackage{layout}

\usepackage{natbib}
\bibpunct{[}{]}{,}{a}{}{;}

\xyoption{poly}
\def\xycompile#1{\POS #1}

\title{Mapping Spaces of $\Gray$-Categories}
\author{Björn Gohla}
\hypersetup{pdfauthor={Björn Gohla},
	pdftitle={Mapping Spaces of Gray-Categories},
	pdfkeywords={Higher gauge theory,  Gray-categories}
	}

\definecolor{Dark}{gray}{.2}
\definecolor{Medium}{gray}{.6}
\definecolor{Light}{gray}{.8}

\usepackage{floatpag}

\thanks{\noindent
  The author was supported by FCT (Portugal) through the doctoral
  grant SFRH/BD/33368/2008. This work was supported by FCT, with
  European Regional Development Fund (COMPETE) and national funds, by
  means of the projects PTDC/MAT/098770/2008 «Invariantes Topológicos
  via Geometria Diferencial» and PTDC/MAT/101503/2008 «Nova Geometria
  e Topologia». The author is a member of  CMUP/Universidade do
  Porto. The hospitality of CMA/Universidade Nova de Lisboa is
  gratefully acknowledged.
}
\address{}
\copyrightyear{2013}
\eaddress{b.gohla@gmx.de}

\begin{document}
\maketitle
\begin{abstract}
  We define a mapping space for $\Gray$-enriched categories
  adapted to higher gauge theory. Our construction differs
  significantly from the canonical mapping space of enriched
  categories in that it is much less rigid. The two essential
  ingredients are a path space construction for $\Gray$-categories and
  a kind of comonadic resolution of the 1-dimensional structure of a
  given $\Gray$-category obtained by lifting the resolution of
  ordinary categories along the canonical fibration of $\Gray\Cat$
  over $\Cat$.
\end{abstract}

\tableofcontents

\section{Introduction}
\label{sec:intro}

It is well known that among algebraic models for homotopy $n$-types
$\Gray$-groupoids model 3-types; \citet{lack2011} gives us a proof
using model category methods. Wanting to study the homotopy 3-type
of the moduli space of 3-connections on a manifold, we thought it apt
to define a mapping space $[\Sc_3(M),\Cc(\Hc)]$ of $\Gray$-groupoids
that could model that moduli space, where $\Sc_3(M)$ is the
fundamental $\Gray$-groupoid and $\Cc(\Hc)$ is the $\Gray$-groupoid
ultimately derived from a 2-crossed Lie-algebra where the
triconnections take their values; see for example \cite{sw11} for
2-connections, to which this is an obvious next step. See
\citep{marpic2011} for the background on the smooth fundamental
$\Gray$-groupoid and triconnections.
The original definition of the $\Gray$-tensor can be found
in \citep{gray}; \citet{gps} give us the definition of tricategories
and show that every tricategory strictifies to a triequivalent
$\Gray$-category. \citet{crans} gives an explicit, elementwise
definition of $\Gray$-categories.

In \citeyear{crans} \citeauthor*{crans} gave a partial solution the mapping space
problem; however, the absence of an interchange law in
$\Gray$-categories prevents lax transformations between
$\Gray$-functors from being composable in general. The slightly
unsatisfactory solution is to restrict to those transformations and
higher cells that can in fact be composed; this does give a mapping
space $\Gray$-category, but a mere stopgap not sufficient for our
purposes. 

Instead, we enlarge the repertoire of maps, and thereby
transformations, in a way that will permit forming all composites of
transformations; specifically we introduce a 2-cocycle that
intermediates coherently between the two possible evaluations of
arrangements of squares shown in \eqref{eq:pasteup} and
\eqref{eq:pastedown}. In analogy with \citet{garner} we introduce a
co-monadic weakening of strict $\Gray$-functors in section
\ref{sec:cofrpl}. The comonad $\fQ^1$ then yields a co-Kleisli
category $\Gray\Cat_{\fQ^1}$. We use in an essential way that
$\Gray\Cat$ is fibered over $\Cat$.

Inspired by \citep{ben1967} we axiomatise lax transformations as maps
into a path-space. In section \ref{sec:pathspc} we introduce a
functorial path-space construction for $\Gray$-categories;
subsequently in section \ref{sec:compaths} we show that this
yields an internal category $\pathspc\H\tto\H$ in $\Gray\Cat_{\fQ^1}$
for a given $\H$ in $\Gray\Cat$.

The $n$-th iterate of $\pathspc{(\_)}$ yields an $n$-truncated internal
cubical object in $\Gray\Cat$. In section \ref{sec:hcells} we
construct an internal $\Gray$-category 
$$\overline{\dblbarspc\H}\tto\dblbarspc\H\tto\pathspc\H\tto\H$$
in $\Gray\Cat_{\fQ^1}$ as a subobject of the third iterated
path-space. It is then a trivial consequence in section
\ref{sec:inthomfct} that we obtain a mapping space $\Gray$-category by
applying the hom
functor $$[\G,\H]:=\Gray\Cat_{\fQ^1}(\G,\overline{\dblbarspc\H}\tto\dblbarspc\H\tto\pathspc\H\tto\H).$$
Furthermore we obtain a restricted mapping space $\{\G,\H\}$, where
everything is as before, except only strict $\Gray$-functors are
permitted between $\G$ and $\H$. This leads to a natural
sesquicategory structure on $\Gray\Cat$.

We hope to be able to prove in a later paper that this internal hom is
part of a monoidal closed structure on $\Gray\Cat_{\fQ^1}$ involving a
suitable extension of \citeauthor{crans}' tensor product.

Finally, in section \ref{sec:together} we give explicit details of
functors, transformations and so on in terms of components.  Lastly,
we remark that if $\H$ is a $\Gray$-groupoid then $\pathspc\H$ as well
as $[\G,\H]$ will be $\Gray$-groupoids.

Similar work was done by \citet{goma2013} concerning 2-crossed
modules, which are equivalent to $\Gray$-groupoids with a single
vertex, that is, $\Gray$-groups. 

A version of this article constituted the author's doctoral thesis
defended at the Faculty of Science, University of Porto. Many thanks
are owed to João Faria Martins for plentiful advice and discussion.

\section{Resolution in Dimension One}
\label{sec:cofrpl}

We define a resolution of the 1-dimensional structure of a
$\Gray$-category using a comonad, by lifting the free category comonad
(called ``path'' in \citep{daw2006}) to $\Gray$-categories; but note that
we use the term in a different way in this paper.

The resulting co-Kleisli category
can be seen as the category of $\Gray$-categories with an enlarged
repertoire of maps, that is flexible enough to carry out our path space
construction. After giving an abstract construction of this category
of pseudo maps we proceed to characterize them explicitly. 

\subsection{Basic Fibrations}
\label{sec:basfib}

There are obvious functors
\begin{equation*}
  \xymatrix{
    \Gray\Cat\ar[r]^-{(\_)_2}&\Sesqui\Cat\ar[r]^-{(\_)_1}&\Cat\ar[r]^-{(\_)_0}&\Set
  }
\end{equation*}
that forget the 3-cells, the 2-cells and 1-cells respectively.
By a slight abuse of language we will denote the composite
$(\_)_1(\_)_2$ by $(\_)_1$ also, it is of
course a fibration as well; we will use it in section \ref{sec:comlift}
to construct the monad $\fQ^1$. We will use the fibration
$(\_)_2(\_)_1(\_)_0=(\_)_0$ in section \ref{sec:inthomfct} to
construct the restricted mapping space $\{\G,\H\}$.

Let $\sS$ be a sesquicategory, $\G$ a $\Gray$-category, and
$F\from\sS\to\G_2$ a sesquifunctor. We define
$\overline{F}\from{}F^*\sS\to\G$ as follows:

\begin{align*}
  (F^*\sS)_0&=\sS_0\\
  (F^*\sS)_1&=\sS_1\\
  (F^*\sS)_2&=\sS_2\\
  (F^*\sS)_3&=\left\{(\Gamma;\alpha,\beta)\,\vline\, \Gamma\from F\alpha\to F\beta\right\}
\end{align*}

Note that the interchange of two 2-cells $\alpha, \beta$ in $F^*\sS$
incident on a 0-cell is given essentially by the interchange of their
images under $F$:
\begin{equation*}
  \beta\ten\alpha=(F\beta\ten{}F\beta;\beta\rhc\alpha,\beta\lhc\alpha)\,.
\end{equation*}

Let us take note of the following useful fact that helps to
characterize the Cartesian maps:
\begin{rem}
  For a functor $p\from\E\to\B$ that preserves co-limits, let
  $D\from\D\to\E$ a diagram in $\E$ with co-limit $(C,k_i)$
  \begin{equation*}
    \begin{xy}
      \xyboxmatrix{
        D_i\ar[r]^{k_i}&C\ar[rrd]^g&{}&{}\\
        {}&{}&A\ar[r]_f&B
      }
    \end{xy}\,,
  \end{equation*}
  assume $p(g)$ factors below as $p(f)u=p(g)$. Furthermore, assume
  that the induced sink $(u_i)=up(k_i)$ has fillers
  $\left<u_i\right>$ above with $f\left<u_i\right>=gk_i$, then the
  co-universally induced map $\left<u\right>\from{}C\to{}A$ is a
filler over $u$. 

This means that to check whether a map $f$ is Cartesian we don't need
to give the filler $u$ directly, but we can define it on presumably
simpler parts of $C$. These then combine into a valid filler.
\end{rem}

\begin{rem}
  Maps Cartesian with respect to $(\_)_2$ are exactly the
  $\Gray$-functors, that are 2-locally isomorphisms of sets. That is, given
  two parallel 2-cells on the intervening 3-cells, the map is bijective.
\end{rem}

\begin{lem}
  $F^*\sS$ is a $\Gray$-category, $\overline{F}$ is a
  $\Gray$-functor and Cartesian with respect to
  $(\_)_2$. \qed
\end{lem}

Similarly, let $\sS$ be a sesquicategory, $\C$ a category, and
$F\from\C\to\sS_1$ a functor, then we define a sesquicategory:

\begin{align*}
  (F^*\C)_0&=\C_0\\
  (F^*\C)_1&=\C_1\\
  (F^*\C)_2&=\left\{(\alpha;f,g)\,\vline\, \alpha\from Ff\to Fg\right\}
\end{align*}

\begin{lem}
  $F^*\C$ is a sesquicategory, $\overline{F}$ is a sesquifunctor, and
  Cartesian with respect to $(\_)_1$.\qed
\end{lem}

\begin{rem}
  Maps Cartesian with respect to $(\_)_1$ are exactly the
  sesquifunctors, that are 1-locally isomorphisms of sets. That is,
  given two parallel 1-cells on the intervening 2-cells, the map is
  bijective.
\end{rem}

For later reference we describe the 
Cartesian liftings of $(\_)_1$ explicitly as well. Let $\G$ be a $\Gray$-category,
$\G_1$ its underlying category. Let $\C$ be an ordinary category and
$F\from\C\to\G_1$ a functor. Then $F^*\G$ is given by:
\begin{align*}
  (F^*\G)_0&=\C_0\\
  (F^*\G)_1&=\C_1\\
  (F^*\G)_2&=\left\{(\alpha;f,g) \,\vline\, f, g\from x\to y,\, \alpha\from
    Ff\to Fg \right\}\\
  (F^*\G)_3&=\left\{(\Gamma;\alpha,\beta;f,g)\,\vline\, f, g\from x\to y,\,
    \Gamma\from F\alpha\to F\beta\right\}
\end{align*}
Source and target maps are as follows:
\begin{align*}
  s_2(\Gamma;\alpha,\beta;f,g)&=(\alpha; f,g) &\quad
  t_2(\Gamma;\alpha,\beta;f,g)&=(\beta;f,g)\\
  s_1(\alpha;f,g)&=f &\quad t_1(\alpha;f,g)&=g\,.
\end{align*}
and $s_0, t_0$ are as given by $\C$. As identities we take:
\begin{align*}
  i_1(f)=(\id_{Ff};f,f)&\quad i_2(\alpha;f,g)=(\id_\alpha;\alpha,\alpha,f,g)\,.
\end{align*}

The tensor in $F^*\G$ of two 2-cells is 
\begin{equation}
  \label{eq:tenpbdef}
  (\beta; g, g')\ten(\alpha; f, f')=
  \left(
    \beta\ten\alpha;
    \beta\lhc\alpha, \beta\rhc\alpha;
    g\#_0f, g'\#_0f'
  \right)
\end{equation}
where 
\begin{equation*}
  \beta\lhc\alpha=(\beta\#_0Ff')\#_1(Fg\#_0\alpha),\quad
  \beta\rhc\alpha=(Fg'\#_0\alpha)\#_1(\beta\#_1Ff)\,.
\end{equation*}

There is an obvious map $\overline{F}\from{}F^*\G\to\G$ over $F$ that
acts like $F$ on 0- and 1-cells, and on 2- and 3-cells as a projection
to $\G$.

\begin{rem} 
  The globular set $F^*\G$ is a $\Gray$-category.
  The composition operations of $F^*\G$ are given by those of $\C$ and
  $\G$ and it is easy to see that they fulfill the axioms of a
  $\Gray$-category. 
\end{rem}

Obviously $G^*F^*\G\iso(FG)^*\G$ and $\id_\C^*\iso\id_{\Gray\Cat_\C}$
coherently. Also, we can always choose $\id_\C^*=\id_{\Gray\Cat_\C}$,
but this is not necessary in what follows.

\begin{lem}
  A map of $\Gray$-categories is Cartesian with respect to
  $\G\mapsto\G_1$  iff it is 1-locally an isomorphism of categories,
  i.e.\@ given two parallel 1-cells the map is bijective on the
  intervening 2-cells and in turn bijective on the 3-cells between
  parallel such.\qed
\end{lem}

\begin{defn}\label{defn:niso}
  We define a map of $\Gray$-categories to be an
  \defterm{$n$-isomorphism} if it is Cartesian with respect to
  $(\_)_n$. It is \defterm{$n$-faithful} if fillers of factorizations
  under $(\_)_n$ are unique, and \defterm{$n$-full} is there (not
  necessarily unique) fillers for all factorizations under $(\_)_n$.  
\end{defn}

By this definition $0$-fidelity is ordinary fidelity of functors,
$1$-fidelity is local fidelity, and so on. 

\begin{rem}\label{rem:pmono}
  One property of Cartesian maps in a fibration $p$ that we are
  going to exploit in the proof of the following theorem is that for
  three arrows upstairs,
  \begin{equation*}
    \begin{xy}\xycompile{
        \xymatrix{
          {}\ar@<+.5ex>[r]^{r}\ar@<-.5ex>[r]_s &{}\ar[r]^f&{}
        }
      }\end{xy}
  \end{equation*}
  with $f$ Cartesian,
  $p(r)=p(s)$ downstairs and $fr=fs$ upstairs imply $r=s$, on account
  of $f$ being $p$-faithful.
\end{rem}

\begin{lem}
  \label{lem:pffcancel}
  If $fg$ is Cartesian with respect to a given fibration $p$ and $f$
  is $p$-faithful, then $g$ is $p$-Cartesian.
\end{lem}
\begin{prf}
  Take $k$ and $u$ such that $p(g)u=p(k)$, then $p(fg)u=p(fk)$ and
  hence by $fg$ being $p$-full there is a filler $\left<u\right>$ such
  that $fg\left<u\right>=fk$. Then by $f$ being $p$-faithful
  $g\left<u\right>=k$.

  By $fg$ being $p$-faithful $\left<u\right>$ is the unique such
  filler. 
\end{prf}

\pagebreak

\subsection{Comonad Liftings}
\label{sec:comlift}

In this section we show that comonads can be lifted along fibrations
of categories.

\begin{defn}
  \label{defn:comonad}
  In an arbitrary 2-category a \defterm{comonad} on an object $A$ is
  given by an endomorphism
  \begin{equation*}
    \begin{xy}
      \xyboxmatrix{
        A\ar[r]^T&A
      } 
    \end{xy}
  \end{equation*}
  and  2-cells 
  \begin{equation*}
    \begin{xy}
      \xyboxmatrix{
        A\ar@/^1.5pc/[r]^T="x"\ar@/_1.5pc/[r]_A="y"&A
        \ar@2"x";"y"**{}?(.3);?(.7)^{\eps}
      }
    \end{xy}
  \end{equation*}
  and
  \begin{equation*}
    \begin{xy}
      \xyboxmatrix{
        A\ar[r]_{T}\ar@/^2pc/[rr]^{T}="x"&A\POS[]="y"\ar[r]_{T}&A
        \ar@2"x";"y"**{}?(.3);?(.7)^{\delta}
    }
    \end{xy}
  \end{equation*}
  such that
  \begin{equation*}
    \begin{xy}
      \xyboxmatrix{
        A\ar[r]^{T}\ar@/^4pc/[rr]^{T}="x1"&A\POS[]="y1"\ar@/^1pc/[r]^T="x"\ar@/_1pc/[r]_A="y"&A
        \ar@2"x";"y"**{}?(.3);?(.7)^{\eps}
        \ar@2"x1";"y1"**{}?(.3);?(.7)^{\delta}
      }
    \end{xy}=
    \begin{xy}
      \xyboxmatrix{
        A\ar@/^1.5pc/[r]^T="x"\ar@/_1.5pc/[r]_T="y"&A
        \ar@2"x";"y"**{}?(.3);?(.7)^{T}
      }
    \end{xy}=
    \begin{xy}
      \xyboxmatrix{
        A\ar@/^4pc/[rr]^{T}="x1"\ar@/^1pc/[r]^T="x"\ar@/_1pc/[r]_A="y"&A\POS[]="y1"\ar[r]^{T}&A
        \ar@2"x";"y"**{}?(.3);?(.7)^{\eps}
        \ar@2"x1";"y1"**{}?(.3);?(.7)^{\delta}
      }
    \end{xy}
  \end{equation*}and
  \begin{equation*}
    \begin{xy}
      \xyboxmatrix{
        A\ar[r]_{T}\ar@/^4pc/[rrr]^{T}="x"&A\POS[]="y"\ar[r]_{T}\ar@/^2pc/[rr]^{T}="x1"&A\POS[]="y1"\ar[r]_{T}&A
        \ar@2"x";"y"**{}?(.3);?(.7)^{\delta}
        \ar@2"x1";"y1"**{}?(.3);?(.7)^{\delta}
      }
    \end{xy}= \begin{xy}
      \xyboxmatrix{
        A\ar[r]_{T}\ar@/^4pc/[rrr]^{T}="x"\ar@/^2pc/[rr]^{T}="x1"&A\POS[]="y1"\ar[r]_{T}&A\POS[]="y"\ar[r]_{T}&A
        \ar@2"x";"y"**{}?(.3);?(.7)^{\delta}
        \ar@2"x1";"y1"**{}?(.3);?(.7)^{\delta}
      }
    \end{xy}\,.
  \end{equation*}
  See, for example, \citet{maclane}.
\end{defn}

If $A$ is a category, $T$ a functor and $\eps$ and $\delta$ natural
transformations, then these equations of course amount to the usual
equations objectwise in $A$:
\begin{equation*}
  \begin{xy}
    \xymatrix{
      {}&Tx\ar[dl]_{Tx}\ar[dr]^{Tx}\ar[d]|{\delta_x}&{}\\
      Tx&TTx\ar[r]_{T\eps_{x}}\ar[l]^{\eps_{Tx}}&Tx
    }
  \end{xy}
\end{equation*}
and 
\begin{equation*}
  \begin{xy}
    \xyboxmatrix{
      Tx\ar[r]^-{\delta_{x}}\ar[d]_-{\delta_{x}}&TTx\ar[d]^{T\delta_x}\\
      TTx\ar[r]_{\delta_{Tx}}&TTTx
    }
  \end{xy}\,.
\end{equation*}

\begin{thm}\label{thm:comonadlift}
  Given a fibration of categories $p\from \E\to \B$,
  a comonad $(Q,\delta,\eps)$ on $\B$ can be lifted to a comonad $(K,
  d, e)$ on $\E$ such that $(K, Q)\from p\to p$ is a comonad in the
  2-category of all fibrations.
\end{thm}

\begin{prf}
  Let $(\_)^*\from\B^\op\to\Cat$ be a chosen cleavage. For every
  $A\in \E_x$ we let $e_A\from(KA=\eps^*_xA)\to A$ be the chosen
  Cartesian lift of $\eps_x\from Qx\to x$. For a morphism $f$ over $j$ in
  \begin{equation*}
    \begin{xy}\xycompile{
      \xymatrix{
        KA\ar[r]^{e_A}\ar@{..>}[dr]_{Kf}& A\ar[dr]^f & \\
        {} & KB\ar[r]_{e_B} & B \\
        Qx\ar[r]^{\eps_x}\ar[dr]_{Qj} & x\ar[dr]^{j} & \\
        {} & Qy\ar[r]_{\eps_y} & y
      }
    }\end{xy}
  \end{equation*}
  the dotted arrow is the unique filler induced by the factorization
  below. This makes $K$ a functor and $e\from K \to \id_\E$ a natural
  transformation.

  We define a family of co-multiplication maps $d_A$ as the unique
  fillers in
  \begin{equation*}
    \begin{xy}\xycompile{
      \xymatrix{
        KA \ar[drr]^{{KA}} \ar@{..>}[dr]_{d_A}&&\\
        {} & KKA \ar[r]_-{e_{KA}}& KA\\
        Qx \ar[drr]^{{Qx}} \ar[dr]_{\delta_x}&& \\
        {} & QQx \ar[r]_{\eps_{Qx}}& Qx
      }
    }\end{xy}
  \end{equation*}
  where the triangle below commutes because is $Q$ co-unital.

  In the diagram 
  \begin{equation*}
    \begin{xy}\xycompile{
      \xymatrix{
        KA\ar@{..>}[drr]\ar[dr]_{d_A}\ar@(r,u)[drr]_{KA}\ar@(r,u)[drrr]^{e_A}&&&\\
        {} &KKA\ar@<.5ex>[r]^{Ke_A}\ar@<-.5ex>[r]_{e_{KA}}
        &KA\ar[r]_{e_A} &A\\
        Qx\ar[dr]_{\delta_x}\ar@(r,u)[drr]_{Qx}\ar@(r,u)[drrr]^{\eps_x}&&&\\
        {} &QQx\ar@<.5ex>[r]^{Q\eps_x}\ar@<-.5ex>[r]_{\eps_{Qx}}
        &Qx\ar[r]_{\eps_x} &x\\
      }
    }\end{xy}
  \end{equation*}
  we see that $e_Ae_{KA}d_A=e_AKe_Ad_A$ by the naturality of $e$, and
  $p(e_{KA}d_A)=p(Ke_Ad_A)$ by $Q$ being a comonad. Hence by remark
  \ref{rem:pmono} the three endomorphisms of $KA$ above have to
  coincide, meaning $d$ is co-unital component wise.

  The naturality of $d$, that is, that $d_BKf=KKfd_A$ is the unique
  filler making the left-hand upstairs square commute
  \begin{equation*}
    \begin{xy}\xycompile{
      \xymatrix{
        KA\ar@{..>}[drr]\ar[r]^-{d_A}\ar[dr]_{Kf}& KKA\ar[dr]^{KKf}& & \\
        {}& KB\ar[r]_-{d_B} & KKB\ar[r]_{e_{KB}} & KB \\
        Qx\ar[r]^-{\delta_x}\ar[dr]_{Qj}& QQx\ar[dr]^{QQj}& & \\
        {}& Qy\ar[r]_-{\delta_y} & QQy\ar[r]_{\eps_{Qy}} & Qy \\
      }
    }\end{xy}
  \end{equation*}
  is obtained by observing that
  $e_{KB}d_BKf=KF=Kfe_{KA}d_A=e_{KB}KKfd_A$, from $e$ being natural
  and a retraction. Also, $p(d_BKf)=p(KKfd_a)$ by naturality of
  $\delta$. We apply \ref{rem:pmono} again.

  Finally, we show that $d$ is co-associative: Consider the diagram
  \begin{equation*}
    \begin{xy}\xycompile{
      \xymatrix{
        KA\ar@{..>}[drr]\ar[r]^-{d_A}\ar[dr]_{d_A}& KKA\ar[dr]^{d_{KA}}&& \\
        {}& KKA\ar[r]_-{Kd_{A}} & KKKA\ar[r]_{e_{KKA}} & KKA \\
        Qx\ar[r]^-{\delta_x}\ar[dr]_{\delta_x}&
        QQx\ar[dr]^{\delta_{Qx}}&& \\
        {}& QQx\ar[r]_-{Q\delta_{x}} & QQQx\ar[r]_{\eps_{QQx}} & QQx
        \,. \\
      }
    }\end{xy}
  \end{equation*}
  We calculate that
  $e_{KKA}Kd_Ad_A=d_Ae_{KA}d_A=d_A=e_{KKA}d_{KA}d_A$, again by
  naturality of $e$ and its retractiveness. Moreover, $\delta$ is
  co-associative, hence we can apply remark \ref{rem:pmono} once more.
\end{prf}

We observe that $K$ preserves Cartesianness of maps, thus in
particular $Ke$ is Cartesian component wise.

Finally we can define our resolution comonad. Let $(Q,\delta,
\epsilon)=(FU,F\eta U, \epsilon)$ be the comonad that arises from the
adjunction
\begin{equation*}
  \begin{xy}
    \xyboxmatrix{
      \RGrph\ar@/^/[r]^-{F}="x"&\Cat\ar@/^/[l]^-{U}="y"
      \ar@{|- }"y";"x"**{}?(.4);?(.6)
    }
  \end{xy}\,.
\end{equation*}
Then, according to theorem \ref{thm:comonadlift}, we obtain the comonad
$(\fQ^1,d,e)$ on $\Gray\Cat$ induced by lifting $Q$ along $(\_)_1$. The
exponent reminds us that this provides a resolution of the
1-dimensional structure of $\Gray$-categories. See section
\ref{sec:adjunctions} for a more abstract point of view on this
construction. In section \ref{sec:spectcel} we will show explicitly
how this comonad acts.

\begin{cor}
  By the above theorem there is a comonad $\fQ^1$ on $\Gray\Cat$ that pulls back the
  $\Gray$-structure onto the free category on the underlying 1-graph.
\end{cor}

If a category $\C$ is already the free category $\C=F\g$ over a
reflexive graph with injection of generators $\eta\from\g\to{}U\C$,
then by adjointness the counit is split
\begin{equation*}
  \begin{xy}
    \xyboxmatrix{
      \C\ar[r]^-{F\eta}\ar@/_2pc/[rr]_{\C}&Q\C\ar[r]^-{\eps}&\C
    }
  \end{xy}\,.
\end{equation*}

\begin{defn}
  \label{defn:freeorder1}
  If a $\Gray$-category $\G$ has an underlying category $\G_1$ of the
  form $F\g$ for some reflexive graph $\g$ we say that $\G$ is
  \defterm{free up to order 1} with generating 1-cells $\g$.
\end{defn}
Let $k\from\G\to\fQ^1\G$ be the filler along $(\_)_1$ for the
factorization $e_1F\eta=(\id_\G)_1$ for the
given generating reflexive graph. This of course gives a splitting
\begin{equation}
  \label{eq:freeorder1split}
  \begin{xy}
    \xyboxmatrix{
      \G\ar[r]^-{k}\ar@/_2pc/[rr]_{\G}&\fQ^1\G\ar[r]^-{e}&\G
    }
  \end{xy}\,.
\end{equation}

If a $\Gray$-category is free up to order 1 we may look at the 1-cells
as follows: every 1-cell $f$ can be written as $[f_1,\ldots, f_n]$,
where the $[f_i]$ are generating 1-cells unique up to insertion and
deletion of units. Now, the action of $k\from\G\to\fQ^1\G$ can be
described as follows:
\begin{enumerate}
\item 0-cells: $k\from{}x\mapsto{}x$
\item 1-cells: $k\from{}f=[f_1,\ldots,f_n]\mapsto[[f_1],\ldots,[f_n]]$
\item 2-cells:
  $k\from{}(\alpha\from{}f\To{}f')\mapsto(\alpha;[[f_1],\ldots,[f_n]],[[f'_1],\ldots,[f'_{n'}]])$
\item 3-cells:
  $k\from{}(\Gamma\from{}\alpha{}\Rrightarrow\alpha')\mapsto(\Gamma;\alpha,
  \alpha';[[f_1],\ldots,[f_n]],[[f'_1],\ldots,[f'_{n'}]])$
\end{enumerate}
This is obviously a section of $e_\G$.

\begin{defn}
  \label{defn:graykleisli}
  The category of $\Gray$-categories and \defterm[pseudo
  $\Gray$-map]{pseudo $\Gray$-maps} is the co-Kleisli-category
  $\Gray\Cat_{\fQ^1}$ of the comonad $\fQ^1$.
\end{defn}

\begin{lem}\label{lem:freeorder1q11}
  The map $k$ for a $\G$ free up to order 1 has the following nice
  behaviour with respect to $\fQ^1$:
  \begin{equation}
    \label{eq:freeorder1q1}
    \begin{xy}
      \xymatrix{
        \G\ar[r]^-{k}\ar[d]_k&\fQ^1\G\ar[d]^{d}\\
        \fQ^1\G\ar[r]_-{\fQ^1k}&\fQ^1\fQ^1\G
      }
    \end{xy}\,.
  \end{equation}
  commutes.
\end{lem}
\begin{prf}
  We apply remark \ref{rem:pmono}: The diagram
  \begin{equation*}
    \begin{xy}
      \xymatrix{
        \G\ar[r]^-{k}\ar[d]_k&\fQ^1\G\ar[d]^{d}\\
        \fQ^1\G\ar[r]_-{\fQ^1k}\ar[d]_{e}&\fQ^1\fQ^1\G\ar[d]^{e}\\
        \G\ar[r]_k&\fQ^1\G
      }
    \end{xy}
  \end{equation*}
  commutes by co-unitality and the definition of $k$.
  Also under $(\_)_1$ the diagram (\ref{eq:freeorder1q1}) becomes 
  \begin{equation*}
    \begin{xy}
      \xymatrix{
        F\g\ar[r]^-{F\eta}\ar[d]_{F\eta}&FUF\g\ar[d]^{F\eta{}UF}\\
        FUF\g\ar[r]_-{FUFU\eta}&FUFUF\g
      }
    \end{xy}
  \end{equation*}
  which commutes by naturality of $\eta$.
\end{prf}

This category has $\Gray$-categories as objects, and morphisms
\begin{equation*}
  \xymatrix{
    \G \ar[r]|*\dir{/}^f & \H
  }
  \text{\qquad are morphisms\qquad}
  \xymatrix{
    \fQ^1\G \ar[r]^f & \H
  }
\end{equation*}
in $\Gray\Cat$. Composition of two maps 
\begin{equation*}
  \xymatrix{
    \G \ar[r]|*\dir{/}^f & \H \ar[r]|*\dir{/}^g& \K
  }
\end{equation*}
is defined by 
\begin{equation*}
  \xymatrix{
    \fQ^1\G\ar[r]^-{d_\G}& \fQ^1\fQ^1\G \ar[r]^-{\fQ^1f} & 
    \fQ^1\H \ar[r]^-g& \K\,.
  }
\end{equation*}
Identities are of the form
\begin{equation*}
  \xymatrix{
    \G \ar[r]|*\dir{/}^{\id_\G} & \G
  }
  =
  \xymatrix{
    \fQ^1\G \ar[r]^{e_\G} & \G\,.
  }
\end{equation*}
By way of notational convenience in diagrams in $\Gray\Cat_{\fQ^1}$ we
use unslashed arrows $f\from\G\to\H$ to denote a strict arrow that is
included in $\Gray\Cat_{\fQ^1}$ as $fe\from\G\laxto\H$.   

The comonad axioms make sure this is a category; c.f.~e.g.~
\citep{maclane}.

There is an adjunction 
\begin{equation*}
  \xymatrix{
    \Gray\Cat\ar@/^/[r]^-{R}="x"&\Gray\Cat_{\fQ^1}\ar@/^/[l]^-{L}="y"
    \ar@{|- }"x";"y"**{}?(.4);?(.6)
  }
\end{equation*}
The functor $R$ takes a strict map $f\from\G\to\H$ to a pseudo map
$fe\from\G\laxto\H$ where $e$ is the co-unit of $\fQ^1$. Moreover, since
$e$ is an epimorphism, $R$ is faithful, and it is bijective on objects, hence $R$
is actually an inclusion; in particular, we have injective maps
\begin{equation}
  \label{eq:strictinclude}
  \begin{xy}
    \xyboxmatrix{
      \GC(\G,\H)\ar[r]^-{e^*}&\GC_{\fQ^1}(\G,\H)
    }
  \end{xy}
\end{equation}
for all $\G$ and $\H$.

We note that the composite of a strict map after a pseudo
map is particularly simple:

\begin{multline}
  \label{eq:strpseudocomp}
  \begin{xy}\xycompile{
    \xymatrix{
      \G\ar[r]|*\dir{/}^{f} &\H\ar[r]|*\dir{/}^{ge} & \K
    }
  }\end{xy}
  =
  \begin{xy}\xycompile{
    \xymatrix{
      \fQ^1\G\ar[r]^-{d_{\fQ^1\G}}\ar@{=}[dr] &\fQ^1\fQ^1\G\ar[r]^{\fQ^1f}\ar[d]|{e_{\fQ^1\G}} & \fQ^1\H\ar[r]^{ge}\ar[d]|{e_{\H}} &\K\\
      {} & \fQ^1\G\ar[r]_f & \H\ar[ur]_{g} &{}
    }
  }\end{xy}\,.
\end{multline}

If $\G$ is free up to order 1 we also get an idempotent function 
\begin{equation}
  \label{eq:freeorderstrictify}
  \begin{xy}
    \xyboxmatrix{
      \Gray\Cat_{\fQ^1}(\G,\H)\ar[r]^{(ke)^*}&\Gray\Cat_{\fQ^1}(\G,\H)
    }
  \end{xy}
\end{equation}
from (\ref{eq:freeorder1split}) we might call strictification (note
the reverse order of $k$ and $e$). It preserves the image of the
functor $R$, that is, strict $\Gray$-functors are preserved.

\begin{lem}
  \label{lem:grcatstrlims}
  The category $\Gray\Cat_{\fQ^1}$ has all limits of diagrams of
  strict maps, that is, those in the subcategory $\Gray\Cat$, that is,
  $\Gray\Cat$ is complete and the inclusion
  $\Gray\Cat\to\Gray\Cat_{\fQ^1}$ preserves all limits.
\end{lem}

\begin{prf}
  Let $D$ be a diagram in $\Gray\Cat$, let $(\ell_i\from L \to D_i)_i$
  be a limiting source in $\Gray\Cat$, we claim its embedding into
  $\Gray\Cat_{\fQ^1}$ is a limiting source there as well.

  Let $(c_i\from C \laxto D_i)_i$ be a source over $D$ in
  $\Gray\Cat_{\fQ^1}$. Thus there is a source $(c_i\from \fQ^1C \to
  D_i)_i$ in $\Gray\Cat$, which induces a map $\pair{c}\from\fQ^1C\to
  L$ and this is of course a map $\pair{c}\from C\laxto L$. The
  diagram 
  \begin{equation*}
    \xymatrix{
      C \ar[dr]^{c_i}|*\dir{/}\ar[d]_{\pair{c}}|*\dir{/}& {} \\
      L \ar[r]_{\ell_i}& D_i
    }
  \end{equation*}
  commutes for all $i$ by the co-unit axiom of $\fQ^1$ and the
  naturality of $e$; c.~f.\@ also \eqref{eq:strpseudocomp}. Because $e$
  is an epimorphism  $\pair{c}$ is the unique filler. 
\end{prf}

In particular, the pullback of two strict maps in $\Gray\Cat_{\fQ^1}$
is the same as its pullback in $\Gray\Cat$. Products are
obviously simply the same in both categories since their diagrams do
not include any nontrivial morphisms.

\begin{rem}  
  \label{rem:indmap}
  For two diagrams $\{a_k\from\G_i\to\G_j\}$,
  $\{b_k\from\H_i\to\H_j\}$ of strict maps of the same type in
  $\Gray\Cat_{Q^1}$ and a natural transformation
  ${f_i\from\G_i\laxto\H_i}$ between them there is an induced map
  $\dot\lim\{f_i\}$ such that:
  \begin{equation}
    \label{eq:indmap1}
    \begin{xy}
      \xyboxmatrix{
        \lim\{\G_i,a_k\}\ar[r]|*\dir{/}^-{\dot\lim{f_i}}\ar[d]_{p_i}&\lim\{\H_i,b_k\}\ar[d]^{p'_i}\\
        \G_i\ar[r]|*\dir{/}_{f_i}&\H_i
      }
    \end{xy}\,.
  \end{equation}
  We unravel this diagram in terms of maps in $\Gray\Cat$ and obtain
  \begin{equation*}
    \begin{xy}
      \xyboxmatrix{
        \fQ^1\lim\{\G_i,a_k\}
        \ar@/^2pc/[rr]^-{\dot\lim{f_i}}\ar[dr]_{\fQ^1p_i}
        \ar[r]^-{\pair{\fQ^1p_i}}&
        \lim\{\fQ^1\G_i,\fQ^1a_k\}\ar[d]^{r_i}\ar[r]^{\lim{f_i}}&
        \lim\{\H_i,b_k\}\ar[d]^{p'_i}\\
        {}&\fQ^1\G_i\ar[r]_{f_i}&\H_i
      }
    \end{xy}
  \end{equation*}
  where the map $\dot\lim{f_i}$ is induced by the universal property of
  the source $\{f_i\fQ^1p_i\}$ in $\Gray\Cat$, that is,
  $\dot\lim\{f_i\}=\pair{f_i\fQ^1p_i}$, which then is the appropriate map in
  $\Gray\Cat_{\fQ^1}$. On the other hand, $\lim{f_i}$ is induced by the cone
$f_ir_i$. By universality $\dot\lim{f_i}=\lim{f_i}\pair{\fQ^1p_i}$. 
\end{rem}

In particular this applies to pullbacks, that is, there is a
canonical map 
\begin{equation*}
f\dot\times{}g\from\G\times_\K\H\laxto\G'\times_{\K'}\H'
\end{equation*}
determined by $f, g, h$ in
\begin{equation}
  \label{eq:pbprodpseud}
  \begin{xy}
    \xyboxmatrix{
      {}&{}&\H\ar[ddl]|!{"2,1";"2,4"}{\hole}^(.7){a}\ar[rrr]^{g}|*\dir{/}&{}&{}&\H'\ar[ddl]^{a'}\\
      \G\ar[dr]_{b}\ar[rrr]^(.3){f}|(.3)*\dir{/}&{}&{}&\G'\ar[dr]_{b'}&{}&{}\\
      {}&\K\ar[rrr]_{h}|*\dir{/}&{}&{}&\K'&{}
    }
\end{xy}\,.
\end{equation}

\begin{rem}
  If in (\ref{eq:indmap1}) the maps $f_i$ are of the form $g_ie$,
  i.e.\@ the $f_i$ come from strict maps, then we have 
  \begin{equation*}
    \dot\lim(g_ie)=(\lim{g_i})e\,.
  \end{equation*}
  In particular in a situation analogous to (\ref{eq:pbprodpseud}) we
  have
  \begin{equation}
    \label{eq:pbprodstr}
    (fe)\dot\times(ge)=(f\times{}g)e
  \end{equation}
\end{rem}

\subsection{Special Cells in the Resolved Space}
\label{sec:spectcel}

We now take a closer look at the structure of $\fQ^1\G$. By definition
1-cells here are non-empty lists $[f_1,\ldots, f_n]$ of composable
$\G$-1-cells modulo insertion or removal of identity 1-cells of $\G$;
composition is concatenation. For composable 1-cells in $\G$, say,
$f_1,\ldots, f_n$ we have several 1-cells in $\fQ^1\G$, in particular
$[f_1,\ldots, f_n]=[f_1]\#_0\cdots\#_0[f_n]$ and $[f_1\#_0\cdots\#_0
f_n]$ and $e_\G$ maps all of these to $f_1\#_0\cdots\#_0
f_n$. Between $[f_1,\ldots, f_n]$ and
$[f_1\#_0\cdots\#_0 f_n]$ we have a 2-cell
\begin{equation*}
  \kappa_{f_1,\ldots, f_n}=(\id_{f_1\#_0\cdots\#_0 f_n}; [f_1,\ldots, f_n], [f_1\#_0\cdots\#_0 f_n])
\end{equation*}
that is the pulled back identity 2-cell of $f_1\#_0\cdots\#_0 f_n$. In
particular we have 
\begin{equation*}
  \begin{xy}\xycompile{
    \xymatrix@+1cm{
      \ar[r]^{[f_2]} \ar[dr]_{[f_1\#_0 f_2]}|{}="y"& \ar[d]^{[f_1]}\\
      {} & {}
      \ar@2"1,2";"y"**{}?(.3);?(.9)|{\kappa_{f_1, f_2}}
    }
  }\end{xy}
\end{equation*}
for all for all pairs $f_1, f_2$ of 1-cells of $\G$. Whiskers and composites of higher cells in
$\fQ^1\G$ are simply carried out in $\G$, hence for example 
\begin{align*}
  \kappa_{f_1, f_2}\#_0[f_3] &=\left(\id_{f_1\#_0 f_2}\#_0 f_3; [f_1,
    f_2]\#_0[f_3], [f_1\#_0 f_2]\#_0[f_3]\right)\\
  &=\left(\id_{f_1\#_0 f_2\#_0 f_3}; [f_1,f_2,f_3], [f_1\#_0 f_2, f_3]\right)
\end{align*}
and
\begin{align*}
  \kappa_{f_1\#_0 f_2, f_3}\#_1\left(\kappa_{f_1, f_2}\#_0[f_3]\right)
  =\left(\id_{f_1\#_0 f_2\#_0 f_3}; [f_1,f_2,f_3], [f_1\#_0 f_2\#_0 f_3]\right)=\kappa_{f_1,f_2,f_3}\,.
\end{align*}
Hence we obtain that
\begin{equation}
  \label{eq:kappaass}
  \begin{xy}\xycompile{
    \xymatrix@+2cm{
      [f_1]\#_0[f_2]\#_0[f_3] \ar@2[r]^-{[f_1]\#_0\kappa_{f_2,f_3}}
      \ar@2[d]_-{\kappa_{f_1, f_2}\#_0[f_3]}\ar@2[dr]|{\kappa_{f_1,f_2,f_3}}& 
      [f_1]\#_0[f_2\#_0 f_3]\ar@2[d]^-{\kappa_{f_1, f_2\#_0f_3}}\\
      [f_1\#_0f_2]\#_0[f_3]\ar@2[r]_-{\kappa_{f_1\#_0f_2,f_3}} & [f_1\#_0f_2\#_2f_3]
    }
  }\end{xy}
\end{equation}
commutes.

We consider the possible
horizontal composites of $\kappa_{f_1,f_2}$ and $\kappa_{f_3,f_4}$ and
their tensor:
\begin{equation*}
  \begin{xy}\xycompile{
    \xymatrix@+1cm"A"{\ar@/^1pc/[r]^{[f_3, f_4]}="x" \ar@/_1pc/[r]_{[f_3\#_0f_4]}="y" 
      \POS \ar@2 "x";"y" **{} ?(.3); ?(.7) _{\kappa_{f_3,f_4}}="X" &
      \ar@/^1pc/[r]^{[f_1,f_2]}="x" \ar@/_1pc/[r]_{[f_1\#_0f_2]}="y" 
      \POS \ar@2 "x";"y" **{} ?(.3); ?(.7) ^{\kappa_{f_1,f_2}}="Y" & {}
      \tria "X"; "Y"}
    \POS +(70,0)
    \xymatrix"B"{\ar@/^1pc/[r]^{[f_3, f_4]}="x" \ar@/_1pc/[r]_{[f_3\#_0f_4]}="y" 
      \POS \ar@2 "x";"y" **{} ?(.3); ?(.7) _{\kappa_{f_3,f_4}}="X" &
      \ar@/^1pc/[r]^{[f_1,f_2]}="x" \ar@/_1pc/[r]_{[f_1\#_0f_2]}="y" 
      \POS \ar@2 "x";"y" **{} ?(.3); ?(.7) ^{\kappa_{f_1,f_2}}="Y" & {}
      \tria "Y"; "X"}
    \POS \ar@3 "A1,3";"B1,1" **{} ?(.3); ?(.7) ^{\kappa_{f_1,f_2}\ten\kappa_{f_3,f_4}}
  }\end{xy}\,.
\end{equation*}

By \eqref{eq:tenpbdef} we obtain
\begin{multline*}
  \kappa_{f_1,f_2}\ten\kappa_{f_3,f_4}= 
  (\id_{f_1\#_0f_2}; [f_1, f_2],[f_1\#_0f_2])\ten(\id_{f_3\#_0f_4}; [f_3, f_4], [f_3\#_0f_4])\\
  =\begin{pmatrix} 
    \id_{f_1\#_0f_2}\ten\id_{f_3\#_0f_4};\\
    (\id_{f_1\#_0f_2}\#_0e[f_3\#_0f_4])\#_1(e[f_1,f_2]\#_0\id_{f_3\#_0f_4}),\\
    (e[f_1\#_0f_2]\#_0\id_{f_3\#_0f_4})\#_1(\id_{f_1\#_0f_2}\#_0e[f_3,f_4]);\\
    [f_1,f_2,f_3,f_4],[f_1\#_0f_2,f_3\#_0f_4]
  \end{pmatrix}\\  
  =\begin{pmatrix} 
    \id_{\id_{f_1\#_0f_2\#_0f_3\#_0f_4}};\\
    (\id_{f_1\#_0f_2}\#_0f_3\#_0f_4)\#_1(f_1\#_0f_2\#_0\id_{f_3\#_0f_4}),\\
    (f_1\#_0f_2\#_0\id_{f_3\#_0f_4})\#_1(\id_{f_1\#_0f_2}\#_0f_3\#_0f_4);\\
    [f_1,f_2,f_3,f_4],[f_1\#_0f_2,f_3\#_0f_4]
  \end{pmatrix}\\  
  =\begin{pmatrix} 
    \id_{\id_{f_1\#_0f_2\#_0f_3\#_0f_4}};\\
    (\id_{f_1\#_0f_2\#_0f_3\#_0f_4})\#_1(\id_{f_1\#_0f_2\#_0f_3\#_0f_4}),\\
    (\id_{f_1\#_0f_2\#_0f_3\#_0f_4})\#_1(\id_{f_1\#_0f_2\#_0f_3\#_0f_4});\\
    [f_1,f_2,f_3,f_4],[f_1\#_0f_2,f_3\#_0f_4]
  \end{pmatrix}\\
  =\begin{pmatrix} 
    \id_{\id_{f_1\#_0f_2\#_0f_3\#_0f_4}};\\
    \id_{f_1\#_0f_2\#_0f_3\#_0f_4},\\
    \id_{f_1\#_0f_2\#_0f_3\#_0f_4};\\
    [f_1,f_2,f_3,f_4],[f_1\#_0f_2,f_3\#_0f_4]
  \end{pmatrix}\,,
\end{multline*}
meaning that this tensor is the identity of the two possible
horizontal composites of $\kappa_{f_1,f_2}$ and $\kappa_{f_3,f_4}$. 

Finally, note that by construction the $\kappa_{f_1,\ldots, f_n}$ are all
invertible.

\subsection{Pseudo Maps Explicitly}

We provide an elementary characterization of pseudo $\Gray$-functors.

\begin{defn}\label{defn:psgrmap}
  A \defterm{pseudo $\fQ^1$ graph map} $F\from \G\to\H$ between
  $\Gray$-categories is a map of 3-globular sets, together with a function
  $F^2\from \G_1\times_{\G_0}\G_1\to \H_2$, such that the
  following conditions hold: 
  \begin{enumerate}
  \item \label{item:psgrmaplocal} the restriction of $F$ to $\G(x, y)$ is a sesquifunctor
    for all 0-cells $x, y$ of $\G$,
  \item $F^2$ is a normalized 2-cocycle, that is, the $F^2_{f_1, f_2}$ are invertible 2-cells 
    $F^2_{f_1,f_2}\from
    F(f_1)\#_0F(f_2)\To F(f_1\#_0f_2)$ with
    \begin{equation}\label{eq:cocycle}
      F^2_{f_1, f_2\#_0f_3}\#_1(F(f_1)\#_0F^2_{f_2,f_3})=F^2_{f_1\#_0f_2,f_3}\#_1(F^2_{f_1,f_2}\#_0F(f_3)),
    \end{equation}
    and for $f_1$ or $f_2$ an identity 1-cell we have
    \begin{equation*}
      F^2_{f_1,f_2}=\id_{Ff_1\#_0Ff_2}, 
    \end{equation*}
  \item left and right whiskers of 2-cells by 1-cells along 0-cells are coherently preserved:
    \begin{align}
      F(\alpha\#_0f)\#_1F^2_{g, f}&=F^2_{g',f}\#_1(F\alpha\#_0Ff)\label{eq:coccoh}\\
      F(g\#_0\beta)\#_1F^2_{g, f}&=F^2_{g, f'}\#_1(Fg\#_0F\beta)\notag
    \end{align}  
  \item left and right whiskers of 3-cells by 1-cells along 0-cells are coherently preserved:
    \begin{align}
      F(\Gamma\#_0f)\#_1F^2_{g, f}&=F^2_{g',f}\#_1(F\Gamma\#_0Ff)\label{eq:coccoh3whsk}\\
      F(g\#_0\Delta)\#_1F^2_{g, f}&=F^2_{g, f'}\#_1(Fg\#_0F\Delta)\notag
    \end{align}
  \item the tensor is coherently preserved:
    \begin{equation}\label{eq:coccohten}
      F(\beta\ten\alpha)\#_1F^2_{g, f}=F^2_{g',f'}\#_1(F\beta\ten F\alpha)
    \end{equation}
  \item the tensors of compositors are trivial:
    \begin{equation}\label{eq:tentriv}
      \left(\begin{xy}\xycompile{
            \xymatrix@+1cm{
              F^2_{f_1,f_2}\lhc F^2_{f_3, f_4}
              \ar@3[r]^-{F^2_{f_1,f_2}\ten F^2_{f_3, f_4}}&   F^2_{f_1,f_2}\rhc F^2_{fe_3, f_4}
            }
          }\end{xy}\right) =\id
    \end{equation}
  \item tensors of 2-co-cycle elements with images of 2-cells vanish:
    \begin{align}
      \label{eq:ten2coccel-l}
      \left(\begin{xy}\xycompile{
            \xymatrix@+1cm{
              F\alpha\lhc F^2_{g,f}
              \ar@3[r]^-{F\alpha\ten F^2_{g,f}}&F\alpha\rhc F^2_{g,f}
            }
          }\end{xy}\right) &=\id  \\
      \label{eq:ten2coccel-r}
      \left(\begin{xy}\xycompile{
            \xymatrix@+1cm{
              F^2_{h,g}\lhc F\beta
              \ar@3[r]^-{F^2_{h,g}\ten F\beta}&F^2_{h,g}\rhc F\beta
            }
          }\end{xy}\right) &=\id
    \end{align}
  \end{enumerate}
  for all suitably incident cells.
  Denote the set of all pseudo $\fQ^1$-graph maps from $\G$ to $\H$ by
  $M(\G,\H)$. 
\end{defn}
Note also how the identity 1-cells of a 0-cells are preserved
strictly, this is part of the globularity condition.

Note furthermore how this definition implies that the horizontal composites are
also coherently preserved as a consequence of (\ref{eq:coccoh}):
\begin{align*}
  F(\alpha\lhc\beta)\#_1F^2_{g,f}&=F^2_{g',f'}\#_1(F\alpha\lhc{}F\beta)\\
  F(\alpha\rhc\beta)\#_1F^2_{g,f}&=F^2_{g',f'}\#_1(F\alpha\rhc{}F\beta)\,.\notag
\end{align*}

\begin{lem}\label{lem:pseudofunchar}
  There is a canonical correspondence between the set of pseudo $\fQ^1$ graph
  maps $M(\G,\H)$ and co-Kleisli maps $\Gray\Cat_{\fQ^1}(\G,\H)$.
  \begin{equation*}
    \begin{xy}
      \xymatrix{
        M(\G,\H)\ar@/^1.5pc/[r]^{(\_)^{\tilde{}}}&
        \Gray\Cat_{\fQ^1}(\G,\H)\ar@/^1.5pc/[l]^{(\_)^\vee}
      }
    \end{xy}
  \end{equation*}
\end{lem}
\begin{prf}
  Given a $\fQ^1$ graph map $F\from\G\to\H$ we define a
  $\Gray$-functor $\tilde F\from\fQ^1\G\to\H$ as follows
  \begin{enumerate}
  \item 0-cells:
    \begin{equation*}
      \tilde F(x)=F(x),
    \end{equation*}
  \item 1-cells: 
    \begin{equation*}
      \tilde F[f_1,\ldots,f_n]=Ff_1\#_0\cdots \#_0Ff_n,  
    \end{equation*}
  \item 2-cells: 
    \begin{equation}\label{eq:twocdef}
      \tilde F(\alpha;[f_1,\ldots, f_n],[g_1,\ldots, g_m])
      =\overline{\tilde F\kappa_{g_1,\ldots,
          g_m}}\#_1F\alpha\#_1\tilde F\kappa_{f_1,\ldots, f_n}
    \end{equation}
    where for $n=2$ the 2-cell $\tilde F\kappa_{f_1,\ldots, f_n}$ is defined
    as $F^2_{f_1,f_2}$ and for $n\geq 3$ as the unique extension due
    to \eqref{eq:cocycle}, \eqref{eq:tentriv},  
  \item 3-cells:
    \begin{equation*}
      \tilde F(\Gamma;\alpha, \beta; [f_1,\ldots, f_n],[g_1,\ldots, g_m])
      =\overline{\tilde F\kappa_{g_1,\ldots,
          g_m}}\#_1F\Gamma\#_1\tilde F\kappa_{f_1,\ldots, f_n}\,.
    \end{equation*}
  \end{enumerate}
  To elucidate, we show that 1-2-whiskers are preserved by $\tilde
  F$. For whiskerable cells  
  \begin{equation*}
    \begin{xy}\xycompile{
      \xymatrix@+1cm"A"{\ar[r]^{[f_1,\ldots, f_n]} &
        \ar@/^1pc/[r]^{[g_1,\ldots, g_m]}="x" \ar@/_1pc/[r]_{[g'_1,\ldots, g'_{m'}]}="y" 
        \POS \ar@2 "x";"y" **{} ?(.3); ?(.7) ^{(\beta;\ldots)}="Y" & {}
      }
    }\end{xy}
  \end{equation*}
  the equation
  \begin{multline*}
    \begin{xy}\xycompile{
      \xymatrix@+1cm"A"{\ar[r]^{\tilde F[f_1,\ldots, f_n]} &
        \ar@/^1pc/[r]^{\tilde F[g_1,\ldots, g_m]}="x"
        \ar@/_1pc/[r]_{\tilde F[g'_1,\ldots, g'_{m'}]}="y" 
        \POS \ar@2@<-1em> "x";"y" **{} ?(.3); ?(.7) ^{\tilde F(\beta;\ldots)}="Y" & {}
      }
    }\end{xy}=
    \begin{xy}\xycompile{
      \xymatrix@+2cm"A"{\ar[r]^{Ff_1\#_0\cdots\#_0Ff_n} &
        \ar@/^4pc/[r]^{Fg_1\#_0\cdots\#_0Fg_m}="x"
        \ar@/^1pc/[r]^{F(g_1\#_0\cdots\#_0g_m)}="y"
        \ar@/_1pc/[r]_{F(g'_1\#_0\cdots\#_0g'_{m'})}="z" 
        \ar@/_4pc/[r]_{Fg'_1\#_0\cdots\#_0Fg'_{m'}}="w" 
        \POS 
        \ar@2 "x";"y" **{} ?(.3); ?(.7) ^{F\kappa_{g_1,\ldots, g_m}}
        \ar@2 "y";"z" **{} ?(.3); ?(.7) ^{F\beta}
        \ar@2 "z";"w" **{} ?(.3); ?(.7) ^{\overline{F\kappa_{g'_1,\ldots, g'_{m'}}}}& {}
      }
    }\end{xy}\\=
    \begin{xy}\xycompile{
      \xymatrix@+5cm"A"{
        {}
        \ar@/^4pc/[r]^{Fg_1\#_0\cdots\#_0Fg_m\#_0Ff_1\#_0\cdots\#_0Ff_n}="x"
        \ar@/^1pc/[r]^{F(g_1\#_0\cdots\#_0g_m\#_0f_1\#_0\cdots\#_0f_n)}="y"
        \ar@/_1pc/[r]_{F(g'_1\#_0\cdots\#_0g'_{m'}\#_0f_1\#_0\cdots\#_0f_n)}="z" 
        \ar@/_4pc/[r]_{Fg'_1\#_0\cdots\#_0Fg'_{m'}\#_0Ff_1\#_0\cdots\#_0Ff_n}="w" 
        \POS 
        \ar@2 "x";"y" **{} ?(.3); ?(.7) ^{F\kappa_{g_1,\ldots, g_m,f_1,\ldots,f_n}}
        \ar@2@<-1em> "y";"z" **{} ?(.3); ?(.7) ^{F(\beta\#_0f_1\#_0\cdots\#_0f_n)}
        \ar@2 "z";"w" **{} ?(.3); ?(.7) ^{\overline{F\kappa_{g'_1,\ldots, g'_{m'},f_1,\ldots,f_n}}}& {}
      }
    }\end{xy}=
    \begin{xy}\xycompile{
      \xymatrix@+4cm"A"{
        \ar@/^2pc/[r]^{\tilde F([g_1,\ldots, g_m]\#_0[f_1,\ldots, f_n])}="x"
        \ar@/_2pc/[r]_{\tilde F([g'_1,\ldots, g'_{m'}]\#_0[f_1,\ldots, f_n]).}="y" 
        \POS \ar@2 "x";"y" **{} ?(.3); ?(.7) |{\tilde F((\beta;\ldots)\#_0[f_1,\ldots, f_n])}="Y" & {}
      }
    }\end{xy}
  \end{multline*}
  is a consequence of \eqref{eq:twocdef}. 
  
  Similarly, we can verify that $\tilde{F}$ preserves tensors: We calculate
  \begin{multline*}
    \tilde{F}\left((\beta; [g_1,\ldots,g_m], [g'_1,\ldots,g'_{m'}])\ten(\alpha; [f_1,\ldots,f_n], [f'_1,\ldots,f'_{n'}])\right)\\=
    \tilde{F}\left(
      \beta\ten\alpha;
      \beta\lhc\alpha, \beta\rhc\alpha;
      [g_1,\ldots,g_m,f_1,\ldots,f_n], [g'_1,\ldots,g'_{m'},f'_1,\ldots,f'_{n'}]
    \right)\\
    =\overline{\tilde{F}\kappa_{g'_1,\ldots,g'_{m'},f'_1,\ldots,f'_{n'}}}\#_1F(\beta\ten\alpha)\#_1\tilde{F}_{g_1,\ldots,g_m,f_1,\ldots,f_n}\\
    =(\overline{\tilde{F}\kappa_{g'_1,\ldots,g'_{m'}}}\ten\overline{\tilde{F}\kappa_{f'_1,\ldots,f'_{n'}}})
    \#_1(F\beta\ten{}F\alpha)\#_1(\tilde{F}_{g_1,\ldots,g_m}\ten\tilde{F}_{f_1,\ldots,f_n})\\
    =(\overline{\tilde{F}\kappa_{g'_1,\ldots,g'_{m'}}}\#_1F\beta\#_1\tilde{F}_{g_1,\ldots,g_m})
    \ten(\overline{\tilde{F}\kappa_{f'_1,\ldots,f'_{n'}}}\#_1F\alpha\#_1\tilde{F}_{f_1,\ldots,f_n})\\
    \tilde{F}(\beta; [g_1,\ldots,g_m], [g'_1,\ldots,g'_{m'}])\ten\tilde{F}(\alpha; [f_1,\ldots,f_n], [f'_1,\ldots,f'_{n'}])
  \end{multline*}
  using (\ref{eq:coccohten}) and (\ref{eq:tentriv}).
  Preservation of the
  remaining operations is equally simple to verify. 

  Conversely, given a $\Gray$-functor $G\from\fQ^1\G\to\H$ we define a
  pseudo $\fQ^1$ graph map $\check{G}\from\G\to\H$ as follows:
  \begin{enumerate}
  \item 0-cells: $\check{G}(x)=G(x)$
  \item 1-cells: $\check{G}(f)=G[f]$
  \item 2-cells:
    $\check{G}(\alpha)=G(\alpha;[f],[f'])$
  \item 3-cells:
    $\check{G}(\Gamma)=G(\Gamma;\alpha,\beta;[f],[f'])$
  \item 2-co-cycle:
    $\check{G}^2_{f_1,f_2}=G\kappa_{f_1,f_2}=G(\id_{f_1\#_0f_2};[f_1\#_0f_2],[f_1,f_2])$
  \end{enumerate}
  This is obviously locally a sesquifunctor. We check the co-cycle
  condition:
  \begin{multline*}
    \check{G}^2_{f_1,f_2\#_0f_3}\#_1(\check{G}f_1\#_0\check{G}^2_{f_2,f_3})\\
    =G(\id_{f_1\#_0f_2\#_0f_3};[f_1,f_2\#_0f_3],[f_1\#_0f_2\#_0f_3])
    \#_1(G[f_1]\#_0G(\id_{f_2\#_0f_3};[f_2,f_3],[f_2\#_0f_3]))\\
    =G(\id_{f_1\#_0f_2\#_0f_3};[f_1,f_2\#_0f_3],[f_1\#_0f_2\#_0f_3])
    \#_1G(\id_{f_1\#_0f_2\#_0f_3};[f_1,f_2,f_3],[f_1,f_2\#_0f_3])\\
    =G(id_{f_1\#_0f_2\#_0f_3};[f_1,f_2,f_3],[f_1\#_0f_2\#_0f_3])\\    
    =G(\id_{f_1\#_0f_2\#_0f_3};[f_1\#_0f_2,f_3],[f_1\#_0f_2\#_0f_3])
    \#_1G(\id_{f_1\#_0f_2\#_0f_3};[f_1,f_2,f_3],[f_1\#_0f_2,f_3])\\
    =G(\id_{f_1\#_0f_2\#_0f_3};[f_1\#_0f_2,f_3],[f_1\#_0f_2\#_0f_3])
    \#_1(G(\id_{f_1\#_0f_2};[f_1,f_2],[f_1\#_0f_2])\#_0G[f_3])\\
    =\check{G}^2_{f_1\#_0f_2,f_3}\#_1(\check{G}^2_{f_1,f_2}\#_0\check{G}f_3)
  \end{multline*}
  Furthermore, we check the coherent preservation of whiskers:
  \begin{multline*}
    \check{G}(\alpha\#_0f)\#_1\check{G}^2_{g,f}\\
    =G(\alpha\#_0f;[g\#_0f],[g'\#_0f])\#_1G(\id_{g\#_0f};[g,f],[g\#_0f])\\
    =G(\alpha\#_0f;[g,f],[g'\#_0f])\\
    =G(\id_{g'\#_0f};[g',f],[g'\#_0f])\#_1G(\alpha\#_0;[g,f],[g',f])\\
    =G(\id_{g'\#_0f};[g',f],[g'\#_0f])\#_1(G(\alpha;[g],[g'])\#_0G[f])\\
    =\check{G}^2_{g',f}\#_1(\check{G}\alpha\#_0\check{G}f)
  \end{multline*}
  The remaining axioms are verified just as easily.

  We verify briefly that $\tilde{\check{G}}=G$, for 1-cells we have
  \begin{equation*}
    \tilde{\check{G}}[f_1,\ldots,f_n]=\check{G}f_1\#_0\ldots\#_0\check{G}f_n
    =G[f_1]\#_0\ldots\#_0G[f_n]=G[f_1,\ldots,f_n]
  \end{equation*}
  and for 2-cells:
  \begin{multline*}
    \tilde{\check{G}}(\alpha;[f_1,\ldots,f_n],[f'_1,\ldots,f'_{n'}])
    =\overline{\check{G}\kappa_{f'_1,\ldots,f'_{n'}}}\#_1\check{G}\#_1\check{G}\kappa_{f_1,\ldots,f_{n}}\\
    =
    \begin{pmatrix}
      G(\id_{f'_1\#_0\cdots\#_0f'_{n'}};[f'_1\#_0\cdots\#_0f'_{n'}],[f'_1,\ldots,f'_{n'}])
      \\\#_1G(\alpha;[f'_1\#_0\cdots\#_0f'_{n'}],[f_1\#_0\cdots\#_0f_{n}])
      \\\#_1G(\id_{f_1\#_0\cdots\#_0f_{n}};[f_1,\ldots,f_{n}],[f_1\#_0\cdots\#_0f_{n}])
    \end{pmatrix}\\
    G(\alpha;[f_1,\ldots,f_n],[f'_1,\ldots,f'_{n'}])
  \end{multline*}

  Finally, $\check{\tilde{F}}=F$.
\end{prf}

\begin{rem}
  \label{rem:psq1grmap}
  Given two pseudo $\fQ^1$ graph maps $F\from\G\to\H$ and
  $G\from\H\to\K$ their composite $GF$ is simply the composite of the
  underlying globular maps with cocycle
  \begin{equation*}
    (GF)^2_{f_1,f_2}=GF^2_{f_1,f_2}\#_1G^2_{Ff_1,Ff_2}\,.
  \end{equation*}
\end{rem}

\begin{lem}
  \label{lem:psq1strictify}
  Under the correspondence in lemma \ref{lem:pseudofunchar} a pseudo
  $\fQ^1$-graph map $F$ has trivial cocycle $F^2$ iff the
  corresponding $\Gray$-functor $\tilde{F}$ is of the form $Ge$. 
\end{lem}
\begin{prf}
  Considering definition \ref{defn:psgrmap} we see that
  $F\in{}M(\G,\H)$ is an ordinary $\Gray$-functor iff $F^2$ is
  trivial, in which case $Fe$ is the embedding of $F$ in
  $\Gray\Cat_{\fQ^1}$ with
  ${(Fe)^\vee}^2_{f_1,f_2}=Fe\kappa_{f_1,f_2}=Fe(\id_{f_1\#_0f_2};[f_1\#_0f_2],[f_1,f_2])
  =F\id_{f_1\#_0f_2}=\id_{F(f_1\#_0f_2)}$. That is actually $G=F$. 

  In turn, if we are given a co-Kleisli map $Ge$ with $G$ a
  $\Gray$-functor we obtain
  ${(Ge)^\vee}^2_{f_1,f_2}=Ge\kappa_{f_1,f_2}=\id_{G(f_1\#_0f_2)}$.
\end{prf}

In particular for $\G$ free up to order 1 with section $k$
(\ref{eq:freeorderstrictify}) induces an idempotent map
\begin{equation}
  \label{eq:q1strictify}
  \begin{xy}
    \xyboxmatrix{
      M(\G,\H)\ar[r]^{(\tilde{(\_)}ke)^\vee}&M(\G,\H)
    }
  \end{xy}
\end{equation}
with image $\Gray\Cat(\G,\H)$.

We spell out the action of this map on an arbitrary pseudo $\fQ^1$
graph map $F\from\G\to\H$ for $\G$, free up to order 1, at the level of
1- and 2-cells. Let $f_1=g_{1,1}\#_0\cdots\#_0{}g_{1,n_{1}}$ and
$f_2=g_{2,1}\#_0\cdots{}\#_0g_{2,n_{2}}$ be unique decompositions up to
units in $\G$ of the 1-cells $f_1, f_2$. This means that
$k(f_1)=[g_{1,1},\ldots,g_{1,n_1}]$,
$k(f_2)=[g_{2,1},\ldots,g_{2,n_2}]$. Furthermore, for a 2-cell
$\alpha\from{}f\To{}f'$ we have
$k(\alpha)=(\alpha;[g_{1,1},\ldots,g_{1,n}],[g'_{1,1},\ldots,g'_{1,n'}])$,
in particular
$k(\id_f)=(\id_f;[g_{1,1},\ldots,g_{1,n}],[g_{1,1},\ldots,g_{1,n}])$. Hence
for a composite we get
\begin{multline}\label{eq:q1strictify1cell}
  (\tilde{F}ke)^\vee(f_1\#_0f_2)\\
  =(\tilde{F}ke)[f_1\#_0f_2]=\tilde{F}k(f_1\#_0f_2)=\tilde{F}[g_{1,1},\ldots,g_{1,n_1},g_{2,1},\ldots,g_{2,n_2}]\\
  =Fg_{1,1}\#_0\cdots,Fg_{1,n_1}\#_0Fg_{2,1}\#_0\cdots\#_0g_{2,n_2}\\
  =(\tilde{F}ke)^\vee(f_1)\#_0(\tilde{F}ke)^\vee(f_2)\,.
\end{multline}
For the 2-cocycle we get 
\begin{multline}\label{eq:q1strictify2coc}
  {(\tilde{F}ke)^\vee}^2_{f_1,f_2}=
  {(\tilde{F}ke)^\vee}(\kappa_{f_1,f_2})=\tilde{F}k(\id_{f_1\#_0f_2})\\
  =\tilde{F}(\id_{f_1\#_0f_2};[g_{1,1},\ldots,g_{1,n_1},g_{2,1},\ldots,g_{2,n_2}],[g_{1,1},\ldots,g_{1,n_1},g_{2,1},\ldots,g_{2,n_2}])\\
  =\overline{\tilde{F}\kappa_{g_{1,1},\ldots,g_{1,n_1},g_{2,1},\ldots,g_{2,n_2}}}
  \#_1F\id_{f_1\#_0f_2}
  \#_1\tilde{F}\kappa_{g_{1,1},\ldots,g_{1,n_1},g_{2,1},\ldots,g_{2,n_2}}\\
    \overline{\tilde{F}\kappa_{g_{1,1},\ldots,g_{1,n_1},g_{2,1},\ldots,g_{2,n_2}}}
  \#_1\id_{F(f_1\#_0f_2)}
  \#_1\tilde{F}\kappa_{g_{1,1},\ldots,g_{1,n_1},g_{2,1},\ldots,g_{2,n_2}}\\
  \overline{\tilde{F}\kappa_{g_{1,1},\ldots,g_{1,n_1},g_{2,1},\ldots,g_{2,n_2}}}
  \#_1\tilde{F}\kappa_{g_{1,1},\ldots,g_{1,n_1},g_{2,1},\ldots,g_{2,n_2}}\\
  =\id_{\tilde{F}[g_{1,1},\ldots,g_{1,n_1},g_{2,1},\ldots,g_{2,n_2}]}\\
  =\id_{Fg_{1,1}\#_0\cdots\#_0Fg_{1,n_1}\#_0Fg_{2,1}\#_0\cdots\#_0Fg_{2,n_2}}\,.
\end{multline}
These equations (\ref{eq:q1strictify1cell}) and
(\ref{eq:q1strictify2coc}) make it palpable how the operation
(\ref{eq:q1strictify}) yiels a strict $\Gray$-functor. 

We will see in section \ref{sec:inthomfct} how $F$ and it's
strictification $Fke$ are related.

\section{Path Spaces}
\label{sec:pathspc}

We construct a path space for $\Gray$-categories and prove some
essential properties. We derived the idea for this construction from
\citet{ben1967}. Maps into this space can be viewed as right
homotopies between functors and are our axiomatization of
transformation for morphisms in $\Gray\Cat_{\fQ^1}$. In section
\ref{sec:compaths} we will introduce an internal category structure
for this path space; its composition operation will allows us to
compose transformations.

\begin{defn}
  \label{defn:pathspc}
  Given a $\Gray$-category $\H$ we define the \defterm{path space}
  $\pathspc\H$ where the cells in each dimension are diagrams in $\H$:
  \begin{align}
    \pathspc\H_0=& \{\xymatrix{{} \ar[r]^f& {}}\}\\
    \pathspc\H_1=& \left\{(g_2;g_0,g_1,f,f')\vline
      \begin{matrix}
        \xymatrix{{} \ar[r]^f \ar[d]_{g_0}&
          {}\ar[d]^{g_1} \ar@2 [dl] **{} ?<(.3); ?>(.7) ^{g_2}\\{} \ar[r]_{f'}& {}}
      \end{matrix}
    \right\}\\
    \pathspc\H_2=& \left\{
      \begin{pmatrix}
        \alpha_3;\alpha_1,\alpha_2,g_2,h_2;\\g_0,g_1,h_0,h_1,f,f'
      \end{pmatrix}\vline\begin{matrix}
        \begin{xy}
          \xymatrix"A"{{} \ar[r]^f \ar[d]|{g_0}="x"\ar@/_3pc/[d]_{h_0}="y" &
            {}\ar[d]^{g_1} \ar@2 [dl] **{} ?(.2); ?(.8) |{g_2}\\{}
            \ar[r]_{f'}& {} \POS \ar@2 "x"; "y" **{} ?<(.3); ?>(.7) ^{\alpha_1} }
          \POS + (30,0)
          \xymatrix"B"{{} \ar[r]^f \ar[d]_{h_0}&
            {}\ar[d]|{h_1}="x" \ar@2 [dl] **{} ?(.2); ?(.8) |{h_2}
            \ar@/^3pc/[d]^{g_1}="y"\\{} \ar[r]_{f'}& {} \POS \ar@2
            "y"; "x" **{} ?<(.3); ?>(.7) ^{\alpha_2}}
          \POS \ar@3 {"A1,1" ; "A2,2" **{} ?(.5)}; {"B1,1" ; "B2,2" **{} ?(.5)} **{}
          ?<(.4);  ?>(.6) ^ {\alpha_3} 
        \end{xy}
      \end{matrix}
    \right\}\\
    \label{eq:3celleq}
    \pathspc\H_3=&\left\{
      \begin{pmatrix}
        \Gamma_1,\Gamma_2,\alpha_3,\beta_3; 
        g_2,h_2,\\
        \alpha_1,\alpha_2,\beta_1,\beta_2;\\
        g_0,g_1,h_0,h_1,
        f,f'
      \end{pmatrix}\vline
      \begin{pmatrix}
        \Gamma_1\from \alpha_1\Rrightarrow\beta_1,\\
        \Gamma_2\from\alpha_2\Rrightarrow\beta_2
      \end{pmatrix} \text{ such that }
      \begin{matrix}
        \beta_3\#_2((f'\#_0\Gamma_1)\#_1g_2)\\=(h'_2\#_1(\Gamma_2\#_0f))\#_2\alpha_3 
      \end{matrix}
    \right\}
  \end{align}
  Compositions and identities arise canonically from pasting of
  diagrams in $\H$, as detailed below.
 \end{defn}

The condition in \eqref{eq:3celleq} on the 3-cells is the commutativity of the following
diagram
\begin{equation}\label{eq:3cellcond}
  \begin{xy}\xycompile{
    \xyboxmatrix"A"{{} \ar[r]^f \ar[d]|{g_0}="x"\ar@/_3pc/[d]|{h_0}="y" &
      {}\ar[d]^{g_1} \ar@2 [dl] **{} ?(.2); ?(.8) |{g_2}\\{}
      \ar[r]_{f'}& {} \POS \ar@2 "x"; "y" **{} ?(.2); ?(.8) |{\alpha_1} }
    \POS + (40,0)
    \xyboxmatrix"B"{{} \ar[r]^f \ar[d]_{h_0}&
      {}\ar[d]|{h_1}="x" \ar@2 [dl] **{} ?(.2); ?(.8) |{h_2}
      \ar@/^3pc/[d]|{g_1}="y"\\{} \ar[r]_{f'}& {} \POS \ar@2
      "y"; "x" **{} ?(.2); ?(.8) |{\alpha_2}}
    \POS + (-40,-40)
    \xyboxmatrix"C"{{} \ar[r]^f \ar[d]|{g_0}="x"\ar@/_3pc/[d]|{h_0}="y" &
      {}\ar[d]^{g_1} \ar@2 [dl] **{} ?(.2); ?(.8) |{g_2}\\{}
      \ar[r]_{f'}& {} \POS \ar@2 "x"; "y" **{} ?(.2); ?(.8) |{\beta_1} }
    \POS + (40,0)
    \xyboxmatrix"D"{{} \ar[r]^f \ar[d]_{h_0}&
      {}\ar[d]|{h_1}="x" \ar@2 [dl] **{} ?(.2); ?(.8) |{h_2}
      \ar@/^3pc/[d]|{g_1}="y"\\{} \ar[r]_{f'}& {} \POS \ar@2
      "y"; "x" **{} ?(.2); ?(.8) |{\beta_2}}
    \ar@3 "A";"B" ^{\alpha_3} 
    \ar@3 "C";"D" _{\beta_3} 
    \ar@3 "A";"C" _{(f'\#_0\Gamma_1)\#_1g_2}
    \ar@3 "B";"D" ^{h_2\#_1(\Gamma_2\#_0f)} 
  }\end{xy}
\end{equation}

The identities in each dimension are obviously the ones
consisting of identity cells. 

\begin{rem}
  By construction the map $(d_0,d_1)\from\pathspc\H\to\H\times\H$ is
  2-faithful in the sense of definition \ref{defn:niso}, but in
  general not full.
\end{rem}

\begin{rem}
  The map $i\from\H\to\pathspc\H$ is 2-Cartesian and 1-faithful, but
  not in general 1-full.
\end{rem}

\subsection{Path Spaces and Cartesian Maps}
\label{sec:pathcart}

\begin{lem}
  The path space construction $\pathspc{(\_)}$ of $\Gray$-categories
  preserves 1-Cartesianness of maps. 
\end{lem}
\begin{prf}
  Let's as assume we have a situation
  \begin{equation*}
    \begin{xy}
      \xyboxmatrix{
        \pathspc\G\ar@<-.5ex>[d]_{d_0}\ar@<+.5ex>[d]^{d_1}\ar[r]^{\pathspc{F}}&
        \pathspc\H\ar@<-.5ex>[d]_{d_0}\ar@<+.5ex>[d]^{d_1}\\
        \G\ar[r]_{{F}}&\H
      }
    \end{xy}\,,
  \end{equation*} 
  take a pair of parallel 1-cells in $\pathspc\G$
  \begin{align*}
    &\xymatrix{
      {} \ar[r]^f \ar[d]_{g_0}&
      {}\ar[d]^{g_1} \ar@2 [dl] **{} ?<(.3); ?>(.7) ^{g_2}\\{}
      \ar[r]_{f'}& {}
    }
    &
    &\xymatrix{
      {} \ar[r]^f \ar[d]_{h_0}&
      {}\ar[d]^{h_1} \ar@2 [dl] **{} ?<(.3); ?>(.7) ^{h_2}\\{}
      \ar[r]_{f'}& {}
    }
  \end{align*}
  we need to show that $\pathspc{F}$ is bijective on the intervening
  2-cells. That means given
  \begin{align*}
    \beta_1&\from F(g_0)\To f(h_0) &   \beta_2&\from F(g_1)\To F(h_1) 
    &\beta_3&\from F(g_2\#_1(\beta_2\#_0f))\Tto
    F((f'\#_0\beta_1)\#_1g_2)
  \end{align*}
  there are unique
  \begin{align*}
    \alpha_1&\from g_0\To h_0 &   \alpha_2&\from g_1\To h_1 
    &\alpha_3&\from g_2\#_1(\alpha_2\#_0f)\Tto
    (f'\#_0\alpha_1)\#_1g_2
  \end{align*}
  with $F(\alpha_i)=\beta_i$. But these exist uniquely by the
  1-Cartesianness of $F$.

  The same kind of argument can be applied to parallel 2-cells in
  $\pathspc\G$. 
\end{prf}
 
\begin{rem} 
  The functor $\pathspc{(\_)}$ preserves 2-Cartesian maps.
\end{rem}

\begin{lem}
  A pullback of a Cartesian map is Cartesian if $p$ preserves pullbacks.
\end{lem}
\begin{prf}
  Let $F$ be $p$-Cartesian, and $G^*F$ the pullback of $F$ along $G$. 
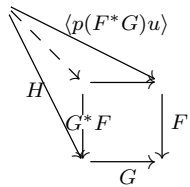
  \begin{equation*}
    \begin{xy}
      \xyboxmatrix{
        {}\ar@{-->}[dr]\ar[ddr]_H\ar[drr]^{\pair{p(F^*G)u}}&&\\
        {}&\ar[d]|{G^*F}\ar[r]&\ar[d]^F\\
        {}&\ar[r]_G&\\
      }
    \end{xy}
  \end{equation*}
  Let $H$ factor through $G$ below as $p(H)=p(G^*F)u$, then $GH$ factors
  through $F$ below as $p(GH)=p(GG^*F)u=p(F)p(F^*G)u$, hence there is
  a unique lift $\pair{p(F^*G)u}$. Hence there is a universally
  induced $\pair{u}$ with $G^*F\pair{u}=H$. 

  The functor $p$ preserving pullbacks ensures that $p\pair{u}=u$.
\end{prf}

\subsection{Vertical Composition Operations in the Path Space}

We need to describe the vertical composition of 1-, 2-, 3-cells along
0-, 1-, 2-cells respectively.

We designate the composition in $\H$ by
$\#_i$ and the interchange by $\ten$, in $\pathspc\H$ we define the
respective operations $\square_i$ and $\bten$ as follows:
\begin{equation*}
  h\square_0 g= (h_2;h_0,h_1,f'',f')\square_0(g_2;g_0,g_1,f,f')= 
  \begin{pmatrix}
    (h_2\#_0g_0)\#_1(h_1\#_0g_2);\\
    h_0\#_0g_0, h_1\#_0g_1,f,f''
  \end{pmatrix}
\end{equation*}
This is just the vertical pasting
\begin{equation}
  \label{eq:vertsquarepaste}
  \begin{matrix}
    \xymatrix{
      {}\ar[r]^f \ar[d]_{g_0}& \ar[d]^{g_1} \ar@2 [dl] **{} ?(.2); ?(.8) |{g_2}\\
      \ar[r]|{f'} \ar[d]_{h_0} & \ar[d]^{h_1} \ar@2 [dl] **{} ?(.2); ?(.8) |{h_2}\\
      \ar[r]_{f''} &{}
    }
  \end{matrix}\,.
\end{equation} Obviously this composition is associative and unital.
 
\begin{rem}
  \label{rem:sqhorcomp} Considering (\ref{eq:vertsquarepaste}) we note
  that if the 1-cells in $\H$ are invertible, with inverse
  $\overline{(\_)}$, then the 2-cell
  \begin{equation*}
    (h_2\#_0g_0)\#_1(h_1\#_0g_2)
  \end{equation*}
  in \eqref{eq:vertsquarepaste} can
  also be written as a horizontal composite in two different ways:
  \begin{equation}
    \label{eq:sqhorcomp}
    (h_2\#_0\overline{f'})\lhc{}g_2=h_2\lhc(\overline{f'}\#_0g_2)
  \end{equation}
  There is of course also the opposite horizontal composite 
  \begin{equation}
    \label{eq:sqhorcompopp}
    (h_2\#_0\overline{f'})\rhc{}g_2=h_2\rhc(\overline{f'}\#_0g_2)
  \end{equation}
  and a 3-cell
  \begin{equation*}
    (h_2\#_0\overline{f'})\ten{}g_2=h_2\ten(\overline{f'}\#_0g_2)
  \end{equation*}
  going from \eqref{eq:sqhorcomp} to \eqref{eq:sqhorcompopp}. The
  picture \eqref{eq:vertsquarepaste}, however,  always means
  \eqref{eq:sqhorcomp}.
\end{rem}

The vertical composite of two 2-cells is
\begin{multline}\label{eq:vert2compdef}
  \beta\square_1\alpha=
  \begin{pmatrix}
    \beta_3;\beta_1,\beta_2,h_2,k_2;\\
    h_0,h_1,k_0,k_1,f,f'
  \end{pmatrix}\square_1
  \begin{pmatrix}
    \alpha_3;\alpha_1,\alpha_2,g_2,h_2;\\g_0,g_1,h_0,h_1,f,f'
  \end{pmatrix}\\=
  \begin{pmatrix}
    (\beta_3\#_1(\alpha_2\#_0f))\#_2((f'\#_0\beta_1)\#_1\alpha_3);\\
    \beta_1\#_1\alpha_1, \beta_2\#_1\alpha_2,g_2,h_2;g_0,g_1,k_0,k_1,f,f'      
  \end{pmatrix}
\end{multline}
which has as its first component  the following composite of $\H$-3-cells
\begin{equation*}
  \begin{xy}\xycompile{
    \xyboxmatrix"A"{
      {} \ar[r]^f \ar[d]|{g_0}="x"
      \ar@/_2.5pc/[d]|{h_0}="y"
      \ar@/_5pc/[d]|{k_0}="z" &
      {}\ar[d]^{g_1} \ar@2 [dl] **{} ?(.2); ?(.8) |{g_2}\\{} \ar[r]_{f'}& {} 
      \POS \ar@2 "x"; "y" **{} ?(.2); ?(.8)|{\alpha_1} 
      \POS \ar@2 "y"; "z" **{} ?(.2); ?(.8)|{\beta_1} 
    } 
    \POS + (55, 0)
    \xyboxmatrix"B"{
      {}\ar[r]^f \ar[d]|{h_0}="x" \ar@/_2.5pc/[d]|{k_0}="y" &
      {}\ar[d]|{h_1}="w" \ar@/^2.5pc/[d]|{g_1}="z" 
      \ar@2 [dl] **{} ?(.2); ?(.8) |{h_2}\\{}
      \ar[r]_{f'}& {} 
      \POS \ar@2 "x"; "y" **{} ?(.2); ?(.8) |{\beta_1} 
      \POS \ar@2 "z"; "w" **{} ?(.2); ?(.8) |{\alpha_2} 
    } 
    \POS + (55, 0)
    \xyboxmatrix"C"{
      {} \ar[r]^f \ar[d]_{k_0}  &
      {}\ar[d]|{k_1}="z" \ar@2 [dl] **{} ?(.2); ?(.8) |{k_2}
      \ar@/^2.5pc/[d]|{h_1}="y" \ar@/^5pc/[d]|{g_1}="x" \\{}
      \ar[r]_{f'}& {} 
      \POS \ar@2 "x"; "y" **{} ?(.2); ?(.8)|{\alpha_2} 
      \POS \ar@2 "y"; "z" **{} ?(.2); ?(.8)|{\beta_2} 
    }
    \ar@3 "A";"B" ^{(f'\#_0\beta_1)\#_1\alpha_3}
    \ar@3 "B";"C" ^{\beta_3\#_1(\alpha_2\#_0f)}
  }\end{xy}\,.
\end{equation*}
We shall henceforth argue mostly diagrammatically in terms of such 3-cell
diagrams, as it is fairly obvious what the lower dimensional
components are.

Vertical composition of $\pathspc\H$-3-cells is particularly simple:
\begin{equation}
  \label{eq:vert3compdef}
  \Delta\square_2\Gamma
  =\begin{pmatrix}
    \Delta_1\from\beta_1\Tto\gamma_1,\\
    \Delta_2\from\beta_2\Tto\gamma_2
  \end{pmatrix}\square_2
  \begin{pmatrix}
    \Gamma_1\from \alpha_1\Tto\beta_1,\\
    \Gamma_2\from\alpha_2\Tto\beta_2
  \end{pmatrix}
  =\begin{pmatrix}
    \Delta_1\#_2\Gamma_1\from \alpha_1\Tto\gamma_1,\\
    \Delta_2\#_2\Gamma_2\from\alpha_2\Tto\gamma_2
  \end{pmatrix}
\end{equation}
The condition \eqref{eq:3cellcond} is obviously satisfied, since we just
paste two instances of the commuting square vertically.

\subsection{Whiskers}

We need to define three whiskering operations, ${}^1\square_0^2$,
${}^1\square_0^3$, ${}^2\square_1^3$, where the raised indices
indicate the dimension of the operands, the lower one the dimension of the
incidence cell. Their symmetry partners are then obvious.

We define right whiskering of a 2-cell by a 1-cell as:
\begin{multline}\label{eq:sq12def}
  k{}^1\square_0^2\alpha=  (k_2;k_0,k_1,f',f''){}^1\square_0^2
  \begin{pmatrix}
    \alpha_3;\alpha_1,\alpha_2;\\g_0,g_1,h_0,h_1,f,f'   
  \end{pmatrix}\\= \left(
    \begin{matrix}
      ((k_2\#_0h_0)\#_1(k_1\#_0\alpha_3))\\
      \#_2(\overline{(k_2\ten\alpha_1)}\#_1(k_1\#_0g_2));\\
      k_0\#_0\alpha_1, k_1\#_0\alpha_2;\\
      k_0\#_0g_0,k_1\#_1g_1,k_0\#_0h_0,k_1\#_0h_1,f,f''
    \end{matrix}
  \right)\,.
\end{multline}
Diagrammatically this is the following composite:
\begin{equation*}
  \begin{xy}\xycompile{ 
      (45,0):
    (0,0)
    \xyboxmatrix"A"{{} \ar[r]^f \ar[d]|{g_0}="x"\ar@/_2.5pc/[d]|{h_0}="y" &
      {}\ar[d]^{g_1} \ar@2 [dl] **{} ?(.2); ?(.8) |{g_2}\\{}
      \ar[r]|{f'} \ar[d]_{k_0}& {}\ar[d]^{k_1} \ar@2 [dl] **{} ?(.2); ?(.8) |{k_2}="Y" 
      \POS \ar@2 "x"; "y" **{} ?(.2); ?(.8)|{\alpha_1}="X" \\
      {}\ar[r]_{f''} & {}
      \tria "Y"; "X"
    }
    ,(1,0)
    \xyboxmatrix"B"{{} \ar[r]^f \ar[d]|{g_0}="x"\ar@/_2.5pc/[d]|{h_0}="y" &
      {}\ar[d]^{g_1} \ar@2 [dl] **{} ?(.2); ?(.8) |{g_2}\\{}
      \ar[r]|{f'} \ar[d]_{k_0}& {}\ar[d]^{k_1} \ar@2 [dl] **{} ?(.2); ?(.8) |{k_2}="Y" 
      \POS \ar@2 "x"; "y" **{} ?(.2); ?(.8)|{\alpha_1}="X" \\
      {}\ar[r]_{f''} & {}
       \tria "X"; "Y"
    }
    +(1,0)
    \xyboxmatrix"C"{{} \ar[r]^f \ar[d]_{h_0}&
      {}\ar[d]|{h_1}="x" \ar@2 [dl] **{} ?(.2); ?(.8) |{h_2}
      \ar@/^2.5pc/[d]|{g_1}="y"\\
      {} \ar[r]|{f'} \ar[d]_{k_0}& 
      {} \ar[d]^{k_1} \ar@2 [dl] **{} ?(.2); ?(.8) |{k_2} 
      \POS \ar@2 "y"; "x" **{} ?(.2); ?(.8) |{\alpha_2}\\
      {} \ar[r]_{f''}&{}
    }
    \ar@3 "A";"B" ^{\overline{(k_2\ten\alpha_1)}\\\;\#_1(k_1\#_0g_2)}    
    \ar@3 "B";"C" ^{(k_2\#_0h_0)\\\;\#_1(k_1\#_0\alpha_3)}    
  }\end{xy}
\end{equation*}

 For
reference  $(\beta_1,\beta_2,\beta_3)\square_0(h_0,h_1,h_2)$ is
\begin{equation*} 
  \begin{xy}\xycompile{
    (45,0) : (0,0)
    \xyboxmatrix"A"{
      {} \ar[r]^f \ar[d]_{h_0} & \ar[d]^{h_1}
      \ar@2[dl]**{}?(.2);?(.8)|{h_2}\\
      \ar[r]|{f'} \ar[d]|{k_0}="x" \ar@/_2.5pc/[d]|{m_0}="y"& \ar[d]
      ^{k_1} \ar@2[dl]**{}?(.2);?(.8)|{k_2}\\
      \ar[r]_{f''} & {}
      \POS\ar@2"x";"y" **{} ?(.2); ?(.8)|{\beta_1}
    }
    +(1,0)
    \xyboxmatrix"B"{
      {} \ar[r]^f \ar[d]_{h_0} & \ar[d]^{h_1}
      \ar@2[dl]**{}?(.2);?(.8)|{h_2}="X"\\
      \ar[r]|{f'} \ar[d]_{m_0} & 
      \ar[d]|{k_1}="y" \ar@2[dl]**{}?(.2);?(.8)|{k_2} 
      \ar@/^2.5pc/[d]|{k_1}="x"\\
      \ar[r]_{f''} & {}
      \POS\ar@2"x";"y" **{} ?(.2); ?(.8)|{\beta_2}="Y"
      \tria"X";"Y"
    }
    +(1,0)
    \xyboxmatrix"C"{
      {} \ar[r]^f \ar[d]_{h_0} & \ar[d]^{h_1}
      \ar@2[dl]**{}?(.2);?(.8)|{h_2}="X"\\
      \ar[r]|{f'} \ar[d]_{m_0} & 
      \ar[d]|{k_1}="y" \ar@2[dl]**{}?(.2);?(.8)|{k_2} 
      \ar@/^2.5pc/[d]|{k_1}="x"\\
      \ar[r]_{f''} & {}
      \POS\ar@2"x";"y" **{} ?(.2); ?(.8)|{\beta_2}="Y"
      \tria"Y";"X"
    }
    \ar@3"A";"B"^{(h_2\#_0k_1)\\\;\#_1(h_0\#_0\beta_3)}    
    \ar@3"B";"C"^{(m_2\#_0h_0)\\\;\#_1(h_2\ten\beta_2)}
  }\end{xy}
\end{equation*}

The action of 1-cells on 3-cells is as follows:
\begin{multline*}
  m{}^1\square_0^3\Gamma=(m_2;m_1,m_2,f',f''){}^1\square_0^3
  \begin{pmatrix}
    \Gamma_1, \Gamma_2,\alpha_3,\beta_3;\\
    \alpha_1,\alpha_2,\beta_1\beta_2,g_2,h_2;\\
    g_0,g_1,h_0,h_1,f,f'
  \end{pmatrix}\\=
  \begin{pmatrix}
    m_0\#_0\Gamma_1, m_1\#_0\Gamma_2,\\
    ((m_2\#_0h_0)\#_1(m_1\#_0\alpha_3))\#_2(\overline{(m_2\ten\alpha_1)}\#_1(m_1\#_0g_2)),\\
    ((m_2\#_0h_0)\#_1(m_1\#_0\beta_3))\#_2(\overline{((m_2\ten\beta_1))}\#_1(m_1\#_0g_2));\\
    m_0\#_0\alpha_1,m_0\#_1\alpha_2,m_0\#_0\beta_1,m_1\#_0\beta_2,\\
    (m_2\#_0g_0)\#_1(m_1\#_0g_2),(m_2\#_0h_0)\#_1(m_1\#_0h_2);\\
    m_0\#_0g_0,m_1\#_0g_1,m_0\#_0h_0,m_1\#_0h_1,f,f''
  \end{pmatrix}
\end{multline*}
We claim this is again a proper 3-cell in $\pathspc\H$, that is,
the whisker satisfies \eqref{eq:3cellcond}, as can be easily seen:
\begin{equation*}
  \begin{xy}\xycompile{
    (50,0):(0,0)
     \xyboxmatrix"A"{{} \ar[r]^f \ar[d]|{g_0}="x"\ar@/_2.5pc/[d]|{h_0}="y" &
      {}\ar[d]^{g_1} \ar@2 [dl] **{} ?(.2); ?(.8) |{g_2}\\{}
      \ar[r]|{f'} \ar[d]_{m_0}& {}\ar[d]^{m_1} \ar@2 [dl] **{} ?(.2); ?(.8) |{m_2}="Y" 
      \POS \ar@2 "x"; "y" **{} ?(.2); ?(.8)|{\alpha_1}="X" \\
      {}\ar[r]_{f''} & {}
      \tria "Y"; "X"
    }
    +(1,0)
    \xyboxmatrix"B"{{} \ar[r]^f \ar[d]|{g_0}="x"\ar@/_2.5pc/[d]|{h_0}="y" &
      {}\ar[d]^{g_1} \ar@2 [dl] **{} ?(.2); ?(.8) |{g_2}\\{}
      \ar[r]|{f'} \ar[d]_{m_0}& {}\ar[d]^{m_1} \ar@2 [dl] **{} ?(.2); ?(.8) |{m_2}="Y" 
      \POS \ar@2 "x"; "y" **{} ?(.2); ?(.8)|{\alpha_1}="X" \\
      {}\ar[r]_{f''} & {}
       \tria "X"; "Y"
    }
    +(1,0)
    \xyboxmatrix"C"{{} \ar[r]^f \ar[d]_{h_0}&
      {}\ar[d]|{h_1}="x" \ar@2 [dl] **{} ?(.2); ?(.8) |{h_2}
      \ar@/^2.5pc/[d]|{g_1}="y"\\
      {} \ar[r]|{f'} \ar[d]_{m_0}& 
      {} \ar[d]^{m_1} \ar@2 [dl] **{} ?(.2); ?(.8) |{m_2} 
      \POS \ar@2 "y"; "x" **{} ?(.2); ?(.8) |{\alpha_2}\\
      {} \ar[r]_{f''}&{}
    }
    \POS 0+(0,-1)
    \xyboxmatrix"D"{{} \ar[r]^f \ar[d]|{g_0}="x"\ar@/_2.5pc/[d]|{h_0}="y" &
      {}\ar[d]^{g_1} \ar@2 [dl] **{} ?(.2); ?(.8) |{g_2}\\{}
      \ar[r]|{f'} \ar[d]_{m_0}& {}\ar[d]^{m_1} \ar@2 [dl] **{} ?(.2); ?(.8) |{m_2}="Y" 
      \POS \ar@2 "x"; "y" **{} ?(.2); ?(.8)|{\beta_1}="X" \\
      {}\ar[r]_{f''} & {}
      \tria "Y"; "X"
    }
    +(1,0)
    \xyboxmatrix"E"{{} \ar[r]^f \ar[d]|{g_0}="x"\ar@/_2.5pc/[d]|{h_0}="y" &
      {}\ar[d]^{g_1} \ar@2 [dl] **{} ?(.2); ?(.8) |{g_2}\\{}
      \ar[r]|{f'} \ar[d]_{m_0}& {}\ar[d]^{m_1} \ar@2 [dl] **{} ?(.2); ?(.8) |{m_2}="Y" 
      \POS \ar@2 "x"; "y" **{} ?(.2); ?(.8)|{\beta_1}="X" \\
      {}\ar[r]_{f''} & {}
       \tria "X"; "Y"
    }
    +(1,0)
    \xyboxmatrix"F"{{} \ar[r]^f \ar[d]_{h_0}&
      {}\ar[d]|{h_1}="x" \ar@2 [dl] **{} ?(.2); ?(.8) |{h_2}
      \ar@/^2.5pc/[d]|{g_1}="y"\\
      {} \ar[r]|{f'} \ar[d]_{m_0}& 
      {} \ar[d]^{m_1} \ar@2 [dl] **{} ?(.2); ?(.8) |{m_2} 
      \POS \ar@2 "y"; "x" **{} ?(.2); ?(.8) |{\beta_2}\\
      {} \ar[r]_{f''}&{}
    }
    \ar@3 "A";"B" ^{\overline{(m_2\ten\alpha_1)}\\\;\#_1(m_1\#_0g_2)}    
    \ar@3 "B";"C" ^{(m_2\#_0h_0)\\\;\#_1(m_1\#_0\alpha_3)}    
    \ar@3 "D";"E" _{\overline{(m_2\ten\beta_1)}\\\;\#_1(m_1\#_0g_2)}    
    \ar@3 "E";"F" _{(m_2\#_0h_0)\\\;\#_1(m_1\#_0\beta_3)}    
    \ar@3 "A";"D" |{
      (f''\#_0m_0\#_0\Gamma_1)
      \\\;\#_1(m_2\#_0g_0)
      \\\;\#_1(m_1\#_0g_2)}
    \ar@3 "B";"E" |{
      (m_2\#_0h_0)
      \\\;\#_1(m_1\#_0f'\#_0\Gamma_1)
      \\\;\#_1(m_1\#_0g_2)}
    \ar@3 "C";"F" |{
      (m_2\#_0h_0)
      \\\;\#_1(m_1\#_0h_2)
      \\\;\#_1(m_1\#_0\Gamma_2\#_0f)}
    \ar@{} "A";"E" |{\text{naturality}}
    \ar@{} "B";"F" |{\text{\eqref{eq:3cellcond}}}
  }\end{xy}
\end{equation*}

Finally, we define 3-2-whiskering:
\begin{multline}
  \label{eq:sq23def}
  \gamma{}^2\square_1^3\Gamma=
  \begin{pmatrix}
    \gamma_3;\gamma_1,\gamma_2,h_2,k_2;\\
    h_0,h_1,k_0,k_1,f,f'
  \end{pmatrix}{}^2\square_1^3
  \begin{pmatrix}
    \Gamma_1,\Gamma_2,\alpha_3,\beta_3;g_2,h_2,\\
    \alpha_1,\alpha_2,\beta_1,\beta_2;\\
    g_0,g_1,h_0,h_1,f,f'
  \end{pmatrix}\\=
  \begin{pmatrix}
    \gamma_1\#_1\Gamma_1,\gamma_2\#_1\Gamma_2,\\
    (\gamma_3\#_1(\alpha_2\#_0f))\#_2((f'\#_0\gamma_1)\#_1\alpha_3),\\
    (\gamma_3\#_1(\beta_2\#_0f))\#_2((f'\#_0\gamma_1)\#_1\beta_3);\\g_2,k_2,
    \gamma_1\#_1\alpha_1,\gamma_2\#_1\alpha_2,\gamma_1\beta_1,\gamma_2\beta_2;\\
    g_0,g_1,k_0,k_1,f,f'
  \end{pmatrix}
\end{multline}
It yeilds a 3-cell in $\pathspc\H$:
\begin{equation}
  \label{eq:sq23defexp}
  \begin{xy}\xycompile{ (50,0):(0,0)
    \xyboxmatrix"A"{
      {} \ar[r]^f \ar[d]|{g_0}="x"
      \ar@/_2.5pc/[d]|{h_0}="y"
      \ar@/_5pc/[d]|{k_0}="z" &
      {}\ar[d]^{g_1} \ar@2 [dl] **{} ?(.2); ?(.8) |{g_2}\\{} \ar[r]_{f'}& {} 
      \POS \ar@2 "x"; "y" **{} ?(.2); ?(.8)|{\alpha_1} 
      \POS \ar@2 "y"; "z" **{} ?(.2); ?(.8)|{\gamma_1} 
    } 
    +(1, 0)
    \xyboxmatrix"B"{
      {}\ar[r]^f \ar[d]|{h_0}="x" \ar@/_2.5pc/[d]|{k_0}="y" &
      {}\ar[d]|{h_1}="w" \ar@/^2.5pc/[d]|{g_1}="z" 
      \ar@2 [dl] **{} ?(.2); ?(.8) |{h_2}\\{}
      \ar[r]_{f'}& {} 
      \POS \ar@2 "x"; "y" **{} ?(.2); ?(.8) |{\gamma_1} 
      \POS \ar@2 "z"; "w" **{} ?(.2); ?(.8) |{\alpha_2} 
    } 
    +(1, 0)
    \xyboxmatrix"C"{
      {} \ar[r]^f \ar[d]_{k_0}  &
      {}\ar[d]|{k_1}="z" \ar@2 [dl] **{} ?(.2); ?(.8) |{k_2}
      \ar@/^2.5pc/[d]|{h_1}="y" \ar@/^5pc/[d]|{g_1}="x" \\{}
      \ar[r]_{f'}& {} 
      \POS \ar@2 "x"; "y" **{} ?(.2); ?(.8)|{\alpha_2} 
      \POS \ar@2 "y"; "z" **{} ?(.2); ?(.8)|{\gamma_2} 
    }
    ,(0,-1)
    \xyboxmatrix"D"{
      {} \ar[r]^f \ar[d]|{g_0}="x"
      \ar@/_2.5pc/[d]|{h_0}="y"
      \ar@/_5pc/[d]|{k_0}="z" &
      {}\ar[d]^{g_1} \ar@2 [dl] **{} ?(.2); ?(.8) |{g_2}\\{} \ar[r]_{f'}& {} 
      \POS \ar@2 "x"; "y" **{} ?(.2); ?(.8)|{\beta_1} 
      \POS \ar@2 "y"; "z" **{} ?(.2); ?(.8)|{\gamma_1} 
    } 
    +(1, 0)
    \xyboxmatrix"E"{
      {}\ar[r]^f \ar[d]|{h_0}="x" \ar@/_2.5pc/[d]|{k_0}="y" &
      {}\ar[d]|{h_1}="w" \ar@/^2.5pc/[d]|{g_1}="z" 
      \ar@2 [dl] **{} ?(.2); ?(.8) |{h_2}\\{}
      \ar[r]_{f'}& {} 
      \POS \ar@2 "x"; "y" **{} ?(.2); ?(.8) |{\gamma_1} 
      \POS \ar@2 "z"; "w" **{} ?(.2); ?(.8) |{\beta_2} 
    } 
    +(1, 0)
    \xyboxmatrix"F"{
      {} \ar[r]^f \ar[d]_{k_0}  &
      {}\ar[d]|{k_1}="z" \ar@2 [dl] **{} ?(.2); ?(.8) |{k_2}
      \ar@/^2.5pc/[d]|{h_1}="y" \ar@/^5pc/[d]|{g_1}="x" \\{}
      \ar[r]_{f'}& {} 
      \POS \ar@2 "x"; "y" **{} ?(.2); ?(.8)|{\beta_2} 
      \POS \ar@2 "y"; "z" **{} ?(.2); ?(.8)|{\gamma_2} 
    }
    \ar@3 "A";"B" ^{(f'\#_0\gamma_1)\\\;\#_1\alpha_3}
    \ar@3 "B";"C" ^{\gamma_3\\\;\#_1(\alpha_2\#_0f)}
    \ar@3 "D";"E" _{(f'\#_0\gamma_1)\\\;\#_1\beta_3}
    \ar@3 "E";"F" _{\gamma_3\\\;\#_1(\beta_2\#_0f)}
    \ar@3 "A";"D" |{(f'\#_0\gamma_1)\\\;\#_1(f'\#_0\Gamma_1)\\\;\#_1g_2}
    \ar@3 "B";"E" |{(f'\#_0\Gamma_1)\\\;\#_1h_2\\\;\#_1(\alpha_2\#_0f)}
    \ar@3 "C";"F" |{k_2\\\;\#_1(\Gamma_2\#_0f)\\\;\#_1(\alpha_2\#_0f)}
    \ar@{} "A";"E" |{\text{\eqref{eq:3cellcond}}}
    \ar@{} "B";"F" |{\text{func.}}
  }\end{xy}
\end{equation} 

\subsection{Horizontal Composition of 2-Cells}
\label{sec:hrz2ccomp}

We shall use the following slightly abbreviated notation for the
higher cells of the mapping space, for example writing
\eqref{eq:sq12def} as:
\newcommand{\sep}{\,\vline\,}
\begin{multline*}
  \xymatrix"A"{
    \ar@/^1pc/[r]^g="x" \ar@/_1pc/[r]_{n}="y" 
    \POS \ar@2 "x";"y" **{} ?(.3); ?(.7) _{\alpha}="X" &
    \ar[r]^k & {}
  }=
  k{}^1\square_0^2\alpha=  (k_2;k_0,k_1,f',f''){}^1\square_0^2
  \begin{pmatrix}
    \alpha_3;\alpha_1,\alpha_2\sep g,n   
  \end{pmatrix}\\= \left(
    \begin{matrix}
      ((k_2\#_0n_0)\#_1(k_1\#_0\alpha_3))
      \#_2(\overline{(k_2\ten\alpha_1)}\#_1(k_1\#_0g_2));\\
      k_0\#_0\alpha_1, k_1\#_0\alpha_2\sep
      k\square_0g,k\square_0n
    \end{matrix}
  \right)\,.
\end{multline*}
In the same spirit we write the opposite whiskering:
\begin{multline*}
  \xymatrix"B"{
    \ar[r]^{n} &
    \ar@/^1pc/[r]^k="x" \ar@/_1pc/[r]_{m}="y" 
    \POS \ar@2 "x";"y" **{} ?(.3); ?(.7) ^{\beta}="Y" & {}
  }=
  \beta{}^2\square_0^1n=
  \begin{pmatrix}
    \beta_3;\beta_1,\beta_2\sep k,m
  \end{pmatrix}\\=
    \begin{pmatrix}
      ((m_2\#_0n_0)\#_1(\beta_2\ten n_2))\#_2(\beta_3\#_1(k_1\#_0n_2));\\
      \beta_1\#_0n_0,\beta_2\#_0n_1\sep k\square_0n, m\square_0n  
    \end{pmatrix}\,.
\end{multline*}

So now we can define the left horizontal composite:

\begin{multline*}
    \xymatrix{
      \ar@/^1pc/[r]^g="x" \ar@/_1pc/[r]_{n}="y" 
      \POS \ar@2 "x";"y" **{} ?(.3); ?(.7) _{\alpha}="X" &
      \ar@/^1pc/[r]^k="x" \ar@/_1pc/[r]_{m}="y" 
      \POS \ar@2 "x";"y" **{} ?(.3); ?(.7) ^{\beta}="Y" & {}
      \tria "X"; "Y"
    }= \beta\blhc\alpha
    =\begin{pmatrix}
      ((m_2\#_0n_0)\#_1(\beta_2\ten n_2))\\
      \#_2(\beta_3\#_1(k_1\#_0n_2));\\
      \beta_1\#_0n_0,\beta_2\#_0n_1\sep k\square_0n, m\square_0n  
    \end{pmatrix}
    \square_1
    \begin{pmatrix}
      ((k_2\#_0n_0)\#_1(k_1\#_0\alpha_3))\\
      \#_2(\overline{(k_2\ten\alpha_1)}\#_1(k_1\#_0g_2));\\
      k_0\#_0\alpha_1, k_1\#_0\alpha_2\sep
      k\square_0g,k\square_0n
    \end{pmatrix}\\
    =
    \begin{pmatrix}
      \begin{pmatrix}
        \begin{pmatrix}
          ((m_2\#_0n_0)\#_1(\beta_2\ten n_2))\\
          \#_2(\beta_3\#_1(k_1\#_0n_2))
        \end{pmatrix}\#_1(k_1\#_0\alpha_2\#_0f)
      \end{pmatrix}\\
      \#_2
      \begin{pmatrix}
        (f''\#_0\beta_1\#_0n_0)\#_1 
        \begin{pmatrix}
          ((k_2\#_0n_0)\#_1(k_1\#_0\alpha_3))\\
          \#_2(\overline{(k_2\ten\alpha_1)}\#_1(k_1\#_0g_2))
        \end{pmatrix}
      \end{pmatrix};\\
      \alpha_1\lhc\beta_1,\alpha_2\lhc\beta_2\sep k\square_0g,m\square_0n
    \end{pmatrix}
\end{multline*}
and conversely, 
\begin{multline*}
    \xymatrix{
      \ar@/^1pc/[r]^g="x" \ar@/_1pc/[r]_{n}="y" 
      \POS \ar@2 "x";"y" **{} ?(.3); ?(.7) _{\alpha}="X" &
      \ar@/^1pc/[r]^k="x" \ar@/_1pc/[r]_{m}="y" 
      \POS \ar@2 "x";"y" **{} ?(.3); ?(.7) ^{\beta}="Y" & {}
      \tria "Y"; "X"
    }= \beta\blhc\alpha
    =\begin{pmatrix}
      ((m_2\#_0n_0)\#_1(\beta_2\ten n_2))\\
      \#_2(\beta_3\#_1(k_1\#_0n_2));\\
      \beta_1\#_0n_0,\beta_2\#_0n_1\sep k\square_0n, m\square_0n  
    \end{pmatrix}
    \square_1
    \begin{pmatrix}
      ((k_2\#_0n_0)\#_1(k_1\#_0\alpha_3))\\
      \#_2(\overline{(k_2\ten\alpha_1)}\#_1(k_1\#_0g_2));\\
      k_0\#_0\alpha_1, k_1\#_0\alpha_2\sep
      k\square_0g,k\square_0n
    \end{pmatrix}\\
    =
    \begin{pmatrix}
      \begin{pmatrix}
        \begin{pmatrix}
          ((m_2\#_0n_0)\#_1(\beta_2\ten n_2))\\
          \#_2(\beta_3\#_1(k_1\#_0n_2))
        \end{pmatrix}\#_1(k_1\#_0\alpha_2\#_0f)
      \end{pmatrix}\\
      \#_2
      \begin{pmatrix}
        (f''\#_0\beta_1\#_0n_0)\#_1 
        \begin{pmatrix}
          ((k_2\#_0n_0)\#_1(k_1\#_0\alpha_3))\\
          \#_2(\overline{(k_2\ten\alpha_1)}\#_1(k_1\#_0g_2))
        \end{pmatrix}
      \end{pmatrix};\\
      \alpha_1\lhc\beta_1,\alpha_2\lhc\beta_2\sep k\square_0g,m\square_0n
    \end{pmatrix}\,.
\end{multline*}

\subsection{Tensors}
\label{sec:pthspcten}

Finally, in
\begin{equation*}
  \begin{xy}\xycompile{
    \xyboxmatrix"A"{\ar@/^1pc/[r]^g="x" \ar@/_1pc/[r]_{n}="y" 
      \POS \ar@2 "x";"y" **{} ?(.3); ?(.7) _{\alpha}="X" &
      \ar@/^1pc/[r]^k="x" \ar@/_1pc/[r]_{m}="y" 
      \POS \ar@2 "x";"y" **{} ?(.3); ?(.7) ^{\beta}="Y" & {}
      \tria "X"; "Y"}
    +(30,0)
    \xyboxmatrix"B"{\ar@/^1pc/[r]^g="x" \ar@/_1pc/[r]_{n}="y" 
      \POS \ar@2 "x";"y" **{} ?(.3); ?(.7) _{\alpha}="X" &
      \ar@/^1pc/[r]^k="x" \ar@/_1pc/[r]_{m}="y" 
      \POS \ar@2 "x";"y" **{} ?(.3); ?(.7) ^{\beta}="Y" & {}
      \tria "Y"; "X"}
    \ar@3 "A";"B" ^{\beta\bten\alpha}
  }\end{xy}
\end{equation*}
letting $\beta\bten\alpha=(\beta_1\ten\alpha_1, \beta_2\ten\alpha_2)$
makes $\pathspc\H$ a $\Gray$-category. This is a well defined
3-cell. 

\subsection{Inverses}
\label{sec:pathspcinv}

If $\H$ has invertible 1- and 2-cells the inverse of a 1-cell
\begin{equation*}
  \begin{xy}
    \xymatrix{{} \ar[r]^f \ar[d]_{g_0}&
      {}\ar[d]^{g_1} \ar@2 [dl] **{} ?<(.3); ?>(.7) ^{g_2}\\{} \ar[r]_{f'}& {}}
  \end{xy}
\end{equation*}
 in
$\pathspc\H$ is given by 
\begin{equation*}
  \begin{matrix}
    \begin{xy}
      \xymatrix{
        {}\ar[rrr]^{f'}\ar[ddd]_{\overline{g_0}}\ar[dr]^{\overline{g_0}}&&&\ar[ddd]^{\overline{g_1}}\\
        {}&\ar[r]^{g_0}\ar[d]_{f}&\ar[d]^{f'}\ar@2 [dl] **{} ?<(.3); ?>(.7) ^{\overline{g_2}}&\\
        {}&\ar[r]_{g_1}&\ar[dr]^{\overline{g_1}}&\\
        \ar[rrr]_{f}&&& }
    \end{xy}
  \end{matrix}
\;.
\end{equation*}

\subsection{Axioms}
\label{sec:pthspcgrax}

This composition of $\pathspc\H$-2-cells is associative: Given three
2-cells 
\begin{align*}
  \alpha=\qquad&\begin{matrix}\begin{xy}\xycompile{
    \xymatrix"A"{{} \ar[r]^f \ar[d]|{g_0}="x"\ar@/_3pc/[d]|{h_0}="y" &
      {}\ar[d]^{g_1} \ar@2 [dl] **{} ?(.2); ?(.8) |{g_2}\\{}
      \ar[r]_{f'}& {} \POS \ar@2 "x"; "y" **{} ?(.2); ?(.8) |{\alpha_1} }
    \POS + (30,0)
    \xymatrix"B"{{} \ar[r]^f \ar[d]_{h_0}&
      {}\ar[d]|{h_1}="x" \ar@2 [dl] **{} ?(.2); ?(.8) |{h_2}
      \ar@/^3pc/[d]|{g_1}="y"\\{} \ar[r]_{f'}& {} \POS \ar@2
      "y"; "x" **{} ?(.2); ?(.8) |{\alpha_2}}
    \POS \ar@3 {"A1,1" ; "A2,2" **{} ?(.5)}; {"B1,1" ; "B2,2" **{} ?(.5)} **{}
    ?<(.4);  ?>(.6) ^ {\alpha_3}
  }\end{xy}\end{matrix}\\
  \beta=\qquad&\begin{matrix}\begin{xy}\xycompile{
    \xymatrix"A"{{} \ar[r]^f \ar[d]|{h_0}="x"\ar@/_3pc/[d]|{k_0}="y" &
      {}\ar[d]^{h_1} \ar@2 [dl] **{} ?(.2); ?(.8) |{h_2}\\{}
      \ar[r]_{f'}& {} \POS \ar@2 "x"; "y" **{} ?(.2); ?(.8) |{\beta_1} }
    \POS + (30,0)
    \xymatrix"B"{{} \ar[r]^f \ar[d]_{k_0}&
      {}\ar[d]|{k_1}="x" \ar@2 [dl] **{} ?(.2); ?(.8) |{k_2}
      \ar@/^3pc/[d]|{h_1}="y"\\{} \ar[r]_{f'}& {} \POS \ar@2
      "y"; "x" **{} ?(.2); ?(.8) |{k_2}}
    \POS \ar@3 {"A1,1" ; "A2,2" **{} ?(.5)}; {"B1,1" ; "B2,2" **{} ?(.5)} **{}
    ?<(.4);  ?>(.6) ^ {\beta_3}
  }\end{xy}\end{matrix}\\
  \gamma=\qquad&\begin{matrix}\begin{xy}\xycompile{
    \xymatrix"A"{{} \ar[r]^f \ar[d]|{k_0}="x"\ar@/_3pc/[d]|{m_0}="y" &
      {}\ar[d]^{k_1} \ar@2 [dl] **{} ?(.2); ?(.8) |{k_2}\\{}
      \ar[r]_{f'}& {} \POS \ar@2 "x"; "y" **{} ?(.2); ?(.8) |{\gamma_1} }
    \POS + (30,0)
    \xymatrix"B"{{} \ar[r]^f \ar[d]_{m_0}&
      {}\ar[d]|{m_1}="x" \ar@2 [dl] **{} ?(.2); ?(.8) |{m_2}
      \ar@/^3pc/[d]|{k_1}="y"\\{} \ar[r]_{f'}& {} \POS \ar@2
      "y"; "x" **{} ?(.2); ?(.8) |{\gamma_2}}
    \POS \ar@3 {"A1,1" ; "A2,2" **{} ?(.5)}; {"B1,1" ; "B2,2" **{} ?(.5)} **{}
    ?<(.4);  ?>(.6) ^ {\gamma_3}
  }\end{xy}\end{matrix}
\end{align*}
we use \eqref{eq:vert2compdef} and the functoriality of the
whiskerings in $\H$ to compute:
\begin{multline*}
  (\gamma\square_1\beta)\square_1\alpha =
  \left(\begin{matrix}
            \underbrace{(\gamma_3\#_1(\beta_2\#_0f))\#_2((f'\#_0\gamma_1)\#_1\beta_3)}_{\omega_3};\\
            \gamma_1\#_1\beta_1, \gamma_2\#_1\beta_2,h_2,m_2;h_0,h_1,m_0,m_1,f,f'
    \end{matrix}\right)\square_1\alpha\\=\left(
    \begin{matrix}
            (\omega_3\#_1(\alpha_2\#_0f))\\\#_2((f'\#_0(\gamma_1\#_1\beta_1))\#_1\alpha_3);\\
            \gamma_1\#_1\beta_1\#_1\alpha_1,\gamma_2\#_1\beta_2\#_1\alpha_2,g_2,m_2;g_0,g_1,m_0,m_1,f,f'
    \end{matrix}
  \right)\\=\left(
    \begin{matrix}
      (((\gamma_3\#_1(\beta_2\#_0f))\#_2((f'\#_0\gamma_1)\#_1\beta_3))\\
      \#_1(\alpha_2\#_0f))\#_2((f'\#_0(\gamma_1\#_1\beta_1))\#_1\alpha_3);\\
      \gamma_1\#_1\beta_1\#_1\alpha_1,\gamma_2\#_1\beta_2\#_1\alpha_2,\\
      g_2,m_2;g_0,g_1,m_0,m_1,f,f'
    \end{matrix}
  \right)=\left(
    \begin{matrix}
      (\gamma_3\#_1(\beta_2\#_0f)\#_1(\alpha_2\#_0f))\\
      \#_2(((f'\#_0\gamma_1)\#_1\beta_3)\#_1(\alpha_2\#_0f))\\
      \#_2((f'\#_0(\gamma_1\#_1\beta_1))\#_1\alpha_3);\\
      \gamma_1\#_1\beta_1\#_1\alpha_1,\gamma_2\#_1\beta_2\#_1\alpha_2,g_2,m_2;\\
      g_0,g_1,m_0,m_1,f,f'
    \end{matrix}
  \right)\\=\left(
    \begin{matrix}
      (\gamma_3\#_1((\beta_2\#_1\alpha_2)\#_0f))\#_2((f'\#_0\gamma_1)\#_1\beta_3\#_1(\alpha_2\#_0f))\\
      \#_2((f'\#_0\gamma_1)\#_1(f'\#_0\beta_1)\#_1\alpha_3);\\
      \gamma_1\#_1\beta_1\#_1\alpha_1,\gamma_2\#_1\beta_2\#_1\alpha_2,g_2,m_2;g_0,g_1,m_0,m_1,f,f'
    \end{matrix}
  \right)\\=\left(
    \begin{matrix}
      (\gamma_3\#_1((\beta_2\#_1\alpha_2)\#_0f))\#_2((f'\#_0\gamma_1)\\
      \#_1\underbrace{((\beta_3\#_1(\alpha_2\#_0f))\#_2((f'\#_0\beta_1)\#_1\alpha_3))}_{\zeta_3});\\
      \gamma_1\#_1\beta_1\#_1\alpha_1,\gamma_2\#_1\beta_2\#_1\alpha_2,g_2,m_2;g_0,g_1,m_0,m_1,f,f'
    \end{matrix}\right)
  \\=\left(
    \begin{matrix}
      (\gamma_3\#_1((\beta_2\#_1\alpha_2)\#_0f))\\\#_2((f'\#_0\gamma_1)\#_1\zeta_3);\\
      \gamma_1\#_1\beta_1\#_1\alpha_1,\gamma_2\#_1\beta_2\#_1\alpha_2,g_2,m_2;g_0,g_1,m_0,m_1,f,f'
    \end{matrix}
  \right)
  \\=\gamma\square_1\left(
    \begin{matrix}
      \zeta_3;\beta_1\#_1\alpha_1,
      \beta_2\#_1\alpha_2,\\
      g_2,k_2;g_0,g_1,k_0,k_1,f,f'
    \end{matrix}
  \right)=\gamma\square_1(\beta\square_1\alpha)\,.
\end{multline*}

We check that 2-1-whiskering in $\pathspc\H$ is functorial, that is,
$m\square_0(\beta\square_1\alpha)=(m\square_0\beta)\square_1(m\square_0\alpha)$. In
diagram \eqref{eq:whisker21func} the diagonal is
$m\square_0(\beta\square_1\alpha)$ and left and down is
$(m\square_0\beta)\square_1(m\square_0\alpha)$. 1-2-whiskering in
$\pathspc\H$ is functorial by duality.
\begin{sidewaysfigure}
  \begin{equation}
    \label{eq:whisker21func}
    \begin{xy}
      (60,0) :
      (0,0)
      \xyboxmatrix"M1"{{} 
        \ar[r]^f
        \ar[d]|{g_0}="x"
        \ar@/_2.5pc/[d]|{h_0}="y"
        \ar@/_5pc/[d]|{k_0}="z"&
        {}\ar[d]^{g_1} \ar@2 [dl] **{} ?(.2); ?(.8) |{g_2}\\{}
        \ar[r]|{f'} \ar[d]_{m_0}& {}\ar[d]^{m_1} \ar@2 [dl] **{} ?(.2); ?(.8) |{m_2}="Y" 
        \POS \ar@2 "x"; "y" **{} ?(.2); ?(.8)|{\alpha_1}="X" 
        \POS \ar@2 "y"; "z" **{} ?(.2); ?(.8)|{\beta_1}\\
        {}\ar[r]_{f''} & {}
        \tria  "Y"; "X"
      }
      ,(1,0)
      \xyboxmatrix"M2"{{} 
        \ar[r]^f \ar[d]|{g_0}="x"
        \ar@/_2.5pc/[d]|{h_0}="y"
        \ar@/_5pc/[d]|{k_0}="z" &
        {}\ar[d]^{g_1} \ar@2 [dl] **{} ?(.2); ?(.8) |{g_2}\\{}
        \ar[r]|{f'} \ar[d]_{m_0}& {}\ar[d]^{m_1} \ar@2 [dl] **{} ?(.2); ?(.8) |{m_2}="Y" 
        \POS \ar@2 "x"; "y" **{} ?(.2); ?(.8)|{\alpha_1}="X" 
        \POS \ar@2 "y"; "z" **{} ?(.2); ?(.8)|{\beta_1}="Z"\\
        {}\ar[r]_{f''} & {}
        \tria "X"; "Y"
      }
      ,(2,0)
      \xyboxmatrix"M3"{{} 
        \ar[r]^f 
        \ar[d]|{h_0}="w"
        \ar@/_2.5pc/[d]|{k_0}="z"&
        {} \ar[d]|{h_1}="x" \ar@2 [dl] **{} ?(.2); ?(.8) |{h_2}
        \ar@/^2.5pc/[d]|{g_1}="y"\\
        {} \ar[r]|{f'} \ar[d]_{m_0}& 
        {} \ar[d]^{m_1} \ar@2 [dl] **{} ?(.2); ?(.8) |{m_2}="X" 
        \POS \ar@2 "y"; "x" **{} ?(.2); ?(.8) |{\alpha_2}
        \POS \ar@2 "w"; "z" **{} ?(.2); ?(.8) |{\beta_1}="Y" \\
        {} \ar[r]_{f''}&{}
        \tria "X"; "Y"
      }
      ,(1,-1)       
      \xyboxmatrix"M4"{{} 
        \ar[r]^f
        \ar[d]|{g_0}="x"
        \ar@/_2.5pc/[d]|{h_0}="y"
        \ar@/_5pc/[d]|{k_0}="z"&
        {}\ar[d]^{g_1} \ar@2 [dl] **{} ?(.2); ?(.8) |{g_2}\\{}
        \ar[r]|{f'} \ar[d]_{m_0}& {}\ar[d]^{m_1} \ar@2 [dl] **{} ?(.2); ?(.8) |{m_2}="Y" 
        \POS \ar@2 "x"; "y" **{} ?(.2); ?(.8)|{\alpha_1} 
        \POS \ar@2 "y"; "z" **{} ?(.2); ?(.8)|{\beta_1}="X"\\
        {}\ar[r]_{f''} & {}
        \tria  "X"; "Y"
      }
      ,(2,-1) 
      \xyboxmatrix"M5"{{} 
        \ar[r]^f 
        \ar[d]|{h_0}="w"
        \ar@/_2.5pc/[d]|{k_0}="z"&
        {} \ar[d]|{h_1}="x" \ar@2 [dl] **{} ?(.2); ?(.8) |{h_2}
        \ar@/^2.5pc/[d]|{g_1}="y"\\
        {} \ar[r]|{f'} \ar[d]_{m_0}& 
        {} \ar[d]^{m_1} \ar@2 [dl] **{} ?(.2); ?(.8) |{m_2}="X" 
        \POS \ar@2 "y"; "x" **{} ?(.2); ?(.8) |{\alpha_2}
        \POS \ar@2 "w"; "z" **{} ?(.2); ?(.8) |{\beta_1}="Y" \\
        {} \ar[r]_{f''}&{}
        \tria "Y"; "X"
      }
      ,(2,-2)
      \xyboxmatrix"M6"{{} 
        \ar[r]^f 
        \ar[d]|{k_0}&
        {} \ar[d]|{k_1}="z" \ar@2 [dl] **{} ?(.2); ?(.8) |{k_2}
        \ar@/^5pc/[d]|{g_1}="x"  \ar@/^2.5pc/[d]|{h_1}="y"\\
        {} \ar[r]|{f'} \ar[d]_{m_0}& 
        {} \ar[d]^{m_1} \ar@2 [dl] **{} ?(.2); ?(.8) |{m_2}="X" 
        \POS \ar@2 "x"; "y" **{} ?(.2); ?(.8) |{\alpha_2}
        \POS \ar@2 "y"; "z" **{} ?(.2); ?(.8) |{\beta_2}="Y" \\
        {} \ar[r]_{f''}&{}
    }
    \ar@3 "M1"; "M2" ^{
      (f''\#_0 m_0\#_0 \beta_1)\\
      \;\#_1\overline{(\alpha_1\ten m_2)}\\
      \;\#_1(m_1\#_0 g_2)
    }
    \ar@3 "M1"; "M4" _{
      \overline{m_2\ten(\beta_1\#_1\alpha_1)}\\
      \;\#_1(m_1\#_0g_2)
    }
    \ar@3 "M2"; "M3" ^{
      (f''\#_0m_0\#_0\beta_1)\\
      \;\#_1(m_2\#_0h_0)\\
      \;\#_1(m_1\#_0\alpha_3)
    }
    \ar@3 "M2"; "M4"|{
      \overline{\beta_1\ten m_2}\\
      \;\#_1(m_1\#_0f'\#_0\alpha_1)\\
      \;\#_1(m_1\#_0g_2)
    }    
    \ar@3 "M3"; "M5"^{
      \overline{\beta_1\ten m_2}\\
      \;\#_1(m_1\#_0h_2)\\
      \;\#_1(m_1\#_0\alpha_2\#_0f)
    }
    \ar@3 "M4"; "M6"_{
      (m_2\#_0k_0)\\
      \;\#_1(m_1\#_0((\beta_3\#_1(\alpha_2\#_0f))\\
      \;\;\#_2((f'\#_0\beta_1)\#_1\alpha_3)))
    }    
    \ar@3 "M4"; "M5"^{
      (m_2\#_0k_0)\\
      \;\#_1(m_1\#_0f'\#_0\beta_1)\\
      \;\#_1(m_1\#_0\alpha_3)
    }
    \ar@3 "M5"; "M6"^{
      (m_2\#_0k_0)\\
      \;\#_1(m_1\#_0\beta_3)\\
      \;\#_1(m_1\#_0\alpha_2\#_0f)
    }
    \POS "M1";{"M2";"M4"**{}?}**{}?(.66)*\labelbox{\text{func.}}
    \POS "M2";"M5"**{}?*\labelbox{\text{func.}} 
  \end{xy}
\end{equation}
\end{sidewaysfigure}
  
It is obvious that 3-1-whiskering is 2-functorial, that is, 
\begin{multline*}
  (m_0,m_1,m_2)\square_0((\Delta_1,\Delta_2)\square_2(\Gamma_1, \Gamma_2))\\=
  (m_0,m_1,m_2)\square_0(\Delta_1\#_2\Gamma_1,\Delta_2\#_2\Gamma_2)\\=
  (m_0\#_0(\Delta_1\#_2\Gamma_1),m_1\#_0(\Delta_2\#_2\Gamma_2))\\=
  (((m_0\#_0\Delta_1)\#_2(m_0\#_0\Gamma_1)),((m_1\#_0\Delta_2)\#_2(m_1\#_0\Gamma_2)))\\=
  ((m_0\#_0\Delta_1),(m_1\#_0\Delta_2))\square_2((m_0\#_0\Gamma_1),(m_1\#_0\Gamma_2))\\=
  ((m_0,m_1,m_2)\square_0(\Delta_1,\Delta_2))\square_2((m_0,m_1,m_2)\square_0(\Gamma_1, \Gamma_2))\,.
\end{multline*}

By duality, 1-2-whiskering in $\pathspc\H$ is functorial as well.
And the 3-2-whiskering thus defined is functorial with
respect to vertical composition of 3-cells, that is,
$\gamma\square_1(\Gamma\square_2\Delta)=(\gamma\square_1\Gamma)\square_2(\gamma\square_1\Delta)$,
as can seen by inspecting 
the following diagram:
\begin{equation*}
  \begin{xy}
    (50,0):(0,0)
    \xyboxmatrix"A"{
      {} \ar[r]^f \ar[d]|{g_0}="x"
      \ar@/_2.5pc/[d]|{h_0}="y"
      \ar@/_5pc/[d]|{k_0}="z" &
      {}\ar[d]^{g_1} \ar@2 [dl] **{} ?(.2); ?(.8) |{g_2}\\{} \ar[r]_{f'}& {} 
      \POS \ar@2 "x"; "y" **{} ?(.2); ?(.8)|{\omega_1} 
      \POS \ar@2 "y"; "z" **{} ?(.2); ?(.8)|{\gamma_1} 
    } 
    +(1, 0)
    \xyboxmatrix"B"{
      {}\ar[r]^f \ar[d]|{h_0}="x" \ar@/_2.5pc/[d]|{k_0}="y" &
      {}\ar[d]|{h_1}="w" \ar@/^2.5pc/[d]|{g_1}="z" 
      \ar@2 [dl] **{} ?(.2); ?(.8) |{h_2}\\{}
      \ar[r]_{f'}& {} 
      \POS \ar@2 "x"; "y" **{} ?(.2); ?(.8) |{\gamma_1} 
      \POS \ar@2 "z"; "w" **{} ?(.2); ?(.8) |{\omega_2} 
    } 
    +(1, 0)
    \xyboxmatrix"C"{
      {} \ar[r]^f \ar[d]_{k_0}  &
      {}\ar[d]|{k_1}="z" \ar@2 [dl] **{} ?(.2); ?(.8) |{k_2}
      \ar@/^2.5pc/[d]|{h_1}="y" \ar@/^5pc/[d]|{g_1}="x" \\{}
      \ar[r]_{f'}& {} 
      \POS \ar@2 "x"; "y" **{} ?(.2); ?(.8)|{\omega_2} 
      \POS \ar@2 "y"; "z" **{} ?(.2); ?(.8)|{\gamma_2} 
    }
    ,(0,-1)
    \xyboxmatrix"A1"{
      {} \ar[r]^f \ar[d]|{g_0}="x"
      \ar@/_2.5pc/[d]|{h_0}="y"
      \ar@/_5pc/[d]|{k_0}="z" &
      {}\ar[d]^{g_1} \ar@2 [dl] **{} ?(.2); ?(.8) |{g_2}\\{} \ar[r]_{f'}& {} 
      \POS \ar@2 "x"; "y" **{} ?(.2); ?(.8)|{\alpha_1} 
      \POS \ar@2 "y"; "z" **{} ?(.2); ?(.8)|{\gamma_1} 
    } 
    +(1, 0)
    \xyboxmatrix"B1"{
      {}\ar[r]^f \ar[d]|{h_0}="x" \ar@/_2.5pc/[d]|{k_0}="y" &
      {}\ar[d]|{h_1}="w" \ar@/^2.5pc/[d]|{g_1}="z" 
      \ar@2 [dl] **{} ?(.2); ?(.8) |{h_2}\\{}
      \ar[r]_{f'}& {} 
      \POS \ar@2 "x"; "y" **{} ?(.2); ?(.8) |{\gamma_1} 
      \POS \ar@2 "z"; "w" **{} ?(.2); ?(.8) |{\alpha_2} 
    } 
    +(1, 0)
    \xyboxmatrix"C1"{
      {} \ar[r]^f \ar[d]_{k_0}  &
      {}\ar[d]|{k_1}="z" \ar@2 [dl] **{} ?(.2); ?(.8) |{k_2}
      \ar@/^2.5pc/[d]|{h_1}="y" \ar@/^5pc/[d]|{g_1}="x" \\{}
      \ar[r]_{f'}& {} 
      \POS \ar@2 "x"; "y" **{} ?(.2); ?(.8)|{\alpha_2} 
      \POS \ar@2 "y"; "z" **{} ?(.2); ?(.8)|{\gamma_2} 
    }
    ,(0,-2)
    \xyboxmatrix"D"{
      {} \ar[r]^f \ar[d]|{g_0}="x"
      \ar@/_2.5pc/[d]|{h_0}="y"
      \ar@/_5pc/[d]|{k_0}="z" &
      {}\ar[d]^{g_1} \ar@2 [dl] **{} ?(.2); ?(.8) |{g_2}\\{} \ar[r]_{f'}& {} 
      \POS \ar@2 "x"; "y" **{} ?(.2); ?(.8)|{\beta_1} 
      \POS \ar@2 "y"; "z" **{} ?(.2); ?(.8)|{\gamma_1} 
    } 
    +(1, 0)
    \xyboxmatrix"E"{
      {}\ar[r]^f \ar[d]|{h_0}="x" \ar@/_2.5pc/[d]|{k_0}="y" &
      {}\ar[d]|{h_1}="w" \ar@/^2.5pc/[d]|{g_1}="z" 
      \ar@2 [dl] **{} ?(.2); ?(.8) |{h_2}\\{}
      \ar[r]_{f'}& {} 
      \POS \ar@2 "x"; "y" **{} ?(.2); ?(.8) |{\gamma_1} 
      \POS \ar@2 "z"; "w" **{} ?(.2); ?(.8) |{\beta_2} 
    } 
    +(1, 0)
    \xyboxmatrix"F"{
      {} \ar[r]^f \ar[d]_{k_0}  &
      {}\ar[d]|{k_1}="z" \ar@2 [dl] **{} ?(.2); ?(.8) |{k_2}
      \ar@/^2.5pc/[d]|{h_1}="y" \ar@/^5pc/[d]|{g_1}="x" \\{}
      \ar[r]_{f'}& {} 
      \POS \ar@2 "x"; "y" **{} ?(.2); ?(.8)|{\beta_2} 
      \POS \ar@2 "y"; "z" **{} ?(.2); ?(.8)|{\gamma_2} 
    }
    \ar@3 "A";"B" ^{(f'\#_0\gamma_1)\\\;\#_1\omega_3}
    \ar@3 "B";"C" ^{\gamma_3\\\;\#_1(\omega_2\#_0f)}
    \ar@3 "A1";"B1" ^{(f'\#_0\gamma_1)\\\;\#_1\alpha_3}
    \ar@3 "B1";"C1" ^{\gamma_3\\\;\#_1(\alpha_2\#_0f)}
    \ar@3 "D";"E" _{(f'\#_0\gamma_1)\\\;\#_1\beta_3}
    \ar@3 "E";"F" _{\gamma_3\\\;\#_1(\beta_2\#_0f)}
    \ar@3 "A";"A1" |{(f'\#_0\gamma_1)\\\;\#_1(f'\#_0\Delta_1)\\\;\#_1g_2}
    \ar@3 "B";"B1" |{(f'\#_0\Delta_1)\\\;\#_1h_2\\\;\#_1(\omega_2\#_0f)}
    \ar@3 "C";"C1" |{k_2\\\;\#_1(\Delta_2\#_0f)\\\;\#_1(\omega_2\#_0f)}
    \ar@3 "A1";"D" |{(f'\#_0\gamma_1)\\\;\#_1(f'\#_0\Gamma_1)\\\;\#_1g_2}
    \ar@3 "B1";"E" |{(f'\#_0\Gamma_1)\\\;\#_1h_2\\\;\#_1(\alpha_2\#_0f)}
    \ar@3 "C1";"F" |{k_2\\\;\#_1(\Gamma_2\#_0f)\\\;\#_1(\alpha_2\#_0f)}
    \ar@{} "A";"B1" |{\text{\eqref{eq:3cellcond}}}
    \ar@{} "B";"C1" |{\text{func.}}
    \ar@{} "A1";"E" |{\text{\eqref{eq:3cellcond}}}
    \ar@{} "B1";"F" |{\text{func.}}
  \end{xy}
\end{equation*}

We see that 2-3-whiskering is functorial:
\begin{multline*}
  (\Delta\square_1\beta)\square_2(\gamma\square_1\Gamma)\\=
  (\Delta_1\#_1\beta_1,\Delta_2\#_1\beta_2)
  \square_2(\gamma_1\#_1\Gamma_1, \gamma_2\#_1\Gamma_2)\\=
  ((\Delta_1\#_1\beta_1)\#_2 (\gamma_1\#_1\Gamma_1),
  ((\Delta_2\#_1\beta_2)\#_2(\gamma_2\#_1\Gamma_2))\\=
  ((\delta_1\#_1\Gamma_1)\#_2(\Delta_1\#_1\alpha_1),
  (\delta_2\#_1\Gamma_2)\#_2(\Delta_2\#_2\alpha_2))\\=
  (\delta_1\#_1\Gamma_1,\delta_2\#_1\Gamma_2)\square_2
  (\Delta_1\#_1\alpha_1, \Delta_2\#_1\alpha_2)\\=
  (\delta\square_1\Gamma)\square_2(\Delta\square_1\alpha)\,.
\end{multline*}
So we can conclude that $\pathspc\H$ is locally a 2-category.

That interchange $\bten$ is natural and functorial in both
arguments follows immediately from the respective properties of $\ten$
in $\H$. Thus we have:

\begin{lem}
  \label{lem:pathspcgrcat}
  The path space $\pathspc\H$ for a $\Gray$-category $\H$ is again a
  $\Gray$-category. \qed
\end{lem}

\begin{lem}
  Given a $\Gray$-functor $F\from\G\to\H$ there is a canonical
  $\Gray$-functor $\pathspc{F}\from\pathspc\G\to\pathspc\H$.
\end{lem}

\begin{prf}
  The $\Gray$-functor $\pathspc{F}$ acts by applying $F$ to all
  components of the cells of $\pathspc\G$:
  \begin{align*}
    \left(\xymatrix{x \ar[r]^f& y}\right)&\mapsto\left(\xymatrix{Fx \ar[r]^{Ff}& Fy}\right) \\
    \left(
      \begin{matrix}
        \xymatrix{{} \ar[r]^f \ar[d]_{g_0}&
          {}\ar[d]^{g_1} \ar@2 [dl] **{} ?<(.3); ?>(.7) ^{g_2}\\{} \ar[r]_{f'}& {}}
      \end{matrix}
    \right)&\mapsto     
    \left(
      \begin{matrix}
        \xymatrix{{} \ar[r]^{Ff} \ar[d]_{Fg_0}&
          {}\ar[d]^{Fg_1} \ar@2 [dl] **{} ?<(.3); ?>(.7) ^{Fg_2}\\{} \ar[r]_{Ff'}& {}}
      \end{matrix}
    \right)\\
    \left(
      \begin{matrix}
        \begin{xy}
          \xymatrix"A"{{} \ar[r]^f \ar[d]|{g_0}="x"\ar@/_3pc/[d]_{h_0}="y" &
            {}\ar[d]^{g_1} \ar@2 [dl] **{} ?(.2); ?(.8) |{g_2}\\{}
            \ar[r]_{f'}& {} \POS \ar@2 "x"; "y" **{} ?<(.3); ?>(.7) ^{\alpha_1} }
          \POS + (30,0)
          \xymatrix"B"{{} \ar[r]^f \ar[d]_{h_0}&
            {}\ar[d]|{h_1}="x" \ar@2 [dl] **{} ?(.2); ?(.8) |{h_2}
            \ar@/^3pc/[d]^{g_1}="y"\\{} \ar[r]_{f'}& {} \POS \ar@2
            "y"; "x" **{} ?<(.3); ?>(.7) ^{\alpha_2}}
          \POS \ar@3 {"A1,1" ; "A2,2" **{} ?(.5)}; {"B1,1" ; "B2,2" **{} ?(.5)} **{}
          ?<(.4);  ?>(.6) ^ {\alpha_3} 
        \end{xy}
      \end{matrix}
    \right)&\mapsto
    \left(
      \begin{matrix}
        \begin{xy}
          \xymatrix"A"{{} \ar[r]^{Ff} \ar[d]|{Fg_0}="x"\ar@/_3pc/[d]_{Fh_0}="y" &
            {}\ar[d]^{Fg_1} \ar@2 [dl] **{} ?(.2); ?(.8) |{Fg_2}\\{}
            \ar[r]_{Ff'}& {} \POS \ar@2 "x"; "y" **{} ?<(.3); ?>(.7) ^{F\alpha_1} }
          \POS + (30,0)
          \xymatrix"B"{{} \ar[r]^{Ff} \ar[d]_{Fh_0}&
            {}\ar[d]|{Fh_1}="x" \ar@2 [dl] **{} ?(.2); ?(.8) |{Fh_2}
            \ar@/^3pc/[d]^{Fg_1}="y"\\{} \ar[r]_{Ff'}& {} \POS \ar@2
            "y"; "x" **{} ?<(.3); ?>(.7) ^{F\alpha_2}}
          \POS \ar@3 {"A1,1" ; "A2,2" **{} ?(.5)}; {"B1,1" ; "B2,2" **{} ?(.5)} **{}
          ?<(.4);  ?>(.6) ^ {F\alpha_3} 
        \end{xy}
      \end{matrix}
    \right)\\
    \begin{pmatrix}
      \begin{xy}\xycompile{ \xyboxmatrix"A"{{} \ar[r]^f
            \ar[d]|{g_0}="x"\ar@/_3pc/[d]|{h_0}="y" & {}\ar[d]^{g_1}
            \ar@2 [dl] **{} ?(.2); ?(.8) |{g_2}\\{} \ar[r]_{f'}& {}
            \POS \ar@2 "x"; "y" **{} ?(.2); ?(.8) |{\alpha_1} } \POS +
          (40,0) \xyboxmatrix"B"{{} \ar[r]^f \ar[d]_{h_0}&
            {}\ar[d]|{h_1}="x" \ar@2 [dl] **{} ?(.2); ?(.8) |{h_2}
            \ar@/^3pc/[d]|{g_1}="y"\\{} \ar[r]_{f'}& {} \POS \ar@2
            "y"; "x" **{} ?(.2); ?(.8) |{\alpha_2}} \POS + (-40,-40)
          \xyboxmatrix"C"{{} \ar[r]^f
            \ar[d]|{g_0}="x"\ar@/_3pc/[d]|{h_0}="y" & {}\ar[d]^{g_1}
            \ar@2 [dl] **{} ?(.2); ?(.8) |{g_2}\\{} \ar[r]_{f'}& {}
            \POS \ar@2 "x"; "y" **{} ?(.2); ?(.8) |{\beta_1} } \POS +
          (40,0) \xyboxmatrix"D"{{} \ar[r]^f \ar[d]_{h_0}&
            {}\ar[d]|{h_1}="x" \ar@2 [dl] **{} ?(.2); ?(.8) |{h_2}
            \ar@/^3pc/[d]|{g_1}="y"\\{} \ar[r]_{f'}& {} \POS \ar@2
            "y"; "x" **{} ?(.2); ?(.8) |{\beta_2}} \ar@3 "A";"B"
          ^{\underline{\alpha_3}} \ar@3 "C";"D"
          _{\underline{\beta_3}} \ar@3 "A";"C"
          |{(f'\#_0\underline{\Gamma_1})\#_1g_2} \ar@3 "B";"D"
          |{h_2\#_1(\underline{\Gamma_2}\#_0f)} }\end{xy}
    \end{pmatrix}
    &\mapsto
    \begin{pmatrix}
      \begin{xy}\xycompile{ \xyboxmatrix"A"{{} \ar[r]^{Ff}
            \ar[d]|{Fg_0}="x"\ar@/_3pc/[d]|{Fh_0}="y" &
            {}\ar[d]^{Fg_1} \ar@2 [dl] **{} ?(.2); ?(.8) |{Fg_2}\\{}
            \ar[r]_{Ff'}& {} \POS \ar@2 "x"; "y" **{} ?(.2); ?(.8)
            |{F\alpha_1} } \POS + (40,0) \xyboxmatrix"B"{{}
            \ar[r]^{Ff} \ar[d]_{Fh_0}& {}\ar[d]|{Fh_1}="x" \ar@2 [dl]
            **{} ?(.2); ?(.8) |{Fh_2} \ar@/^3pc/[d]|{Fg_1}="y"\\{}
            \ar[r]_{Ff'}& {} \POS \ar@2 "y"; "x" **{} ?(.2); ?(.8)
            |{F\alpha_2}} \POS + (-40,-40) \xyboxmatrix"C"{{}
            \ar[r]^{Ff} \ar[d]|{Fg_0}="x"\ar@/_3pc/[d]|{Fh_0}="y" &
            {}\ar[d]^{Fg_1} \ar@2 [dl] **{} ?(.2); ?(.8) |{Fg_2}\\{}
            \ar[r]_{Ff'}& {} \POS \ar@2 "x"; "y" **{} ?(.2); ?(.8)
            |{F\beta_1} } \POS + (40,0) \xyboxmatrix"D"{{} \ar[r]^{Ff}
            \ar[d]_{Fh_0}& {}\ar[d]|{Fh_1}="x" \ar@2 [dl] **{} ?(.2);
            ?(.8) |{Fh_2} \ar@/^3pc/[d]|{Fg_1}="y"\\{} \ar[r]_{Ff'}&
            {} \POS \ar@2 "y"; "x" **{} ?(.2); ?(.8) |{F\beta_2}}
          \ar@3 "A";"B" ^{\underline{F\alpha_3}} \ar@3 "C";"D"
          _{\underline{F\beta_3}} \ar@3 "A";"C"
          |{(Ff'\#_0\underline{F\Gamma_1})\#_1Fg_2} \ar@3 "B";"D"
          |{Fh_2\#_1(\underline{F\Gamma_2}\#_0Ff)} }\end{xy}
    \end{pmatrix}
  \end{align*}
  This preserves the structure of $\pathspc\G$ since $F$ preserves all
  commuting diagrams on the nose. 
\end{prf}

\begin{thm}
  \label{thm:pathspcfnctr}
  Furthermore, $\pathspc{(-)}$ is canonically an
  endofunctor of $\Gray\Cat$.
\end{thm}

\begin{prf}
  Obviously $\pathspc{GF}=\pathspc{G}\pathspc{F}$.
\end{prf}

We finally note the following:
\begin{lem}
\label{lem:pathlim}
  The functor $\pathspc{(-)}\from\Gray\Cat\to\Gray\Cat$  preserves limits.
\end{lem}
\begin{prf}
  This is obviously true for products. 

  For the equalizer $\Bbb E$ of two strict maps $F,G$ we remember that
  the action of $\pathspc F$ and $\pathspc G$ is defined by the
  component wise action of $F$ and $G$, that is, a cell of
  $\pathspc{\Bbb E}$ is equal under $\pathspc F$ and $\pathspc G$ iff
  its components are so under $F$ and $G$.
\end{prf}

A straightforward calculation shows how this forms part of an adjunction
\begin{equation*}
  \xymatrix{
    \Gray\Cat\ar@/^/[r]^-{\pathspc{(\_)}}="x"&\Gray\Cat\ar@/^/[l]^-{\_\ten\I}="y"
    \ar@{|- }"x";"y"**{}?(.4);?(.6)
  }
\end{equation*}
where $\I$ is the free $\Gray$-category on a single 1-cell
$(01)\from{}0\to{}1$ and $\ten$ is Crans' tensor of
$\Gray$-categories.

\section{Composition of Paths}
\label{sec:compaths}

We want to turn the path space that we constructed in the previous
section into the arrow part of an internal category, which requires us
to define a composition map as follows:  

\begin{defn}
  \label{defn:pathcomp}
  We define the \defterm{composite of paths} as a pseudo $\fQ^1$ graph
  map $m\from\pathspc\H\times_{\H}\pathspc\H\laxto\pathspc\H$ by
  horizontal pasting in the following fashion:
  \begin{enumerate}
  \item \label{item:pathcomp0} 0-cells
    \begin{equation*}
      \begin{pmatrix}
        \begin{xy}\xycompile{ \xymatrix{y\ar[r]^{\hat{f}}&z}
          }\end{xy},
        \begin{xy}\xycompile{ \xymatrix{x\ar[r]^f&y} }\end{xy}
      \end{pmatrix}
      \mapsto
      \begin{pmatrix}
        \begin{xy}\xycompile{ \xymatrix{ x\ar[rr]^{\hat{f}\#_0f}&{}& z
            } }\end{xy}
      \end{pmatrix}
    \end{equation*}
  \item \label{item:pathcomp1} 1-cells
    \begin{multline*}
      \begin{pmatrix}
        \begin{matrix}
          \xymatrix{ {}\ar[r]^{\hat f}\ar[d]_{\hat{g_0}=g_1}&
            {}\ar[d]^{\hat{g_1}} \ar@2
            [dl]**{}?(.2);?(.8)|{\hat{g_2}}\\
            {}\ar[r]_{\hat{f'}}&{} }
        \end{matrix}\;,
        \begin{matrix}
          \xymatrix{
            {}\ar[r]^f\ar[d]_{g_0}&{}\ar[d]^{g_1}\ar@2[dl]**{}?(.2);?(.8)|{g_2}\\
            {}\ar[r]_{f'}&{} }
        \end{matrix}
      \end{pmatrix}
      \mapsto
      \begin{pmatrix}
        \begin{matrix}
          \xymatrix{ {}\ar[r]^f\ar[d]_{g_0}&
            {}\ar[d]|{g_1}\ar@2[dl]**{}?(.2);?(.8)|{g_2} \ar[r]^{\hat
              f}&
            {}\ar[d]^{\hat{g_1}}\ar@2[dl]**{}?(.2);?(.8)|{\hat{g_2}}\\
            {}\ar[r]_{f'}&{}\ar[r]_{\hat{f'}}&{} }
        \end{matrix}
      \end{pmatrix}
      \\=\begin{pmatrix}
        \begin{matrix}
          \xymatrix@+1cm{ {}\ar[rr]^{\hat{f}\#_0f}\ar[d]_{g_0}& {}&
            {}\ar[d]^{\hat{g_1}} \ar@2[dll]**{}?(.2);?(.8)|{
              (\hat{f'}\#_0g_2)\\\#_1(\hat{g_2}\#_0{f})
            }\\
            {}\ar[rr]_{\hat{f'}\#_0f'}&&{} }
        \end{matrix}
      \end{pmatrix}
    \end{multline*}
  \item \label{item:pathcomp2} 2-cells
    \begin{multline*}
      \begin{pmatrix}
        \begin{matrix}
          \begin{xy}\xycompile{ \xymatrix"A"{{} \ar[r]^{\hat{f}}
                \ar[d]|{\hat{g_0}}="x"\ar@/_3pc/[d]_{\hat{h_0}}="y" &
                {}\ar[d]^{\hat{g_1}} \ar@2 [dl] **{} ?(.2); ?(.8)
                |{\hat{g_2}}\\{} \ar[r]_{\hat{f'}}& {} \POS \ar@2 "x";
                "y" **{} ?<(.3); ?>(.7) ^{\hat{\alpha_1}=\alpha_2}}
              \POS + (30,0) \xymatrix"B"{{} \ar[r]^{\hat{f}}
                \ar[d]_{\hat{h_0}}& {}\ar[d]|{\hat{h_1}}="x" \ar@2
                [dl] **{} ?(.2); ?(.8) |{\hat{h_2}}
                \ar@/^3pc/[d]^{\hat{g_1}}="y"\\{} \ar[r]_{\hat{f'}}&
                {} \POS \ar@2 "y"; "x" **{} ?<(.3); ?>(.7)
                ^{\hat{\alpha_2}}} \POS \ar@3 {"A1,1" ; "A2,2" **{}
                ?(.5)}; {"B1,1" ; "B2,2" **{} ?(.5)} **{} ?<(.4);
              ?>(.6) ^ {\hat{\alpha_3}} }\end{xy}
        \end{matrix}\;,
        \begin{matrix}
          \begin{xy}\xycompile{ \xymatrix"A"{{} \ar[r]^f
                \ar[d]|{g_0}="x"\ar@/_3pc/[d]_{h_0}="y" &
                {}\ar[d]^{g_1} \ar@2 [dl] **{} ?(.2); ?(.8) |{g_2}\\{}
                \ar[r]_{f'}& {} \POS \ar@2 "x"; "y" **{} ?<(.3);
                ?>(.7) ^{\alpha_1} } \POS + (30,0) \xymatrix"B"{{}
                \ar[r]^f \ar[d]_{h_0}& {}\ar[d]|{h_1}="x" \ar@2 [dl]
                **{} ?(.2); ?(.8) |{h_2} \ar@/^3pc/[d]^{g_1}="y"\\{}
                \ar[r]_{f'}& {} \POS \ar@2 "y"; "x" **{} ?<(.3);
                ?>(.7) ^{\alpha_2}} \POS \ar@3 {"A1,1" ; "A2,2" **{}
                ?(.5)}; {"B1,1" ; "B2,2" **{} ?(.5)} **{} ?<(.4);
              ?>(.6) ^ {\alpha_3} }\end{xy}
        \end{matrix}
      \end{pmatrix}\\\mapsto
      \begin{pmatrix}
        \begin{matrix}
          \begin{xy}\xycompile{ \xymatrix"A"{ {} \ar[r]^f
                \ar[d]|{g_0}="x"\ar@/_3pc/[d]_{h_0}="y" &
                {}\ar[d]|{g_1}\ar@2[dl]**{}?(.2);?(.8)|{g_2}{}\ar[r]^{\hat
                  f}&
                {}\ar[d]^{\hat{g_1}}\ar@2[dl]**{}?(.2);?(.8)|{\hat{g_2}}\\
                {}\ar[r]_{f'}&{}
                \POS\ar@2"x";"y"**{}?<(.3);?>(.7)^{\alpha_1}{}\ar[r]_{\hat{f'}}&{}
              } \POS + (45,0) \xymatrix"B"{ {}\ar[r]^f\ar[d]_{h_0}&
                {}\ar@/_1pc/[d]|{h_1}="x" \ar@2
                {[l]**{}?};[dl]**{}?(.2);?(.8)|{h_2}
                \ar@/^1pc/[d]|{g_1}="y"\ar[r]^{\hat{f}} &
                {}\ar[d]^{\hat{g_1}}
                \ar@2 {[d];[dl]**{}?}**{}?(.2);?(.8)|{\hat{g_2}}\\
                {} \ar[r]_{f'}& \ar[r]_{\hat{f'}}&{}
                \POS\ar@2"y";"x"**{}?<(.3);?>(.7)^{\alpha_2} } \POS +
              (45,0) \xymatrix"C"{ {} \ar[r]^f \ar[d]_{h_0} &
                {}\ar[d]|{h_1}\ar@2[dl]**{}?(.2);?(.8)|{h_2}{}\ar[r]^{\hat
                  f}&
                {}\ar[d]|{\hat{h_1}}="x"\ar@2[dl]**{}?(.2);?(.8)|{\hat{h_2}}\ar@/^3pc/[d]^{h_0}="y"\\
                {}\ar[r]_{f'}&{}
                \POS\ar@2"y";"x"**{}?<(.3);?>(.7)^{\hat{\alpha_2}}{}\ar[r]_{\hat{f'}}&{}
              } \POS \ar@3 {"A1,1" ; "A2,3" **{} ?(.5)}; {"B1,1" ;
                "B2,3" **{} ?(.5)} **{} ?<(.4); ?>(.6) ^
              {(\hat{f'}\#_0\alpha_3)\\\;\#_1(\hat{g_2}\#_0f)} \POS
              \ar@3 {"B1,1" ; "B2,3" **{} ?(.5)}; {"C1,1" ; "C2,3"
                **{} ?(.5)} **{} ?<(.4); ?>(.6) ^
              {(\hat{f'}\#_0h_2)\\\;\#_1(\hat{\alpha_3}\#_0f)} }
          \end{xy}
        \end{matrix}
      \end{pmatrix}
    \end{multline*}
  \item \label{item:pathcomp3} 3-cells
    \begin{multline*}
      \begin{pmatrix}
        \begin{matrix}
          \begin{xy}\xycompile{ \xymatrix"A"{{} \ar[r]^{\hat{f}}
                \ar[d]|{\hat{g_0}}="x"\ar@/_3pc/[d]_{\hat{h_0}}="y" &
                {}\ar[d]^{\hat{g_1}} \ar@2 [dl] **{} ?(.2); ?(.8)
                |{\hat{g_2}}\\{} \ar[r]_{\hat{f'}}& {} \POS \ar@2 "x";
                "y" **{} ?<(.3); ?>(.7) ^{\hat{\alpha_1}=\alpha_2}}
              \POS + (30,0) \xymatrix"B"{{} \ar[r]^{\hat{f}}
                \ar[d]_{\hat{h_0}}& {}\ar[d]|{\hat{h_1}}="x" \ar@2
                [dl] **{} ?(.2); ?(.8) |{\hat{h_2}}
                \ar@/^3pc/[d]^{\hat{g_1}}="y"\\{} \ar[r]_{\hat{f'}}&
                {} \POS \ar@2 "y"; "x" **{} ?<(.3); ?>(.7)
                ^{\hat{\alpha_2}}} \POS (0,-30) \xymatrix"C"{{}
                \ar[r]^{\hat{f}}
                \ar[d]|{\hat{g_0}}="x"\ar@/_3pc/[d]_{\hat{h_0}}="y" &
                {}\ar[d]^{\hat{g_1}} \ar@2 [dl] **{} ?(.2); ?(.8)
                |{\hat{g_2}}\\{} \ar[r]_{\hat{f'}}& {} \POS \ar@2 "x";
                "y" **{} ?<(.3); ?>(.7) ^{\hat{\beta_1}=\beta_2}} \POS
              + (30,0) \xymatrix"D"{{} \ar[r]^{\hat{f}}
                \ar[d]_{\hat{h_0}}& {}\ar[d]|{\hat{h_1}}="x" \ar@2
                [dl] **{} ?(.2); ?(.8) |{\hat{h_2}}
                \ar@/^3pc/[d]^{\hat{g_1}}="y"\\{} \ar[r]_{\hat{f'}}&
                {} \POS \ar@2 "y"; "x" **{} ?<(.3); ?>(.7)
                ^{\hat{\beta_2}}} \POS \ar@3 {"A1,1" ; "A2,2" **{}
                ?(.5)}; {"B1,1" ; "B2,2" **{} ?(.5)} **{} ?<(.4);
              ?>(.6) ^ {\hat{\alpha_3}} \POS \ar@3 {"C1,1" ; "C2,2"
                **{} ?(.5)}; {"D1,1" ; "D2,2" **{} ?(.5)} **{} ?<(.4);
              ?>(.6) _ {\hat{\beta_3}} \POS \ar@3 {"A1,1" ; "A2,2"
                **{} ?(.5)}; {"C1,1" ; "C2,2" **{} ?(.5)} **{} ?<(.4);
              ?>(.6) _
              {(\hat{f'}\#_0\hat{\Gamma_1})\#_1\hat{g_2}\\
                =(\hat{f'}\#_0{\Gamma_2})\#_1\hat{g_2}} \POS \ar@3
              {"B1,1" ; "B2,2" **{} ?(.5)}; {"D1,1" ; "D2,2" **{}
                ?(.5)} **{} ?<(.4); ?>(.6) ^
              {\hat{h_2}\#_1(\hat{\Gamma_2}\#_0\hat{f})} }\end{xy}
        \end{matrix}\;,
        \begin{matrix}
          \begin{xy}\xycompile{ \xymatrix"A"{{} \ar[r]^{{f}}
                \ar[d]|{{g_0}}="x"\ar@/_3pc/[d]_{{h_0}}="y" &
                {}\ar[d]^{{g_1}} \ar@2 [dl] **{} ?(.2); ?(.8)
                |{{g_2}}\\{} \ar[r]_{{f'}}& {} \POS \ar@2 "x"; "y"
                **{} ?<(.3); ?>(.7) ^{{\alpha_1}=\alpha_2}} \POS +
              (30,0) \xymatrix"B"{{} \ar[r]^{{f}} \ar[d]_{{h_0}}&
                {}\ar[d]|{{h_1}}="x" \ar@2 [dl] **{} ?(.2); ?(.8)
                |{{h_2}} \ar@/^3pc/[d]^{{g_1}}="y"\\{} \ar[r]_{{f'}}&
                {} \POS \ar@2 "y"; "x" **{} ?<(.3); ?>(.7)
                ^{{\alpha_2}}} \POS (0,-30) \xymatrix"C"{{}
                \ar[r]^{{f}}
                \ar[d]|{{g_0}}="x"\ar@/_3pc/[d]_{{h_0}}="y" &
                {}\ar[d]^{{g_1}} \ar@2 [dl] **{} ?(.2); ?(.8)
                |{{g_2}}\\{} \ar[r]_{{f'}}& {} \POS \ar@2 "x"; "y"
                **{} ?<(.3); ?>(.7) ^{{\beta_1}=\beta_2}} \POS +
              (30,0) \xymatrix"D"{{} \ar[r]^{{f}} \ar[d]_{{h_0}}&
                {}\ar[d]|{{h_1}}="x" \ar@2 [dl] **{} ?(.2); ?(.8)
                |{{h_2}} \ar@/^3pc/[d]^{{g_1}}="y"\\{} \ar[r]_{{f'}}&
                {} \POS \ar@2 "y"; "x" **{} ?<(.3); ?>(.7)
                ^{{\beta_2}}} \POS \ar@3 {"A1,1" ; "A2,2" **{} ?(.5)};
              {"B1,1" ; "B2,2" **{} ?(.5)} **{} ?<(.4); ?>(.6) ^
              {{\alpha_3}} \POS \ar@3 {"C1,1" ; "C2,2" **{} ?(.5)};
              {"D1,1" ; "D2,2" **{} ?(.5)} **{} ?<(.4); ?>(.6) _
              {{\beta_3}} \POS \ar@3 {"A1,1" ; "A2,2" **{} ?(.5)};
              {"C1,1" ; "C2,2" **{} ?(.5)} **{} ?<(.4); ?>(.6) _
              {({f'}\#_0{\Gamma_1})\#_1{g_2}} \POS \ar@3 {"B1,1" ;
                "B2,2" **{} ?(.5)}; {"D1,1" ; "D2,2" **{} ?(.5)} **{}
              ?<(.4); ?>(.6) ^ {{h_2}\#_1({\Gamma_2}\#_0{f})}
            }\end{xy}
        \end{matrix}
      \end{pmatrix}\\\mapsto
      \begin{pmatrix}
        \begin{matrix}
          \begin{xy}\xycompile{ \xymatrix"A"{ {} \ar[r]^f
                \ar[d]|{g_0}="x"\ar@/_3pc/[d]_{h_0}="y" &
                {}\ar[d]|{g_1}\ar@2[dl]**{}?(.2);?(.8)|{g_2}{}\ar[r]^{\hat
                  f}&
                {}\ar[d]^{\hat{g_1}}\ar@2[dl]**{}?(.2);?(.8)|{\hat{g_2}}\\
                {}\ar[r]_{f'}&{}
                \POS\ar@2"x";"y"**{}?<(.3);?>(.7)^{\alpha_1}{}\ar[r]_{\hat{f'}}&{}
              } \POS + (45,0) \xymatrix"B"{ {}\ar[r]^f\ar[d]_{h_0}&
                {}\ar@/_1pc/[d]|{h_1}="x" \ar@2
                {[l]**{}?};[dl]**{}?(.2);?(.8)|{h_2}
                \ar@/^1pc/[d]|{g_1}="y"\ar[r]^{\hat{f}} &
                {}\ar[d]^{\hat{g_1}}
                \ar@2 {[d];[dl]**{}?}**{}?(.2);?(.8)|{\hat{g_2}}\\
                {} \ar[r]_{f'}& \ar[r]_{\hat{f'}}&{}
                \POS\ar@2"y";"x"**{}?<(.3);?>(.7)^{\alpha_2} } \POS +
              (45,0) \xymatrix"C"{ {} \ar[r]^f \ar[d]_{h_0} &
                {}\ar[d]|{h_1}\ar@2[dl]**{}?(.2);?(.8)|{h_2}{}\ar[r]^{\hat
                  f}&
                {}\ar[d]|{\hat{h_1}}="x"\ar@2[dl]**{}?(.2);?(.8)|{\hat{h_2}}\ar@/^3pc/[d]^{h_0}="y"\\
                {}\ar[r]_{f'}&{}
                \POS\ar@2"y";"x"**{}?<(.3);?>(.7)^{\hat{\alpha_2}}{}\ar[r]_{\hat{f'}}&{}
              } \POS(0,-45) \xymatrix"D"{ {} \ar[r]^f
                \ar[d]|{g_0}="x"\ar@/_3pc/[d]_{h_0}="y" &
                {}\ar[d]|{g_1}\ar@2[dl]**{}?(.2);?(.8)|{g_2}{}\ar[r]^{\hat
                  f}&
                {}\ar[d]^{\hat{g_1}}\ar@2[dl]**{}?(.2);?(.8)|{\hat{g_2}}\\
                {}\ar[r]_{f'}&{}
                \POS\ar@2"x";"y"**{}?<(.3);?>(.7)^{\beta_1}{}\ar[r]_{\hat{f'}}&{}
              } \POS + (45,0) \xymatrix"E"{ {}\ar[r]^f\ar[d]_{h_0}&
                {}\ar@/_1pc/[d]|{h_1}="x" \ar@2
                {[l]**{}?};[dl]**{}?(.2);?(.8)|{h_2}
                \ar@/^1pc/[d]|{g_1}="y"\ar[r]^{\hat{f}} &
                {}\ar[d]^{\hat{g_1}}
                \ar@2 {[d];[dl]**{}?}**{}?(.2);?(.8)|{\hat{g_2}}\\
                {} \ar[r]_{f'}& \ar[r]_{\hat{f'}}&{}
                \POS\ar@2"y";"x"**{}?<(.3);?>(.7)^{\beta_2} } \POS +
              (45,0) \xymatrix"F"{ {} \ar[r]^f \ar[d]_{h_0} &
                {}\ar[d]|{h_1}\ar@2[dl]**{}?(.2);?(.8)|{h_2}{}\ar[r]^{\hat
                  f}&
                {}\ar[d]|{\hat{h_1}}="x"\ar@2[dl]**{}?(.2);?(.8)|{\hat{h_2}}\ar@/^3pc/[d]^{h_0}="y"\\
                {}\ar[r]_{f'}&{}
                \POS\ar@2"y";"x"**{}?<(.3);?>(.7)^{\hat{\beta_2}}{}\ar[r]_{\hat{f'}}&{}
              } \POS \ar@3 {"A1,1" ; "A2,3" **{} ?(.5)}; {"B1,1" ;
                "B2,3" **{} ?(.5)} **{} ?<(.4); ?>(.6) ^
              {(\hat{f'}\#_0\alpha_3)\\\;\#_1(\hat{g_2}\#_0f)} \POS
              \ar@3 {"B1,1" ; "B2,3" **{} ?(.5)}; {"C1,1" ; "C2,3"
                **{} ?(.5)} **{} ?<(.4); ?>(.6) ^
              {(\hat{f'}\#_0h_2)\\\;\#_1(\hat{\alpha_3}\#_0f)} \POS
              \ar@3 {"D1,1" ; "D2,3" **{} ?(.5)}; {"E1,1" ; "E2,3"
                **{} ?(.5)} **{} ?<(.4); ?>(.6) _
              {(\hat{f'}\#_0\beta_3)\\\;\#_1(\hat{g_2}\#_0f)} \POS
              \ar@3 {"E1,1" ; "E2,3" **{} ?(.5)}; {"F1,1" ; "F2,3"
                **{} ?(.5)} **{} ?<(.4); ?>(.6) _
              {(\hat{f'}\#_0h_2)\\\;\#_1(\hat{\beta_3}\#_0f)} \POS
              \ar@3 {"A1,1" ; "A2,3" **{} ?(.5)}; {"D1,1" ; "D2,3"
                **{} ?(.5)} **{} ?<(.4); ?>(.6) _
              {(\hat{f'}\#_0f'\#_0\Gamma_1)\\
                \;\#_1(\hat{f'}\#_0g_2)\\
                \;\#_1(\hat{g_2}\#_0f)} \POS \ar@3 {"C1,1" ; "C2,3"
                **{} ?(.5)}; {"F1,1" ; "F2,3" **{} ?(.5)} **{} ?<(.4);
              ?>(.6) ^
              {(\hat{f'}\#_0h_2)\\
                \;\#_1(\hat{h_2}\#_0f)\\
                \;\#_1(\hat{\Gamma_2}\#_0\hat{f}\#_0f)} \POS \ar@3
              {"B1,1" ; "B2,3" **{} ?(.5)}; {"E1,1" ; "E2,3" **{}
                ?(.5)} **{} ?<(.4); ?>(.6) |
              {(\hat{f'}\#_0h_2)\#_1(\hat{f'}\#_0\Gamma_2\#_0f)\#_1(\hat{g_2}\#_0f)}
            }
          \end{xy}
        \end{matrix}
      \end{pmatrix}
    \end{multline*}
  \item \label{item:pathcomp2coc} the 2-cocycle: for a (vertically) composable pair in
    $\pathspc\H\times_\H\pathspc\H$ we have the composite of the
    images and the image of the composites under $m$:
    \begin{equation}
      \label{eq:pasteup}
      \begin{matrix}
        m\begin{pmatrix}
          \begin{matrix}
            \xymatrix{ {}\ar[r]^{\hat f}\ar[d]_{\hat{g_0}\\=g_1}&
              {}\ar[d]^{\hat{g_1}} \ar@2
              [dl]**{}?(.2);?(.8)|{\hat{g_2}}\\
              {}\ar[r]_{\hat{f'}}&{} }
          \end{matrix}\;,
          \begin{matrix}
            \xymatrix{
              {}\ar[r]^f\ar[d]_{g_0}&{}\ar[d]^{g_1}\ar@2[dl]**{}?(.2);?(.8)|{g_2}\\
              {}\ar[r]_{f'}&{} }
          \end{matrix}
        \end{pmatrix}\\\square_0\\
        m\begin{pmatrix}
          \begin{matrix}
            \xymatrix{ {}\ar[r]^{\hat f'}\ar[d]_{\hat{g'_0}\\=g'_1}&
              {}\ar[d]^{\hat{g'_1}} \ar@2
              [dl]**{}?(.2);?(.8)|{\hat{g'_2}}\\
              {}\ar[r]_{\hat{f''}}&{} }
          \end{matrix}\;,
          \begin{matrix}
            \xymatrix{
              {}\ar[r]^f\ar[d]_{g'_0}&{}\ar[d]^{g'_1}\ar@2[dl]**{}?(.2);?(.8)|{g'_2}\\
              {}\ar[r]_{f''}&{} }
          \end{matrix}
        \end{pmatrix}
      \end{matrix}
      =
      \begin{pmatrix}
        \begin{matrix}
          \xymatrix{ {}\ar[r]^f\ar[d]_{g_0}&
            {}\ar[d]|{g_1}\ar@2[dl]**{}?(.2);?(.8)|{g_2}="x"
            \ar[r]^{\hat f}&
            {}\ar[d]^{\hat{g_1}}\ar@2[dl]**{}?(.2);?(.8)|{\hat{g_2}}\\
            {}\ar[r]|{f'}\ar[d]_{g'_0}&
            {}\ar[r]|{\hat{f'}}\ar[d]|{\hat{g'_1}}\ar@2[dl]**{}?(.2);?(.8)|{{g'_2}}&{}
            {}\ar[d]^{\hat{g'_1}}\ar@2[dl]**{}?(.2);?(.8)|{\hat{g'_2}}="y"\\
            {}\ar[r]_{f''}&{}\ar[r]_{\hat{f''}}&{} \POS"x"\tria"y" }
        \end{matrix}
      \end{pmatrix}
    \end{equation}
    \begin{equation}
      \label{eq:pastedown}
      m\begin{pmatrix}
        \begin{matrix}
          \xymatrix{ {}\ar[r]^{\hat f}\ar[d]_{\hat{g_0}\\=g_1}&
            {}\ar[d]^{\hat{g_1}} \ar@2
            [dl]**{}?(.2);?(.8)|{\hat{g_2}}\\
            {}\ar[r]_{\hat{f'}}&{} }\\\square_0\\
          \xymatrix{ {}\ar[r]^{\hat f'}\ar[d]_{\hat{g'_0}\\=g'_1}&
            {}\ar[d]^{\hat{g'_1}} \ar@2
            [dl]**{}?(.2);?(.8)|{\hat{g'_2}}\\
            {}\ar[r]_{\hat{f''}}&{} }
        \end{matrix}\;,
        \begin{matrix}
          \xymatrix{
            {}\ar[r]^f\ar[d]_{g_0}&{}\ar[d]^{g_1}\ar@2[dl]**{}?(.2);?(.8)|{g_2}\\
            {}\ar[r]_{f'}&{} }\\\square_0\\
          \xymatrix{
            {}\ar[r]^f\ar[d]_{g'_0}&{}\ar[d]^{g'_1}\ar@2[dl]**{}?(.2);?(.8)|{g'_2}\\
            {}\ar[r]_{f''}&{} }
        \end{matrix}
      \end{pmatrix}
      =
      \begin{pmatrix}
        \begin{matrix}
          \xymatrix{ {}\ar[r]^f\ar[d]_{g_0}&
            {}\ar[d]|{g_1}\ar@2[dl]**{}?(.2);?(.8)|{g_2}="x"
            \ar[r]^{\hat f}&
            {}\ar[d]^{\hat{g_1}}\ar@2[dl]**{}?(.2);?(.8)|{\hat{g_2}}\\
            {}\ar[r]|{f'}\ar[d]_{g'_0}&
            {}\ar[r]|{\hat{f'}}\ar[d]|{\hat{g'_1}}\ar@2[dl]**{}?(.2);?(.8)|{{g'_2}}&{}
            {}\ar[d]^{\hat{g'_1}}\ar@2[dl]**{}?(.2);?(.8)|{\hat{g'_2}}="y"\\
            {}\ar[r]_{f''}&{}\ar[r]_{\hat{f''}}&{} \POS"y"\tria"x" }
        \end{matrix}
      \end{pmatrix}
    \end{equation}
    And the 2-cocycle going between them is:
    \begin{equation}\label{eq:pastecocdef}
      m^2\begin{pmatrix}
        \begin{pmatrix}
          \begin{matrix}
            \xymatrix{ {}\ar[r]^{\hat f}\ar[d]_{\hat{g_0}\\=g_1}&
              {}\ar[d]^{\hat{g_1}} \ar@2
              [dl]**{}?(.2);?(.8)|{\hat{g_2}}\\
              {}\ar[r]_{\hat{f'}}&{} }
          \end{matrix}\;,
          \begin{matrix}
            \xymatrix{
              {}\ar[r]^f\ar[d]_{g_0}&{}\ar[d]^{g_1}\ar@2[dl]**{}?(.2);?(.8)|{g_2}\\
              {}\ar[r]_{f'}&{} }
          \end{matrix}
        \end{pmatrix},\\
        \begin{pmatrix}
          \begin{matrix}
            \xymatrix{ {}\ar[r]^{\hat f'}\ar[d]_{\hat{g'_0}\\=g'_1}&
              {}\ar[d]^{\hat{g'_1}} \ar@2
              [dl]**{}?(.2);?(.8)|{\hat{g'_2}}\\
              {}\ar[r]_{\hat{f''}}&{} }
          \end{matrix}\;,
          \begin{matrix}
            \xymatrix{
              {}\ar[r]^f\ar[d]_{g'_0}&{}\ar[d]^{g'_1}\ar@2[dl]**{}?(.2);?(.8)|{g'_2}\\
              {}\ar[r]_{f''}&{} }
          \end{matrix}
        \end{pmatrix}
      \end{pmatrix}\from
      \begin{matrix}
        \begin{xy}
          \xycompile{ \xymatrix"A"{ {}\ar[r]^f\ar[d]_{g_0}&
              {}\ar[d]|{g_1}\ar@2[dl]**{}?(.2);?(.8)|{g_2}="x"
              \ar[r]^{\hat f}&
              {}\ar[d]^{\hat{g_1}}\ar@2[dl]**{}?(.2);?(.8)|{\hat{g_2}}\\
              {}\ar[r]|{f'}\ar[d]_{g'_0}&
              {}\ar[r]|{\hat{f'}}\ar[d]|{\hat{g'_1}}\ar@2[dl]**{}?(.2);?(.8)|{{g'_2}}&{}
              {}\ar[d]^{\hat{g'_1}}\ar@2[dl]**{}?(.2);?(.8)|{\hat{g'_2}}="y"\\
              {}\ar[r]_{f''}&{}\ar[r]_{\hat{f''}}&{} \POS"x"\tria"y" }
            \POS+(60,0) \xymatrix"B"{ {}\ar[r]^f\ar[d]_{g_0}&
              {}\ar[d]|{g_1}\ar@2[dl]**{}?(.2);?(.8)|{g_2}="x"
              \ar[r]^{\hat f}&
              {}\ar[d]^{\hat{g_1}}\ar@2[dl]**{}?(.2);?(.8)|{\hat{g_2}}\\
              {}\ar[r]|{f'}\ar[d]_{g'_0}&
              {}\ar[r]|{\hat{f'}}\ar[d]|{\hat{g'_1}}\ar@2[dl]**{}?(.2);?(.8)|{{g'_2}}&{}
              {}\ar[d]^{\hat{g'_1}}\ar@2[dl]**{}?(.2);?(.8)|{\hat{g'_2}}="y"\\
              {}\ar[r]_{f''}&{}\ar[r]_{\hat{f''}}&{} \POS"y"\tria"x" }
            \POS\ar@3
            {"A1,1";"A3,3"**{}?(.5)};{"B1,1";"B3,3"**{}?(.5)}**{}?<(.4);?>(.6)
            ^{(\hat{f''}\#_0{g'_2}\#_0{g_0})\\\#_1(\hat{g'_2}\ten{g_2})\\\#_1(\hat{g'_1}\#_0\hat{g_2}\#_0{f})}
          }
        \end{xy}
      \end{matrix}
    \end{equation} For completeness' sake we give it in the algebraic
    notation:
    \begin{equation*}
      \begin{pmatrix}
        (\hat{f''}\#_0{g'_2}\#_0{g_0})\#_1(\hat{g'_2}\ten{g_2})\#_1(\hat{g'_1}\#_0\hat{g_2}\#_0{f});\\
        \id_{g'_0\#_0g_0}, \id_{\hat{g'_1}\#_0\hat{g_1}},\\
        (\hat{f''}\#_0g'_2\#_0g_0)\#_1(\hat{g'_2}\lhc
        g_2)\#_1(\hat{g'_1}\#_0\hat{g_2}\#_0f),\\
        (\hat{f''}\#_0g'_2\#_0g_0)\#_1(\hat{g'_2}\rhc
        g_2)\#_1(\hat{g'_1}\#_0\hat{g_2}\#_0f);\\
        g'_0\#_0g_0, \hat{g'_1}\#_0\hat{g_1}, g'_0\#_0g_0,
        \hat{g'_1}\#_0\hat{g_1}, \hat{f}\#_0f, \hat{f''}\#_0f''
      \end{pmatrix}
    \end{equation*}
  \end{enumerate}
\end{defn}

\begin{lem}
  \label{lem:mqmapprf}
  The map $m\from\pathspc\H\times_\H\pathspc\H\laxto\pathspc\H$ is a
  pseudo $\fQ^1$ graph map and hence by lemma \ref{lem:pseudofunchar}
  uniquely defines a pseudo $\Gray$-functor.
\end{lem}
\begin{prf}
  As defined above, $m$ is obviously a 3-globular map. We verify that
  it is locally a sesquifunctor: Let $(\beta^1,\beta^2)$ and
  $(\alpha^1,\alpha^2)$ be two pairs of 2-cells in
  $\pathspc\H\times_\H\pathspc\H$ composable along a pair of
  1-cells. Then
  \begin{equation*}
    m((\beta^1,\beta^2)\square_1(\alpha^1,\alpha^2))=m((\beta^1\square_1\alpha^1),(\beta^2\square_1\alpha^2))
    =m(\beta^1,\beta^2)\square_1m(\alpha^1,\alpha^2)
  \end{equation*}
  follows obviously from the fact that in $\H$ 3-cells compose along a
  2-cells interchangeably.  Let $(\Delta^1,\Delta^2)$ and
  $(\Gamma^1,\Gamma^2)$ be two pairs of 3-cells in
  $\pathspc\H\times_\H\pathspc\H$ composable along a pair of
  2-cells. Then
  \begin{multline*}
    m((\Delta^1,\Delta^2)\square_2(\Gamma^1,\Gamma^2))=m((\Delta^1\square_2\Gamma^1),(\Delta^2\square_2\Gamma^2))\\
    =m((\Delta^1_1\#_2\Gamma^1_1,\Delta^1_2\#_2\Gamma^1_2),(\Delta^2_1\#_2\Gamma^2_1,\Delta^2_2\#_2\Gamma^2_2))
    =(\Delta^1_1\#_2\Gamma^1_1,\Delta^2_2\#_2\Gamma^2_2)\\
    =(\Delta^1_1,\Delta^2_2)\square_2(\Gamma^1_1,\Gamma^2_2)
    =m((\Delta^1_1,\Delta^1_2),(\Delta^2_1,\Delta^2_2))\square_2m((\Gamma^1_1,\Gamma^1_2),(\Gamma^2_1,\Gamma^2_2))\\
    =m(\Delta^1,\Delta^2)\square_2m(\Gamma^1,\Gamma^2)\,.
  \end{multline*}
  For the vertical composition of 3-cells see
  (\ref{eq:vert3compdef}), their images under $m$ are pastings of
  commuting diagrams, so preservation is immediate. Preservation of whiskers of 3-cells by
  2-cells given for each component of $\pathspc\H\times_\H\pathspc\H$
  in (\ref{eq:sq23defexp}), again according to definition
  \ref{defn:pathcomp}.\ref{item:pathcomp3}
  $m$ pastes two such commuting diagrams horizontally. Preservation of
  units is trivially satisfied. This concludes verification of
  \ref{defn:psgrmap}.\ref{item:psgrmaplocal}. 

  We verify that $m^2$ is a 2-cocycle in (\ref{eq:mqmapprf2coc}).
  \begin{sidewaysfigure}
    \begin{equation}
      \label{eq:mqmapprf2coc}
      \begin{xy}
        \save
        (70,0):(0,-1)::
        ,(0,0)="A11"
        ,(1,0)="A21"
        ,(2,0)="A31"
        ,(3,0)="A41"
        ,(0,1)="A12"
        ,(1,1)="A22"
        ,(2,1)="A32"
        ,(3,1)="A42"
        \restore
        \POS 
        ,"A11"
        \xyboxmatrix"A11"{
          \ar[d]_{g_0^1}\ar[r]^{f^1}&\ar[d]|{g_1^1}\ar[r]^{f^2}\ar@2[dl]**{}?(.2);?(.8)|{g_2^1}="X1"&
          \ar[d]^{g_1^2}\ar@2[dl]**{}?(.2);?(.8)|{g_2^2}="Y1"\\
          \ar[d]_{h_0^1}\ar[r]|{f'^1}&\ar[d]|{h_1^1}\ar[r]|{f'^2}\ar@2[dl]**{}?(.2);?(.8)|{h_2^1}="X2"&
          \ar[d]^{h_1^2}\ar@2[dl]**{}?(.2);?(.8)|{h_2^2}="Y2"\\
          \ar[d]_{k_0^1}\ar[r]|{f''^1}&\ar[d]|{k_1^1}\ar[r]|{f''^2}\ar@2[dl]**{}?(.2);?(.8)|{k_2^1}="X3"&
          \ar[d]^{k_1^2}\ar@2[dl]**{}?(.2);?(.8)|{k_2^2}="Y3"\\
          \ar[r]_{f'''^1}&\ar[r]_{f'''^2}&{}          
          \POS
          \tria"X1";"Y2"
          \tria"X2";"Y3"
        }
        ,"A21"
        \xyboxmatrix"A21"{
          \ar[d]_{g_0^1}\ar[r]^{f^1}&\ar[d]|{g_1^1}\ar[r]^{f^2}\ar@2[dl]**{}?(.2);?(.8)|{g_2^1}="X1"&
          \ar[d]^{g_1^2}\ar@2[dl]**{}?(.2);?(.8)|{g_2^2}="Y1"\\
          \ar[d]_{h_0^1}\ar[r]|{f'^1}&\ar[d]|{h_1^1}\ar[r]|{f'^2}\ar@2[dl]**{}?(.2);?(.8)|{h_2^1}="X2"&
          \ar[d]^{h_1^2}\ar@2[dl]**{}?(.2);?(.8)|{h_2^2}="Y2"\\
          \ar[d]_{k_0^1}\ar[r]|{f''^1}&\ar[d]|{k_1^1}\ar[r]|{f''^2}\ar@2[dl]**{}?(.2);?(.8)|{k_2^1}="X3"&
          \ar[d]^{k_1^2}\ar@2[dl]**{}?(.2);?(.8)|{k_2^2}="Y3"\\
          \ar[r]_{f'''^1}&\ar[r]_{f'''^2}&{}          
          \POS
          \tria"Y2";"X1"
          \tria"X2";"Y3"
        }
        ,"A31"
        \xyboxmatrix"A31"{
          \ar[d]_{g_0^1}\ar[r]^{f^1}&\ar[d]|{g_1^1}\ar[r]^{f^2}\ar@2[dl]**{}?(.2);?(.8)|{g_2^1}="X1"&
          \ar[d]^{g_1^2}\ar@2[dl]**{}?(.2);?(.8)|{g_2^2}="Y1"\\
          \ar[d]_{h_0^1}\ar[r]|{f'^1}&\ar[d]|{h_1^1}\ar[r]|{f'^2}\ar@2[dl]**{}?(.2);?(.8)|{h_2^1}="X2"&
          \ar[d]^{h_1^2}\ar@2[dl]**{}?(.2);?(.8)|{h_2^2}="Y2"\\
          \ar[d]_{k_0^1}\ar[r]|{f''^1}&\ar[d]|{k_1^1}\ar[r]|{f''^2}\ar@2[dl]**{}?(.2);?(.8)|{k_2^1}="X3"&
          \ar[d]^{k_1^2}\ar@2[dl]**{}?(.2);?(.8)|{k_2^2}="Y3"\\
          \ar[r]_{f'''^1}&\ar[r]_{f'''^2}&{}          
          \POS
          \tria"Y2";"X1"
          \tria"X1";"Y3"
          \tria"Y3";"X2"
        }
        ,"A41"
        \xyboxmatrix"A41"{
          \ar[d]_{g_0^1}\ar[r]^{f^1}&\ar[d]|{g_1^1}\ar[r]^{f^2}\ar@2[dl]**{}?(.2);?(.8)|{g_2^1}="X1"&
          \ar[d]^{g_1^2}\ar@2[dl]**{}?(.2);?(.8)|{g_2^2}="Y1"\\
          \ar[d]_{h_0^1}\ar[r]|{f'^1}&\ar[d]|{h_1^1}\ar[r]|{f'^2}\ar@2[dl]**{}?(.2);?(.8)|{h_2^1}="X2"&
          \ar[d]^{h_1^2}\ar@2[dl]**{}?(.2);?(.8)|{h_2^2}="Y2"\\
          \ar[d]_{k_0^1}\ar[r]|{f''^1}&\ar[d]|{k_1^1}\ar[r]|{f''^2}\ar@2[dl]**{}?(.2);?(.8)|{k_2^1}="X3"&
          \ar[d]^{k_1^2}\ar@2[dl]**{}?(.2);?(.8)|{k_2^2}="Y3"\\
          \ar[r]_{f'''^1}&\ar[r]_{f'''^2}&{}          
          \POS
          \tria"Y3";"X1"
        }
        ,"A12"
        \xyboxmatrix"A12"{
          \ar[d]_{g_0^1}\ar[r]^{f^1}&\ar[d]|{g_1^1}\ar[r]^{f^2}\ar@2[dl]**{}?(.2);?(.8)|{g_2^1}="X1"&
          \ar[d]^{g_1^2}\ar@2[dl]**{}?(.2);?(.8)|{g_2^2}="Y1"\\
          \ar[d]_{h_0^1}\ar[r]|{f'^1}&\ar[d]|{h_1^1}\ar[r]|{f'^2}\ar@2[dl]**{}?(.2);?(.8)|{h_2^1}="X2"&
          \ar[d]^{h_1^2}\ar@2[dl]**{}?(.2);?(.8)|{h_2^2}="Y2"\\
          \ar[d]_{k_0^1}\ar[r]|{f''^1}&\ar[d]|{k_1^1}\ar[r]|{f''^2}\ar@2[dl]**{}?(.2);?(.8)|{k_2^1}="X3"&
          \ar[d]^{k_1^2}\ar@2[dl]**{}?(.2);?(.8)|{k_2^2}="Y3"\\
          \ar[r]_{f'''^1}&\ar[r]_{f'''^2}&{}          
          \POS
          \tria"X1";"Y2"
          \tria"X2";"Y3"
        }
        ,"A22"
        \xyboxmatrix"A22"{
          \ar[d]_{g_0^1}\ar[r]^{f^1}&\ar[d]|{g_1^1}\ar[r]^{f^2}\ar@2[dl]**{}?(.2);?(.8)|{g_2^1}="X1"&
          \ar[d]^{g_1^2}\ar@2[dl]**{}?(.2);?(.8)|{g_2^2}="Y1"\\
          \ar[d]_{h_0^1}\ar[r]|{f'^1}&\ar[d]|{h_1^1}\ar[r]|{f'^2}\ar@2[dl]**{}?(.2);?(.8)|{h_2^1}="X2"&
          \ar[d]^{h_1^2}\ar@2[dl]**{}?(.2);?(.8)|{h_2^2}="Y2"\\
          \ar[d]_{k_0^1}\ar[r]|{f''^1}&\ar[d]|{k_1^1}\ar[r]|{f''^2}\ar@2[dl]**{}?(.2);?(.8)|{k_2^1}="X3"&
          \ar[d]^{k_1^2}\ar@2[dl]**{}?(.2);?(.8)|{k_2^2}="Y3"\\
          \ar[r]_{f'''^1}&\ar[r]_{f'''^2}&{}          
          \POS
          \tria"X1";"Y2"
          \tria"Y3";"X2"
        }
        ,"A32"
        \xyboxmatrix"A32"{
          \ar[d]_{g_0^1}\ar[r]^{f^1}&\ar[d]|{g_1^1}\ar[r]^{f^2}\ar@2[dl]**{}?(.2);?(.8)|{g_2^1}="X1"&
          \ar[d]^{g_1^2}\ar@2[dl]**{}?(.2);?(.8)|{g_2^2}="Y1"\\
          \ar[d]_{h_0^1}\ar[r]|{f'^1}&\ar[d]|{h_1^1}\ar[r]|{f'^2}\ar@2[dl]**{}?(.2);?(.8)|{h_2^1}="X2"&
          \ar[d]^{h_1^2}\ar@2[dl]**{}?(.2);?(.8)|{h_2^2}="Y2"\\
          \ar[d]_{k_0^1}\ar[r]|{f''^1}&\ar[d]|{k_1^1}\ar[r]|{f''^2}\ar@2[dl]**{}?(.2);?(.8)|{k_2^1}="X3"&
          \ar[d]^{k_1^2}\ar@2[dl]**{}?(.2);?(.8)|{k_2^2}="Y3"\\
          \ar[r]_{f'''^1}&\ar[r]_{f'''^2}&{}          
          \POS
          \tria"Y2";"X1"
          \tria"X1";"Y3"
          \tria"Y3";"X2"
        }
        ,"A42"
        \xyboxmatrix"A42"{
          \ar[d]_{g_0^1}\ar[r]^{f^1}&\ar[d]|{g_1^1}\ar[r]^{f^2}\ar@2[dl]**{}?(.2);?(.8)|{g_2^1}="X1"&
          \ar[d]^{g_1^2}\ar@2[dl]**{}?(.2);?(.8)|{g_2^2}="Y1"\\
          \ar[d]_{h_0^1}\ar[r]|{f'^1}&\ar[d]|{h_1^1}\ar[r]|{f'^2}\ar@2[dl]**{}?(.2);?(.8)|{h_2^1}="X2"&
          \ar[d]^{h_1^2}\ar@2[dl]**{}?(.2);?(.8)|{h_2^2}="Y2"\\
          \ar[d]_{k_0^1}\ar[r]|{f''^1}&\ar[d]|{k_1^1}\ar[r]|{f''^2}\ar@2[dl]**{}?(.2);?(.8)|{k_2^1}="X3"&
          \ar[d]^{k_1^2}\ar@2[dl]**{}?(.2);?(.8)|{k_2^2}="Y3"\\
          \ar[r]_{f'''^1}&\ar[r]_{f'''^2}&{}          
          \POS
          \tria"Y3";"X1"
        }
        \POS 
        \ar@3"A11";"A21"^{m(k)\square_0m^2_{h,g}}_{=(f'''^2\#_0k_2^1\#_0h_0^1\#_0g_0^1)\\
          \#_1(k_2^2\lhc{}h_2^1\#_0g_0^1)\\
          \#_1(k_1^2\#_0h_2^2\ten{}g_2^1)\\
          \#_1(k_1^2\#_0h_1^2\#_0g_2^2\#_0f^1)}
        \ar@3@/^7pc/"A21";"A41"^{m^2_{k,h\square_0g}}
        \ar@3"A21";"A31"_{(f'''^2\#_0k_2^1\#_0h_0^1\#_0g_0^1)\\\#_1(k_2^2\ten{}h_2^1\#_0g_0^1)\\
          \#_1(k_1^2\#_0h_2^2\rhc{}g_2^1)\\\#_1(k_1^2\#_0h_1^2\#_0g_2^2\#_0f^1)}
        \ar@3"A31";"A41"_{\big(f'''^2\#_0((k_2^1\#_0h_0^1)\\
          \#_1(k_1^1\#_0h_2^1))\#_0g_0^1\big)\\
          \#_1(k_2^2\ten{}h_1^1\#_0g_2^1)\\
          \#_1\big(k_1^2\#_0((h_2^2\#_0g_1^1)\\
          \#_1(h_1^2\#_0g_2^2))\#_0f^1\big)}
        \ar@3"A12";"A22"_{=m^2_{k,h}\square_0m(g)}^{(f'''^2\#_0k_2^1\#_0h_0^1\#_0g_0^1)\\
          \#_1(k_2^2\ten{}h_2^1\#_0g_0^1)\\
          \#_1(k_1^2\#_0h_2^2\lhc{}g_2^1)\\
          \#_1(k_1^2\#_0h_1^2\#_0g_2^2\#_0f^1)}
        \ar@3@/^7pc/"A21";"A41"^{m^2_{k,h\square_0g}}
        \ar@3@/_7pc/"A22";"A42"_{m^2_{k\square_0h,g}}
        \ar@3"A22";"A32"^{(f'''^2\#_0k_2^1\#_0h_0^1\#_0g_0^1)\\\#_1(k_2^2\rhc{}h_2^1\#_0g_0^1)\\
          \#_1(k_1^2\#_0h_2^2\ten{}g_2^1)\\\#_1(k_1^2\#_0h_1^2\#_0g_2^2\#_0f^1)}
        \ar@3"A32";"A42"
        \ar@{=}"A11";"A12"
        \ar@{=}"A31";"A32"
        \ar@{=}"A41";"A42"
      \end{xy}
    \end{equation}
  \end{sidewaysfigure}
  Note that in the last column of (\ref{eq:mqmapprf2coc}) 
  \begin{equation*}
    \begin{pmatrix}
      (f'''^2\#_0\underline{k_2^1}\#_0h_0^1\#_0g_0^1)\\\\
      \#_1(\underline{k_2^2}\rhc{}\underline{h_2^1}\#_0g_0^1)\\
      \#_1(k_1^2\#_0\underline{h_2^2}\rhc{}\underline{g_2^1})\\\\
      \#_1(k_1^2\#_0h_1^2\#_0\underline{g_2^2}\#_0f^1)
    \end{pmatrix}
    =
    \begin{pmatrix}
      (f'''^2\#_0\underline{k_2^1}\#_0h_0^1\#_0g_0^1)\\
      \#_1(f'''^2\#_0k_1^1\#_0\underline{h_2^1}\#_0g_0^1)\\
      \#_1(\underline{k_2^2}\#_0h_1^1\#_0f'^1\#_0g_0^1)\\
      \#_1(k_1^2\#_0f''^2\#_0h_1^1\#_0\underline{g_2^1})\\
      \#_1(k_1^2\#_0\underline{h_2^2}\#_0g_1^1\#_0f^1)\\
      \#_1(k_1^2\#_0h_1^2\#_0\underline{g_2^2}\#_0f^1)
    \end{pmatrix}
    =
    \begin{pmatrix}
      \big(f'''^2\#_0((\underline{k_2^1}\#_0h_0^1)\\\#_1(k_1^1\#_0\underline{h_2^1}))\#_0g_0^1\big)\\
      \#_1(\underline{k_2^2}\lhc{}h_1^1\#_0\underline{g_2^1})\\
      \#_1\big(k_1^2\#_0((\underline{h_2^2}\#_0g_1^1)\\\#_1(h_1^2\#_0\underline{g_2^2}))\#_0f^1\big)
    \end{pmatrix}\,,
  \end{equation*} showing how the multiple horizontal composites of
  squares can be simplified.
  And the left hand rectangle in (\ref{eq:mqmapprf2coc}) commutes by
  local interchange. Also, $m^2$ is normalized by the unitality of the
  tensor in $\H$.

  We check the coherent preservation of whiskers of
  2-cells by 1-cells on the left, that is,
  \begin{equation*}
    m^2_{\tilde{h},g}\square_1(m(\alpha)\square_0m(g))=m(\alpha\square_0g)\square_1m^2_{h,g}
  \end{equation*}
  in (\ref{eq:mqmapprf21whskpres}), where the parts
  commute by the naturality of the tensor and the local
  interchange. The corresponding condition for right whiskers is
  verified similarly. Coherent preservation of whiskers of 3-cells by
  1-cells is checked in the same way using in addition the naturality
  of the horizontal composition of a 3-cell by a 2-cell along a
  0-cell. This proves conditions (\ref{eq:coccoh}) and
  (\ref{eq:coccoh3whsk}).  
  \begin{sidewaysfigure}
    \begin{equation}
      \label{eq:mqmapprf21whskpres}
      \begin{xy}
        \save
        (48,0):(0,-1)::
        ,(0,0)="A11"
        ,(1,0)="A21"
        ,(2,0)="A31"
        ,(3,0)="A41"
        ,(4,0)="A51"
        ,(0,1)="A12"
        ,(1,1)="A22"
        ,(2,1)="A32"
        ,(3,1)="A42"
        ,(4,1)="A52"
        \restore
        \POS 
        ,"A11"
        \xyboxmatrix"A11"{
          \ar[d]_{g_0^1}\ar[r]^{f^1}&\ar[d]|{g_1^1}\ar[r]^{f^2}\ar@2[dl]**{}?(.2);?(.8)|{g_2^1}="X1"&
          \ar[d]^{g_1^2}\ar@2[dl]**{}?(.2);?(.8)|{g_2^2}="Y1"\\
          \ar[d]|{h_0^1}="W1"\ar[r]|{f'^1}\ar@/_2pc/[d]_{\tilde{h}_0^1}="W2"&
          \ar[d]|{h_1^1}\ar[r]|{f'^2}\ar@2[dl]**{}?(.2);?(.8)|{h_2^1}="X2"&
          \ar[d]^{h_1^2}\ar@2[dl]**{}?(.2);?(.8)|{h_2^2}="Y2"\\
          \ar[r]_{f''^1}&\ar[r]_{f''^2}&{}          
          \POS
          \tria"X1";"Y2"
          \ar@2"W1";"W2"**{}?(.2);?(.8)|{\alpha_1^1}
        },"A21"
        \xyboxmatrix"A21"@W+.6cm{
          \ar[d]_{g_0^1}\ar[r]^{f^1}&\ar[d]|{g_1^1}\ar[r]^{f^2}\ar@2[dl]**{}?(.2);?(.8)|{g_2^1}="X1"&
          \ar[d]^{g_1^2}\ar@2[dl]**{}?(.2);?(.8)|{g_2^2}="Y1"\\
          \ar[r]|{f'^1}\ar[d]_{\tilde{h}_0^1}="W2"&
          \ar@/_1pc/[d]|{\tilde{h}_1^1}="Z2"\ar@/^1pc/[d]|{h_1^1}="Z1"\ar[r]|{f'^2}\ar@2[dl]**{}?(.2);?(.8)|{\tilde{h}_2^1}="X2"&
          \ar[d]^{h_1^2}\ar@2[dl]**{}?(.2);?(.8)|{h_2^2}="Y2"\\
          \ar[r]_{f''^1}&\ar[r]_{f''^2}&{}          
          \POS
          \tria"X1";"Y2"
          \ar@2"Z1";"Z2"**{}?(.2);?(.8)|{\alpha_2^1}
        },"A31"
        \xyboxmatrix"A31"{
          \ar[d]_{g_0^1}\ar[r]^{f^1}&\ar[d]|{g_1^1}\ar[r]^{f^2}\ar@2[dl]**{}?(.2);?(.8)|{g_2^1}="X1"&
          \ar[d]^{g_1^2}\ar@2[dl]**{}?(.2);?(.8)|{g_2^2}="Y1"\\
          \ar[d]_{\tilde{h}_0^1}\ar[r]|{f'^1}&
          \ar[d]|{\tilde{h}_1^1}\ar[r]|{f'^2}\ar@2[dl]**{}?(.2);?(.8)|{\tilde{h}_2^1}="X2"&
          \ar[d]|{\tilde{h}_1^2}="W2"\ar@/^2pc/[d]^{h_1^2}="W1"\ar@2[dl]**{}?(.2);?(.8)|{\tilde{h}_2^2}\\
          \ar[r]_{f''^1}&\ar[r]_{f''^2}&{}          
          \POS
          \ar@2"W1";"W2"**{}?(.2);?(.8)|{\alpha_2^2}="Y2"
          \tria"X1";"Y2"
        },"A51"
        \xyboxmatrix"A51"{
          \ar[d]_{g_0^1}\ar[r]^{f^1}&\ar[d]|{g_1^1}\ar[r]^{f^2}\ar@2[dl]**{}?(.2);?(.8)|{g_2^1}="X1"&
          \ar[d]^{g_1^2}\ar@2[dl]**{}?(.2);?(.8)|{g_2^2}="Y1"\\
          \ar[d]_{\tilde{h}_0^1}\ar[r]|{f'^1}&
          \ar[d]|{\tilde{h}_1^1}\ar[r]|{f'^2}\ar@2[dl]**{}?(.2);?(.8)|{\tilde{h}_2^1}="X2"&
          \ar[d]|{\tilde{h}_1^2}="W2"\ar@/^2pc/[d]^{h_1^2}="W1"\ar@2[dl]**{}?(.2);?(.8)|{\tilde{h}_2^2}="Y2"\\
          \ar[r]_{f''^1}&\ar[r]_{f''^2}&{}          
          \POS
          \tria"X1";"Y2"
          \ar@2"W1";"W2"**{}?(.2);?(.8)|{\alpha_2^2}="Z1"
          \tria"Z1";"Y1"
        },"A41"
        \xyboxmatrix"A41"{
          \ar[d]_{g_0^1}\ar[r]^{f^1}&\ar[d]|{g_1^1}\ar[r]^{f^2}\ar@2[dl]**{}?(.2);?(.8)|{g_2^1}="X1"&
          \ar[d]^{g_1^2}\ar@2[dl]**{}?(.2);?(.8)|{g_2^2}="Y1"\\
          \ar[d]_{\tilde{h}_0^1}\ar[r]|{f'^1}&
          \ar[d]|{\tilde{h}_1^1}\ar[r]|{f'^2}\ar@2[dl]**{}?(.2);?(.8)|{\tilde{h}_2^1}="X2"&
          \ar[d]|{\tilde{h}_1^2}="W2"\ar@/^2pc/[d]^{h_1^2}="W1"\ar@2[dl]**{}?(.2);?(.8)|{\tilde{h}_2^2}="Y2"\\
          \ar[r]_{f''^1}&\ar[r]_{f''^2}&{}          
          \POS
          \tria"X1";"Y2"
          \ar@2"W1";"W2"**{}?(.2);?(.8)|{\alpha_2^2}="Z1"
          \tria"Y1";"Z1"
        },"A12"
        \xyboxmatrix"A12"{
          \ar[d]_{g_0^1}\ar[r]^{f^1}&\ar[d]|{g_1^1}\ar[r]^{f^2}\ar@2[dl]**{}?(.2);?(.8)|{g_2^1}="X1"&
          \ar[d]^{g_1^2}\ar@2[dl]**{}?(.2);?(.8)|{g_2^2}="Y1"\\
          \ar[d]|{h_0^1}="W1"\ar[r]|{f'^1}\ar@/_2pc/[d]_{\tilde{h}_0^1}="W2"&
          \ar[d]|{h_1^1}\ar[r]|{f'^2}\ar@2[dl]**{}?(.2);?(.8)|{h_2^1}="X2"&
          \ar[d]^{h_1^2}\ar@2[dl]**{}?(.2);?(.8)|{h_2^2}="Y2"\\
          \ar[r]_{f''^1}&\ar[r]_{f''^2}&{}          
          \POS
          \tria"Y2";"X1"
          \ar@2"W1";"W2"**{}?(.2);?(.8)|{\alpha_1^1}
        },"A22"
        \xyboxmatrix"A22"@W+.6cm{
          \ar[d]_{g_0^1}\ar[r]^{f^1}&\ar[d]|{g_1^1}\ar[r]^{f^2}\ar@2[dl]**{}?(.2);?(.8)|{g_2^1}="X1"&
          \ar[d]^{g_1^2}\ar@2[dl]**{}?(.2);?(.8)|{g_2^2}="Y1"\\
          \ar[r]|{f'^1}\ar[d]_{\tilde{h}_0^1}="W2"&
          \ar@/_1pc/[d]|{\tilde{h}_1^1}="Z2"\ar@/^1pc/[d]|{h_1^1}="Z1"\ar[r]|{f'^2}\ar@2[dl]**{}?(.2);?(.8)|{\tilde{h}_2^1}="X2"&
          \ar[d]^{h_1^2}\ar@2[dl]**{}?(.2);?(.8)|{h_2^2}="Y2"\\
          \ar[r]_{f''^1}&\ar[r]_{f''^2}&{}          
          \POS
          \ar@2"Z1";"Z2"**{}?(.2);?(.8)|{\alpha_2^1}="R"
          \tria"X1";"R"
          \tria"Y2";"X1"
        },"A32"
        \xyboxmatrix"A32"@W+.6cm{
          \ar[d]_{g_0^1}\ar[r]^{f^1}&\ar[d]|{g_1^1}\ar[r]^{f^2}\ar@2[dl]**{}?(.2);?(.8)|{g_2^1}="X1"&
          \ar[d]^{g_1^2}\ar@2[dl]**{}?(.2);?(.8)|{g_2^2}="Y1"\\
          \ar[r]|{f'^1}\ar[d]_{\tilde{h}_0^1}="W2"&
          \ar@/_1pc/[d]|{\tilde{h}_1^1}="Z2"\ar@/^1pc/[d]|{h_1^1}="Z1"\ar[r]|{f'^2}\ar@2[dl]**{}?(.2);?(.8)|{\tilde{h}_2^1}="X2"&
          \ar[d]^{h_1^2}\ar@2[dl]**{}?(.2);?(.8)|{h_2^2}="Y2"\\
          \ar[r]_{f''^1}&\ar[r]_{f''^2}&{}          
          \POS
          \ar@2"Z1";"Z2"**{}?(.2);?(.8)|{\alpha_2^1}="R"
          \tria"R";"X1"
        },"A42"
        \xyboxmatrix"A42"{
          \ar[d]_{g_0^1}\ar[r]^{f^1}&\ar[d]|{g_1^1}\ar[r]^{f^2}\ar@2[dl]**{}?(.2);?(.8)|{g_2^1}="X1"&
          \ar[d]^{g_1^2}\ar@2[dl]**{}?(.2);?(.8)|{g_2^2}="Y1"\\
          \ar[d]_{\tilde{h}_0^1}\ar[r]|{f'^1}&
          \ar[d]|{\tilde{h}_1^1}\ar[r]|{f'^2}\ar@2[dl]**{}?(.2);?(.8)|{\tilde{h}_2^1}="X2"&
          \ar[d]|{\tilde{h}_1^2}="W2"\ar@/^2pc/[d]^{h_1^2}="W1"\ar@2[dl]**{}?(.2);?(.8)|{\tilde{h}_2^2}="Y2"\\
          \ar[r]_{f''^1}&\ar[r]_{f''^2}&{}          
          \POS
          \tria"Y2";"X1"
          \ar@2"W1";"W2"**{}?(.2);?(.8)|{\alpha_2^2}="Z1"
          \tria"Y1";"Z1"
        },"A52"
        \xyboxmatrix"A52"{
          \ar[d]_{g_0^1}\ar[r]^{f^1}&\ar[d]|{g_1^1}\ar[r]^{f^2}\ar@2[dl]**{}?(.2);?(.8)|{g_2^1}="X1"&
          \ar[d]^{g_1^2}\ar@2[dl]**{}?(.2);?(.8)|{g_2^2}="Y1"\\
          \ar[d]_{\tilde{h}_0^1}\ar[r]|{f'^1}&
          \ar[d]|{\tilde{h}_1^1}\ar[r]|{f'^2}\ar@2[dl]**{}?(.2);?(.8)|{\tilde{h}_2^1}="X2"&
          \ar[d]|{\tilde{h}_1^2}="W2"\ar@/^2pc/[d]^{h_1^2}="W1"\ar@2[dl]**{}?(.2);?(.8)|{\tilde{h}_2^2}="Y2"\\
          \ar[r]_{f''^1}&\ar[r]_{f''^2}&{}          
          \POS
          \tria"Y2";"X1"
          \ar@2"W1";"W2"**{}?(.2);?(.8)|{\alpha_2^2}="Z1"
          \tria"Z1";"Y1"
        }
        \ar@3@/^4.5pc/"A11";"A51"^{m(\alpha)\square_0m(g)}
        \ar@3@/_4.5pc/"A12";"A52"_{m(\alpha\square_0g)}
        \ar@3"A11";"A12"_{m^2_{h,g}}^{=(f''^2\#_0f''^1\#_0\alpha_1^1\#_0g_0^1)\\\#_1(f''^2\#_0h_2^1\#_0g_0^1)\\\#_1(h_2^2\ten{}g_2^1)\\\#_1(h_1^2\#_0g_2^2\#_0f^1)}
        \ar@3"A51";"A52"_{m^2_{\tilde{h},g}}^{=(f''^2\#_0\tilde{h}_2^1\#_0g_0^1)\\\#_1(\tilde{h}_2^2\ten{}g_2^1)\\\#_1(\alpha_2^2\rhc{}g_2^2\#_0f^1)}
        \ar@3"A11";"A21"|{}="x"\save\POS"x"+(0,30)*!C\labelbox{(f''^2\#_0\underline{\alpha_3^1}\#_0g_0^1)\\\#_1(h_2^2\lhc{}g_2^1)\\\#_1(h_1^2\#_0g_2^2\#_0f^1)}\ar@{.}"x";c\restore
        \ar@3"A21";"A31"|{}="x"\save\POS"x"+(0,30)*!C\labelbox{(f''^2\#_0\tilde{h}_2^1\#_0g_0^1)\\\#_1(\underline{\alpha_3^2}\#_0f'^1\#_0g_0^1)\\\#_1(h_1^2\#_0f'^2\#_0g_2^1)\\\#_1(h_2^1\#_0g_2^2\#_0f^1)}\ar@{.}"x";c\restore
        \ar@3"A31";"A41"|{}="x"\save\POS"x"+(0,30)*!C\labelbox{(f''^2\#_0\tilde{h}_2^1\#_0g_0^1)\\\#_1(\tilde{h}_2^2\#_0f'^1\#_0g_0^1)\\\#_1(\alpha_2^2\#_0f'^2\ten{}g_2^1)\\\#_1(h_2^1\#_0g_2^2\#_0f^1)}\ar@{.}"x";c\restore
        \ar@3"A41";"A51"|{}="x"\save\POS"x"+(0,30)*!C\labelbox{(f''^2\#_0\tilde{h}_2^1\#_0g_0^1)\\\#_1(\tilde{h}_2^2\lhc{}g_2^1)\\\#_1(\alpha_2^2\ten{}g_2^2\#_0f^1)}\ar@{.}"x";c\restore
        \ar@3"A12";"A22"|{}="x"\save\POS"x"+(0,-30)*!C\labelbox{(f''^2\#_0\underline{\alpha_3^1}\#_0g_0^1)\\\#_1(h_2^2\rhc{}g_2^1)\\\#_1(h_1^2\#_0g_2^2\#_0f^1)}\ar@{.}"x";c\restore
        \ar@3"A22";"A32"|{}="x"\save\POS"x"+(0,-30)*!C\labelbox{(f''^2\#_0\tilde{h}_2^1\#_0g_0^1)\\\#_1(f''^2\#_0\alpha_2^1\ten{}g_2^1)\\\#_1(h_2^2\#_0g_1^1\#_0f^1)\\\#_1(h_1^2\#_0g_2^2\#_0f^1)}\ar@{.}"x";c\restore
        \ar@3"A32";"A42"|{}="x"\save\POS"x"+(0,-30)*!C\labelbox{(f''^2\#_0\tilde{h}_2^1\#_0g_0^1)\\\#_1(f''^2\#_0\tilde{h}_1^1\#_0g_2^1)\\\#_1(\underline{\alpha_3^2}\#_0g_1^1\#_0f^1)\\\#_1(h_1^2\#_0g_2^2\#_0f^1)}\ar@{.}"x";c\restore
        \ar@3"A42";"A52"|{}="x"\save\POS"x"+(0,-30)*!C\labelbox{(f''^2\#_0\tilde{h}_2^1\#_0g_0^1)\\\#_1(\tilde{h}_2^2\rhc{}g_2^1)\\\#_1(\alpha_2^2\ten{}g_2^2\#_0f^1)}\ar@{.}"x";c\restore
        \ar@3"A41";"A42"|{(f''^2\#_0\tilde{h}_2^1\#_0g_0^1)\#_1(\tilde{h}_2^2\ten{}g_2^1)\\\#_1(\alpha_2^2\rhc{}g_2^1\#_0f^1)}
        \ar@3"A21";"A22"|{}="x"\save\POS"x"+(-60,-47)*!C\labelbox{(f''^2\#_0\tilde{h}_2^1\#_0g_0^1)\\\#_1(f''^2\#_0\alpha_2^1\#_0f'^1\#_0g_0^1)\\\#_1(h_2^2\ten{}g_2^1)\\\#_1(h_1^2\#_0g_2^2\#_0f^1)}\ar@{.}"x";c\restore
        \ar@3"A21";"A32"|{}="x"\save\POS"x"+(-31,-75)*!C\labelbox{(f''^2\#_0\tilde{h}_2^2\#_0g_0^1)\\\#_1(((f''^2\#_0\alpha_2^1)\#_1h_2^2)\ten{g_2^1})\\\#_1(h_1^2\#_0g_2^2\#_0f^1)}\ar@{.}"x";c\restore
        \ar@3"A31";"A42"|{}="x"\save\POS"x"+(-31,-75)*!C\labelbox{(f''^2\#_0\tilde{h}_2^2\#_0g_0^1)\\\#_1((\tilde{h}_2^2\#_1(\alpha_2^2\#_0f'^2))\ten{g_2^1})\\\#_1(h_1^2\#_0g_2^2\#_0f^1)}\ar@{.}"x";c\restore
      \end{xy}
    \end{equation}
  \end{sidewaysfigure}
  \begin{sidewaysfigure}
    \begin{equation}
      \label{eq:mqmapprf21whskprescubetopsub}
      \begin{xy}
        \save
        (50,0):(0,-1)::
        ,(0,0)="A11"
        ,(1,0)="A21"
        ,(2,0)="A31"
        ,(3,0)="A41"
        ,(4,0)="A51"
        ,(0,1)="A12"
        ,(1,1)="A22"
        ,(2,1)="A32"
        ,(3,1)="A42"
        ,(4,1)="A52"
        ,(0,2)="A13"
        ,(1,2)="A23"
        ,(2,2)="A33"
        ,(3,2)="A43"
        ,(4,2)="A53"
        ,(0,3)="A14"
        ,(1,3)="A24"
        ,(2,3)="A34"
        ,(3,3)="A44"
        ,(4,3)="A54"
        \restore
        \POS,"A11"
        \xyboxmatrix"A11"{
          \ar[r]^{f^1}\ar[d]_{\tilde{h}_0^1}&
          \ar@/^1pc/[d]|{h_1^1}="Z1"\ar@/_1pc/[d]|{\tilde{h}_1^0}="Z2"
          \ar[r]^{f^2}\ar@2 {[l]**{}?(.5)};[dl]**{}?(.2);?(.8)|{\tilde{h}_2^1}="X1"&
          \ar[d]^{h_1^2}\ar@2 {[d];[dl]**{}?(.5)}**{}?(.2);?(.8)|{h_2^2}="Y1"\\
          \ar[d]|{k_0^1}="W1"\ar[r]|{f'^1}\ar@/_2pc/[d]_{\tilde{k}_0^1}="W2"&
          \ar[d]|{k_1^1}\ar[r]|{f'^2}\ar@2[dl]**{}?(.2);?(.8)|{k_2^1}="X2"&
          \ar[d]^{k_1^2}\ar@2[dl]**{}?(.2);?(.8)|{k_2^2}="Y2"\\
          \ar[r]_{f''^1}&\ar[r]_{f''^2}&{}          
          \POS
          \ar@2"W1";"W2"**{}?(.2);?(.8)|{\beta_1^1}
          \ar@2"Z1";"Z2"**{}?(.2);?(.8)|{\alpha_2^1}="K"
          \tria"X1";"Y2"
        },"A21"
        \xyboxmatrix"A21"{
          \ar[d]_{\tilde{h}_0^1}\ar[r]^{f^1}&
          \ar[d]|{\tilde{h}_1^1}\ar[r]^{f^2}\ar@2[dl]**{}?(.2);?(.8)|{\tilde{h}_2^1}="X1"&
          \ar[d]|{\tilde{h}_1^2}="Z2"\ar@2[dl]**{}?(.2);?(.8)|{\tilde{h}_2^2}="Y1"\ar@/^2pc/[d]^{{h}_1^2}="Z1"\\
          \ar[d]|{k_0^1}="W1"\ar[r]|{f'^1}\ar@/_2pc/[d]_{\tilde{k}_0^1}="W2"&
          \ar[d]|{k_1^1}\ar[r]|{f'^2}\ar@2[dl]**{}?(.2);?(.8)|{k_2^1}="X2"&
          \ar[d]^{k_1^2}\ar@2[dl]**{}?(.2);?(.8)|{k_2^2}="Y2"\\
          \ar[r]_{f''^1}&\ar[r]_{f''^2}&{}          
          \POS
          \ar@2"W1";"W2"**{}?(.2);?(.8)|{\beta_1^1}
          \ar@2"Z1";"Z2"**{}?(.2);?(.8)|{\alpha_2^2}
          \tria"X1";"Y2"
        }
        ,"A22"
        \xyboxmatrix"A22"{
          \ar[d]_{\tilde{h}_0^1}\ar[r]^{f^1}&
          \ar[d]|{\tilde{h}_1^1}\ar[r]^{f^2}\ar@2[dl]**{}?(.2);?(.8)|{\tilde{h}_2^1}="X1"&
          \ar[d]|{\tilde{h}_1^2}="Z2"\ar@2[dl]**{}?(.2);?(.8)|{\tilde{h}_2^2}="Y1"\ar@/^2pc/[d]^{{h}_1^2}="Z1"\\
          \ar[d]|{k_0^1}="W1"\ar[r]|{f'^1}\ar@/_2pc/[d]_{\tilde{k}_0^1}="W2"&
          \ar[d]|{k_1^1}\ar[r]|{f'^2}\ar@2[dl]**{}?(.2);?(.8)|{k_2^1}="X2"&
          \ar[d]^{k_1^2}\ar@2[dl]**{}?(.2);?(.8)|{k_2^2}="Y2"\\
          \ar[r]_{f''^1}&\ar[r]_{f''^2}&{}          
          \POS
          \ar@2"W1";"W2"**{}?(.2);?(.8)|{\beta_1^1}
          \ar@2"Z1";"Z2"**{}?(.2);?(.8)|{\alpha_2^2}
          \tria"Y2";"X1"
        }
        ,"A31"
        \xyboxmatrix"A31"{
          \ar[d]_{\tilde{h}_0^1}\ar[r]^{f^1}&
          \ar[d]|{\tilde{h}_1^1}\ar[r]^{f^2}\ar@2[dl]**{}?(.2);?(.8)|{\tilde{h}_2^1}="X1"&
          \ar[d]|{\tilde{h}_1^2}="Z2"\ar@2[dl]**{}?(.2);?(.8)|{\tilde{h}_2^2}="Y1"\ar@/^2pc/[d]^{{h}_1^2}="Z1"\\
          \ar[r]|{f'^1}\ar[d]_{\tilde{k}_1^1}&\ar@/^1pc/[d]|{k_0^2}="W1"\ar@/_1pc/[d]|{\tilde{k}_1^1}="W2"
          \ar[r]|{f'^2}\ar@2 {[l]**{}?(.5)};[dl]**{}?(.2);?(.8)|{\tilde{k}_2^1}="X2"&
          \ar[d]^{k_1^2}\ar@2 {[d];[dl]**{}?(.5)}**{}?(.2);?(.8)|{k_2^2}="Y2"\\
          \ar[r]_{f''^1}&\ar[r]_{f''^2}&{}          
          \POS
          \ar@2"W1";"W2"**{}?(.2);?(.8)|{\beta_2^1}
          \ar@2"Z1";"Z2"**{}?(.2);?(.8)|{\alpha_2^2}
          \tria"X1";"Y2"
        }
        ,"A24"
        \xyboxmatrix"A24"{
          \ar[d]_{\tilde{h}_0^1}\ar[rr]^{f^1}&&\ar[dl]|{\tilde{h}_1^1}\ar[dr]|{h_1^1}\ar[rr]^{f^2}
          \ar@2[dll]**{}?(.2);?(.8)|{\tilde{h}_2^1}="X1"&&
          \ar[d]^{h_1^2}\ar@2[dl]**{}?(.2);?(.8)|{h_2^2}="Y1"\\
          \ar[r]|{f'^1}\ar[d]_{\tilde{k}_0^1}="W2"&\ar[dr]|{\tilde{k}_1^1}="Z2"\ar@2[dl]**{}?(.2);?(.8)|{\tilde{k}_2^1}="X2"
          &&\ar@2[ll]**{}?(.2);?(.8)|{\beta_2^1\lhc\alpha_2^1}="R"
          \ar[r]|{f'^2}\ar[dl]|{k_1^1}="Z1"&
          \ar[d]^{k_1^2}\ar@2[dll]**{}?(.2);?(.8)|{k_2^2}="Y2"\\
          \ar[rr]_{f''^1}&&\ar[rr]_{f''^2}&&{}          
        }
        ,"A23"
        \xyboxmatrix"A23"{
          \ar[d]_{\tilde{h}_0^1}\ar[rr]^{f^1}&&\ar[dl]|{\tilde{h}_1^1}\ar[dr]|{h_1^1}\ar[rr]^{f^2}
          \ar@2[dll]**{}?(.2);?(.8)|{\tilde{h}_2^1}="X1"&&
          \ar[d]^{h_1^2}\ar@2[dl]**{}?(.2);?(.8)|{h_2^2}="Y1"\\
          \ar[r]|{f'^1}\ar[d]_{\tilde{k}_0^1}="W2"&\ar[dr]|{\tilde{k}_1^1}="Z2"\ar@2[dl]**{}?(.2);?(.8)|{\tilde{k}_2^1}="X2"
          &&\ar@2[ll]**{}?(.2);?(.8)|{\beta_2^1\rhc\alpha_2^1}="R"
          \ar[r]|{f'^2}\ar[dl]|{k_1^1}="Z1"&
          \ar[d]^{k_1^2}\ar@2[dll]**{}?(.2);?(.8)|{k_2^2}="Y2"\\
          \ar[rr]_{f''^1}&&\ar[rr]_{f''^2}&&{}          
        }
        ,"A14"
        \xyboxmatrix"A14"{
          \ar[r]^{f^1}\ar[d]_{\tilde{h}_0^1}&
          \ar@/^1pc/[d]|{h_1^1}="Z1"\ar@/_1pc/[d]|{\tilde{h}_1^0}="Z2"
          \ar[r]^{f^2}\ar@2 {[l]**{}?(.5)};[dl]**{}?(.2);?(.8)|{\tilde{h}_2^1}="X1"&
          \ar[d]^{h_1^2}\ar@2 {[d];[dl]**{}?(.5)}**{}?(.2);?(.8)|{h_2^2}="Y1"\\
          \ar[d]|{\tilde{k}_0^1}="W1"\ar[r]|{f'^1}&
          \ar@/^1pc/[d]|{k_1^1}="W1"\ar@/_1pc/[d]|{\tilde{k}_1^0}="W2"
          \ar@2 {[l]**{}?(.5)};[dl]**{}?(.2);?(.8)|{\tilde{k}_2^1}="X2"\ar[r]|{f'^2}&
          \ar[d]^{k_1^2}\ar@2 {[d];[dl]**{}?(.5)}**{}?(.2);?(.8)|{k_2^2}="Y2"\\
          \ar[r]_{f''^1}&\ar[r]_{f''^2}&{}          
          \POS
          \ar@2"W1";"W2"**{}?(.2);?(.8)|{\beta_2^1}="K2"
          \ar@2"Z1";"Z2"**{}?(.2);?(.8)|{\alpha_2^1}="K1"
          \tria"X1";"K2"
          \tria"Y2";"K1"
        }
        ,"A12"
        \xyboxmatrix"A12"{
          \ar[r]^{f^1}\ar[d]_{\tilde{h}_0^1}&
          \ar@/^1pc/[d]|{h_1^1}="Z1"\ar@/_1pc/[d]|{\tilde{h}_1^0}="Z2"
          \ar[r]^{f^2}\ar@2 {[l]**{}?(.5)};[dl]**{}?(.2);?(.8)|{\tilde{h}_2^1}="X1"&
          \ar[d]^{h_1^2}\ar@2 {[d];[dl]**{}?(.5)}**{}?(.2);?(.8)|{h_2^2}="Y1"\\
          \ar[d]|{k_0^1}="W1"\ar[r]|{f'^1}\ar@/_2pc/[d]_{\tilde{k}_0^1}="W2"&
          \ar[d]|{k_1^1}\ar[r]|{f'^2}\ar@2[dl]**{}?(.2);?(.8)|{k_2^1}="X2"&
          \ar[d]^{k_1^2}\ar@2[dl]**{}?(.2);?(.8)|{k_2^2}="Y2"\\
          \ar[r]_{f''^1}&\ar[r]_{f''^2}&{}          
          \POS
          \ar@2"W1";"W2"**{}?(.2);?(.8)|{\beta_1^1}
          \ar@2"Z1";"Z2"**{}?(.2);?(.8)|{\alpha_2^1}="K"
          \tria"K";"Y2"
          \tria"Y2";"X1"
        },"A34"
        \xyboxmatrix"A34"{
          \ar[r]^{f^1}\ar[d]_{\tilde{h}_0^1}&
          \ar@/^1pc/[d]|{h_1^1}="Z1"\ar@/_1pc/[d]|{\tilde{h}_1^0}="Z2"
          \ar[r]^{f^2}\ar@2 {[l]**{}?(.5)};[dl]**{}?(.2);?(.8)|{\tilde{h}_2^1}="X1"&
          \ar[d]^{h_1^2}\ar@2 {[d];[dl]**{}?(.5)}**{}?(.2);?(.8)|{h_2^2}="Y1"\\
          \ar[d]|{\tilde{k}_0^1}="W1"\ar[r]|{f'^1}&
          \ar@/^1pc/[d]|{k_1^1}="W1"\ar@/_1pc/[d]|{\tilde{k}_1^0}="W2"
          \ar@2 {[l]**{}?(.5)};[dl]**{}?(.2);?(.8)|{\tilde{k}_2^1}="X2"\ar[r]|{f'^2}&
          \ar[d]^{k_1^2}\ar@2 {[d];[dl]**{}?(.5)}**{}?(.2);?(.8)|{k_2^2}="Y2"\\
          \ar[r]_{f''^1}&\ar[r]_{f''^2}&{}          
          \POS
          \ar@2"W1";"W2"**{}?(.2);?(.8)|{\beta_2^1}="K2"
          \ar@2"Z1";"Z2"**{}?(.2);?(.8)|{\alpha_2^1}="K1"
          \tria"K2";"X1"
          \tria"K1";"Y2"
        },"A32"
        \xyboxmatrix"A32"{
          \ar[d]_{\tilde{h}_0^1}\ar[r]^{f^1}&
          \ar[d]|{\tilde{h}_1^1}\ar[r]^{f^2}\ar@2[dl]**{}?(.2);?(.8)|{\tilde{h}_2^1}="X1"&
          \ar[d]|{\tilde{h}_1^2}="Z2"\ar@2[dl]**{}?(.2);?(.8)|{\tilde{h}_2^2}="Y1"\ar@/^2pc/[d]^{{h}_1^2}="Z1"\\
          \ar[r]|{f'^1}\ar[d]_{\tilde{k}_1^1}&\ar@/^1pc/[d]|{k_0^2}="W1"\ar@/_1pc/[d]|{\tilde{k}_1^1}="W2"
          \ar[r]|{f'^2}\ar@2 {[l]**{}?(.5)};[dl]**{}?(.2);?(.8)|{\tilde{k}_2^1}="X2"&
          \ar[d]^{k_1^2}\ar@2 {[d];[dl]**{}?(.5)}**{}?(.2);?(.8)|{k_2^2}="Y2"\\
          \ar[r]_{f''^1}&\ar[r]_{f''^2}&{}          
          \POS
          \ar@2"W1";"W2"**{}?(.2);?(.8)|{\beta_2^1}="K"
          \ar@2"Z1";"Z2"**{}?(.2);?(.8)|{\alpha_2^2}
          \tria"X1";"K"
          \tria"Y2";"X1"
        }
        \ar@3"A11";"A12"
        \ar@3"A11";"A21"
        \ar@3"A12";"A14"
        \ar@3"A12";"A22"
        \ar@3"A12";"A23"
        \ar@3"A14";"A24"
        \ar@3"A21";"A22"
        \ar@3"A21";"A31"
        \ar@3"A22";"A32"
        \ar@3"A23";"A14"
        \ar@3"A23";"A34"
        \ar@3"A24";"A34"
        \ar@3"A31";"A32"
        \ar@3"A32";"A23"
        \ar@3"A34";"A32"
      \end{xy}
    \end{equation}
  \end{sidewaysfigure}

  We verify the coherent preservation of tensors, i.~e.\@ 
  that
  \begin{equation}
    \label{eq:mqmapprftenpres1}
    m(\beta\bten\alpha)\square_1m^2_{k,h}=m^2_{\tilde{k},\tilde{h}}\square_1(m(\beta)\bten m(\alpha))\,,
  \end{equation}
  where $\alpha, \beta, k,h,\tilde{k},\tilde{h}$ are 2- and 1-cells
  respectively in $\pathspc\H\times_\H\pathspc\H$. In terms of
  constituent cells (\ref{eq:mqmapprftenpres1}) can be drawn as
  (\ref{eq:mqmapprftenpres1cube}), where the pasting of the center and
  right squares corresponds to the right hand side of the
  equation~(\ref{eq:mqmapprftenpres1}), and the pasting of the left
  and outer squares corresponds to the left hand side. Equality in
  (\ref{eq:mqmapprftenpres1}) is equivalent to the top and bottom
  squares commuting, since the aforementioned ones do so by assumption.
  \begin{sidewaysfigure}
    \begin{equation}
      \label{eq:mqmapprftenpres1cube}
      \begin{xy}
        \save
        (70,0):(0,-1)::
        ,(0,0)="A11"
        ,(1,0)="A21"
        ,(0,1)="A12"
        ,(1,1)="A22"
        \restore
        \save
        (-50,50);p+(170,0):(0,-1)::
        ,(0,0)="B11"
        ,(1,0)="B21"
        ,(0,1)="B12"
        ,(1,1)="B22"
        \restore
        \POS 
        ,"A11"
        \xyboxmatrix"A11"{
          \ar[d]|{h_0^1}\ar[r]^{f^1}\ar@/_4pc/[dd]|{\tilde{k}_0^1\#_0\tilde{h}_0^1}="W2"&
          \ar[d]|{h_1^1}\ar[r]^{f^2}\ar@2[dl]**{}?(.2);?(.8)|{h_2^1}="X1"&
          \ar[d]^{h_1^2}\ar@2[dl]**{}?(.2);?(.8)|{h_2^2}="Y1"\\
          {}\POS*+{}="W1"\ar[d]|{k_0^1}\ar[r]|{f'^1}&
          \ar[d]|{k_1^1}\ar[r]|{f'^2}\ar@2[dl]**{}?(.2);?(.8)|{k_2^1}="X2"&
          {}\POS*+{}="Z1"\ar[d]^{k_1^2}\ar@2[dl]**{}?(.2);?(.8)|{k_2^2}="Y2"\\
          \ar[r]_{f''^1}&\ar[r]_{f''^2}&{}          
          \POS
          \tria"X1";"Y2"
          \ar@2"W1";"W2"**{}?<(.2);?>(.8)_{\beta_1^1\lhc\alpha_1^1}
        }
        ,"A21"
        \xyboxmatrix"A21"{
          \ar[d]_{\tilde{h}_0^1}\ar[r]^{f^1}&
          \ar[d]|{\tilde{h}_1^1}\ar[r]^{f^2}\ar@2[dl]**{}?(.2);?(.8)|{\tilde{h}_2^1}="X1"&
          \ar[d]|{\tilde{h}_1^2}\ar@2[dl]**{}?(.2);?(.8)|{\tilde{h}_2^2}="Y1"\ar@/^4pc/[dd]|{k_1^2\#_0h_1^2}="Z2"\\
          {}\POS*+{}="W1"\ar[d]_{\tilde{k}_0^1}\ar[r]|{f'^1}&
          \ar[d]|{\tilde{k}_1^1}\ar[r]|{f'^2}\ar@2[dl]**{}?(.2);?(.8)|{\tilde{k}_2^1}="X2"&
          {}\POS*+{}="Z1"\ar[d]|{\tilde{k}_1^2}\ar@2[dl]**{}?(.2);?(.8)|{\tilde{k}_2^2}="Y2"\\
          \ar[r]_{f''^1}&\ar[r]_{f''^2}&{}          
          \POS
          \tria"X1";"Y2"
          \ar@2"Z2";"Z1"**{}?<(.2);?>(.8)_{\beta_2^2\lhc\alpha_2^2}
        } 
        ,"A12"
        \xyboxmatrix"A12"{
          \ar[d]|{h_0^1}\ar[r]^{f^1}\ar@/_4pc/[dd]|{\tilde{k}_0^1\#_0\tilde{h}_0^1}="W2"&
          \ar[d]|{h_1^1}\ar[r]^{f^2}\ar@2[dl]**{}?(.2);?(.8)|{h_2^1}="X1"&
          \ar[d]^{h_1^2}\ar@2[dl]**{}?(.2);?(.8)|{h_2^2}="Y1"\\
          {}\POS*+{}="W1"\ar[d]|{k_0^1}\ar[r]|{f'^1}&
          \ar[d]|{k_1^1}\ar[r]|{f'^2}\ar@2[dl]**{}?(.2);?(.8)|{k_2^1}="X2"&
          {}\POS*+{}="Z1"\ar[d]^{k_1^2}\ar@2[dl]**{}?(.2);?(.8)|{k_2^2}="Y2"\\
          \ar[r]_{f''^1}&\ar[r]_{f''^2}&{}          
          \POS
          \tria"X1";"Y2"
          \ar@2"W1";"W2"**{}?<(.2);?>(.8)_{\beta_1^1\rhc\alpha_1^1}
        }
        ,"A22"
        \xyboxmatrix"A22"{
          \ar[d]_{\tilde{h}_0^1}\ar[r]^{f^1}&
          \ar[d]|{\tilde{h}_1^1}\ar[r]^{f^2}\ar@2[dl]**{}?(.2);?(.8)|{\tilde{h}_2^1}="X1"&
          \ar[d]|{\tilde{h}_1^2}\ar@2[dl]**{}?(.2);?(.8)|{\tilde{h}_2^2}="Y1"\ar@/^4pc/[dd]|{k_1^2\#_0h_1^2}="Z2"\\
          {}\POS*+{}="W1"\ar[d]_{\tilde{k}_0^1}\ar[r]|{f'^1}&
          \ar[d]|{\tilde{k}_1^1}\ar[r]|{f'^2}\ar@2[dl]**{}?(.2);?(.8)|{\tilde{k}_2^1}="X2"&
          {}\POS*+{}="Z1"\ar[d]|{\tilde{k}_1^2}\ar@2[dl]**{}?(.2);?(.8)|{\tilde{k}_2^2}="Y2"\\
          \ar[r]_{f''^1}&\ar[r]_{f''^2}&{}          
          \POS
          \tria"X1";"Y2"
          \ar@2"Z2";"Z1"**{}?<(.2);?>(.8)_{\beta_2^2\rhc\alpha_2^2}
        }
        ,"B11"
        \xyboxmatrix"B11"{
          \ar[d]|{h_0^1}\ar[r]^{f^1}\ar@/_4pc/[dd]|{\tilde{k}_0^1\#_0\tilde{h}_0^1}="W2"&
          \ar[d]|{h_1^1}\ar[r]^{f^2}\ar@2[dl]**{}?(.2);?(.8)|{h_2^1}="X1"&
          \ar[d]^{h_1^2}\ar@2[dl]**{}?(.2);?(.8)|{h_2^2}="Y1"\\
          {}\POS*+{}="W1"\ar[d]|{k_0^1}\ar[r]|{f'^1}&
          \ar[d]|{k_1^1}\ar[r]|{f'^2}\ar@2[dl]**{}?(.2);?(.8)|{k_2^1}="X2"&
          {}\POS*+{}="Z1"\ar[d]^{k_1^2}\ar@2[dl]**{}?(.2);?(.8)|{k_2^2}="Y2"\\
          \ar[r]_{f''^1}&\ar[r]_{f''^2}&{}          
          \POS
          \tria"Y2";"X1"
          \ar@2"W1";"W2"**{}?<(.2);?>(.8)_{\beta_1^1\lhc\alpha_1^1}
        }
        ,"B21"
        \xyboxmatrix"B21"{
          \ar[d]_{\tilde{h}_0^1}\ar[r]^{f^1}&
          \ar[d]|{\tilde{h}_1^1}\ar[r]^{f^2}\ar@2[dl]**{}?(.2);?(.8)|{\tilde{h}_2^1}="X1"&
          \ar[d]|{\tilde{h}_1^2}\ar@2[dl]**{}?(.2);?(.8)|{\tilde{h}_2^2}="Y1"\ar@/^4pc/[dd]|{k_1^2\#_0h_1^2}="Z2"\\
          {}\POS*+{}="W1"\ar[d]_{\tilde{k}_0^1}\ar[r]|{f'^1}&
          \ar[d]|{\tilde{k}_1^1}\ar[r]|{f'^2}\ar@2[dl]**{}?(.2);?(.8)|{\tilde{k}_2^1}="X2"&
          {}\POS*+{}="Z1"\ar[d]|{\tilde{k}_1^2}\ar@2[dl]**{}?(.2);?(.8)|{\tilde{k}_2^2}="Y2"\\
          \ar[r]_{f''^1}&\ar[r]_{f''^2}&{}          
          \POS
          \tria"Y2";"X1"
          \ar@2"Z2";"Z1"**{}?<(.2);?>(.8)_{\beta_2^2\lhc\alpha_2^2}
        } 
        ,"B12"
        \xyboxmatrix"B12"{
          \ar[d]|{h_0^1}\ar[r]^{f^1}\ar@/_4pc/[dd]|{\tilde{k}_0^1\#_0\tilde{h}_0^1}="W2"&
          \ar[d]|{h_1^1}\ar[r]^{f^2}\ar@2[dl]**{}?(.2);?(.8)|{h_2^1}="X1"&
          \ar[d]^{h_1^2}\ar@2[dl]**{}?(.2);?(.8)|{h_2^2}="Y1"\\
          {}\POS*+{}="W1"\ar[d]|{k_0^1}\ar[r]|{f'^1}&
          \ar[d]|{k_1^1}\ar[r]|{f'^2}\ar@2[dl]**{}?(.2);?(.8)|{k_2^1}="X2"&
          {}\POS*+{}="Z1"\ar[d]^{k_1^2}\ar@2[dl]**{}?(.2);?(.8)|{k_2^2}="Y2"\\
          \ar[r]_{f''^1}&\ar[r]_{f''^2}&{}          
          \POS
          \tria"Y2";"X1"
          \ar@2"W1";"W2"**{}?<(.2);?>(.8)_{\beta_1^1\rhc\alpha_1^1}
        }
        ,"B22"
        \xyboxmatrix"B22"{
          \ar[d]_{\tilde{h}_0^1}\ar[r]^{f^1}&
          \ar[d]|{\tilde{h}_1^1}\ar[r]^{f^2}\ar@2[dl]**{}?(.2);?(.8)|{\tilde{h}_2^1}="X1"&
          \ar[d]|{\tilde{h}_1^2}\ar@2[dl]**{}?(.2);?(.8)|{\tilde{h}_2^2}="Y1"\ar@/^4pc/[dd]|{k_1^2\#_0h_1^2}="Z2"\\
          {}\POS*+{}="W1"\ar[d]_{\tilde{k}_0^1}\ar[r]|{f'^1}&
          \ar[d]|{\tilde{k}_1^1}\ar[r]|{f'^2}\ar@2[dl]**{}?(.2);?(.8)|{\tilde{k}_2^1}="X2"&
          {}\POS*+{}="Z1"\ar[d]|{\tilde{k}_1^2}\ar@2[dl]**{}?(.2);?(.8)|{\tilde{k}_2^2}="Y2"\\
          \ar[r]_{f''^1}&\ar[r]_{f''^2}&{}          
          \POS
          \tria"Y2";"X1"
          \ar@2"Z2";"Z1"**{}?<(.2);?>(.8)_{\beta_2^2\rhc\alpha_2^2}
        }
        \ar@3"A11";"A21"**{}?<;?>^{(m(\beta)\blhc{}m(\alpha))_3}
        \ar@3"A11";"A12"**{}?<;?>_{(f''^2\#_0f''^1\#_0\underline{(\beta\bten\alpha)_1})\\\#_1(m(k)\square_0m(h))}
        \ar@3"A21";"A22"**{}?<;?>^{(m(k)\square_0m(h))\\\#_1(\underline{(\beta\bten\alpha)_2}\#_0f''^2\#_0f''^1)}
        \ar@3"A12";"A22"**{}?<;?>_{(m(\beta)\brhc{}m(\alpha))_3}
        \ar@3"B11";"B21"**{}?<;?>^{(m(\beta\blhc{}\alpha))_3}
        \ar@3"B11";"B12"**{}?<;?>_{(f''^2\#_0f''^1\#_0\underline{(\beta\bten\alpha)_1})\\\#_1(m(k\square_0h))}
        \ar@3"B21";"B22"**{}?<;?>^{(m(k\square_0h))\\\#_1(\underline{(\beta\bten\alpha)_2}\#_0f''^2\#_0f''^1)}
        \ar@3"B12";"B22"**{}?<;?>_{(m(\beta\brhc{}\alpha))_3}
        \ar@3"A11";"B11"**{}?<;?>_{(f''^2\#_0f''^1\#_0{(\beta\blhc\alpha)_1})\\\#_1\underline{m^2_{k,h}}}
        \ar@3"A12";"B12"**{}?<;?>^{(f''^2\#_0f''^1\#_0{(\beta\brhc\alpha)_1})\\\#_1\underline{m^2_{k,h}}}
        \ar@3"A21";"B21"**{}?<;?>^{\underline{m^2_{\tilde{k},\tilde{h}}}\\\#_1({(\beta\blhc\alpha)_2}\#_0f''^2\#_0f''^1)}
        \ar@3"A22";"B22"**{}?<;?>_{\underline{m^2_{\tilde{k},\tilde{h}}}\\\#_1({(\beta\brhc\alpha)_2}\#_0f''^2\#_0f''^1)}
      \end{xy}
    \end{equation}
  \end{sidewaysfigure}

  We thus spell out the details of the top and bottom squares in
  (\ref{eq:mqmapprftenpres1cube}): The diagram
  (\ref{eq:mqmapprf21whskprescubetop}) shows the details of the top
  square of (\ref{eq:mqmapprftenpres1cube}). The central octagon of
  (\ref{eq:mqmapprf21whskprescubetop}) is broken down in
  (\ref{eq:mqmapprf21whskprescubetopsub}). The parts of these two
  diagrams commute essentially by the $\Gray$-category axioms and the
  definitions of 2- and 3-cells in the path space. The bottom square on
  (\ref{eq:mqmapprftenpres1cube}) is analogous.
  
  \begin{sidewaysfigure}
    \begin{equation}
      \label{eq:mqmapprf21whskprescubetop}
      \begin{xy}
        \save
        (50,0):(0,-3)::
        ,(0,0)="A11"
        ,(1,0)="A21"
        ,(2,0)="A31"
        ,(3,0)="A41"
        ,(4,0)="A51"
        ,(0,1)="A12"
        ,(1,1)="A22"
        ,(2,1)="A32"
        ,(3,1)="A42"
        ,(4,1)="A52"
        \save
        ,(.5,.25);p+(1,0):(0,-1)::
        ,(0,0)="B11"  
        ,(1,0)="B21"
        ,(2,0)="B31"
        ,(3,0)="B41"
        ,(4,0)="B51"
        \restore
        \save
        ,"A12",{c;"A52"**{}?(0)}+(0,-.25);"A12",{c;"A52"**{}?(1)}+(0,-.25):
        **{},?(.2):
        ,(0,0)="C11"  
        ,(1,0)="C21"
        ,(2,0)="C31"
        ,(3,0)="C41"
        ,(4,0)="C51"
        ,(5,0)="C61"
        \restore
        \restore
        \xyboxmatrix"A11"{
          \ar[d]|{h_0^1}\ar[r]^{f^1}\ar@/_4pc/[dd]|{\tilde{k}_0^1\#_0\tilde{h}_0^1}="W2"&
          \ar[d]|{h_1^1}\ar[r]^{f^2}\ar@2[dl]**{}?(.2);?(.8)|{h_2^1}="X1"&
          \ar[d]^{h_1^2}\ar@2[dl]**{}?(.2);?(.8)|{h_2^2}="Y1"\\
          {}\POS*+{}="W1"\ar[d]|{k_0^1}\ar[r]|{f'^1}&
          \ar[d]|{k_1^1}\ar[r]|{f'^2}\ar@2[dl]**{}?(.2);?(.8)|{k_2^1}="X2"&
          {}\POS*+{}="Z1"\ar[d]^{k_1^2}\ar@2[dl]**{}?(.2);?(.8)|{k_2^2}="Y2"\\
          \ar[r]_{f''^1}&\ar[r]_{f''^2}&{}          
          \POS
          \tria"X1";"Y2"
          \ar@2"W1";"W2"**{}?<(.2);?>(.8)_{\beta_1^1\lhc\alpha_1^1}
        },"A21"
        \xyboxmatrix"A21"{
          \ar[d]|{h_0^1}="Z1"\ar[r]^{f^1}\ar@/_2pc/[d]_{\tilde{h}_0^1}="Z2"&
          \ar[d]|{h_1^1}\ar[r]^{f^2}\ar@2[dl]**{}?(.2);?(.8)|{h_2^1}="X1"&
          \ar[d]^{h_1^2}\ar@2[dl]**{}?(.2);?(.8)|{h_2^2}="Y1"\\
          \ar[d]|{k_0^1}="W1"\ar[r]|{f'^1}\ar@/_2pc/[d]_{\tilde{k}_0^1}="W2"&
          \ar[d]|{k_1^1}\ar[r]|{f'^2}\ar@2[dl]**{}?(.2);?(.8)|{k_2^1}="X2"&
          \ar[d]^{k_1^2}\ar@2[dl]**{}?(.2);?(.8)|{k_2^2}="Y2"\\
          \ar[r]_{f''^1}&\ar[r]_{f''^2}&{}          
          \POS
          \ar@2"W1";"W2"**{}?(.2);?(.8)|{\beta_1^1}
          \ar@2"Z1";"Z2"**{}?(.2);?(.8)|{\alpha_1^1}="K"
          \tria"K";"Y2"
        },"B11"
        \xyboxmatrix"B11"{
          \ar[d]|{h_0^1}="Z1"\ar[r]^{f^1}\ar@/_2pc/[d]_{\tilde{h}_0^1}="Z2"&
          \ar[d]|{h_1^1}\ar[r]^{f^2}\ar@2[dl]**{}?(.2);?(.8)|{h_2^1}="X1"&
          \ar[d]^{h_1^2}\ar@2[dl]**{}?(.2);?(.8)|{h_2^2}="Y1"\\
          \ar[d]|{k_0^1}="W1"\ar[r]|{f'^1}\ar@/_2pc/[d]_{\tilde{k}_0^1}="W2"&
          \ar[d]|{k_1^1}\ar[r]|{f'^2}\ar@2[dl]**{}?(.2);?(.8)|{k_2^1}="X2"&
          \ar[d]^{k_1^2}\ar@2[dl]**{}?(.2);?(.8)|{k_2^2}="Y2"\\
          \ar[r]_{f''^1}&\ar[r]_{f''^2}&{}          
          \POS
          \ar@2"W1";"W2"**{}?(.2);?(.8)|{\beta_1^1}
          \ar@2"Z1";"Z2"**{}?(.2);?(.8)|{\alpha_1^1}="K"
          \tria"Y2";"K"
          \tria"X1";"Y2"
          \tria"K";"X2"
        },"B21"
        \xyboxmatrix"B21"{
          \ar[r]^{f^1}\ar[d]_{\tilde{h}_0^1}&
          \ar@/^1pc/[d]|{h_1^1}="Z1"\ar@/_1pc/[d]|{\tilde{h}_1^0}="Z2"
          \ar[r]^{f^2}\ar@2 {[l]**{}?(.5)};[dl]**{}?(.2);?(.8)|{\tilde{h}_2^1}="X1"&
          \ar[d]^{h_1^2}\ar@2 {[d];[dl]**{}?(.5)}**{}?(.2);?(.8)|{h_2^2}="Y1"\\
          \ar[d]|{k_0^1}="W1"\ar[r]|{f'^1}\ar@/_2pc/[d]_{\tilde{k}_0^1}="W2"&
          \ar[d]|{k_1^1}\ar[r]|{f'^2}\ar@2[dl]**{}?(.2);?(.8)|{k_2^1}="X2"&
          \ar[d]^{k_1^2}\ar@2[dl]**{}?(.2);?(.8)|{k_2^2}="Y2"\\
          \ar[r]_{f''^1}&\ar[r]_{f''^2}&{}          
          \POS
          \ar@2"W1";"W2"**{}?(.2);?(.8)|{\beta_1^1}
          \ar@2"Z1";"Z2"**{}?(.2);?(.8)|{\alpha_2^1}="K"
          \tria"X1";"Y2"
        },"A31"
        \xyboxmatrix"A31"{
          \ar[d]_{\tilde{h}_0^1}\ar[r]^{f^1}&
          \ar[d]|{\tilde{h}_1^1}\ar[r]^{f^2}\ar@2[dl]**{}?(.2);?(.8)|{\tilde{h}_2^1}="X1"&
          \ar[d]|{\tilde{h}_1^2}="Z2"\ar@2[dl]**{}?(.2);?(.8)|{\tilde{h}_2^2}="Y1"\ar@/^2pc/[d]^{{h}_1^2}="Z1"\\
          \ar[d]|{k_0^1}="W1"\ar[r]|{f'^1}\ar@/_2pc/[d]_{\tilde{k}_0^1}="W2"&
          \ar[d]|{k_1^1}\ar[r]|{f'^2}\ar@2[dl]**{}?(.2);?(.8)|{k_2^1}="X2"&
          \ar[d]^{k_1^2}\ar@2[dl]**{}?(.2);?(.8)|{k_2^2}="Y2"\\
          \ar[r]_{f''^1}&\ar[r]_{f''^2}&{}          
          \POS
          \ar@2"W1";"W2"**{}?(.2);?(.8)|{\beta_1^1}
          \ar@2"Z1";"Z2"**{}?(.2);?(.8)|{\alpha_2^2}
          \tria"X1";"Y2"
        }
        ,"B31"
        \xyboxmatrix"B31"{
          \ar[d]_{\tilde{h}_0^1}\ar[r]^{f^1}&
          \ar[d]|{\tilde{h}_1^1}\ar[r]^{f^2}\ar@2[dl]**{}?(.2);?(.8)|{\tilde{h}_2^1}="X1"&
          \ar[d]|{\tilde{h}_1^2}="Z2"\ar@2[dl]**{}?(.2);?(.8)|{\tilde{h}_2^2}="Y1"\ar@/^2pc/[d]^{{h}_1^2}="Z1"\\
          \ar[r]|{f'^1}\ar[d]_{\tilde{k}_1^1}&\ar@/^1pc/[d]|{k_0^2}="W1"\ar@/_1pc/[d]|{\tilde{k}_1^1}="W2"
          \ar[r]|{f'^2}\ar@2 {[l]**{}?(.5)};[dl]**{}?(.2);?(.8)|{\tilde{k}_2^1}="X2"&
          \ar[d]^{k_1^2}\ar@2 {[d];[dl]**{}?(.5)}**{}?(.2);?(.8)|{k_2^2}="Y2"\\
          \ar[r]_{f''^1}&\ar[r]_{f''^2}&{}          
          \POS
          \ar@2"W1";"W2"**{}?(.2);?(.8)|{\beta_2^1}
          \ar@2"Z1";"Z2"**{}?(.2);?(.8)|{\alpha_2^2}
          \tria"X1";"Y2"
        }
        ,"A51"
        \xyboxmatrix"A51"{
          \ar[d]_{\tilde{h}_0^1}\ar[r]^{f^1}&
          \ar[d]|{\tilde{h}_1^1}\ar[r]^{f^2}\ar@2[dl]**{}?(.2);?(.8)|{\tilde{h}_2^1}="X1"&
          \ar[d]|{\tilde{h}_1^2}="Z2"\ar@2[dl]**{}?(.2);?(.8)|{\tilde{h}_2^2}="Y1"\ar@/^4pc/[dd]|{k_1^2\#_0h_1^2}="Z1"\\
          \ar[d]_{\tilde{k}_0^1}\ar[r]|{f'^1}&
          \ar[d]|{\tilde{k}_1^1}\ar[r]|{f'^2}\ar@2[dl]**{}?(.2);?(.8)|{\tilde{k}_2^1}="X2"&
          {}\POS*+{}="Z2"\ar[d]|{\tilde{k}_1^2}="W2"\ar@2[dl]**{}?(.2);?(.8)|{\tilde{k}_2^2}="Y2"\\
          \ar[r]_{f''^1}&\ar[r]_{f''^2}&{}          
          \POS
          \ar@2"Z1";"Z2"**{}?<(.2);?>(.8)_{\beta_2^2\lhc\alpha_2^2}
          \tria"X1";"Y2"
        },"A41"
        \xyboxmatrix"A41"{
          \ar[d]_{\tilde{h}_0^1}\ar[r]^{f^1}&
          \ar[d]|{\tilde{h}_1^1}\ar[r]^{f^2}\ar@2[dl]**{}?(.2);?(.8)|{\tilde{h}_2^1}="X1"&
          \ar[d]|{\tilde{h}_1^2}="Z2"\ar@2[dl]**{}?(.2);?(.8)|{\tilde{h}_2^2}="Y1"\ar@/^2pc/[d]^{{h}_1^2}="Z1"\\
          \ar[d]_{\tilde{k}_0^1}\ar[r]|{f'^1}&
          \ar[d]|{\tilde{k}_1^1}\ar[r]|{f'^2}\ar@2[dl]**{}?(.2);?(.8)|{\tilde{k}_2^1}="X2"&
          \ar[d]|{\tilde{k}_1^2}="W2"\ar@2[dl]**{}?(.2);?(.8)|{\tilde{k}_2^2}="Y2"\ar@/^2pc/[d]^{\tilde{k}_1^2}="W1"\\
          \ar[r]_{f''^1}&\ar[r]_{f''^2}&{}          
          \POS
          \ar@2"W1";"W2"**{}?(.2);?(.8)|{\beta_2^2}="K"
          \ar@2"Z1";"Z2"**{}?(.2);?(.8)|{\alpha_2^2}
          \tria"X1";"K"
        },"B41"
        \xyboxmatrix"B41"{
          \ar[d]_{\tilde{h}_0^1}\ar[r]^{f^1}&
          \ar[d]|{\tilde{h}_1^1}\ar[r]^{f^2}\ar@2[dl]**{}?(.2);?(.8)|{\tilde{h}_2^1}="X1"&
          \ar[d]|{\tilde{h}_1^2}="Z2"\ar@2[dl]**{}?(.2);?(.8)|{\tilde{h}_2^2}="Y1"\ar@/^2pc/[d]^{{h}_1^2}="Z1"\\
          \ar[d]_{\tilde{k}_0^1}\ar[r]|{f'^1}&
          \ar[d]|{\tilde{k}_1^1}\ar[r]|{f'^2}\ar@2[dl]**{}?(.2);?(.8)|{\tilde{k}_2^1}="X2"&
          \ar[d]|{\tilde{k}_1^2}="W2"\ar@2[dl]**{}?(.2);?(.8)|{\tilde{k}_2^2}="Y2"\ar@/^2pc/[d]^{\tilde{k}_1^2}="W1"\\
          \ar[r]_{f''^1}&\ar[r]_{f''^2}&{}          
          \POS
          \ar@2"W1";"W2"**{}?(.2);?(.8)|{\beta_2^2}="K"
          \ar@2"Z1";"Z2"**{}?(.2);?(.8)|{\alpha_2^2}
          \tria"K";"X1"
          \tria"Y1";"K"
          \tria"X1";"Y2"
        },"A12"
        \xyboxmatrix"A12"{
          \ar[d]|{h_0^1}\ar[r]^{f^1}\ar@/_4pc/[dd]|{\tilde{k}_0^1\#_0\tilde{h}_0^1}="W2"&
          \ar[d]|{h_1^1}\ar[r]^{f^2}\ar@2[dl]**{}?(.2);?(.8)|{h_2^1}="X1"&
          \ar[d]^{h_1^2}\ar@2[dl]**{}?(.2);?(.8)|{h_2^2}="Y1"\\
          {}\POS*+{}="W1"\ar[d]|{k_0^1}\ar[r]|{f'^1}&
          \ar[d]|{k_1^1}\ar[r]|{f'^2}\ar@2[dl]**{}?(.2);?(.8)|{k_2^1}="X2"&
          {}\POS*+{}="Z1"\ar[d]^{k_1^2}\ar@2[dl]**{}?(.2);?(.8)|{k_2^2}="Y2"\\
          \ar[r]_{f''^1}&\ar[r]_{f''^2}&{}          
          \POS
          \tria"Y2";"X1"
          \ar@2"W1";"W2"**{}?<(.2);?>(.8)_{\beta_1^1\lhc\alpha_1^1}
        };"B11"**{}?(.68)
        \xyboxmatrix"C11"{
          \ar[d]|{h_0^1}="Z1"\ar[r]^{f^1}\ar@/_2pc/[d]_{\tilde{h}_0^1}="Z2"&
          \ar[d]|{h_1^1}\ar[r]^{f^2}\ar@2[dl]**{}?(.2);?(.8)|{h_2^1}="X1"&
          \ar[d]^{h_1^2}\ar@2[dl]**{}?(.2);?(.8)|{h_2^2}="Y1"\\
          \ar[d]|{k_0^1}="W1"\ar[r]|{f'^1}\ar@/_2pc/[d]_{\tilde{k}_0^1}="W2"&
          \ar[d]|{k_1^1}\ar[r]|{f'^2}\ar@2[dl]**{}?(.2);?(.8)|{k_2^1}="X2"&
          \ar[d]^{k_1^2}\ar@2[dl]**{}?(.2);?(.8)|{k_2^2}="Y2"\\
          \ar[r]_{f''^1}&\ar[r]_{f''^2}&{}          
          \POS
          \ar@2"W1";"W2"**{}?(.2);?(.8)|{\beta_1^1}
          \ar@2"Z1";"Z2"**{}?(.2);?(.8)|{\alpha_1^1}="K"
          \tria"Y2";"X1"
          \tria"K";"X2"
        },"A22"
        \xyboxmatrix"A22"@W+.6cm{
          {}
        },"A32"
        \xyboxmatrix"A32"{
          \ar[d]_{\tilde{h}_0^1}\ar[rr]^{f^1}&&\ar[dl]|{\tilde{h}_1^1}\ar[dr]|{h_1^1}\ar[rr]^{f^2}
          \ar@2[dll]**{}?(.2);?(.8)|{\tilde{h}_2^1}="X1"&&
          \ar[d]^{h_1^2}\ar@2[dl]**{}?(.2);?(.8)|{h_2^2}="Y1"\\
          \ar[r]|{f'^1}\ar[d]_{\tilde{k}_0^1}="W2"&\ar[dr]|{\tilde{k}_1^1}="Z2"\ar@2[dl]**{}?(.2);?(.8)|{\tilde{k}_2^1}="X2"
          &&\ar@2[ll]**{}?(.2);?(.8)|{\beta_2^1\lhc\alpha_2^1}="R"
          \ar[r]|{f'^2}\ar[dl]|{k_1^1}="Z1"&
          \ar[d]^{k_1^2}\ar@2[dll]**{}?(.2);?(.8)|{k_2^2}="Y2"\\
          \ar[rr]_{f''^1}&&\ar[rr]_{f''^2}&&{}          
        },"A42"
        \xyboxmatrix"A42"{
          {}
        },"A52"
        \xyboxmatrix"A52"{
          \ar[d]_{\tilde{h}_0^1}\ar[r]^{f^1}&
          \ar[d]|{\tilde{h}_1^1}\ar[r]^{f^2}\ar@2[dl]**{}?(.2);?(.8)|{\tilde{h}_2^1}="X1"&
          \ar[d]|{\tilde{h}_1^2}\ar@2[dl]**{}?(.2);?(.8)|{\tilde{h}_2^2}="Y1"\ar@/^4pc/[dd]|{k_1^2\#_0h_1^2}="Z2"\\
          {}\POS*+{}="W1"\ar[d]_{\tilde{k}_0^1}\ar[r]|{f'^1}&
          \ar[d]|{\tilde{k}_1^1}\ar[r]|{f'^2}\ar@2[dl]**{}?(.2);?(.8)|{\tilde{k}_2^1}="X2"&
          {}\POS*+{}="Z1"\ar[d]|{\tilde{k}_1^2}\ar@2[dl]**{}?(.2);?(.8)|{\tilde{k}_2^2}="Y2"\\
          \ar[r]_{f''^1}&\ar[r]_{f''^2}&{}          
          \POS
          \tria"Y2";"X1"
          \ar@2"Z2";"Z1"**{}?<(.2);?>(.8)_{\beta_2^2\lhc\alpha_2^2}
        },"C21"
        \xyboxmatrix"C21"{
          \ar[r]^{f^1}\ar[d]_{\tilde{h}_0^1}&
          \ar@/^1pc/[d]|{h_1^1}="Z1"\ar@/_1pc/[d]|{\tilde{h}_1^0}="Z2"
          \ar[r]^{f^2}\ar@2 {[l]**{}?(.5)};[dl]**{}?(.2);?(.8)|{\tilde{h}_2^1}="X1"&
          \ar[d]^{h_1^2}\ar@2 {[d];[dl]**{}?(.5)}**{}?(.2);?(.8)|{h_2^2}="Y1"\\
          \ar[d]|{k_0^1}="W1"\ar[r]|{f'^1}\ar@/_2pc/[d]_{\tilde{k}_0^1}="W2"&
          \ar[d]|{k_1^1}\ar[r]|{f'^2}\ar@2[dl]**{}?(.2);?(.8)|{k_2^1}="X2"&
          \ar[d]^{k_1^2}\ar@2[dl]**{}?(.2);?(.8)|{k_2^2}="Y2"\\
          \ar[r]_{f''^1}&\ar[r]_{f''^2}&{}          
          \POS
          \ar@2"W1";"W2"**{}?(.2);?(.8)|{\beta_1^1}
          \ar@2"Z1";"Z2"**{}?(.2);?(.8)|{\alpha_2^1}="K"
          \tria"Y2";"K"
        },"C31"
        \xyboxmatrix"C31"{
          \ar[r]^{f^1}\ar[d]_{\tilde{h}_0^1}&
          \ar@/^1pc/[d]|{h_1^1}="Z1"\ar@/_1pc/[d]|{\tilde{h}_1^0}="Z2"
          \ar[r]^{f^2}\ar@2 {[l]**{}?(.5)};[dl]**{}?(.2);?(.8)|{\tilde{h}_2^1}="X1"&
          \ar[d]^{h_1^2}\ar@2 {[d];[dl]**{}?(.5)}**{}?(.2);?(.8)|{h_2^2}="Y1"\\
          \ar[d]|{\tilde{k}_0^1}="W1"\ar[r]|{f'^1}&
          \ar@/^1pc/[d]|{k_1^1}="W1"\ar@/_1pc/[d]|{\tilde{k}_1^0}="W2"
          \ar@2 {[l]**{}?(.5)};[dl]**{}?(.2);?(.8)|{\tilde{k}_2^1}="X2"\ar[r]|{f'^2}&
          \ar[d]^{k_1^2}\ar@2 {[d];[dl]**{}?(.5)}**{}?(.2);?(.8)|{k_2^2}="Y2"\\
          \ar[r]_{f''^1}&\ar[r]_{f''^2}&{}          
          \POS
          \ar@2"W1";"W2"**{}?(.2);?(.8)|{\beta_2^1}="K2"
          \ar@2"Z1";"Z2"**{}?(.2);?(.8)|{\alpha_2^1}="K1"
          \tria"X1";"K2"
          \tria"Y2";"K1"
        }
        ,"C31";"B21"**{}?(.5)
        \xyboxmatrix"D1"{
          \ar[r]^{f^1}\ar[d]_{\tilde{h}_0^1}&
          \ar@/^1pc/[d]|{h_1^1}="Z1"\ar@/_1pc/[d]|{\tilde{h}_1^0}="Z2"
          \ar[r]^{f^2}\ar@2 {[l]**{}?(.5)};[dl]**{}?(.2);?(.8)|{\tilde{h}_2^1}="X1"&
          \ar[d]^{h_1^2}\ar@2 {[d];[dl]**{}?(.5)}**{}?(.2);?(.8)|{h_2^2}="Y1"\\
          \ar[d]|{k_0^1}="W1"\ar[r]|{f'^1}\ar@/_2pc/[d]_{\tilde{k}_0^1}="W2"&
          \ar[d]|{k_1^1}\ar[r]|{f'^2}\ar@2[dl]**{}?(.2);?(.8)|{k_2^1}="X2"&
          \ar[d]^{k_1^2}\ar@2[dl]**{}?(.2);?(.8)|{k_2^2}="Y2"\\
          \ar[r]_{f''^1}&\ar[r]_{f''^2}&{}          
          \POS
          \ar@2"W1";"W2"**{}?(.2);?(.8)|{\beta_1^1}
          \ar@2"Z1";"Z2"**{}?(.2);?(.8)|{\alpha_2^1}="K"
          \tria"K";"Y2"
          \tria"Y2";"X1"
        },"C41"
        \xyboxmatrix"C41"{
          \ar[r]^{f^1}\ar[d]_{\tilde{h}_0^1}&
          \ar@/^1pc/[d]|{h_1^1}="Z1"\ar@/_1pc/[d]|{\tilde{h}_1^0}="Z2"
          \ar[r]^{f^2}\ar@2 {[l]**{}?(.5)};[dl]**{}?(.2);?(.8)|{\tilde{h}_2^1}="X1"&
          \ar[d]^{h_1^2}\ar@2 {[d];[dl]**{}?(.5)}**{}?(.2);?(.8)|{h_2^2}="Y1"\\
          \ar[d]|{\tilde{k}_0^1}="W1"\ar[r]|{f'^1}&
          \ar@/^1pc/[d]|{k_1^1}="W1"\ar@/_1pc/[d]|{\tilde{k}_1^0}="W2"
          \ar@2 {[l]**{}?(.5)};[dl]**{}?(.2);?(.8)|{\tilde{k}_2^1}="X2"\ar[r]|{f'^2}&
          \ar[d]^{k_1^2}\ar@2 {[d];[dl]**{}?(.5)}**{}?(.2);?(.8)|{k_2^2}="Y2"\\
          \ar[r]_{f''^1}&\ar[r]_{f''^2}&{}          
          \POS
          \ar@2"W1";"W2"**{}?(.2);?(.8)|{\beta_2^1}="K2"
          \ar@2"Z1";"Z2"**{}?(.2);?(.8)|{\alpha_2^1}="K1"
          \tria"K2";"X1"
          \tria"K1";"Y2"
        },"C51"
        \xyboxmatrix"C51"{
          \ar[d]_{\tilde{h}_0^1}\ar[r]^{f^1}&
          \ar[d]|{\tilde{h}_1^1}\ar[r]^{f^2}\ar@2[dl]**{}?(.2);?(.8)|{\tilde{h}_2^1}="X1"&
          \ar[d]|{\tilde{h}_1^2}="Z2"\ar@2[dl]**{}?(.2);?(.8)|{\tilde{h}_2^2}="Y1"\ar@/^2pc/[d]^{{h}_1^2}="Z1"\\
          \ar[r]|{f'^1}\ar[d]_{\tilde{k}_1^1}&\ar@/^1pc/[d]|{k_0^2}="W1"\ar@/_1pc/[d]|{\tilde{k}_1^1}="W2"
          \ar[r]|{f'^2}\ar@2 {[l]**{}?(.5)};[dl]**{}?(.2);?(.8)|{\tilde{k}_2^1}="X2"&
          \ar[d]^{k_1^2}\ar@2 {[d];[dl]**{}?(.5)}**{}?(.2);?(.8)|{k_2^2}="Y2"\\
          \ar[r]_{f''^1}&\ar[r]_{f''^2}&{}          
          \POS
          \ar@2"W1";"W2"**{}?(.2);?(.8)|{\beta_2^1}="K"
          \ar@2"Z1";"Z2"**{}?(.2);?(.8)|{\alpha_2^2}
          \tria"K";"X1"
        },"C41";"B31"**{}?(.5)
        \xyboxmatrix"D2"{
          \ar[d]_{\tilde{h}_0^1}\ar[r]^{f^1}&
          \ar[d]|{\tilde{h}_1^1}\ar[r]^{f^2}\ar@2[dl]**{}?(.2);?(.8)|{\tilde{h}_2^1}="X1"&
          \ar[d]|{\tilde{h}_1^2}="Z2"\ar@2[dl]**{}?(.2);?(.8)|{\tilde{h}_2^2}="Y1"\ar@/^2pc/[d]^{{h}_1^2}="Z1"\\
          \ar[r]|{f'^1}\ar[d]_{\tilde{k}_1^1}&\ar@/^1pc/[d]|{k_0^2}="W1"\ar@/_1pc/[d]|{\tilde{k}_1^1}="W2"
          \ar[r]|{f'^2}\ar@2 {[l]**{}?(.5)};[dl]**{}?(.2);?(.8)|{\tilde{k}_2^1}="X2"&
          \ar[d]^{k_1^2}\ar@2 {[d];[dl]**{}?(.5)}**{}?(.2);?(.8)|{k_2^2}="Y2"\\
          \ar[r]_{f''^1}&\ar[r]_{f''^2}&{}          
          \POS
          \ar@2"W1";"W2"**{}?(.2);?(.8)|{\beta_2^1}="K"
          \ar@2"Z1";"Z2"**{}?(.2);?(.8)|{\alpha_2^2}
          \tria"X1";"K"
          \tria"Y2";"X1"
        },"A52";"B41"**{}?(.68)
        \xyboxmatrix"C61"{
          \ar[d]_{\tilde{h}_0^1}\ar[r]^{f^1}&
          \ar[d]|{\tilde{h}_1^1}\ar[r]^{f^2}\ar@2[dl]**{}?(.2);?(.8)|{\tilde{h}_2^1}="X1"&
          \ar[d]|{\tilde{h}_1^2}="Z2"\ar@2[dl]**{}?(.2);?(.8)|{\tilde{h}_2^2}="Y1"\ar@/^2pc/[d]^{{h}_1^2}="Z1"\\
          \ar[d]_{\tilde{k}_0^1}\ar[r]|{f'^1}&
          \ar[d]|{\tilde{k}_1^1}\ar[r]|{f'^2}\ar@2[dl]**{}?(.2);?(.8)|{\tilde{k}_2^1}="X2"&
          \ar[d]|{\tilde{k}_1^2}="W2"\ar@2[dl]**{}?(.2);?(.8)|{\tilde{k}_2^2}="Y2"\ar@/^2pc/[d]^{\tilde{k}_1^2}="W1"\\
          \ar[r]_{f''^1}&\ar[r]_{f''^2}&{}          
          \POS
          \ar@2"W1";"W2"**{}?(.2);?(.8)|{\beta_2^2}="K"
          \ar@2"Z1";"Z2"**{}?(.2);?(.8)|{\alpha_2^2}
          \tria"Y1";"K"
          \tria"Y2";"X1"
        }
        \ar@3@<-2ex>"A11";"A12"
        \ar@3"A11";"A21"
        \ar@3"A11";"B11"
        \ar@3"A12";"A32"
        \ar@3"A12";"C11"
        \ar@3"A21";"A31"
        \ar@3"A21";"B21"
        \ar@3"A31";"A41"
        \ar@3"A31";"B31"
        \ar@3"A32";"A52"
        \ar@3"A32";"C41"
        \ar@3"A41";"A51"
        \ar@3"A41";"B41"
        \ar@3@<+2ex>"A51";"A52"
        \ar@3"B11";"A21"
        \ar@3"B11";"C11"
        \ar@3"B21";"A31"
        \ar@3"B21";"C21"
        \ar@3"B21";"D1"
        \ar@3"B31";"A41"
        \ar@3"B31";"C51"
        \ar@3"B31";"D2"
        \ar@3"B41";"A51"
        \ar@3"B41";"C61"
        \ar@3"C11";"C21"
        \ar@3"C21";"C31"
        \ar@3"C31";"A32"
        \ar@3"C41";"C51"
        \ar@3"C41";"D2"
        \ar@3"C51";"C61"
        \ar@3"C61";"A52"
        \ar@3"D1";"C21"
        \ar@3"D1";"C31"
        \ar@3"D2";"C51"
        \ar@3@/^4pc/"A21";"C11"
        \ar@3@/_4pc/"A41";"C61"
      \end{xy}
    \end{equation}
  \end{sidewaysfigure}
  This proves (\ref{eq:coccohten}).

  Furthermore, we check that tensors of cocycle elements are trivial: We
  calculate according to section \ref{sec:pthspcten}:
  \begin{equation*}
    m^2_{f_1,f_2}\bten{}m^2_{f_3,f_4}=((m^2_{f_1,f_2})_1\ten{}(m^2_{f_3,f_4})_1,(m^2_{f_1,f_2})_2\ten{}(m^2_{f_3,f_4})_2)\,,
  \end{equation*}
  where according to (\ref{eq:pastecocdef}) all the arguments on the
  right are trivial, hence their tensors are trivial, that is,
  (\ref{eq:tentriv}) holds.

  Lastly, images of 2-cells tensor trivially with co-cycle components
  by the unitality of the tensor in $\H$ and the fact that the 2-cell
  faces of $m^2$ are trivial, hence verifying (\ref{eq:ten2coccel-l})
  and (\ref{eq:ten2coccel-r}).
\end{prf}

\begin{thm}
  There is a 
  pseudo $\Gray$-functor $m$ such that 
  \begin{equation}
    \label{eq:internalgraycat}
    \xymatrix@+.5cm{
      \pathspc\H\times_\H\pathspc\H \ar[r]|-*\dir{/}^-m& 
      \pathspc\H\ar@/^/[r]^{d_1}\ar@/_/[r]_{d_0}& \H \ar[l]|{i}
    }
  \end{equation}
  is an internal category object in $\Gray\Cat_{\fQ^1}$.
\end{thm}

\begin{prf}
  We need to verify that $m$ is an associative and unital
  operation.   We need to check first that
  \begin{equation*}
    \begin{xy}
      \xyboxmatrix{
        \pathspc\H\times_{d_0,d_1}\pathspc\H\times_{d_0,d_1}\pathspc\H
        \ar[r]^-{\pathspc\H\times{}m}|-*\dir{/}\ar[d]_{m\times\pathspc\H}|-*\dir{/}&
        \pathspc\H\times_{d_0,d_1}\pathspc\H\ar[d]^m|-*\dir{/}\\
        \pathspc\H\times_{d_0,d_1}\pathspc\H\ar[r]_-{m}|-*\dir{/}&\pathspc\H
      }
    \end{xy}\,,
  \end{equation*}
  where $m\times\pathspc\H$ and $\pathspc\H\times{}m$ exist by the
  observation in remark \ref{rem:indmap}. On the level of globular maps this
  is obvious, since it is just pasting according to definition
  \ref{defn:pathcomp}. Proving that the cocylces both ways around are
  the same, means drawing a diagram that looks like
  (\ref{eq:mqmapprf2coc}) with each array transposed.

  Unitality is obvious, source and target conditions
  
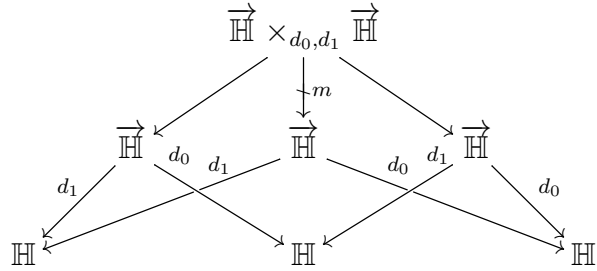
\begin{equation*}
    \begin{xy}
      \xyboxmatrix{
        {}&{}&\pathspc{\H}\times_{d_0,d_1}\pathspc{\H}\ar[dl]\ar[dr]\ar[d]^{m}|*-\dir{/}&{}&{}\\
        {}&\pathspc{\H}\ar[dl]_{d_1}\ar[dr]^(.2){d_0}&
        \pathspc{\H}\ar[dll]_(.3){d_1}|!{"2,2";"3,3"}\hole
        \ar[drr]^(.3){d_0}|!{"2,4";"3,3"}\hole&\pathspc{\H}\ar[dl]_(.2){d_1}\ar[dr]^{d_0}&\\
        \H&&\H&&\H
      }
    \end{xy}
  \end{equation*}
  hold by definition \ref{defn:pathcomp}. In particular, the 2-cell components of
  $m^2$ are trivial, thus $d_0m$ and  $d_1m$ are strict $\Gray$-functors, even
  though $m$ is pseudo.
\end{prf}

\begin{lem}
  \label{lem:multstrictnat}
  For a strict $\Gray$-functor $F$ the multiplication map $m$ is natural,
  that is
  \begin{equation}
    \label{eq:multstrictnat}
    \begin{xy}
      \xyboxmatrix{
        \pathspc\H\times_{d_0,d_1}\pathspc\H\ar[r]|-*\dir{/}^-m
        \ar[d]_{\pathspc{F}\times\pathspc{F}}&
        \pathspc\H\ar[d]^{F}\\
        \pathspc\K\times_{d_0,d_1}\pathspc\K\ar[r]|-*\dir{/}_-m&\pathspc\K
      }
    \end{xy}\,.
  \end{equation}
\end{lem}
Note that by (\ref{eq:pbprodstr}) we have
$(\pathspc{F}e)\dot\times(\pathspc{F}e)=(\pathspc{F}\times\pathspc{F})e.$

\begin{prf}
  Verifying (\ref{eq:multstrictnat}) elementwise is straightforward.
\end{prf}

We can define the 1-cell  inverse to
\begin{equation}\label{eq:2dimcell}
  \begin{xy}
    \xymatrix{{} \ar[r]^f \ar[d]_{g_0}&
      {}\ar[d]^{g_1} \ar@2 [dl] **{} ?<(.3); ?>(.7) ^{g_2}\\{} \ar[r]_{f'}& {}}
  \end{xy}
\end{equation}
 with respect to $m$ as
\begin{equation}\label{eq:2dimcellinv}
  \begin{matrix}
    \begin{xy}
      \xyboxmatrix{
        {}\ar[rrr]^{\overline{f}}\ar[ddd]_{g_1}\ar[dr]^{\overline{f}}&&&\ar[ddd]^{g_0}\\
        {}&\ar[r]^{g_0}\ar[d]_{f}&\ar[d]^{f'}\ar@2 [dl] **{} ?<(.3); ?>(.7) ^{\overline{g_2}}&\\
        {}&\ar[r]_{g_1}&\ar[dr]^{\overline{f'}}&\\
        \ar[rrr]_{\overline{f'}}&&& }
    \end{xy}
  \end{matrix}
\end{equation}
where $\overline{(\_)}$ is the respective vertical inverse in $\H$.

\begin{lem}
  The path space 1-cell in \eqref{eq:2dimcellinv} is a left and right
  inverse to \eqref{eq:2dimcell} with respect to $m$.
\end{lem}
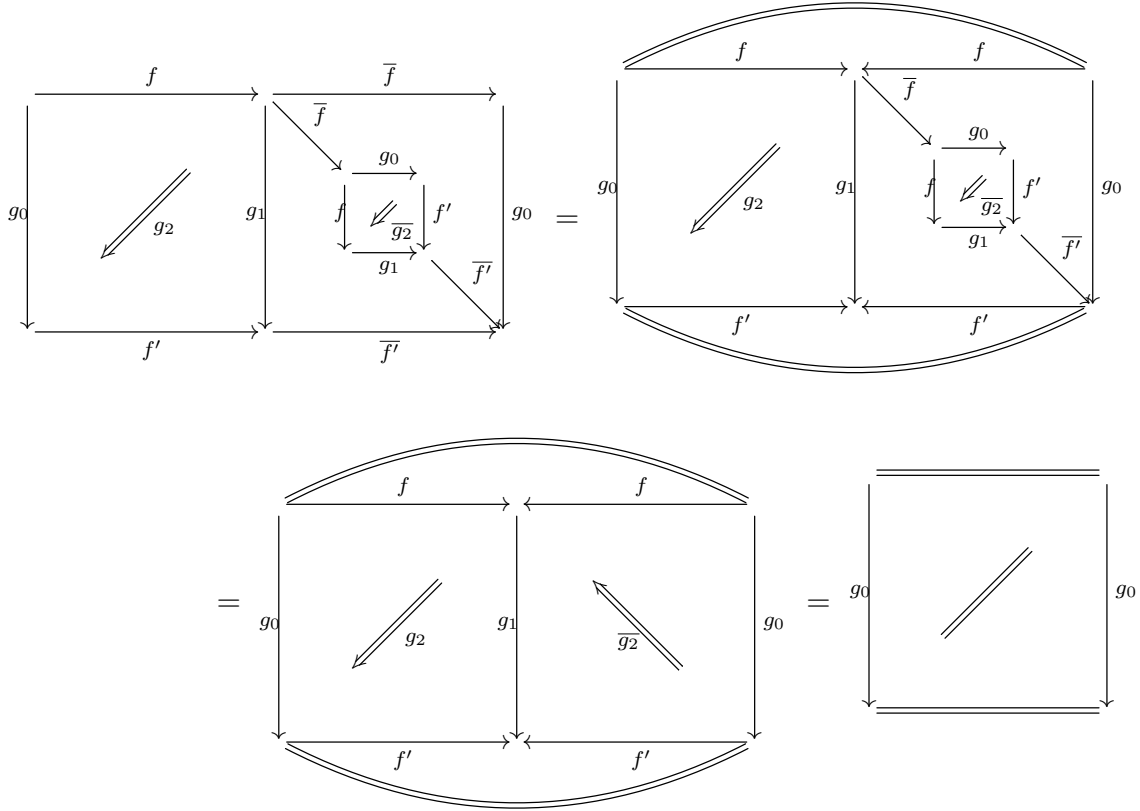
\begin{prf}
  \begin{multline*}
    \begin{xy}
      \xyboxmatrix{
        {}\ar[rrr]^f\ar[ddd]_{g_0}&&&\ar[rrr]^{\overline{f}}
        \ar[ddd]_{g_1}\ar[dr]^{\overline{f}}\ar@2[dddlll]**{}?<(.3);?>(.7)^{g_2}&&&
        \ar[ddd]^{g_0}\\
        {}&&&&\ar[r]^{g_0}\ar[d]_{f}&\ar[d]^{f'}\ar@2 [dl] **{} ?<(.3); ?>(.7)^{\overline{g_2}}&\\
        {}&&&&\ar[r]_{g_1}&\ar[dr]^{\overline{f'}}&\\
        {}\ar[rrr]_{f'}&&&\ar[rrr]_{\overline{f'}}&&&}
    \end{xy}=
    \begin{xy}
      \xyboxmatrix{
        {}\ar[rrr]^f\ar[ddd]_{g_0}\ar@{=}@/^2pc/[rrrrrr]&&&
        \ar[ddd]_{g_1}\ar[dr]^{\overline{f}}\ar@2[dddlll]**{}?<(.3);?>(.7)^{g_2}&&&
        \ar[ddd]^{g_0}\ar[lll]_{f}\\
        {}&&&&\ar[r]^{g_0}\ar[d]_{f}&\ar[d]^{f'}\ar@2 [dl] **{} ?<(.3); ?>(.7)^{\overline{g_2}}&\\
        {}&&&&\ar[r]_{g_1}&\ar[dr]^{\overline{f'}}&\\
        {}\ar[rrr]_{f'}\ar@{=}@/_2pc/[rrrrrr]_{}&&&&&&\ar[lll]^{f'}}
    \end{xy}\\=
    \begin{xy}
      \xyboxmatrix{
        {}\ar[rrr]^f\ar[ddd]_{g_0}\ar@{=}@/^2pc/[rrrrrr]^{}&&&
        \ar[ddd]_{g_1}\ar@2[dddlll]**{}?<(.3);?>(.7)^{g_2}&&&
        \ar[ddd]^{g_0}\ar[lll]_{f}\\
        {}&&&&&&\\
        {}&&&&&&\\
        {}\ar[rrr]_{f'}\ar@{=}@/_2pc/[rrrrrr]&&&&&&\ar[lll]^{f'}
        \ar@2[uuulll]**{}?<(.3);?>(.7)^{\overline{g_2}}}
    \end{xy}=
    \begin{xy}
      \xyboxmatrix{
        \ar[ddd]_{g_0}\ar@{=}[rrr]&&&\ar[ddd]^{g_0}\ar@{=}[dddlll]**{}?<(.3);?>(.7)\\
        {}&&&\\
        {}&&&\\
        \ar@{=}[rrr]&&&
      }
    \end{xy}
  \end{multline*}
  And similarly for the right inverse.
\end{prf}

Furthermore, these inverses behave well with respect to the internal category structure:
\begin{thm}
  Given the situation in \eqref{eq:internalgraycat}, assume $\H$ is a
  $\Gray$-groupoid, then there is a $\fQ^1$-map
  $o\from\pathspc\H\laxto\pathspc\H$ (``opposite'') such that
  \eqref{eq:internalgraycat} becomes an internal groupoid in $\Gray\Cat_{\fQ^1}$.
\end{thm}
\begin{prf}
  The action of $o$ on 0- and 1-cells is already given in
  \eqref{eq:2dimcellinv}, the effect on 2- and 3-cells of $\pathspc\H$
  is analogous.

  Furthermore, we need to give a 2-cocycle $o^2_{h,g}\from
  o(h)\square_0o(g)\to{}o(h\square_0g)$ the non-trivial part of which
  is the following 3-cell:
  {\allowdisplaybreaks
    \begin{multline*}
    o\begin{pmatrix}
      \begin{xy}
        \xyboxmatrix{{}\ar[r]^{f'}\ar[d]_{h_0}&
          {}\ar[d]^{h_1}\ar@2[dl]**{}?<(.3);?>(.7)^{h_2}\\
          {}\ar[r]_{f''}& {}}
      \end{xy}   
    \end{pmatrix}\square_0
    o\begin{pmatrix}
      \begin{xy}
        \xyboxmatrix{{}\ar[r]^f\ar[d]_{g_0}&
          {}\ar[d]^{g_1}\ar@2 [dl]**{}?<(.3);?>(.7)^{g_2}\\
          {}\ar[r]_{f'}&{}}
      \end{xy}      
    \end{pmatrix}=
    \begin{xy}
      \xyboxmatrix{
        {}\ar[rrr]^{\overline{f}}\ar[ddd]_{g_1}\ar[dr]^{\overline{f}}&&&\ar[ddd]^{g_0}\\
        {}&\ar[r]^{g_0}\ar[d]_{f}&\ar[d]^{f'}\ar@2 [dl] **{} ?<(.3); ?>(.7) ^{\overline{g_2}}&\\
        {}&\ar[r]_{g_1}&\ar[dr]^{\overline{f'}}&\\
        \ar[rrr]|{\overline{f'}}\ar[ddd]_{h_1}\ar[dr]^{\overline{f'}}&&&\ar[ddd]^{h_0}\\
        {}&\ar[r]^{h_0}\ar[d]_{f'}&\ar[d]^{f''}\ar@2 [dl] **{} ?<(.3); ?>(.7) ^{\overline{h_2}}&\\
        {}&\ar[r]_{h_1}&\ar[dr]^{\overline{f''}}&\\
        \ar[rrr]_{\overline{f''}}&&& 
      }
    \end{xy}
    \\=\begin{xy}
      \xyboxmatrix"A"{
        {}\ar[rrr]^{\overline{f}}\ar[ddd]_{g_1}\ar[dr]^{\overline{f}}&&&\ar[ddd]^{g_0}\\
        {}&\ar[r]^{g_0}\ar[d]_{f}&\ar[d]^{f'}\ar@2[dl]**{}?<(.3);?>(.7)^{\overline{g_2}}|{}="x"&\\
        {}&\ar[r]_(.1){g_1}&\ar[ddl]^(.3){\overline{f'}}&\\
        \ar[ddd]_{h_1}&&&\ar[ddd]^{h_0}\\
        {}&\ar[r]^(.9){h_0}\ar[d]_{f'}&\ar[d]^{f''}\ar@2[dl]**{}?<(.3);?>(.7)^{\overline{h_2}}|{}="y"&\\
        {}&\ar[r]_{h_1}&\ar[dr]^{\overline{f''}}&\\
        \ar[rrr]_{\overline{f''}}&&& 
        {}\tria "x";"y"
      }
      +(65,0)
      \xyboxmatrix"B"{
        {}\ar[rrr]^{\overline{f}}\ar[ddd]_{g_1}\ar[dr]^{\overline{f}}&&&\ar[ddd]^{g_0}\\
        {}&\ar[r]^{g_0}\ar[d]_{f}&\ar[d]^{f'}\ar@2[dl]**{}?<(.3);?>(.7)^{\overline{g_2}}|{}="x"&\\
        {}&\ar[r]_(.1){g_1}&\ar[ddl]^(.3){\overline{f'}}&\\
        \ar[ddd]_{h_1}&&&\ar[ddd]^{h_0}\\
        {}&\ar[r]^(.9){h_0}\ar[d]_{f'}&\ar[d]^{f''}\ar@2[dl]**{}?<(.3);?>(.7)^{\overline{h_2}}|{}="y"&\\
        {}&\ar[r]_{h_1}&\ar[dr]^{\overline{f''}}&\\
        \ar[rrr]_{\overline{f''}}&&& 
        {}\tria "y";"x"
      }
      \POS\ar@3 "A";"B" ^{\overline{f''}\#_0\\\;\left((\overline{h_2}\#_0\overline{f'})\ten\overline{g_2}\right)\\\;\#_0\overline{f}}
    \end{xy}
    \\=
    \begin{xy}
      \xyboxmatrix{
        {}\ar[rrr]^{\overline{f}}\ar[ddd]_{h_1\#_0g_1}\ar[dr]^{\overline{f}}&&&\ar[ddd]^{h_0\#_0g_0}\\
        {}&\ar[r]^{h_0\#_0g_0}\ar[d]_{f}&\ar[d]^{f''}\ar@2 [dl] **{} ?<(.3); ?>(.7)|{}="x"
        \save\POS+(30,10)*!C\labelbox{(\overline{h_2}\#_0\overline{f'})\rhc{}\overline{g_2}\\
          \;=\overline{(h_2\#_0\overline{f'})\lhc{}g_2}\\
          \;=\overline{h_2\lhc{}(\overline{f'}\#_0g_2)}\\
          \;=\overline{(h_2\#_0g_0)\#_1(h_1\#_0g_2)}}\ar@{.}"x";c\restore&\\
        {}&\ar[r]_{h_1\#_0g_1}&\ar[dr]^{\overline{f''}}&\\
        \ar[rrr]_{\overline{f''}}&&& }
    \end{xy}=
    o\begin{pmatrix}
      \begin{xy}
        \xyboxmatrix{{}\ar[r]^f\ar[d]_{g_0}&
          {}\ar[d]^{g_1}\ar@2 [dl]**{}?<(.3);?>(.7)^{g_2}\\
          {}\ar[r]_{f'}\ar[d]_{h_0}&
          {}\ar[d]^{h_1}\ar@2[dl]**{}?<(.3);?>(.7)^{h_2}\\
          {}\ar[r]_{f''}& {}}
      \end{xy}      
    \end{pmatrix}
  \end{multline*}
}
  For the relationship between horizontal composition and pasting of
  squares see remark \ref{rem:sqhorcomp}. 

  We check that $o^2$ is indeed a 2-cocycle. Given suitably incident
  1-cells of $\H$ we need to verify that the analog of \eqref{eq:cocycle} hold, that is, 
  \begin{equation*}
    o^2_{k, h\square_0g}\square_1(o(k)\square_0o^2_{h,g})=o^2_{k\square_0h,g}\square_1(o^2_{k,h}\square_0o(g)),
  \end{equation*} 
  hence (\ref{eq:oppfunct}) commutes.
  \begin{figure}
    \centering
    \begin{equation}
      \label{eq:oppfunct}
      \begin{xy}
        \xyboxmatrix"A"{
          {}\ar[rrr]^{\overline{f}}\ar[ddd]_{g_1}\ar[dr]^{\overline{f}}&&&\ar[ddd]^{g_0}\\
          {}&\ar[r]^{g_0}\ar[d]_{f}
          &\ar[d]^{f'}\ar@2[dl]**{}?<(.3);?>(.7)^{\overline{g_2}}|{}="x"&\\
          {}&\ar[r]_(.1){g_1}&\ar[ddl]^(.3){\overline{f'}}&\\
          \ar[ddd]_{h_1}&&&\ar[ddd]^{h_0}\\
          {}&\ar[r]^(.9){h_0}\ar[d]_{f'}\POS.[rdddd]!C*+++\frm{.}="yz"&
          \ar[d]^{f''}\ar@2[dl]**{}?<(.3);?>(.7)^{\overline{h_2}}|{}="y"&\\
          {}&\ar[r]_{h_1}&\ar[ddl]^(.3){\overline{f''}}&\\
          \ar[ddd]_{k_1}&&&\ar[ddd]^{k_0}{}\\
          {}&\ar[r]^(.9){k_0}\ar[d]_{f''}&\ar[d]^{f'''}\ar@2[dl]**{}?<(.3);?>(.7)^{\overline{k_2}}|{}="z"&\\
          {}&\ar[r]_{k_1}&\ar[dr]^{\overline{f'''}}&\\
          \ar[rrr]_{\overline{f'''}}&&& {}\tria "y";"z" \tria "x";"yz"
        } ,(60,0) 
        \xyboxmatrix"B"{
          {}\ar[rrr]^{\overline{f}}\ar[ddd]_{g_1}\ar[dr]^{\overline{f}}&&&\ar[ddd]^{g_0}\\
          {}&\ar[r]^{g_0}\ar[d]_{f}\POS .[rdddd]!C*+++\frm{.}="xy"
          &\ar[d]^{f'}\ar@2[dl]**{}?<(.3);?>(.7)^{\overline{g_2}}|{}="x"&\\
          {}&\ar[r]_(.1){g_1}&\ar[ddl]^(.3){\overline{f'}}&\\
          \ar[ddd]_{h_1}&&&\ar[ddd]^{h_0}\\
          {}&\ar[r]^(.9){h_0}\ar[d]_{f'}&\ar[d]^{f''}\ar@2[dl]**{}?<(.3);?>(.7)^{\overline{h_2}}|{}="y"&\\
          {}&\ar[r]_{h_1}&\ar[ddl]^(.3){\overline{f''}}&\\
          \ar[ddd]_{k_1}&&&\ar[ddd]^{k_0}{}\tria
          "x";"y"
          \\
          {}&\ar[r]^(.9){k_0}\ar[d]_{f''}&\ar[d]^{f'''}\ar@2[dl]**{}?<(.3);?>(.7)^{\overline{k_2}}|{}="z"&\\
          {}&\ar[r]_{k_1}&\ar[dr]^{\overline{f'''}}&\\
          \ar[rrr]_{\overline{f'''}}&&& {}\tria "xy";"z" 
        } ,(120,0)
        \xyboxmatrix"C"{
          {}\ar[rrr]^{\overline{f}}\ar[ddd]_{g_1}\ar[dr]^{\overline{f}}&&&\ar[ddd]^{g_0}\\
          {}&\ar[r]^{g_0}\ar[d]_{f}\POS.[rdddd]!C*+++\frm{.}="xy"
          &\ar[d]^{f'}\ar@2[dl]**{}?<(.3);?>(.7)^{\overline{g_2}}|{}="x"&\\
          {}&\ar[r]_(.1){g_1}&\ar[ddl]^(.3){\overline{f'}}&\\
          \ar[ddd]_{h_1}&&&\ar[ddd]^{h_0}\\
          {}&\ar[r]^(.9){h_0}\ar[d]_{f'}&\ar[d]^{f''}\ar@2[dl]**{}?<(.3);?>(.7)^{\overline{h_2}}|{}="y"&\\
          {}&\ar[r]_{h_1}&\ar[ddl]^(.3){\overline{f''}}&\\
          \ar[ddd]_{k_1}&&&\ar[ddd]^{k_0}{}
          \\
          {}&\ar[r]^(.9){k_0}\ar[d]_{f''}&\ar[d]^{f'''}\ar@2[dl]**{}?<(.3);?>(.7)^{\overline{k_2}}|{}="z"&\\
          {}&\ar[r]_{k_1}&\ar[dr]^{\overline{f'''}}&\\
          \ar[rrr]_{\overline{f'''}}&&& {}\tria "xy";"z" \tria "y";"x"
        } 
        ,(0,-110) 
        \xyboxmatrix"D"{
          {}\ar[rrr]^{\overline{f}}\ar[ddd]_{g_1}\ar[dr]^{\overline{f}}&&&\ar[ddd]^{g_0}\\
          {}&\ar[r]^{g_0}\ar[d]_{f}
          &\ar[d]^{f'}\ar@2[dl]**{}?<(.3);?>(.7)^{\overline{g_2}}|{}="x"&\\
          {}&\ar[r]_(.1){g_1}&\ar[ddl]^(.3){\overline{f'}}&\\
          \ar[ddd]_{h_1}&&&\ar[ddd]^{h_0}\\
          {}&\ar[r]^(.9){h_0}\ar[d]_{f'}\POS.[rdddd]!C*+++\frm{.}="yz"
          &\ar[d]^{f''}\ar@2[dl]**{}?<(.3);?>(.7)^{\overline{h_2}}|{}="y"&\\
          {}&\ar[r]_{h_1}&\ar[ddl]^(.3){\overline{f''}}&\\
          \ar[ddd]_{k_1}&&&\ar[ddd]^{k_0}{}\\
          {}&\ar[r]^(.9){k_0}\ar[d]_{f''}&\ar[d]^{f'''}\ar@2[dl]**{}?<(.3);?>(.7)^{\overline{k_2}}|{}="z"&\\
          {}&\ar[r]_{k_1}&\ar[dr]^{\overline{f'''}}&\\
          \ar[rrr]_{\overline{f'''}}&&& {}\tria "x";"yz" \tria "z";"y"
        }*{+} ,(60,-110) 
        \xyboxmatrix"E"{
          {}\ar[rrr]^{\overline{f}}\ar[ddd]_{g_1}\ar[dr]^{\overline{f}}&&&\ar[ddd]^{g_0}\\
          {}&\ar[r]^{g_0}\ar[d]_{f}
          &\ar[d]^{f'}\ar@2[dl]**{}?<(.3);?>(.7)^{\overline{g_2}}|{}="x"&\\
          {}&\ar[r]_(.1){g_1}&\ar[ddl]^(.3){\overline{f'}}&\\
          \ar[ddd]_{h_1}&&&\ar[ddd]^{h_0}\\
          {}&\ar[r]^(.9){h_0}\ar[d]_{f'}\POS.[rdddd]!C*+++\frm{.}="yz"
          &\ar[d]^{f''}\ar@2[dl]**{}?<(.3);?>(.7)^{\overline{h_2}}|{}="y"&\\
          {}&\ar[r]_{h_1}&\ar[ddl]^(.3){\overline{f''}}&\\
          \ar[ddd]_{k_1}&&&\ar[ddd]^{k_0}{}\\
          {}&\ar[r]^(.9){k_0}\ar[d]_{f''}&\ar[d]^{f'''}\ar@2[dl]**{}?<(.3);?>(.7)^{\overline{k_2}}|{}="z"&\\
          {}&\ar[r]_{k_1}&\ar[dr]^{\overline{f'''}}&\\
          \ar[rrr]_{\overline{f'''}}&&& {}\tria "yz";"x" \tria "z";"y"
        } 
        ,(120,-110) 
        \xyboxmatrix"F"{
          {}\ar[rrr]^{\overline{f}}\ar[ddd]_{g_1}\ar[dr]^{\overline{f}}&&&\ar[ddd]^{g_0}\\
          {}&\ar[r]^{g_0}\ar[d]_{f}\POS .[rdddd]!C*+++\frm{.}="xy"
          &\ar[d]^{f'}\ar@2[dl]**{}?<(.3);?>(.7)^{\overline{g_2}}|{}="x"&\\
          {}&\ar[r]_(.1){g_1}&\ar[ddl]^(.3){\overline{f'}}&\\
          \ar[ddd]_{h_1}&&&\ar[ddd]^{h_0}\\
          {}&\ar[r]^(.9){h_0}\ar[d]_{f'}&\ar[d]^{f''}\ar@2[dl]**{}?<(.3);?>(.7)^{\overline{h_2}}|{}="y"&\\
          {}&\ar[r]_{h_1}&\ar[ddl]^(.3){\overline{f''}}&\\
          \ar[ddd]_{k_1}&&&\ar[ddd]^{k_0}{}
          \\
          {}&\ar[r]^(.9){k_0}\ar[d]_{f''}&\ar[d]^{f'''}\ar@2[dl]**{}?<(.3);?>(.7)^{\overline{k_2}}|{}="z"&\\
          {}&\ar[r]_{k_1}&\ar[dr]^{\overline{f'''}}&\\
          \ar[rrr]_{\overline{f'''}}&&& {}\tria "z";"xy" \tria "y";"x"
        } \ar@{=} "A";"B"**{}?<;?> \ar@3 "B";"C"**{}?<;?> \ar@3
        "A";"D"**{}?<;?> \ar@3 "D";"E"**{}?<;?> \ar@{=}
        "E";"F"**{}?<;?> \ar@3 "C";"F"**{}?<;?>
      \end{xy}
    \end{equation}
  \end{figure}
\end{prf}

\section{Higher Cells}
\label{sec:hcells}

In order to describe higher transformations between maps of
$\Gray$-categories we construct an internal $\Gray$-category in
$\Gray\Cat_{\fQ^1}$ as a substructure of the iterated path space.

\subsection{Combining Path Spaces and Resolutions}
\label{sec:semidist}

We begin by describing explicitly the action of 
$\pathspc e\from\pathspc{\fQ^1\G}\to\pathspc\G$ as follows:
{\allowdisplaybreaks
\begin{align*}
  \pathspc e\begin{pmatrix}
    \xymatrix{
      {}\ar[r]^{[f_1,\ldots,f_{n_f}]}&{}
    }
  \end{pmatrix}&=
  \begin{pmatrix}
  \xymatrix{
      {}\ar[r]^{f_1\#_0\cdots\#_0f_{n_f}}&{}
    }
  \end{pmatrix}\\
  \pathspc e\begin{pmatrix}
    \xymatrix@+2cm{{}\ar[r]^{[f_1,\ldots,f_{n_f}]}\ar[d]_{[g_{0,0},\ldots,\\g_{0,n_{g_0}}]}&
      {}\ar[d]^{[g_{1,0},\ldots,\\g_{1,n_{g_1}}]}
      \ar@2[dl]**{}?<(.1);?>(.9)|{(g_2;[g_{1,0},\ldots,g_{1,n_{g_1}},\\\,\,f_1,\ldots,f_{n_f}],\\\,
        [f'_1,\ldots,f'_{n_{f'}},\\\,\,g_{0,0},\ldots,g_{0,n_{g_0}}])}\\
      {}\ar[r]_{[f'_1,\ldots,f'_{n_{f'}}]}&{}}
  \end{pmatrix}&=
  \begin{pmatrix}
    \xymatrix@+2cm{{}\ar[r]^{f_1\#_0\cdots\#_0f_{n_f}}\ar[d]_{g_{0,0}\#_0\cdots\\\#_0g_{0,n_{g_0}}}&
      {}\ar[d]^{g_{1,0}\#_0\cdots\\\#_0g_{1,n_{g_1}}}\ar@2[dl]**{}?<(.1);?>(.9)^{g_2}\\
      {}\ar[r]_{f'_1\#_0\cdots\#_0f'_{n_{f'}}}&{}}
  \end{pmatrix}\\  \pathspc e\begin{pmatrix}
    \begin{pmatrix}
      \alpha_3;[g_{1,1},\ldots,g_{1,n_{g_1}},f_{1,1},\ldots,f_{1,n_{f}}],\\
      \,[f'_{1,1},\ldots,f'_{1,n_{f'}},h_{0,1},\ldots,h_{0,n_{h_0}}]
    \end{pmatrix};\\
    (\alpha_1;[g_{0,1},\ldots,g_{0,n_{g_0}}],[h_{0,1},\ldots,h_{0,n_{h_0}}]),\\
    (\alpha_2;[g_{1,1},\ldots,g_{1,n_{g_1}}],[h_{1,1},\ldots,h_{1,n_{h_1}}]),\\
    \begin{pmatrix}
      g_2;[g_{1,1},\ldots,g_{1,n_{g_1}},f_{1,1},\ldots,f_{1,n_{f}}],\\
      \,[f'_{1,1},\ldots,f'_{1,n_{f'}},g_{0,1},\ldots,g_{0,n_{h_0}}]
    \end{pmatrix},\\
    \begin{pmatrix}
      h_2;[h_{1,1},\ldots,h_{1,n_{h_1}},f_{1,1},\ldots,f_{1,n_{f}}],\\
      \,[f'_{1,1},\ldots,f'_{1,n_{f'}},h_{0,1},\ldots,h_{0,n_{h_0}}]
    \end{pmatrix};\\
    [g_{0,1},\ldots,g_{0,n_{g_0}}],[g_{1,1},\ldots,g_{1,n_{g_1}}],\\
    [h_{0,1},\ldots,h_{0,n_{h_0}}],[h_{1,1},\ldots,h_{1,n_{h_1}}],\\
    [f_{1,1},\ldots,f_{1,n_{f}}],[f'_{1,1},\ldots,f'_{1,n_{f'}}]
  \end{pmatrix}&=
  \begin{pmatrix}
    \alpha_3;
    \alpha_1,\alpha_2,g_2,h_2;\\
    g_{0,1}\#_0\cdots\#_0g_{0,n_{g_0}},g_{1,1}\#_0\cdots\#_0g_{1,n_{g_1}},\\
    h_{0,1}\#_0\cdots\#_0h_{0,n_{h_0}},h_{1,1}\#_0\cdots\#_0h_{1,n_{h_1}},\\
    f_{1,1}\#_0\cdots\#_0f_{1,n_{f}},f'_{1,1}\#_0\cdots\#_0f'_{1,n_{f'}}
  \end{pmatrix}\\  
  \pathspc e\begin{pmatrix}
    (\Gamma_1;\alpha_1,\beta_1,[g_{0,1},\ldots,g_{0,n_{g_0}}],[h_{0,1},\ldots,h_{0,n_{h_0}}]),\\
    (\Gamma_2;\alpha_2,\beta_2,[g_{1,1},\ldots,g_{1,n_{g_1}}],[h_{1,1},\ldots,h_{1,n_{h_1}}])
  \end{pmatrix}&=
  \begin{pmatrix}
    \Gamma_1,\Gamma_2
  \end{pmatrix}
\end{align*}
}
where for the 3-cells we used the abbreviated notation of \eqref{eq:3celleq}.

\begin{lem}
  \label{lem:pathepscart}
  The map $\pathspc e\from \pathspc{\fQ^1\G}\to \pathspc \G$ is
  Cartesian with respect $(\_)_1$. 
\end{lem}
\begin{prf}
  $\pathspc e$ is obviously surjective on 0- and 1-cells and 2-locally an isomorphism.
\end{prf}

Let $F\ladj U \from \Cat \to \RGrph$ be the usual adjunction, then
$(\pathspc e)_1\from\pathspc{\fQ^1\G}_1 \to \pathspc \G _1$ has a
splitting $s\from U(\pathspc\G_1)\to U(\pathspc{\fQ^1\G}_1)$ under
$U$ as follows:
\begin{align*}
  s\left(\xymatrix{{}\ar[r]^f&{}}\right)&=\left(\xymatrix{{}\ar[r]^{[f]}&{}}\right)\\
  s\begin{pmatrix}
    \xymatrix{{}\ar[r]^f\ar[d]_{g_0}&
      {}\ar[d]^{g_1}\ar@2[dl]**{}?<(.3);?>(.7)^{g_2}\\{}\ar[r]_{f'}&{}}
  \end{pmatrix}&=
  \begin{pmatrix}
    \xymatrix@+2cm{{}\ar[r]^{[f]}\ar[d]_{[g_0]}&
      {}\ar[d]^{[g_1]}\ar@2[dl]**{}?<(.3);?>(.7)|{(g_2;[g_1,f],[f',g_0])}\\{}\ar[r]_{[f']}&{}}
  \end{pmatrix}
\end{align*}
Obviously in $\RGrph$ we have $U(\pathspc
e_1)s=\id_{U(\pathspc\G_1)}$, taking the transpose $\overline s$ we
get
\begin{equation}\label{eq:sdist1}
  \begin{xy}
    \xyboxmatrix{ FU(\pathspc \G_1)={\fQ^1\pathspc\G}_1\ar[r]^-{\overline
        s}\ar[dr]_{\eps=e_1}&\pathspc{\fQ^1\G}_1\ar[d]^{\pathspc
        e_1}\\
      {}&\pathspc\G_1 }
  \end{xy}\,,
\end{equation} 
since $\pathspc e$ is Cartesian we can lift $\overline s$ through
$(\_)_1$ to obtain
$\psi\from\fQ^1\pathspc\G\to\pathspc{\fQ^1\G}$ satisfying
\begin{equation}
  \label{eq:sdist1a}
  \begin{xy}
    \xyboxmatrix{
      {\fQ^1\pathspc\G}
      \ar[r]^-{\psi_{\G}}
      \ar[dr]_{e_{\pathspc{\G}}}&
      \pathspc{\fQ^1\G}\ar[d]^{\pathspc{e_{\G}}}\\
      {}&\pathspc\G
    }
  \end{xy}\,.
\end{equation}

Let us consider the action of
$\overline{s}\from{\fQ^1\pathspc\G}_1\to\pathspc{\fQ^1\G}_1$. On
0-cells it acts just like $s$, on 1-cells we have the assignment:
\begin{equation*}
  \overline{s}\begin{pmatrix}
    \xymatrix{{}\ar[r]^{f^n}\ar[d]_{g_0^n}&
      {}\ar[d]^{g_1^n}\ar@2[dl]**{}?<(.3);?>(.7)^{g_2^n}\\{}\ar[r]_{f^{n-1}}&{}}\\
    \vdots\\
    \xymatrix{{}\ar[r]^{f^1}\ar[d]_{g_0^1}&
      {}\ar[d]^{g_1^1}\ar@2[dl]**{}?<(.3);?>(.7)^{g_2^1}\\{}\ar[r]_{f^0}&{}}
  \end{pmatrix}=
  \begin{pmatrix}
    \xymatrix@+4cm{{}\ar[r]^{[f^n]}\ar[d]_{[g_0^1,\ldots,g_0^n]}&
      {}\ar[d]^{[g_1^1,\ldots,g_1^n]}\ar@2[dl]**{}?<(.1);?>(.9)|
      {\big((g_2^1\#_0g_0^2\#_0\cdots\#_0g_0^n)\\
        \#_1\cdots\\
        \#_1(g_1^1\#_0\cdots\#_0g_2^i\#_0\cdots\#_0g_0^n)\#_1\cdots\\
        \#_1(g_1^1\#_0\cdots\#_0g_1^{n-1}\#_0g_2^n);\\
        [g_1^1,\ldots,g_1^n,f^n],
        [f^0,g_0^1,\ldots,g_0^n]\big)
      }\\
      {}\ar[r]_{[f^0]}&{}
    }
  \end{pmatrix}
\end{equation*}

\begin{lem}
  \label{lem:sdistnat}
  The family $\psi$ is natural with respect to maps $F\from\G\to\H$.
\end{lem}
\begin{prf}
  Consider the diagram
  \begin{equation*}
    \begin{xy}
      \xyboxmatrix{
        \fQ^1\pathspc{\G}\ar[r]^{\psi_\G}\ar@/^1.5pc/[rr]^{e_{\pathspc\G}}\ar[d]_{\fQ^1\pathspc{F}}&
        \pathspc{\fQ^1\G}\ar[r]^{\pathspc{e_\G}}\ar[d]_{\pathspc{\fQ^1F}}&
        \pathspc{\G}\ar[d]^{\pathspc{F}}\\
        \fQ^1\pathspc{\H}\ar[r]_{\psi_\H}\ar@/_1.5pc/[rr]_{e_{\pathspc\H}}&
        \pathspc{\fQ^1\H}\ar[r]_{\pathspc{e_\H}}& \pathspc{\H} }
    \end{xy}\,,
  \end{equation*}
  since the top and bottom triangles as well as the right hand square
  commute we obtain
  $\pathspc{e_\H}\psi_\H\fQ^1\pathspc{F}=\pathspc{e_\H}\pathspc{\fQ^1F}\psi_\G$.
  Since $\psi_1=\overline{s}$ we need to only verify that
  $\overline{s}_\H(\fQ^1\pathspc{F})_1=(\pathspc{\fQ^1F})_1\overline{s}_\G$,
  but this is immediate from the action of $\pathspc{(\_)}$ and
  $\fQ^1$. Naturality then follows by remark \ref{rem:pmono}.
\end{prf}

It remains to be verified that $\psi$ is compatible with the
co-multiplication $d\from \fQ^1\to\fQ^1\fQ^1$, that is,
\begin{equation}\label{eq:sdist2}
  \xymatrix{
    \fQ^1\pathspc\G\ar[r]^-{d_{\pathspc\G}}\ar[d]_-{\psi_\G}&\fQ^1\fQ^1\pathspc\G\ar[r]^{\fQ^1\psi_\G}&
    \fQ^1\pathspc{\fQ^1\G}\ar[d]^{\psi_{\fQ^1\G}}\\
    \pathspc{\fQ^1\G}\ar[rr]_{\pathspc{d_\G}}&&\pathspc{\fQ^1\fQ^1\G}
  }
\end{equation}
commutes. We will prove this using, again, remark \ref{rem:pmono} with
$\pathspc e$ and the commutativity of the underlying diagram of
categories
\begin{equation*}
  \xymatrix{
    FU(\pathspc\G_1)\ar[r]^-{F\eta U}\ar[d]_-{\overline s}&FUFU(\pathspc\G_1)\ar[r]^{FU\overline s}&
    FU(\pathspc{\fQ^1\G}_1)\ar[d]^{\overline s}\\
    \pathspc{\fQ^1\G}_1\ar[rr]_{\pathspc{d_\G}_1}&&\pathspc{\fQ^1\fQ^1\G}_1\,.
  }
\end{equation*}
But because the upper left object is free over the reflexive graph
$U(\pathspc \G_1)$ it is sufficient to check for generating 0- and
1-cells.   

For 0-cells we compute:
\begin{multline*}
  \pathspc{d_\G}_1\overline{s}\begin{pmatrix}\xymatrix{{}\ar[r]^f&{}}\end{pmatrix}=
  \pathspc{d_\G}_1\begin{pmatrix}\xymatrix{{}\ar[r]^{[f]}&{}}\end{pmatrix}=
  \begin{pmatrix}\xymatrix{{}\ar[r]^{[[f]]}&{}}\end{pmatrix}\\=
  \overline{s}\begin{pmatrix}\xymatrix{{}\ar[r]^{[f]}&{}}\end{pmatrix}=
  \overline{s}(FU\overline{s})\begin{pmatrix}\xymatrix{{}\ar[r]^{f}&{}}\end{pmatrix}=
  \overline{s}(FU\overline{s})(F\eta{}U)\begin{pmatrix}\xymatrix{{}\ar[r]^{f}&{}}\end{pmatrix}
\end{multline*}
And likewise for 1-cells:
\begin{multline*}
  \pathspc{d_\G}_1\overline{s}\begin{pmatrix}
    \xymatrix{ {}\ar[r]^{f}\ar[d]_{g_0}&
      {}\ar[d]^{g_1}\ar@2[dl]**{}?<(.3);?>(.7)^{g_2}\\{}\ar[r]_{f'}&{}
    }\end{pmatrix}=
  \pathspc{d_\G}_1\begin{pmatrix}
    \xymatrix@+1cm{{}\ar[r]^{[f]}\ar[d]_{[g_0]}&
      {}\ar[d]^{[g_1]}\ar@2[dl]**{}?<(.1);?>(.9)|{(g_2;[g_1,f],\\[f',g_0])}\\{}\ar[r]_{[f']}&{}
    }\end{pmatrix}=
  \begin{pmatrix}
    \xymatrix@+1cm{{}\ar[r]^{[[f]]}\ar[d]_{[[g_0]]}&
      {}\ar[d]^{[[g_1]]}\ar@2[dl]**{}?<(.1);?>(.9)|{(g_2;[[g_1],[f]],\\[[f]',[g_0]])}\\{}\ar[r]_{[[f']]}&{}
    }\end{pmatrix}\\=
  \overline{s}\begin{pmatrix}
    \xymatrix@+1cm{{}\ar[r]^{[f]}\ar[d]_{[g_0]}&
      {}\ar[d]^{[g_1]}\ar@2[dl]**{}?<(.1);?>(.9)|{(g_2;[g_1,f],\\[f',g_0])}\\{}\ar[r]_{[f']}&{}
    }\end{pmatrix}=
  \overline{s}(FU\overline{s})\begin{pmatrix}
    \xymatrix{{}\ar[r]^{f}\ar[d]_{g_0}&
      {}\ar[d]^{g_1}\ar@2[dl]**{}?<(.3);?>(.7)^{g_2}\\{}\ar[r]_{f'}&{}
    }\end{pmatrix}=
  \overline{s}(FU\overline{s})(F\eta{}U)\begin{pmatrix}
    \xymatrix{{}\ar[r]^{f}\ar[d]_{g_0}&
      {}\ar[d]^{g_1}\ar@2[dl]**{}?<(.3);?>(.7)^{g_2}\\{}\ar[r]_{f'}&{}
    }\end{pmatrix}
\end{multline*}
Furthermore, we can check that post-composing \eqref{eq:sdist2} with
$\pathspc e$ gives a commuting diagram:
\begin{equation*}
  \begin{xy}
    \xyboxmatrix{
      \fQ^1\pathspc\G\ar[r]^{d_{\pathspc\G}}\ar[dr]_{\fQ^1\pathspc\G}\ar[d]_{\psi_\G}&
      \fQ^1\fQ^1\pathspc\G\ar[d]|{e_{\fQ^1\pathspc\G}}\ar[r]^{\fQ^1\psi_\G}&
      \fQ^1\pathspc{\fQ^1\G}\ar[r]^{\psi_{\fQ^1\G}}\ar[ddr]^{e_{\pathspc{\fQ^1\G}}}&
      \pathspc{\fQ^1\fQ^1\G}\ar[dd]^{\pathspc{e_{\fQ^1\G}}}\\
      \pathspc{\fQ^1\G}\ar[d]_{\pathspc{d_{\G}}}\ar[rrrd]_{\pathspc{\fQ^1\G}}&\fQ^1\pathspc\G\ar[rrd]^{\psi_\G}&&\\
      \pathspc{\fQ^1\fQ^1\G}\ar[rrr]_{\pathspc{e_{\fQ^1\G}}}&&&\pathspc{\fQ^1\G}
    }
  \end{xy}
\end{equation*}
where we use \eqref{eq:sdist1a}, naturality of $\psi$ in lemma \ref{lem:sdistnat},
and the fact that $\fQ^1$ is a comonad. Hence we can cancel
$\pathspc{e}$ and obtain (\ref{eq:sdist2}).
 
So, we have proved the following
\begin{lem}
  \label{lem:psidef}
  There is a natural transformation
  $\psi\from\fQ^1\pathspc{(\_)}\to\pathspc{\fQ^1(\_)}$ satisfying
  properties \eqref{eq:sdist1a} and  \eqref{eq:sdist2}. We call
  $\psi$ a \defterm{semi-distributive law}. \qed 
\end{lem}
\def\GC{\Gray\Cat}
\begin{rem}
  In terms of formal category theory the pair $(\pathspc{(\_)}, \psi)$
  is an endomorphism of the comonad $(\fQ^1,d,e)$, that is, 
  \begin{equation*}
    \begin{xy}
      \xyboxmatrix{
        \GC\ar[r]^{\pathspc{(\_)}}\ar[d]|{\fQ^1}="x"\ar@/_2pc/[d]_\id="y"&\GC\ar[d]^{\fQ^1}
        \ar@2[dl]**{}?(.3);?(.7)_\psi\\
        \GC\ar[r]_{\pathspc{(\_)}}&\GC     \ar@2"x";"y"**{}?(.3);?(.7)_e
      }
    \end{xy}
    =
    \begin{xy}
      \xyboxmatrix{
        \GC\ar[r]^{\pathspc{(\_)}}\ar[d]_{\id}&\GC\ar[d]|{\id}="y"\ar@/^2pc/[d]^{\fQ^1}="x"\\
        \GC\ar[r]_{\pathspc{(\_)}}&\GC  \ar@2"x";"y"**{}?(.3);?(.7)_e} 
    \end{xy} \end{equation*}
  and
  \begin{equation*}
    \begin{xy}
      \xyboxmatrix{
        \GC\ar[r]^{\pathspc{(\_)}}\ar[dd]|{\fQ^1}="x"\ar@/_4pc/[dd]_{\fQ^1\fQ^1}="y"&\GC\ar[dd]^{\fQ^1}
        \ar@2[ddl]**{}?(.3);?(.7)_\psi\\
        {}&{}\\
        \GC\ar[r]_{\pathspc{(\_)}}&\GC\ar@2"x";"y"**{}?(.3);?(.7)_d \\
      }
    \end{xy}
    =
    \begin{xy}
      \xyboxmatrix{
        \GC\ar[r]^{\pathspc{(\_)}}\ar[d]_{\fQ^1}&\GC\ar[d]|{\fQ^1}
        \ar@/^4pc/[dd]^{\fQ^1}="x"\ar@2[dl]**{}?(.3);?(.7)_\psi\\
        \GC\ar[r]|{\pathspc{(\_)}}\ar[d]_{\fQ^1}&\GC\POS[]="y"\ar[d]|{\fQ^1}\ar@2[dl]**{}?(.3);?(.7)_\psi\\
        \GC\ar[r]_{\pathspc{(\_)}}&\GC\ar@2"x";"y"**{}?(.3);?>(.8)_d\\
      }
    \end{xy}\,.
  \end{equation*}
\end{rem}

\begin{lem}
  \label{lem:pathspcprop1}
  The functor $\pathspc{(\_)}$ extends canonically to an endofunctor $\qopth$
  of $\Gray\Cat_{\fQ^1}$ by
  \begin{equation*}
    \qopth\begin{pmatrix}
      \xymatrix{
        \G\ar[r]|*\dir{/}^{f}&\H
      }
    \end{pmatrix}
    =\begin{pmatrix}
      \xymatrix{
        \fQ^1\pathspc\G\ar[r]^-\psi&\pathspc{\fQ^1\G}\ar[r]^{\pathspc{f}}&\pathspc\H
      }
    \end{pmatrix}
    =\begin{pmatrix}
      \xymatrix{
        \pathspc{\G}\ar[r]|*\dir{/}^{\qopth(f)}&\pathspc\H
      } 
    \end{pmatrix}.
  \end{equation*}
  Furthermore, it preserves strictness of maps.
\end{lem}
\begin{prf}
  We use the properties of $\psi$ to check that this assignment is
  functorial. Given two maps $f\from\G\laxto\H$ and $g\from\H\laxto\K$
  we compare $\qopth(g)\qopth(f)$ at the top and $\qopth(gf)$ at the bottom:
  \begin{equation*}
    \xymatrix{
      \fQ^1\pathspc\G\ar[r]^-d\ar[dr]_-{\psi}&\fQ^1\fQ^1\pathspc\G\ar[r]^-{\fQ^1\psi}&
      \fQ^1\pathspc{\fQ^1\G}\ar[r]^{\fQ^1\pathspc{f}}\ar[d]_-{\psi}&\fQ^1\pathspc\H\ar[r]^\psi&\pathspc{\fQ^1\H}\ar[r]^{\pathspc{g}}&
      \pathspc\K\\
      {}&\pathspc{\fQ^1\G}\ar[r]_-{\pathspc{d}}&\pathspc{\fQ^1\fQ^1\G}\ar[rru]_{\pathspc{\fQ^1f}}&&&
    }\,.
  \end{equation*}
  The naturality of $\psi$ and (\ref{eq:sdist2}) make sure they are
  equal. Preservation of units is exactly \eqref{eq:sdist1a}.

  We remember that a strict map in $\Gray\Cat_{\fQ^1}$ is given by
  $fe_\G$ where $f\from\G\to\H$ is from $\Gray\Cat$ and $e$ is the
  co-unit of $\fQ^1$. Then by (\ref{eq:sdist1a}) we get
  \begin{equation*}
    \qopth({fe_\G})=\pathspc{f}\pathspc{e_\G}\psi_\G=\pathspc{f}e_{\pathspc{\G}}\,,
  \end{equation*} Meaning that $\qopth$ acts on strict maps like
  $\pathspc{(\_)}$, in particular, it takes identities to identities.
\end{prf}

\begin{lem}
  \label{lem:plimpres}
  The functor $\qopth\from\GC_{\fQ^1}\to\GC_{\fQ^1}$ preserves limits
  of diagrams of strict maps.
\end{lem}
\begin{prf}
  Finally, by lemma \ref{lem:pathlim} the restriction $\pathspc{(\_)}$
  of $\qopth$ to $\Gray\Cat$ preserves limits:  
  Let $p_i\from\lim\{\H_i,b_k\}\to\H_i$ be a limit cone in $\GC$, let
  $f_i\from\G\laxto\pathspc\H_i$ be a cone in $\GC_{\fQ^1}$.
  \begin{equation*}
    \begin{xy}
      \xyboxmatrix{
        \fQ^1\G\ar@{-->}[r]^-{\pair{f_i}}\ar[dr]_{f_i}&\pathspc{\lim\{\H_i,b_k\}}\ar[d]^{\pathspc{p_i}}\\
        {}&\pathspc{\H_i}
      }
    \end{xy}
  \end{equation*}
  $\pathspc{p_i}$ is a limit cone, hence there is the unique weak map
  $\pair{f_i}\from\G\laxto\pathspc{\lim\{\H_i,b_k\}}$.
\end{prf}

\begin{lem}
  \label{lem:pindpres}
  The functor $\qopth\from\GC_{\fQ^1}\to\GC_{\fQ^1}$ preserves induced
  maps of limits of strict diagrams, that is, $\qopth(\dot\lim
  f_i)=\dot\lim(\qopth f_i)$. 
\end{lem}

\begin{prf}
  Consider
  \begin{equation*}
    \begin{xy}
      \xyboxmatrix{
        \fQ^1\pathspc{\lim\{\G_i,a_k\}}\ar[rr]^-{\psi}\ar[dd]_{\fQ^1\pathspc{p_i}}
        \ar[dr]^{\fQ^1\pair{\pathspc{p_i}}}&{}&
        \pathspc{\fQ^1\lim\{\G_i,a_k\}}\ar[r]^-{\pathspc{\dot\lim{f_i}}}
        \ar[dd]|!{"2,2";"1,4"}\hole^{\pathspc{\fQ^1p_i}}&
        \pathspc{\lim\{\H_i,b_k\}}\ar[dd]^{\pathspc{p'_i}}\\
        {}&\fQ^1\lim\{\pathspc{\G_i},\pathspc{a_k}\}\ar[urr]_(.7){\dot\lim\qopth{f_i}}\ar[dl]^{\fQ^1p''_i}&{}&{}\\
        \fQ^1\pathspc{\G_i}\ar[rr]_{\psi}&{}&\pathspc{\fQ^1\G_i}\ar[r]_{\pathspc{f_i}}&\pathspc{\H_i}
      }
    \end{xy}
  \end{equation*} using the conventions of remark \ref{rem:indmap}. Also, note
  that $\pathspc{\dot\lim{f_i}}\psi=\qopth(\dot\lim{f_i})$ by
  definition. $\dot\lim{f_i}$ is the induced arrow for the source
  $f_i(\fQ^1p_i)$, $\dot\lim\qopth{f_i}$ is the induced arrow for $\qopth(f_i)\fQ^1(p''_i)$.
  Since
  \begin{equation*}
    \pathspc{p'_i}(\dot\lim\qopth{f_i})\fQ^1\pair{\pathspc{p_i}}=\pathspc{p'_i}\pathspc{\dot\lim{f_i}}\psi
  \end{equation*}
  and $\pathspc{p'_i}$ is a limit cone we obtain
  \begin{equation*}
    (\dot\lim\qopth{f_i})\fQ^1\pair{\pathspc{p_i}}=\pathspc{\dot\lim{f_i}}\psi\,.
  \end{equation*} 
\end{prf}

If the limit is, for example, a product we may now say
that
\begin{equation}
  \label{eq:pprodindpres}
  \qopth(f\dot\times{}g)=\qopth{f}\dot\times\qopth{g}\,.
\end{equation}
From now on however we shall use $\times$ for the product of arrows in
$\GC_{\fQ^1}$.

\begin{lem}
  \label{lem:faceweaknat}
  The face maps are natural with respect to weak maps, that is
  \begin{equation}
    \label{eq:faceweaknat}
    \begin{xy}
      \xyboxmatrix{
        \pathspc\G\ar[d]_{\qopth{f}}|*-\dir{/}\ar@<.5ex>[r]^{d_0}\ar@<-.5ex>[r]_{d_1}&\G\ar[d]^f|*-\dir{/}\\
        \pathspc\H\ar@<.5ex>[r]^{d_0}\ar@<-.5ex>[r]_{d_1}&\H
      }
    \end{xy}
  \end{equation}
  commutes.
\end{lem}
\begin{prf}
  We write (\ref{eq:faceweaknat}) in terms of its underlying maps:
  \begin{equation}
    \label{eq:faceweaknatunder}
    \begin{xy}
      \xyboxmatrix{
        \fQ^1\pathspc\G\ar[r]^-d\ar[d]_-d&\fQ^1\fQ^1\pathspc\G\ar[r]^{\fQ^1e}&
        \fQ^1\pathspc\G\ar@<.5ex>[r]^-{\fQ^1d_0}\ar@<-.5ex>[r]_-{\fQ^1d_1}\ar[d]_\psi&
        \fQ^1\G\ar@{=}[d]^{}\\
        \fQ^1\fQ^1\pathspc\G\ar[r]_{\fQ^1\psi}\ar[rru]_{e}&\fQ^1\pathspc{\fQ^1\G}\ar[r]^-{e}\ar[d]_{\fQ^1\pathspc{f}}&
        \pathspc{\fQ^1\G}\ar[d]_{\pathspc{f}}\ar@<.5ex>[r]^-{d_0}\ar@<-.5ex>[r]_-{d_1}&\fQ^1\G\ar[d]^{f}\\
        {}&\fQ^1\pathspc\H\ar[r]_-{e}&\pathspc\H\ar@<.5ex>[r]^-{d_0}\ar@<-.5ex>[r]_-{d_1}&\H
      }
    \end{xy}\,,
  \end{equation}
  that is, (\ref{eq:faceweaknat}) commuting is equivalent to the outer
  frame in (\ref{eq:faceweaknatunder}) commuting. All parts are given
  by naturality and the co-unit laws of $\fQ^1$, except the upper right
  square.

  We use remark \ref{rem:pmono} to conclude $d_0\psi=\fQ^1d_0$ and
  $d_1\psi=\fQ^1d_1$:    
  By naturality and semi-distributivity we get
  $ed_0\psi=d_0\pathspc{e}\psi=d_0e=e\fQ^1d_0$, furthermore,
  $(d_0\psi)_1=(\fQ^1d_0)_1$ is immediate from the definition of
  $\psi$. The map $d_1$ is obviously treated in the same way.
\end{prf}

\begin{lem}
  \label{lem:idweaknat}
  The degeneracy maps of the path space are natural with respect to
  weak maps: 
  \begin{equation*}
    \begin{xy}
      \xyboxmatrix{
        \G\ar[d]_f|*-\dir{/}\ar[r]^i&\pathspc\G\ar[d]^{\qopth{f}}|*-\dir{/}\\
        \H\ar[r]_i&\pathspc\H
      }
    \end{xy}\,.
  \end{equation*}
\end{lem}

\begin{prf}
  Consider
  \begin{equation*}
    \begin{xy}
      \xyboxmatrix{
        \fQ^1\G\ar[r]^d\ar[d]_d&\fQ^1\fQ^1\G\ar[r]^{\fQ^1e}&
        \fQ^1\G\ar@{=}[d]\ar[r]^{\fQ^1i}&\fQ^1\pathspc\G\ar[d]^{\psi}\\
        \fQ^1\fQ^1\G\ar[rr]^e\ar[d]_{\fQ^1f}&{}&\fQ^1\G\ar[d]^f\ar[r]^i&\pathspc{\fQ^1\G}\ar[d]^{\pathspc{f}}\\
        \fQ^1\H\ar[rr]_e&{}&\H\ar[r]_i&\pathspc\H
      }
    \end{xy}\,.
  \end{equation*}
  We conclude that then top right square commutes by computing
  $\pathspc{e}i=ie=e\fQ^1i=\pathspc{e}\psi\fQ^1i$ and
  checking that $(\psi\fQ^1i)_1=i_1$ and again applying
  remark \ref{rem:pmono} together with lemma \ref{lem:pathepscart}. 
\end{prf}

The functor $\qopth$ can also be applied to $\fQ^1$-graph maps by setting
$\qopth'=(\qopth\tilde{G})^{\vee}$; see lemma \ref{lem:pseudofunchar} for the
notation.  For the sake of completeness we describe briefly the effect
of $\qopth'$ at the level of 1-cells as well as its 2-co-cycle. Let
$G\from\G\to\H$ be a $\fQ^1$-graph map. We take a 1-cell
$g\from{}f\to{}f'$ from $\pathspc\G$ and calculate: 
\begin{multline}
  \label{eq:q1grphmappathspc}
  (\qopth'G)(g)={\left(\pathspc{\tilde{G}}\psi\right)}^{\vee}(g)=
  \pathspc{\tilde{G}}\psi\begin{bmatrix}
    \xymatrix{{}\ar[r]^{f}\ar[d]_{g_0}&
      {}\ar[d]^{g_1}\ar@2[dl]**{}?<(.1);?>(.9)|
      {g_2}\\
      {}\ar[r]_{f'}&{}
    }
  \end{bmatrix}
  \\
  =\begin{pmatrix}
    \xymatrix@+1cm{{}\ar[r]^{\tilde{G}[f]}\ar[d]_{\tilde{G}[g_0]}&
      {}\ar[d]^{\tilde{G}[g_1]}\ar@2[dl]**{}?<(.1);?>(.9)|
      {\tilde{G}(g_2;
        [g_1,f],\\
        [f,g_0])
      }\\
      {}\ar[r]_{\tilde{G}[f']}&{}
    }
  \end{pmatrix}
  =\begin{pmatrix}
    \xymatrix@+1cm{{}\ar[r]^{Gf}\ar[d]_{Gg_0}&
      {}\ar[d]^{Gg_1}\ar@2[dl]**{}?<(.1);?>(.9)|
      {\overline{G^2_{f',g_0}}\\\#_1Gg_2\\\#_1G^2_{g_1,f}}\\
      {}\ar[r]_{Gf'}&{}
    }
  \end{pmatrix}
\end{multline}
Taking two composable 1-cells $g\from{}f\to{}f'$ and
$h\from{}f'\to{}f''$ of $\pathspc\G$ we get a  2-cocycle with
components as shown in (\ref{eq:q1grphmappathspc2coc}), where in the
end the $\tilde{G}\kappa_{\ldots}$ are iterated 2-cocycles of $G$.
\begin{figure}
  \centering
  \begin{multline}
    \label{eq:q1grphmappathspc2coc}
    ((\qopth'G)^\vee)^2_{h,g}=((\pathspc{\tilde{G}}\psi)^\vee)^2_{h,g}=\pathspc{\tilde{G}}\psi(\kappa_{h,g})\\
    =\pathspc{\tilde{G}}\psi\left(
      \begin{matrix}
        \begin{xy}
          \xyboxmatrix"A"@+1cm{
            {}\ar[r]^{f}\ar[d]|{h_0\#_0g_0}="x"\ar@/_3pc/[d]_{h_0\#_0g_0}="y"&
            {}\ar[d]^{h_1\#_0g_1}\ar@2[dl]**{}?<(.1);?>(.9)|
            {(h_2\#_0g_0)\\
              \#_1(h_1\#_0g_2)}\\
            {}\ar[r]_{f''}&{} \POS \ar@2 "x"; "y" **{} ?<(.1); ?>(.9)
            ^{\id} } \POS - (0,40)
          \xyboxmatrix"B"@+1cm{{}\ar[r]^{f}\ar[d]_{h_0\#_0g_0}&
            {}\ar[d]|{h_1\#_0g_1}="x"\ar@/^3pc/[d]^{h_1\#_0g_1}="y"\ar@2[dl]**{}?<(.1);?>(.9)|
            {(h_2\#_0g_0)\\
              \#_1(h_1\#_0g_2)}\\
            {}\ar[r]_{f''}&{} \POS \ar@2 "y"; "x" **{} ?<(.1); ?>(.9)
            ^{\id} } \POS \ar@3 "A";"B" **{}
          ?<(0);?>(1)^{\id_{(h_2\#_0g_0)\#_1(h_1\#_0g_2)}}
        \end{xy}
      \end{matrix}
      ;
      \begin{bmatrix}
        \begin{xy}
          \xyboxmatrix{{}\ar[r]^{f'}\ar[d]_{h_0}&
            {}\ar[d]^{h_1}\ar@2[dl]**{}?<(.1);?>(.9)|
            {h_2}\\
            {}\ar[r]_{f''}&{} }
        \end{xy}
        ,
        \begin{xy}
          \xyboxmatrix{{}\ar[r]^{f}\ar[d]_{g_0}&
            {}\ar[d]^{g_1}\ar@2[dl]**{}?<(.1);?>(.9)|
            {g_2}\\
            {}\ar[r]_{f'}&{} }
        \end{xy}
      \end{bmatrix},
      \begin{bmatrix}
        \begin{xy}
          \xyboxmatrix@+1cm{{}\ar[r]^{f}\ar[d]_{h_0\#_0g_0}&
            {}\ar[d]^{h_1\#_0g_1}\ar@2[dl]**{}?<(.1);?>(.9)|
            {(h_2\#_0g_0)\\
              \#_1(h_1\#_0g_2)}\\
            {}\ar[r]_{f''}&{} }
        \end{xy}
      \end{bmatrix}
    \right)\\
    =\pathspc{\tilde{G}}\left(
      \begin{matrix}
        \begin{xy}
          \xyboxmatrix"A"@+2cm{
            {}\ar[r]^{[f]}\ar[d]|{[h_0,g_0]}="x"\ar@/_3pc/[d]_{[h_0\#_0g_0]}="y"&
            {}\ar[d]^(.7){[h_1,g_1]}\ar@2[dl]**{}?<(.1);?>(.9)|
            {((h_2\#_0g_0)\\\#_1(h_1\#_0g_2);\\[h_1,g_1,f],\\[f'',h_0,g_0])}\\
            {}\ar[r]_{[f'']}&{} \POS \ar@2 "x"; "y" **{} ?<(.1);
            ?>(.9) ^{\kappa_{h_0,g_0}} } \POS + (80,0)
          \xyboxmatrix"B"@+2cm{{}\ar[r]^{[f]}\ar[d]_(.7){[h_0\#_0g_0]}&
            {}\ar[d]|{[h_1\#_0g_1]}="x"\ar@/^3pc/[d]^{[h_1,g_1]}="y"\ar@2[dl]**{}?<(.1);?>(.9)|
            {((h_2\#_0g_0)\\\#_1(h_1\#_0g_2);\\[h_1\#_0g_1,f],\\[f'',h_0\#_0g_0])}\\
            {}\ar[r]_{[f'']}&{} \POS \ar@2 "y"; "x" **{} ?<(.1);
            ?>(.9) ^{\kappa_{h_1,g_1}} } \POS \ar@3 "A";"B" **{}
          ?<(0);?>(1)^{(\id_{(h_2\#_0g_0)\#_1(h_1\#_0g_2)};\\
            (h_2\#_0g_0)\#_1(h_1\#_0g_2),\\
            (h_2\#_0g_0)\#_1(h_1\#_0g_2);\\[h_1,g_1,f],[f'',h_0\#_0g_0])}
        \end{xy}
      \end{matrix}
    \right)\\
    =\left(
      \begin{matrix}
        \begin{xy}
          \xyboxmatrix"A"@+2cm{
            {}\ar[r]^{\tilde{G}[f]}\ar[d]|{\tilde{G}[h_0,g_0]}="x"\ar@/_3pc/[d]_{\tilde{G}[h_0\#_0g_0]}="y"&
            {}\ar[d]^(.7){\tilde{G}[h_1,g_1]}\ar@2[dl]**{}?<(.1);?>(.9)|
            {\tilde{G}((h_2\#_0g_0)\\\#_1(h_1\#_0g_2);\\[h_1,g_1,f],\\[f'',h_0,g_0])}\\
            {}\ar[r]_{\tilde{G}[f'']}&{} \POS \ar@2 "x"; "y" **{}
            ?<(.1); ?>(.9) ^{\tilde{G}\kappa_{h_0,g_0}} } \POS +
          (80,0)
          \xyboxmatrix"B"@+2cm{{}\ar[r]^{\tilde{G}[f]}\ar[d]_(.7){\tilde{G}[h_0\#_0g_0]}&
            {}\ar[d]|{\tilde{G}[h_1\#_0g_1]}="x"\ar@/^3pc/[d]^{\tilde{G}[h_1,g_1]}="y"\ar@2[dl]**{}?<(.1);?>(.9)|
            {\tilde{G}((h_2\#_0g_0)\\\#_1(h_1\#_0g_2);\\[h_1\#_0g_1,f],\\[f'',h_0\#_0g_0])}\\
            {}\ar[r]_{\tilde{G}[f'']}&{} \POS \ar@2 "y"; "x" **{}
            ?<(.1); ?>(.9) ^{\tilde{G}\kappa_{h_1,g_1}} } \POS \ar@3
          "A";"B" **{}
          ?<(0);?>(1)^{\tilde{G}(\id_{(h_2\#_0g_0)\#_1(h_1\#_0g_2)};\\
            (h_2\#_0g_0)\#_1(h_1\#_0g_2),\\
            (h_2\#_0g_0)\#_1(h_1\#_0g_2);\\[h_1,g_1,f],[f'',h_0\#_0g_0])}
        \end{xy}
      \end{matrix}
    \right)\\
    =\left(
      \begin{matrix}
        \begin{xy}
          \xyboxmatrix"A"@+2cm{
            {}\ar[r]^{Gf}\ar[d]|{Gh_0\#_0Gg_0}="x"\ar@/_3pc/[d]_{G(h_0\#_0g_0)}="y"&
            {}\ar[d]^(.7){Gh_1\#_0Gg_1}\ar@2[dl]**{}?<(.1);?>(.9)|
            {\overline{\tilde{G}\kappa_{f'',h_0,g_0}}\\
              \#_1G((h_2\#_0g_0)\\\#_1(h_1\#_0g_2))\\
              \#_1\tilde{G}\kappa_{h_1,g_1,f}}\\
            {}\ar[r]_{Gf''}&{} \POS \ar@2 "x"; "y" **{} ?<(.1); ?>(.9)
            ^{G^2_{h_0,g_0}} } \POS + (80,0)
          \xyboxmatrix"B"@+2cm{{}\ar[r]^{Gf}\ar[d]_(.7){G(h_0\#_0g_0)}&
            {}\ar[d]|{G(h_1\#_0g_1)}="x"\ar@/^3pc/[d]^{Gh_1\#_0Gg_1}="y"\ar@2[dl]**{}?<(.1);?>(.9)|
            {G((h_2\#_0g_0)\\\#_1(h_1\#_0g_2);\\[h_1\#_0g_1,f],\\[f'',h_0\#_0g_0])}\\
            {}\ar[r]_{\tilde{G}[f'']}&{} \POS \ar@2 "y"; "x" **{}
            ?<(.1); ?>(.9) ^{\tilde{G}\kappa_{h_1,g_1}} } \POS \ar@3
          "A";"B" **{} ?<(0);?>(1)^{\id}
        \end{xy}
      \end{matrix}
    \right)
  \end{multline}
\end{figure}

\subsection{Iterating the Path Space Construction}
\label{sec:pathit}

\begin{rem}
  As a consequence of lemma \ref{lem:faceweaknat}, lemma
  \ref{lem:idweaknat}, and lemma \ref{lem:multstrictnat}
  the maps $i, d_0, d_1$ and $m$ for all $\Gray$-categories $\H$
  constitute natural transformations with respect to strict maps. 
\end{rem}

For reference, this means that for all $f\from\H\to\K$ the following diagram commutes sequentially:
\begin{equation*}
  \xymatrix{
    \pathspc\H\times_{\H}\pathspc\H\ar[r]|-*\dir{/}^-m\ar[d]_{\pathspc{f}\times\pathspc{f}}&
    \pathspc\H\ar@<-1ex>[r]_{{d_0}}\ar@<1ex>[r]^{{d_1}}\ar[d]_{\pathspc{f}}&\H\ar[d]^{{f}}\ar[l]|{i}\\
    \pathspc\K\times_{\K}\pathspc\K\ar[r]|-*\dir{/}_-m&\pathspc\K\ar@<-1ex>[r]_{{d_0}}\ar@<1ex>[r]^{{d_1}}&\K\ar[l]|{i}
  }
\end{equation*}

Iterating the arrow construction yields an internal cubical set, so it
allows us to talk about higher cells in the internal language of
$\Gray\Cat$. But since we want to construct an internal
$\Gray$-category we need to restrict to cubical cells with certain
degeneracies.  The general recipe beyond the construction in section
\ref{sec:pathspc} is to apply $\pathspc{(\_)}$ and squash the excess
faces given by $\pathspc{d_{0,1}}$ so that the only non-trivial faces
of each cubical element are the ones given by $d_{0,1}$.

This general procedure will canonically yield an internal reflexive
$n$-graph, we will furthermore have to provide the operations in each
degree to actuallty obtain a $\Gray$-category. We carry out this
construction for the degrees 2 and 3 in sections \ref{sec:twopaths} and
\ref{sec:threepaths}.

\subsubsection{2-Paths}
\label{sec:twopaths}

We construct the space of 2-paths $\overline{\overline{\H}}$ over
$\pathspc{\H}$ and  give the vertical composition of 2-paths and
their whiskers by 1-paths. 

The 0-cells in $\dblpathspc\H$ are squares, and
we want to filter out those square that are actually bigons, that is,
have identity arrows as left and right sides. That is exactly what we
get by forming the double pullback on the left:
\begin{equation}
  \label{eq:dbldef}
  \xymatrix{
    \overline{\overline{\H}}\ar[r]^{j}\ar@<-.5ex>[d]_{\overline{d_0}}\ar@<.5ex>[d]^{\overline{d_1}}&
    \dblpathspc\H\ar@<-.5ex>[d]_{\pathspc{d_0}}\ar@<.5ex>[d]^{\pathspc{d_1}}\ar@<-.5ex>[r]_{{d_0}}\ar@<.5ex>[r]^{{d_1}}& 
    \pathspc\H\ar@<-.5ex>[d]_{{d_0}}\ar@<.5ex>[d]^{{d_1}}\\
    \H\ar[r]_i& \pathspc\H\ar@<-.5ex>[r]_{{d_0}}\ar@<.5ex>[r]^{{d_1}}&\H
  }
\end{equation}
where $\overline{\overline{\H}}$ is the intersection of the
pullbacks of $d_{0}$ and $d_{1}$ along $i$. Let $d_{0}^j=d_{0}j$
and $d_{1}^j=d_{1}j$.
\begin{lem}
  The diagram
  \begin{equation}
    \label{eq:dblglob}
    \xymatrix{
      \overline{\overline{\H}}\ar@<-.5ex>[r]_{{d_0^j}}\ar@<.5ex>[r]^{{d_1^j}}&
      \pathspc\H\ar@<-.5ex>[r]_{{d_0}}\ar@<.5ex>[r]^{{d_1}}&
      \H
    }
  \end{equation}
  is a globular object, i.~e.\@ $d_0d_0^j=d_0d_1^j$ 
  and $d_1d_0^j=d_1d_1^j$.
\end{lem}
\begin{prf}
  Using the naturality of $d_{0}$ and $d_{1}$ we calculate:
  \begin{equation*}
    d_0d_0^j=d_0d_0j=d_0\pathspc{d_0}j=d_0i\overline{d_0}=d_1i\overline{d_0}=d_1\pathspc{d_0}j=d_0d_1j=d_1d_0^j\,,
  \end{equation*}
  and similarly for $d_1$.
\end{prf}

To get a unit for $\overline{\overline{\H}}$, that is, an identity
2-paths for 1-paths, we consider the following
diagram:
\begin{equation*}
  \xymatrix{
    \pathspc\H\ar@{..>}[rd]|{\overline{i}}\ar@<-.5ex>[ddr]_{{d_0}}\ar@<.5ex>[ddr]^{{d_1}}\ar[drr]^{{i}}&&&\\
    &\overline{\overline{\H}}\ar[r]|{j}\ar@<-.5ex>[d]_{\overline{d_0}}\ar@<.5ex>[d]^{\overline{d_1}}&
    \dblpathspc\H\ar@<-.5ex>[d]_{\pathspc{d_0}}\ar@<.5ex>[d]^{\pathspc{d_1}}\ar@<-.5ex>[r]_{{d_0}}\ar@<.5ex>[r]^{{d_1}}& 
    \pathspc\H\ar@<-.5ex>[d]_{{d_0}}\ar@<.5ex>[d]^{{d_1}}\\
    &\H\ar[r]_i& \pathspc\H\ar@<-.5ex>[r]_{{d_0}}\ar@<.5ex>[r]^{{d_1}}&\H
  }
\end{equation*}
The upper left span is a compatible source by the naturality of $i$.
The induced arrow $\overline{i}$ is a joint section of ${d_{0}^j}$ and
${d_{1}^j}$. Hence we get:

\begin{lem}
  The diagram 
  \begin{equation}
    \label{eq:dblref}
    \xymatrix{
      \overline{\overline{\H}}\ar@<-1ex>[r]_{{d_0^j}}\ar@<1ex>[r]^{{d_1^j}}&\pathspc\H\ar[l]|{\overline{i}}
    }
  \end{equation}
  is a reflexive graph.\qed
\end{lem}

\begin{lem}
  \label{lem:dblbarfunct}
  The mapping $\overline{\overline{(-)}}$ extends to a sub-functor of
  $\dblpathspc{(-)}\from\Gray\Cat\to\Gray\Cat$ with natural embedding $j$.
\end{lem}
\begin{prf}
  For each $\H$ the map $j$ is a monomorphism by construction and
  $\overline{\overline{(-)}}$ extends to morphisms by the universal
    property.
\end{prf}

\begin{lem}
  \label{lem:pbmult}
  There is a multiplication
  \begin{equation*}
    \xymatrix{
      \overline{\overline{\H}}\times_{d_0^j, d_1^j}\overline{\overline{\H}}\ar[r]^-{\overline{m}}|-*\dir{/}&
      \overline{\overline{\H}}
    }
  \end{equation*}
  with
  \begin{align}
    \label{eq:vercompstcond}
    d^j_0\overline{m}&=d^j_0p_1\\
    d^j_1\overline{m}&=d^j_1p_0 \notag 
  \end{align}
  uniquely induced by $m_{\pathspc\H}$. 
\end{lem}
\begin{prf}
  All we need to show is that $m(j\times j)$ factors through
  $j$, that is, show that the two outer
  rectangles commute:
  \begin{equation}
    \label{eq:vertcompdef}
    \xymatrix{
      \overline{\overline{\H}}\times_{d_0^j,d_1^j}\overline{\overline{\H}}
      \ar[r]^{j\times j}\ar@<-.5ex>[d]_{{p_0}}\ar@<.5ex>[d]^{{p_1}}
      \ar@<-.5ex>@/_3pc/[dd]_{{d'_0}}\ar@<.5ex>@/_3pc/[dd]^{{d'_1}}
      \ar@{..>}@/^3pc/[d]^{\overline m}|-*\dir{/}&
      \dblpathspc\H\times_{d_0,d_1}\dblpathspc\H\ar[d]^{m}|-*\dir{/}\\
      \overline{\overline{\H}}\ar[r]^{j}\ar@<-.5ex>[d]_{{\overline
          d_0}}\ar@<.5ex>[d]^{{\overline d_1}}&
      \dblpathspc\H\ar@<-.5ex>[d]_{{\pathspc
          d_0}}\ar@<.5ex>[d]^{{\pathspc d_1}}\\
      \H\ar[r]_{i}&\pathspc\H
    }
  \end{equation}
  that is, we shall verify that
  \begin{align*}
    \pathspc{d_0}m(j\times{}j)&=id'_0\\
    \pathspc{d_1}m(j\times{}j)&=id'_1
  \end{align*}
  in order to obain $\overline{m}$ as a universally induced arrow.

  First we prove that $\overline d_0p_0=\overline d_0 p_1$:
  \begin{equation}
    \label{eq:local1}
    \overline d_0p_0=d_0i\overline d_0p_0=d_0\pathspc d_0j p_0=
    d_0d_0jp_0=d_0d_0^jp_0=d_0d_1^jp_1=d_0d_0^jp_1=\overline d_0p_1
  \end{equation}
  which holds by \eqref{eq:dblref}, \eqref{eq:dblglob} and
  \eqref{eq:dbldef}. Similarly $\overline d_1p_0=\overline
  d_1p_1$. Thus we may define $d'_0=\overline d_0p_0$ and
  $d'_1=\overline d_1p_0$. Note that $j\times j$ is universally
  induced by $d_0 j p_0=d_1 j p_1$. 
  
  Furthermore, we need that $(i\overline d_0\times i\overline
  d_0)=(i,i)d'_0$ and $(i\overline d_1\times i\overline
  d_1)=(i,i)d'_1$. Consider
  \begin{equation*}
    \xymatrix{
      \overline{\overline{\H}}\times_{d_0^j,d_1^j}\overline{\overline{\H}}\ar[rrr]^{p_1}\ar[ddd]_{p_0}\ar[dr]^{d'_0}
      \ar@/_3pc/[drdr]_{(i\overline d_0\times i\overline d_0)}&&&
      \overline{\overline{\H}}\ar[d]^{\overline d_0}\\
      {}&\H\ar[rr]^{}\ar[dd]_{}\ar[dr]^{(i,i)}&&\H\ar[d]^{i}\\
      {}&&\pathspc\H\times_{d_0, d_1}\pathspc\H\ar[r]\ar[d]&\pathspc\H\ar[d]^{d_1}\\
      \overline{\overline{\H}}\ar[r]_{\overline
        d_0}&\H\ar[r]_i&\pathspc\H\ar[r]_{d_0}&\H
    }
  \end{equation*}
  The top and left squares commute by \eqref{eq:local1}, and
  \eqref{eq:dblglob} makes the pair $(i\overline d_0p_0,i\overline
  d_0p_1)$ a compatible source for lower right pullback square. The
  universality thus proves our equation. 
  
  Finally, we verify that 
  \begin{equation*}
    \pathspc d_0m(j\times j)=m(\pathspc d_0\times\pathspc d_0)(j\times
    j)=m(\pathspc d_0j\times\pathspc d_0j)=m(i\overline d_0j\times
    i\overline d_0j)=m(i,i)d'_0=id'_0 \,.
  \end{equation*}
  By the same token $d_1m(j\times j)=id'_1$ hence we get the desired
  $\overline m$.

  To check (\ref{eq:vercompstcond}) we calculate:
  \begin{equation*}
    d^j_0\overline{m}=d_0j\overline{m}=d_0m(j\times{}j)=d_0p_1(j\times{}j)=d_0jp_1=d^j_0p_1\,.
  \end{equation*}
\end{prf}

\begin{lem}
  \label{lem:intvertuass}
  The composition $\overline m$ is unital and associative, that is, it
  makes \eqref{eq:dblref} a category.
\end{lem}
\begin{prf}
  Obviously since $m_{\pathspc\H}$ is so: Using the notation of
  (\ref{eq:vertcompdef}) we can formulate the associativity condition
  as the two composites in the left hand column being equal:
  \begin{equation*}
    \begin{xy}
      \xyboxmatrix{
        (\dblbarspc{\H})^3\ar[r]^-{j\times{}j\times{}j}
        \ar@<-.5ex>[d]_-{\dblbarspc{\H}\times{}\overline{m}}|-*\dir{/}
        \ar@<+.5ex>[d]^-{\overline{m}\times{}\dblbarspc{\H}}|-*\dir{/}&
        (\dblpathspc{\H})^3
        \ar@<-.5ex>[d]_-{\dblpathspc{\H}\times{}{m}}|-*\dir{/}
        \ar@<+.5ex>[d]^-{{m}\times{}\dblpathspc{\H}}|-*\dir{/}\\
        \dblbarspc{\H}\times_{d^j_0,d^j_1}\dblbarspc{\H}\ar[d]_{\overline{m}}|-*\dir{/}\ar[r]^-{j\times{}j}&
        \dblpathspc{\H}\times_{d_0,d_1}\dblpathspc{\H}\ar[d]^{m}|-*\dir{/}\\
        \dblbarspc{\H}\ar[r]_{j}&\dblpathspc{\H}
      }
    \end{xy}
  \end{equation*}
  whence we conclude that
  $j\overline{m}(\dblbarspc{\H}\times\overline{m})=j\overline{m}(\overline{m}\times\dblbarspc{\H})$,
  and by $j$ mono we get the desired
  $\overline{m}(\dblbarspc{\H}\times\overline{m})=\overline{m}(\overline{m}\times\dblbarspc{\H})$. 
  
  For the unit we can argue in the same manner:
  \begin{equation*}
    \begin{xy}
      \xyboxmatrix{
        \dblbarspc{\H}\ar[dr]^{\pair{\dblbarspc{\H},\overline{i}}}\ar[rr]^j
        \ar[dd]_{\dblbarspc{\H}}&&\dblpathspc{\H}\ar[dr]^{\pair{\dblpathspc{\H},{i}}}
        \ar[dd]|!{"2,2";"2,4"**{}}{\hole}_(.3){\dblbarspc{\H}}&\\
        {}&\dblbarspc{\H}\times_{d^j_0,d^j_1}\dblbarspc{\H}\ar[rr]^(.3){j\times{}j}
        \ar[dl]^{\overline{m}}|-*\dir{/}&&\dblpathspc{\H}\times_{d_0,d_1}\dblpathspc{\H}\ar[dl]^{{m}}|-*\dir{/}\\
        \dblbarspc{\H}\ar[rr]_j&&\dblpathspc{\H}&
      }
    \end{xy}\,.
  \end{equation*}
\end{prf}

\begin{lem}
  \label{lem:peedintcat} Applying $\qopth$ to an internal category  
  \begin{equation}
    \label{eq:intcat}
    \xymatrix{
      {{\K}}\times_{d_0, d_1}{{\K}}\ar[r]^-{{m}}|-*\dir{/}&\K\ar@<-1ex>[r]_{{d_0}}\ar@<1ex>[r]^{{d_1}}&\H\ar[l]|{{i}}
    }
  \end{equation}
  yields an internal category 
 \begin{equation*}
    \xymatrix{
      \pathspc{\K}\times_{ \pathspc{d_0},  \pathspc{d_1}} \pathspc{\K}\ar@{=}[r] &\pathspc{\K\times_{d_0, d_1}\K}\ar[r]^-{\qopth{m}}|-*\dir{/}&\pathspc\K\ar@<-1ex>[r]_{\pathspc{d_0}}\ar@<1ex>[r]^{\pathspc{d_1}}
    &\pathspc\H\ar[l]|{\pathspc{i}}\,.
  }
  \end{equation*}
\end{lem}

\begin{prf}
  This is true since $\qopth$ is an endofunctor of $\GC_{\fQ^1}$ that
  by lemma \ref{lem:pathlim} preserves pullbacks of strict diagrams. In
  particular 
  \begin{equation*}
    \begin{xy}
      \xyboxmatrix@+.5cm{
        \pathspc\K\times_{\pathspc{d}_0,\pathspc{d}_1}\pathspc\K\times_{\pathspc{d}_0,\pathspc{d}_1}\pathspc\K
        \ar[r]^-{\pathspc\K\dot\times\qopth{m}}|-*\dir{/}
        \ar[d]_-{\qopth{m}\dot\times\pathspc\K}|-*\dir{/}&
        \pathspc\K\times_{\pathspc{d}_0,\pathspc{d}_1}\pathspc\K
        \ar[d]^-{\qopth{m}}|-*\dir{/}\\
        \pathspc\K\times_{\pathspc{d}_0,\pathspc{d}_1}\pathspc\K\ar[r]_-{\qopth{m}}|-*\dir{/}&\pathspc\K
      }
    \end{xy}
  \end{equation*}
  commutes since by (\ref{eq:pprodindpres}) $\qopth(\K\dot\times{}m)=\pathspc\K\dot\times\qopth{m}$.
\end{prf}

\begin{lem}
  \label{lem:whiskdef}
  There are left and right whiskering maps
  \begin{equation*}
    \xymatrix{
      \overline{\overline{\H}}\times_{\overline{d_0}, d_1}\pathspc{\H}\ar[r]^-{w_\ell}|-*\dir{/}&
      \overline{\overline{\H}}
    }
  \end{equation*}
  \begin{equation*}
    \xymatrix{
      \pathspc{\H}\times_{d_0, \overline{d_1}}\overline{\overline{\H}}\ar[r]^-{w_r}|-*\dir{/}&
      \overline{\overline{\H}}
    }
  \end{equation*}
  induced uniquely by $\qopth(m)$. 
\end{lem}
\begin{prf}
  We construct a restricted horizontal composition $m'_r\from \pathspc\H\times_{d_0,\overline{d_1}}\overline{\overline{\H}}
  \laxto\dblpathspc\H$ in the following diagram:
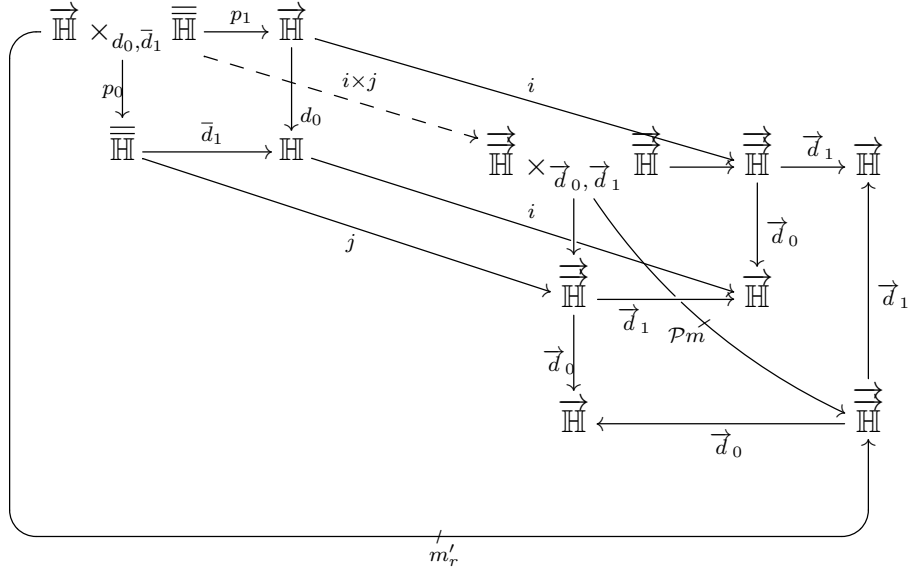
\begin{equation*}
    \begin{xy}
      \xymatrix"B"{
        \dblpathspc\H\times_{\pathspc d_0,\pathspc
          d_1}\dblpathspc\H\ar[r]
        \ar[d]\ar@/_1pc/[rrdd]|!{"B2,1";"B2,2"**{}}\hole_{\qopth{m}}|-*\dir{/}&
        \dblpathspc\H\ar[d]^{\pathspc d_0}\ar[r]^{\pathspc d_1}&\pathspc{\H}\\
        \dblpathspc\H\ar[r]_(.3){\pathspc d_1}\ar[d]_{\pathspc d_0}&\pathspc\H&{}\\
        \pathspc\H&{}&\dblpathspc\H\ar[ll]^{\pathspc d_0}\ar[uu]_{\pathspc d_1}
      }
      \POS+(-60,18)
      \xymatrix"A"{
        \pathspc\H\times_{d_0,\overline{d}_1}\overline{\overline{\H}}\ar[r]^-{p_1}\ar[d]_{p_0}\ar@{-->}"B1,1"^{i\times j}&
        \pathspc\H\ar[d]^(.7){d_0}\ar"B1,2"^i\\
        \overline{\overline{\H}}\ar[r]^{\overline{d}_1}\ar"B2,1"_j&\H\ar"B2,2"^i
        |!{"B1,1";"B2,1"**{}}{\hole}
      }
      \POS"A1,1"\ar`l+(-15,-15)`d"B3,1"+(0,-15)`r"B3,3"_{m'_r}|-*\dir{/}"B3,3"
    \end{xy}
  \end{equation*}  
  where $i\times j$ is universally induced and $m'_r$ is defined as the
  composite $\qopth(m)(i\times j)$. We need to show that $m'_r$ factors
  through $\overline{\overline{\H}}$.

  Consider the defining pullback for $\overline{\overline{\H}}$:
  \begin{equation}
    \label{eq:whiskdef2}
    \begin{xy}
      \xymatrix{
        \pathspc\H\times_{d_0,\overline{d_1}}\overline{\overline{\H}}\ar@{..>}[rd]|(.7){w_r}|-*\dir{/}
        \ar@<-.5ex>[ddr]_{\overline{d_0}p_0}\ar@<.5ex>[ddr]^{{d_1p_1}}\ar[drr]^{m'_r}|-*\dir{/}&&&\\
        &\overline{\overline{\H}}\ar[r]|{j}\ar@<-.5ex>[d]_{\overline{d_0}}\ar@<.5ex>[d]^{\overline{d_1}}&
        \dblpathspc\H\ar@<-.5ex>[d]_{\pathspc{d_0}}\ar@<.5ex>[d]^{\pathspc{d_1}}\ar@<-.5ex>[r]_{{d_0}}\ar@<.5ex>[r]^{{d_1}}& 
        \pathspc\H\ar@<-.5ex>[d]_{{d_0}}\ar@<.5ex>[d]^{{d_1}}\\
        &\H\ar[r]_i& \pathspc\H\ar@<-.5ex>[r]_{{d_0}}\ar@<.5ex>[r]^{{d_1}}&\H\,.
      }
    \end{xy}
  \end{equation}
  We need to show that $\pathspc{d}_0m'_r=i\overline d_0p_0$ and
  $\pathspc{d}_1m'_r=i d_1p_1$ to obtain a universal $w_r$, hence we calculate:
  \begin{align*}
    \pathspc d_0 m'_r&=\pathspc d_0 \qopth(m)(i\times j)= \pathspc
    d_0jp_0=\overline d_0p_0\\
    \pathspc d_1 m'_r&=\pathspc d_1 \qopth(m)(i\times j)= \pathspc
    d_1ip_1=\overline d_1p_1 
  \end{align*}
  using the definitions of $i\times j$ and $j$ as well as the
  naturality of $i$.
  
  For $w_\ell$ there is a corresponding argument.
\end{prf}

\begin{lem}
  \label{lem:whiskasscomp}
  Left and right whiskering are compatible and associative, that is, the diagrams
  \begin{equation*}
    \xymatrix{
      \pathspc\H\times_{d_0,d_1}\pathspc\H\times_{d_0,
        \overline{d}_1}\overline{\overline{\H}}
      \ar[r]^-{\pathspc\H\times w_r}|-*\dir{/}
      \ar[d]_-{m\times\pathspc\H}|-*\dir{/}&
      \pathspc\H\times_{d_0, \overline d_1}\overline{\overline{\H}}
      \ar[d]^{w_r}|-*\dir{/}\\
      \pathspc\H\times_{d_0,\overline d_1}\overline{\overline{\H}}
      \ar[r]_-{w_r}|-*\dir{/}&
      \overline{\overline\H}
    }
  \end{equation*}
  \begin{equation*}
    \xymatrix{
      \overline{\overline{\H}}\times_{\overline d_0, d_1}\pathspc\H\times_{d_0,
        d_1}\pathspc\H
      \ar[r]^-{w_\ell\times \pathspc\H}|-*\dir{/}
      \ar[d]_-{\pathspc\H\times m}|-*\dir{/}&
      \overline{\overline{\H}}\times_{\overline d_0, d_1}\pathspc\H
      \ar[d]^{w_\ell}|-*\dir{/}\\
      \overline{\overline{\H}}\times_{\overline d_0,d_1}\pathspc\H
      \ar[r]_-{w_\ell}|-*\dir{/}&
      \overline{\overline\H}
    }
  \end{equation*}
  \begin{equation*}
    \xymatrix{
      \pathspc\H\times_{d_0,\overline d_1}\overline{\overline{\H}}\times_{\overline d_0,
        d_1}\pathspc\H
      \ar[r]^-{w_r\times \pathspc\H}|-*\dir{/}
      \ar[d]_-{\pathspc\H\times w_\ell}|-*\dir{/}&
      \overline{\overline{\H}}\times_{\overline d_0, d_1}\pathspc\H
      \ar[d]^{w_\ell}|-*\dir{/}\\
      \pathspc\H\times_{d_0,\overline d_1}\overline{\overline{\H}}
      \ar[r]_-{w_r}|-*\dir{/}&
      \overline{\overline\H}
    }
  \end{equation*}
  commute.
\end{lem}
\begin{prf}
  The objects in the above diagram embed into pullbacks of
  $\dblpathspc\H$ by $j$ and these pullbacks being preserved by
  $\qopth$ and the monicity of $j$ yield the desired result.
\end{prf}

\begin{lem}
  \label{lem:whiskextlem}
  $w_\ell$ and $w_r $ extend $m$. That is
  \begin{equation*}
    \begin{xy}
      \xyboxmatrix@+1cm{
        \pathspc{\H}\times_{d_0,\overline{d_1}}\overline{\overline{\H}}\ar[r]^-{w_r}|-*\dir{/}
        \ar@<-1ex>[d]_(.4){{\pathspc\H\times{}d_0^j}}\ar@<1ex>[d]^(.4){{\pathspc\H\times{}d_1^j}}&
        \overline{\overline{\H}}
        \ar@<-1ex>[d]_(.4){{d_0^j}}\ar@<1ex>[d]^(.4){{d_1^j}}\\
        \pathspc{\H}\times_{d_0,d_1}\pathspc\H\ar[r]_-{m}|-*\dir{/}\ar[u]|(.4){\pathspc\H\times\overline{i}}&
        \pathspc\H\ar[u]|(.4){\overline{i}}
      }
    \end{xy}\qquad
    \begin{xy}
      \xyboxmatrix@+1cm{
        \overline{\overline{\H}}\times_{\overline{d_0},d_1}\pathspc{\H}\ar[r]^-{w_\ell}|-*\dir{/}
        \ar@<-1ex>[d]_(.4){{d_0\times{}\pathspc\H}}\ar@<1ex>[d]^(.4){{d_1\times{}\pathspc\H}}&
        \overline{\overline{\H}}
        \ar@<-1ex>[d]_(.4){{d_0^j}}\ar@<1ex>[d]^(.4){{d_1^j}}\\
        \pathspc{\H}\times_{d_0,d_1}\pathspc\H\ar[r]_-{m}|-*\dir{/}\ar[u]|(.4){\pathspc\H\times\overline{i}}&
        \pathspc\H\ar[u]|(.4){\overline{i}}
      }
    \end{xy}
  \end{equation*} commute serially, and the outside 0-faces are
  preserved:
  \begin{align}
    \label{eq:whskextofcpres}
    \overline{d_0}w_r&=\overline{d_0}p_1&\qquad\overline{d_0}w_\ell&={d_0}p_1\\
    \overline{d_1}w_r&={d_1}p_0&\qquad\overline{d_1}w_\ell&=\overline{d_1}p_0\notag
  \end{align}
\end{lem}
\begin{prf}
  Considering the proof of lemma \ref{lem:whiskdef} we calculate:
  \begin{equation*}
    d_0^jw_r=d_0jw_r=d_0m'_r=d_0\qopth{m}(i\times{}j)=m(d_0\times{}d_0)(i\times{}j)=m(\pathspc\H\times{}d_0^j)\,.
  \end{equation*}
 Similarly for $d_1^j$ and $w_\ell$.

 The equations (\ref{eq:whskextofcpres}) hold by the construction as
 given in (\ref{eq:whiskdef2}).
\end{prf}

Lemma \ref{lem:whiskextlem} allows us to define left and right
horizontal composites. Call the composite along the middle in the
following diagram $h_\ell\from\overline{\overline{\H}}\times_{\overline{d_0},\overline{d_1}}\overline{\overline{\H}}\laxto\overline{\overline{\H}}$:
\begin{equation}
  \label{eq:lhcompdef}
  \begin{xy}
    \xyboxmatrix{
      \pathspc{\H}\times_{d_0,\overline{d_1}}\overline{\overline{\H}}\ar[r]^-{w_r}|-*\dir{/}&
      \overline{\overline{\H}}&{}\\
      \overline{\overline{\H}}\times_{\overline{d_0},\overline{d_1}}\overline{\overline{\H}}\ar@{..>}[r]^-{}|-*\dir{/}
      \ar[d]_{\overline{\overline{\H}}\times{}d_1^j}\ar[u]^{d_0^j\times{}\overline{\overline{\H}}}&
      \overline{\overline{\H}}\times_{d_0^j,d_1^j}\overline{\overline{\H}}\ar[d]_{}\ar[u]_{}\ar[r]^-{\overline{m}}|-*\dir{/}&
      \overline{\overline{\H}}\\
      \overline{\overline{\H}}\times_{\overline{d_0},d_1}\pathspc{\H}\ar[r]_-{w_\ell}|-*\dir{/}&\overline{\overline{\H}}&{}
    }
  \end{xy}\,,
\end{equation}
and correspondingly $h_r\from\overline{\overline{\H}}\times_{\overline{d_0},\overline{d_1}}\overline{\overline{\H}}\laxto\overline{\overline{\H}}$:
\begin{equation}
  \label{eq:rhcompdef}
  \begin{xy}
    \xyboxmatrix{
      \overline{\overline{\H}}\times_{\overline{d_0},d_1}\pathspc{\H}\ar[r]^-{w_\ell}|-*\dir{/}&\overline{\overline{\H}}&{}\\
      \overline{\overline{\H}}\times_{\overline{d_0},\overline{d_1}}\overline{\overline{\H}}\ar@{..>}[r]^-{}|-*\dir{/}
      \ar[d]_{d_0^j\times{}\overline{\overline{\H}}}\ar[u]^{\overline{\overline{\H}}\times{}d_1^j}&
      \overline{\overline{\H}}\times_{d_0^j,d_1^j}\overline{\overline{\H}}\ar[d]_{}\ar[u]_{}\ar[r]^-{\overline{m}}|-*\dir{/}&
      \overline{\overline{\H}}\\
      \pathspc{\H}\times_{d_0,\overline{d_1}}\overline{\overline{\H}}\ar[r]_-{w_r}|-*\dir{/}&
      \overline{\overline{\H}}&{}
    }
  \end{xy}\,.
\end{equation}

\begin{lem}
  Left and right horizontal composites give a globular object
  \begin{equation}
    \label{eq:hcompglob}
    \xymatrix{
      \overline{\overline{\H}}\times_{\overline{d_0},\overline{d_1}}\overline{\overline{\H}}
      \ar@<-.5ex>[r]_-{h_r}|-*\dir{/}\ar@<.5ex>[r]^-{h_\ell}|-*\dir{/}&
      \overline{\overline{\H}}\ar@<-.5ex>[r]_{{d_0^j}}\ar@<.5ex>[r]^{{d_1^j}}& 
      \pathspc\H
    }\,.
  \end{equation}
\end{lem}

\begin{prf}
  We calculate:
  \begin{align*}
    d_0^jh_\ell
    &\stackrel{\tqref{eq:lhcompdef}}{=}d_0j\overline{m}\left<w_r(d_0^j\times\overline{\overline\H}),w_\ell(\overline{\overline\H}\times{}d_1^j)\right>\\
    &\stackrel{\tqref{eq:vertcompdef}}{=}d_0m(j\times{}j)\left<w_r(d_0^j\times\overline{\overline\H}),w_\ell(\overline{\overline\H}\times{}d_1^j)\right>\\
    &\stackrel{\tqref{eq:whiskdef2}}{=}d_0p_0\left<m'_r(d_0^j\times\overline{\overline\H}),m'_\ell(\overline{\overline\H}\times{}d_1^j)\right>\\
    &=d_0m'_r(d_0^j\times\overline{\overline\H})\\
    &=d_0\qopth{m}(i\times{}j)(d_0^j\times\overline{\overline\H})\\
    &\stackrel{\tqref{eq:faceweaknat}}{=}m(d_0\times{}d_0)(i\times{}j)(d_0^j\times\overline{\overline\H})\\
    &=m(d_0^j\times{}d_0^j)
  \end{align*}
  and by the same token
  \begin{equation}
    \label{eq:hcompface}
    d_0^jh_r=m(d_0^j\times{}d_0^j)\,.
  \end{equation}
  Analogously for $d_1^j$.
\end{prf}

\subsubsection{3-Paths}
\label{sec:threepaths}

We proceed to construct the internal 3-path object and the operations
involving 3-cells. Note that the $\overline{(\_)}$ and $\tilde{(\_)}$
used in this section are not at all functors.

We apply the construction in (\ref{eq:dbldef}) to
(\ref{eq:dblref}) as follows:
\begin{equation*}
  \xymatrix{
    \overline{\overline{\overline{\H}}}\ar[r]^{j}\ar@<-.5ex>[d]_{\overline{d_0^j}}\ar@<.5ex>[d]^{\overline{d_1^j}}&
    \pathspc{\overline{\overline{\H}}}\ar@<-.5ex>[d]_{\pathspc{d_0^j}}\ar@<.5ex>[d]^{\pathspc{d_1^j}}\ar@<-.5ex>[r]_{{d_0}}\ar@<.5ex>[r]^{{d_1}}& 
    \overline{\overline{\H}}\ar@<-.5ex>[d]_{{d_0^j}}\ar@<.5ex>[d]^{{d_1^j}}\\
    \pathspc{\H}\ar[r]_i& \dblpathspc\H\ar@<-.5ex>[r]_{{d_0}}\ar@<.5ex>[r]^{{d_1}}&\pathspc\H
  }
\end{equation*}
By (\ref{eq:dblref}) we get a reflexive graph
\begin{equation*}
  \xymatrix{
    \overline{\overline{\overline{\H}}}\ar@<-1ex>[r]_{{d_0^j}}\ar@<1ex>[r]^{{d_1^j}}&\overline{\overline{\H}}\ar[l]|{\overline{i}}
  }
\end{equation*}
where by (\ref{eq:dblglob})
\begin{equation*}
  \xymatrix{
    \overline{\overline{\overline{\H}}}\ar@<-.5ex>[r]_{{d_0^j}}\ar@<.5ex>[r]^{{d_1^j}}&
    \overline{\overline{\H}}\ar@<-.5ex>[r]_{{d_0^j}}\ar@<.5ex>[r]^{{d_1^j}}&
    \pathspc\H\ar@<-.5ex>[r]_{{d_0}}\ar@<.5ex>[r]^{{d_1}}&
    \H
  }
\end{equation*}
is a 3-globular object. Furthermore, by applying the reasoning of
lemma \ref{lem:pbmult} we get a vertical multiplication map
\begin{equation*}
    \xymatrix{
      \overline{\overline{\overline{\H}}}\times_{d_0^j, d_1^j}\overline{\overline{\overline{\H}}}\ar[r]^-{\overline{\overline{m}}}|-*\dir{/}&
      \overline{\overline{\overline{\H}}}
    }
  \end{equation*}
arising as a restriction of $m_{\overline{\overline{\H}}}$:
\begin{equation*}
  \xymatrix{
    \overline{\overline{\overline{\H}}}\times_{d_0^j,d_1^j} \overline{\overline{\overline{\H}}}
    \ar[r]^{j\times j}\ar@<-.5ex>[d]_{{p_0}}\ar@<.5ex>[d]^{{p_1}}
    \ar@<-.5ex>@/_3pc/[dd]_{{d'_0}}\ar@<.5ex>@/_3pc/[dd]^{{d'_1}}
    \ar@{..>}@/^3pc/[d]^{\overline{\overline{m}}}|-*\dir{/}&
    \pathspc{\overline{\overline\H}}\times_{d_0,d_1}\pathspc{\overline{\overline\H}}\ar[d]^{m}|-*\dir{/}\\
    \overline{\overline{\overline{\H}}}\ar[r]^{j}\ar@<-.5ex>[d]_{{\overline
        d_0}}\ar@<.5ex>[d]^{{\overline d_1}}&
    \pathspc{\overline{\overline\H}}\ar@<-.5ex>[d]_{{\pathspc
        d_0}}\ar@<.5ex>[d]^{{\pathspc d_1}}\\
    \pathspc\H\ar[r]_{i}&\overline{\overline\H}
  }
\end{equation*}
where $d'_0=\overline{d_0}p_0$ and $d'_1=\overline{d_1}p_1$.

\begin{lem}
  \label{lem:whsk1on3}
   There are left and right whiskering maps
  \begin{equation*}
    \xymatrix{
      \overline{\overline{\overline{\H}}}\times_{d_0\overline{d_0^j}, d_1}\pathspc{\H}\ar[r]^-{\overline{w_\ell}}|-*\dir{/}&
      \overline{\overline{\overline{\H}}}
    }
  \end{equation*}
  \begin{equation*}
    \xymatrix{
      \pathspc{\H}\times_{d_0, d_1\overline{d_1^j}}\overline{\overline{\overline{\H}}}\ar[r]^-{\overline{w_r}}|-*\dir{/}&
      \overline{\overline{\overline{\H}}}
    }
  \end{equation*}
  induced uniquely by $\qopth{w_\ell}$ and $\qopth{w_r}$. 
\end{lem}
\begin{prf}
  We define $\overline{w_\ell}$ as the universally induced arrow in
  the following diagram:
  \begin{equation}
    \label{eq:whsk1on3def}
    \begin{xy}
      \xyboxmatrix{
        \overline{\overline{\overline{\H}}}\times_{d_0\overline{d^j_0},d_1}\pathspc{\H}
        \ar[r]^{j\times{}i}
        \ar@{-->}[d]|(.3)*\dir{/}^{\overline{w_\ell}}
        \ar@<-.5ex>@/_3pc/[dd]|(.3)*\dir{/}_{r_0}
        \ar@<.5ex>@/_3pc/[dd]|(.3)*\dir{/}^{r_1}&
        \pathspc{\overline{\overline{\H}}}\times_{\pathspc{\overline{d_0}},\pathspc{d_1}}\dblpathspc{\H}
        \ar@<-.5ex>[r]_{{d_0}}\ar@<.5ex>[r]^{{d_1}}\ar[d]|(.3)*\dir{/}|{\qopth{w_\ell}}&
        \dblbarspc{\H}\times_{\overline{d_0},d_1}\pathspc{\H}
        \ar@<-.5ex>[rd]_{{d^j_0\times{}\pathspc{\H}}}\ar@<.5ex>[rd]^{d^j_1\times{}\pathspc{\H}}
        \ar[d]|(.3)*\dir{/}_{{w_\ell}}&
        {}\\
        \overline{\overline{\overline{\H}}}
        \ar@<-.5ex>[d]_{\overline{d^j_0}}\ar@<.5ex>[d]^{\overline{d^j_1}}
        \ar[r]^{j}&
        \pathspc{\overline{\overline{\H}}}\ar@<-.5ex>[r]_{{d_0}}\ar@<.5ex>[r]^{{d_1}}\ar@<-.5ex>[d]_{\pathspc{d^j_0}}
        \ar@<.5ex>[d]^{\pathspc{d^j_1}}&
        \overline{\overline{\H}}\ar@<-.5ex>[d]_{{d^j_0}}\ar@<.5ex>[d]^{{d^j_1}}&
        \pathspc{\H}\times_{d_0,d_1}\pathspc{\H}\ar[dl]|-*\dir{/}^{m}\\
        \pathspc{\H}\ar[r]_{i}&\dblpathspc{\H}\ar@<-.5ex>[r]_{{d_0}}\ar@<.5ex>[r]^{{d_1}}&\pathspc{\H}&{}
      }
    \end{xy}
  \end{equation}
  where $r_0=m(\overline{d^j_0}\times\pathspc\H)$ and
  $r_1=m(\overline{d^j_1}\times\pathspc\H)$.
  We calculate
  \begin{multline*}
    ir_0\\
    =im(\overline{d^j_0}\times\pathspc\H)=\qopth{m}(i\times{}i)(\overline{d^j_0}\times\pathspc\H)
    =\qopth{m}(i\overline{d^j_0}\times{}i)=\qopth{m}(\pathspc{d^j_0}j\times{}i)
    =\qopth(d^j_0w_\ell)(j\times{}i)\\
    =\pathspc{d^j_0}\qopth{w_\ell}(j\times{}i)\,,
  \end{multline*}
  and likewise for $r_1$ and $\pathspc{d^j_1}$. And hence we obtain
  $\overline{w_\ell}$, and $\overline{w_r}$ by analogy.
\end{prf}

\begin{lem}
  \label{lem:whsk13extlem}
  $\overline{w_\ell}$ and $\overline{w_r}$ extend $w_\ell$ and $w_r$
  respectively. That is
  \begin{equation}
    \label{eq:whsk13exteq}
    \begin{xy}
      \xyboxmatrix{
        \pathspc{\H}\times_{d_0,d_1\overline{d_1^j}}\overline{\overline{\overline{\H}}}\ar[r]^-{\overline{w_r}}|-*\dir{/}
        \ar@<-.5ex>[d]_{{\pathspc\H\times{}d_0^j}}\ar@<.5ex>[d]^{{\pathspc\H\times{}d_1^j}}&
        \overline{\overline{\overline{\H}}}
        \ar@<-.5ex>[d]_{{d_0^j}}\ar@<.5ex>[d]^{{d_1^j}}\\
        \pathspc{\H}\times_{d_0,d_1^j}\overline{\overline{\H}}\ar[r]_-{w_r}|-*\dir{/}&\overline{\overline{\H}}
      }
    \end{xy}\qquad
    \begin{xy}
      \xyboxmatrix{
        \overline{\overline{\overline{\H}}}\times_{d_0\overline{d_0^j},d_1}\pathspc{\H}\ar[r]^-{\overline{w_\ell}}|-*\dir{/}
        \ar@<-.5ex>[d]_{{d_0^j\times{}\pathspc\H}}\ar@<.5ex>[d]^{{d_1^j\times{}\pathspc\H}}&
        \overline{\overline{\overline{\H}}}
        \ar@<-.5ex>[d]_{{d_0^j}}\ar@<.5ex>[d]^{{d_1^j}}\\
        \overline{\overline{\H}}\times_{d_0^j,d_1}\pathspc\H\ar[r]_-{w_\ell}|-*\dir{/}&\overline{\overline{\H}}
      }
    \end{xy}
  \end{equation} commute serially.
\end{lem}
\begin{prf}
  Inspecting (\ref{eq:whsk1on3def}) we can calculate
  \begin{multline*}
    d_0^j\overline{w_\ell}\\=d_0j\overline{w_\ell}=d_0\qopth(w_\ell)(j\times{}i)
    =w_\ell{}d_0(j\times{}i)=w_\ell{}(d_0\times{}d_0)(j\times{}i)\\=w_\ell{}(d_0^j\times{}\pathspc\H)\,.
  \end{multline*}
  And likewise for the other squares in (\ref{eq:whsk13exteq}).
\end{prf}

Lastly, we need the whiskering of a 3-path by a 2-path along a
1-path. We can reapply the basic scheme of lemma \ref{lem:whiskdef}.
\begin{lem}
  \label{lem:whiskdef23}
  There are left and right whiskering maps
  \begin{equation*}
    \xymatrix{
      \overline{\overline{\overline{\H}}}\times_{\overline{d^j_0}, d^j_1}\overline{\overline{\H}}\ar[r]^-{\tilde{w}_\ell}|-*\dir{/}&
      \overline{\overline{\overline{\H}}}
    }
  \end{equation*}
  \begin{equation*}
    \xymatrix{
      \overline{\overline{\H}}\times_{d^j_0, \overline{d^j_1}}\overline{\overline{\overline{\H}}}\ar[r]^-{\tilde{w}_r}|-*\dir{/}&
      \overline{\overline{\overline{\H}}}
    }
  \end{equation*}
  induced uniquely by $\qopth(\overline{m})$. 
  
  And these extend $\overline{m}$, that is
  \begin{align}
    \label{eq:whiskdef23ext}
    d^j_0\tilde{w}_r=\overline{m}(\overline{\overline{\H}}\times{}d^j_0)&
    \qquad{}d^j_1\tilde{w}_r=\overline{m}(\overline{\overline{\H}}\times{}d^j_1)\\
    d^j_0\tilde{w}_\ell=\overline{m}(d^j_0\times{}\overline{\overline{\H}})&
    \qquad{}d^j_1\tilde{w}_\ell=\overline{m}(d^j_1\times{}\overline{\overline{\H}})\,.
  \end{align}
\end{lem}

\begin{prf}
  The desired map arises as a universal arrow in the following diagram:
  \begin{equation}
    \label{eq:whiskdef23}
    \begin{xy}
      \xyboxmatrix@+.5cm{
        \overline{\overline{\H}}\times_{d^j_0,\overline{d^j_1}}\overline{\overline{\overline{\H}}}
        \ar@{..>}[rd]^{\tilde{w}_r}|-*\dir{/}
        \ar[rr]^-{i\times{}j}
        \ar@<-.5ex>[ddr]_{\overline{d^1_0}p_0}\ar@<.5ex>[ddr]^{{d^j_1p_1}}&&
        \pathspc{\overline{\overline{\H}}}\times_{\pathspc{d^j_0},\pathspc{d^j_1}}\pathspc{\overline{\overline{\H}}}
        \ar[d]^{\qopth{\overline{m}}}|-*\dir{/}\ar@<-.5ex>[r]_{{d_0}}\ar@<.5ex>[r]^{{d_1}}&
        \overline{\overline{\H}}\times_{{d^j_0},d^j_1}\overline{\overline{\H}}
        \ar[d]^{{\overline{m}}}|-*\dir{/}\\
        &\overline{\overline{\overline{\H}}}\ar[r]^{j}\ar@<-.5ex>[d]_{\overline{d^j_0}}\ar@<.5ex>[d]^{\overline{d^j_1}}&
        \pathspc{\overline{\overline{\H}}}
        \ar@<-.5ex>[d]_{\pathspc{d^j_0}}\ar@<.5ex>[d]^{\pathspc{d^j_1}}\ar@<-.5ex>[r]_{{d_0}}\ar@<.5ex>[r]^{{d_1}}& 
        \overline{\overline{\H}}\ar@<-.5ex>[d]_{{d^j_0}}\ar@<.5ex>[d]^{{d^j_1}}\\
        &\pathspc{\H}\ar[r]_i& \dblpathspc\H\ar@<-.5ex>[r]_{{d_0}}\ar@<.5ex>[r]^{{d_1}}&\H\,.
      }
    \end{xy}\,.
  \end{equation}
  Now, we can verify
  $i\overline{d^j_0}p_0=\pathspc{d^j_0}jp_0=\pathspc{d^j_0}p_0(i\times{}j)=\pathspc{d^j_0}
  \qopth{\overline{m}}(i\times{}j)$ and $i{d^j_1}p_1=\pathspc{d^j_1}jp_1=\pathspc{d^j_1}p_1(i\times{}j)=\pathspc{d^j_1}
  \qopth{\overline{m}}(i\times{}j)$.

  The equations (\ref{eq:whiskdef23ext}) are now immediate.
\end{prf}

\subsection{The Space of Parallel Cells}
\label{sec:parcel}

For a $\Gray$-category $\H$ we define the space of parallel 1-cells
$P^1(\H)$ as the following limit:
\begin{equation*}
  \begin{xy}
    \xyboxmatrix{
      {}&P^1(\H)\ar[dl]_{p_0}\ar[dr]^{p_1}&{}\\
      \pathspc\H\ar[d]_{d_0}\ar[drr]^(.3){d_1}&{}&\pathspc\H\ar[dll]_(.3){d_0}|\hole\ar[d]^{d_1}\\
      \H&{}&\H\\
    }
  \end{xy}
\end{equation*}

\begin{equation}
  \label{eq:par2celdef}
  \begin{xy}
    \xyboxmatrix{
      {}&P^2(\H)\ar[dl]_{p_0}\ar[dr]^{p_1}&{}\\
      \overline{\overline{\H}}\ar[d]_{d_0^j}\ar[drr]^(.3){d_1^j}&{}&\overline{\overline{\H}}\ar[dll]_(.3){d_0^j}|\hole\ar[d]^{d_1^j}\\
      \pathspc\H&{}&\pathspc\H\\
    }
  \end{xy}
\end{equation}

\begin{lem}\label{lem:1cart}
  The canonical map
  $\pair{d_0^j,d_1^j}\from\overline{\overline{\overline{\H}}}\to{}P^2(\H)$ is 1-Cartesian.
\end{lem}
\begin{prf}
  Consider the following cells in
  $\overline{\overline{\overline{\H}}}$
  \begin{align*}
    f=(f_4;f_2,f_3;f_0,f_1)&\\
    g=(g_4;g_2,g_3;g_0,g_1)&\\
    h=(h_4,h_5;h_2,h_3;h_0,h_1)\from{}f\to{}g&\\
    k=(k_4,k_5;k_2,k_3;k_0,k_1)\from{}f\to{}g&\\
    \alpha=(\alpha_3;\alpha_1,\alpha_2)\from{}h\To{}k&
  \end{align*}

  By construction the map $\pair{d_0^j,d_1^j}$ acts on this data as
  follows:
  \begin{align*}
    f&\mapsto((f_2;f_0,f_1),(f_3;f_0,f_1))\\
    g&\mapsto((g_2;g_0,g_1),(g_3;g_0,g_1))\\
    h&\mapsto((h_4;h_2,h_3;h_0,h_1),(h_5;h_2,h_3;h_0,h_1))\\
    k&\mapsto((k_4;k_2,k_3;k_0,k_1),(k_5;k_2,k_3;k_0,k_1))\\
    \alpha&\mapsto((\alpha_3;\alpha_1,\alpha_2),(\alpha_3;\alpha_1,\alpha_2))
  \end{align*}
  where on the right we find parallel pairs of cells from
  $\overline{\overline{\H}}$, that is, in (\ref{eq:1cartprf1}) the
  central square, the outer square, and the left and right hand trapezoids
  commute by assumption.

  The requisite compatibility conditions for $f, g, h, k, \alpha$ to
  be cells of $\overline{\overline{\H}}$ are displayed in
  (\ref{eq:1cartprf1}). We obverse that the remaining trapezoids at
  the top and the bottom commute by naturality of $\#_1$ and $\ten$ in
  $\H$. Hence we conclude that given 1-cells $h, k$ in
  $\overline{\overline{\overline{\H}}}$ all higher cells, including
  3-cells, between them are determined by their image under
  $\pair{d_0^j,d_1^j}$.
  \begin{sidewaysfigure}
    \begin{equation}
      \label{eq:1cartprf1}
      \begin{xy}
        \save
        (40,0):(0,-.8)::
        ,(0,0)="A11"
        ,(0,1)="A12"
        ,(0,2)="A13"
        ,(0,3)="A14"
        ,(0,4)="A15"
        ,(0,5)="A16"
        ,(1,0)="A21"
        ,(1,1)="A22"
        ,(1,2)="A23"
        ,(1,3)="A24"
        ,(1,4)="A25"
        ,(1,5)="A26"
        ,(2,0)="A31"
        ,(2,1)="A32"
        ,(2,1.9)="A33"
        ,(2,3.1)="A34"
        ,(2,4)="A35"
        ,(2,5)="A36"
        ,(3,0)="A41"
        ,(3,1)="A42"
        ,(3,1.9)="A43"
        ,(3,3.1)="A44"
        ,(3,4)="A45"
        ,(3,5)="A46"
        ,(4,0)="A51"
        ,(4,1)="A52"
        ,(4,2)="A53"
        ,(4,3)="A54"
        ,(4,4)="A55"
        ,(4,5)="A56"
        ,(5,0)="A61"
        ,(5,1)="A62"
        ,(5,2)="A63"
        ,(5,3)="A64"
        ,(5,4)="A65"
        ,(5,5)="A66"
        \restore
        \POS,"A11"
        \xyboxmatrix"A11"@+.5cm{
          {}\ar@/^1pc/[r]^{f_0}="r"\ar@/_1pc/[r]_{f_1}="s" 
          \ar[d]|{k_0}="x"\ar@/_3pc/[d]_{h_0}="y" 
          &
          {}\ar[d]^{k_1}="x1"
          \ar@<+1.5ex>@2[dl] **{}?(.3);?(.7)|{k_3}
          \\
          \ar@/_1pc/[r]_{g_1}="s1" 
          & 
          {}
          \POS\ar@2 "x" ; "y" **{}?<(.3);?>(.7)^{\alpha_1} 
          \ar@2 "r" ; "s" **{}?(.3);?(.7)_{f_2}
        }
        ,"A22"
        \xyboxmatrix"A22"@+.5cm{
          {}\ar@/^1pc/[r]^{f_0}="r"\ar@/_1pc/[r]_{f_1}="s" 
          \ar[d]_{h_0}="x" 
          &
          {}\ar[d]|{h_1}="x1"\ar@/^3pc/[d]^{k_1}="y1"
          \ar@<+1.5ex>@2[dl] **{}?(.3);?(.7)|{h_3}
          \\
          {}\ar@/_1pc/[r]_{g_1}="s1" 
          & 
          {}
          \POS
          \ar@2 "y1" ; "x1" **{}?<(.3);?>(.7)^{\alpha_2}="X"
          \ar@2 "r" ; "s" **{}?(.3);?(.7)_{f_2}="Y"
          \tria "Y";"X"
        }
        ,"A33"
        \xyboxmatrix"A33"@+.5cm{
          {}\ar@/^1pc/[r]^{f_0}="r"\ar@/_1pc/[r]_{f_1}="s" 
          \ar[d]_{h_0}="x" 
          &
          {}\ar[d]|{h_1}="x1"\ar@/^3pc/[d]^{k_1}="y1"
          \ar@<+1.5ex>@2[dl] **{}?(.3);?(.7)|{h_3}
          \\
          {}\ar@/_1pc/[r]_{g_1}="s1" 
          & 
          {}
          \POS
          \ar@2 "y1" ; "x1" **{}?<(.3);?>(.7)^{\alpha_2}="X"
          \ar@2 "r" ; "s" **{}?(.3);?(.7)_{f_2}="Y"
          \tria "X";"Y"
        }
        ,"A61"
        \xyboxmatrix"A61"@+.5cm{
          {}\ar@/^1pc/[r]^{f_0}="r"\ar@/_1pc/[r]_{f_1}="s" 
          \ar[d]|{k_0}="x"\ar@/_3pc/[d]_{h_0}="y" 
          &
          {}\ar[d]^{k_1}="x1"
          \ar@<+1.5ex>@2[dl] **{}?(.3);?(.7)|{k_3}
          \\
          {}\ar@/_1pc/[r]_{g_1}="s1" 
          & 
          {}
          \POS\ar@2 "x" ; "y" **{}?<(.3);?>(.7)^{\alpha_1} 
          \ar@2 "r" ; "s" **{}?(.3);?(.7)_{f_3}
        }
        ,"A52"
        \xyboxmatrix"A52"@+.5cm{
          {}\ar@/^1pc/[r]^{f_0}="r"
          \ar@/_1pc/[r]_{f_1}="s"
          \ar[d]_{h_0}="x"
          &
          {}\ar[d]|{h_1}="x1"
          \ar@/^3pc/[d]^{k_1}="y1"
          \ar@<+1.5ex>@2[dl] **{}?(.3);?(.7)|{h_3}
          \\
          {}\ar@/_1pc/[r]_{g_1}="s1" 
          & 
          {}
          \POS
          \ar@2 "y1" ; "x1" **{}?<(.3);?>(.7)^{\alpha_2}="Y" 
          \ar@2 "r" ; "s" **{}?(.3);?(.7)_{f_3}="X"
          \tria "X";"Y"
        }
        ,"A43"
        \xyboxmatrix"A43"@+.5cm{
          {}\ar@/^1pc/[r]^{f_0}="r"
          \ar@/_1pc/[r]_{f_1}="s"
          \ar[d]_{h_0}="x"
          &
          {}\ar[d]|{h_1}="x1"
          \ar@/^3pc/[d]^{k_1}="y1"
          \ar@<+1.5ex>@2[dl] **{}?(.3);?(.7)|{h_3}
          \\
          {}\ar@/_1pc/[r]_{g_1}="s1" 
          & 
          {}
          \POS
          \ar@2 "y1" ; "x1" **{}?<(.3);?>(.7)^{\alpha_2}="Y" 
          \ar@2 "r" ; "s" **{}?(.3);?(.7)_{f_3}="X"
          \tria "Y";"X"          
        }
        ,"A16"
        \xyboxmatrix"A16"@+.5cm{
          {}\ar@/^1pc/[r]^{f_0}="r"
          \ar[d]|{k_0}="x"\ar@/_3pc/[d]_{h_0}="y" 
          &
          {}\ar[d]^{k_1}="x1"
          \ar@<-1.5ex>@2[dl] **{}?(.3);?(.7)|{k_2}
          \\
          {}\ar@/^1pc/[r]^{g_0}="r1"\ar@/_1pc/[r]_{g_1}="s1" 
          & 
          {}
          \POS
          \ar@2 "x" ; "y" **{}?<(.3);?>(.7)^{\alpha_1}="X"
          \ar@2 "r1" ; "s1" **{}?(.3);?(.7)_{g_2}="Y"
          \tria "Y";"X"
        }
        ,"A25"
        \xyboxmatrix"A25"@+.5cm{
          {}\ar@/^1pc/[r]^{f_0}="r"
          \ar[d]|{k_0}="x"\ar@/_3pc/[d]_{h_0}="y" 
          &
          {}\ar[d]^{k_1}="x1"
          \ar@<-1.5ex>@2[dl] **{}?(.3);?(.7)|{k_2}
          \\
          {}\ar@/^1pc/[r]^{g_0}="r1"\ar@/_1pc/[r]_{g_1}="s1" 
          & 
          {}
          \POS
          \ar@2 "x" ; "y" **{}?<(.3);?>(.7)^{\alpha_1}="X"
          \ar@2 "r1" ; "s1" **{}?(.3);?(.7)_{g_2}="Y"
          \tria "X";"Y"         
        }
        ,"A34"
        \xyboxmatrix"A34"@+.5cm{
          {}\ar@/^1pc/[r]^{f_0}="r"
          \ar[d]_{h_0}="x"
          &
          {}\ar[d]|{h_1}="x1"\ar@/^3pc/[d]^{k_1}="y1"
          \ar@<-1.5ex>@2[dl] **{}?(.3);?(.7)|{h_2}
          \\
          {}\ar@/^1pc/[r]^{g_0}="r1"\ar@/_1pc/[r]_{g_1}="s1" 
          & 
          {}
          \POS
          \ar@2 "y1" ; "x1" **{}?<(.3);?>(.7)^{\alpha_2} 
          \ar@2 "r1" ; "s1" **{}?(.3);?(.7)_{g_2}
        }
        ,"A66"
        \xyboxmatrix"A66"@+.5cm{
          {}\ar@/^1pc/[r]^{f_0}="r"
          \ar[d]|{k_0}="x"\ar@/_3pc/[d]_{h_0}="y" 
          &
          {}\ar[d]^{k_1}="x1"
          \ar@<-1.5ex>@2[dl] **{}?(.3);?(.7)|{k_2}
          \\
          {}\ar@/^1pc/[r]^{g_0}="r1"\ar@/_1pc/[r]_{g_1}="s1" 
          & 
          {}
          \POS\ar@2 "x" ; "y" **{}?<(.3);?>(.7)^{\alpha_1}="X" 
          \ar@2 "r1" ; "s1" **{}?(.3);?(.7)_{g_3}="Y"
          \tria "Y";"X"
        }
        ,"A55"
        \xyboxmatrix"A55"@+.5cm{
          {}\ar@/^1pc/[r]^{f_0}="r"
          \ar[d]|{k_0}="x"\ar@/_3pc/[d]_{h_0}="y" 
          &
          {}\ar[d]^{k_1}="x1"
          \ar@<-1.5ex>@2[dl] **{}?(.3);?(.7)|{k_2}
          \\
          {}\ar@/^1pc/[r]^{g_0}="r1"\ar@/_1pc/[r]_{g_1}="s1" 
          & 
          {}
          \POS\ar@2 "x" ; "y" **{}?<(.3);?>(.7)^{\alpha_1}="X" 
          \ar@2 "r1" ; "s1" **{}?(.3);?(.7)_{g_3}="Y"
          \tria "X";"Y"
        }
        ,"A44"
        \xyboxmatrix"A44"@+.5cm{
          {}\ar@/^1pc/[r]^{f_0}="r"
          \ar[d]_{h_0}="x"
          &
          {}\ar[d]|{h_1}="x1"\ar@/^3pc/[d]^{k_1}="y1"
          \ar@<-1.5ex>@2[dl] **{}?(.3);?(.7)|{h_2}
          \\
          {}\ar@/^1pc/[r]^{g_0}="r1"\ar@/_1pc/[r]_{g_1}="s1" 
          & 
          {}
          \POS
          \ar@2 "y1" ; "x1" **{}?<(.3);?>(.7)^{\alpha_2} 
          \ar@2 "r1" ; "s1" **{}?(.3);?(.7)_{g_3}
        }
        \ar@3"A11";"A22"^{\underline{\alpha_3}\#_1(k_1\#_0f_2)}
        \ar@3"A22";"A33"^{h_3\#_1(\alpha_2\ten{}f_2)}
        \ar@3"A11";"A61"^{(g_1\#_0\alpha_1)\#_1k_3\#_1(k_1\#_0\underline{f_4})}
        \ar@3"A61";"A52"^{\underline{\alpha_3}\#_1(k_1\#_0f_3)}
        \ar@3"A52";"A43"^{h_3\#_1(\alpha_2\ten{}f_3)}
        \ar@3"A33";"A43"|{}="x"
        \save\POS"x"+(0,17)*!C\labelbox{h_3\#_1(h_1\#_0\underline{f_4})\\\#_1(\alpha_2\#_0f_0)}\ar@{.}"x";c\restore
        \ar@3"A11";"A16"|{(g_1\#_0\alpha_1)\#_1\underline{k_4}}
        \ar@3"A33";"A34"_{\underline{h_4}\#_1(\alpha_2\#_0f_0)}
        \ar@3"A43";"A44"^{\underline{h_5}\#_1(\alpha_2\#_0f_0)}
        \ar@3"A34";"A44"_{}="x"
        \save\POS"x"+(0,-17)*!C\labelbox{(\underline{g_4}\#_0h_0)\\\#_1h_2\#_1(\alpha_2\#_0f_0)}\ar@{.}"x";c\restore
        \ar@3"A61";"A66"|{(g_1\#_0\alpha_1)\#_1\underline{k_5}}
        \ar@3"A16";"A25"_{\overline{g_2\ten{}\alpha_1}\#_1k_2}
        \ar@3"A25";"A34"_{(g_2\#_0h_0)\#_1\underline{\alpha_3}}
        \ar@3"A16";"A66"_{(g_1\#_0\alpha_1)\#_1(\underline{g_4}\#_0k_0)\#_1k_2}
        \ar@3"A55";"A44"_{(g_3\#_0h_0)\#_1\alpha_3}
        \ar@3"A66";"A55"^{\overline{g_3\ten{}\alpha_1}\#_1k_2}
        \ar@3"A22";"A52"^{h_3\#_1(\alpha_2\#_0f_1)\#_1(k_1\#_0\underline{f_4})}
        \ar@3"A25";"A55"_{(\underline{g_4}\#_0h_0)\#_1(g_0\#_0\alpha_1)\#_1k_2}
      \end{xy}
    \end{equation}
  \end{sidewaysfigure}
\end{prf}

\begin{lem}
  The 3-paths compose horizontally along 2-paths, that is,
  \label{lem:whisk23inter}
  \begin{equation*}
    \begin{xy}
      \xyboxmatrix@+2.5cm{
        \overline{\overline{\overline{\H}}}\times_{\overline{d^j_0},\overline{d^j_1}}\overline{\overline{\overline{\H}}}
        \ar[r]^-{\pair{\tilde{w_\ell}(\overline{\overline{\overline{\H}}}\times{}d^j_1),
            \tilde{w_r}(d^j_10\times{}\overline{\overline{\overline{\H}}})}}
        \ar[d]_-{\pair{\tilde{w_r}(d^j_1\times{}\overline{\overline{\overline{\H}}}),
            \tilde{w_\ell}(\overline{\overline{\overline{\H}}}\times{}d^j_0)}}&
        \overline{\overline{\overline{\H}}}\times_{{d^j_0},{d^j_1}}\overline{\overline{\overline{\H}}}
        \ar[d]^-{\overline{\overline{m}}}\\
        \overline{\overline{\overline{\H}}}\times_{{d^j_0},{d^j_1}}\overline{\overline{\overline{\H}}}
        \ar[r]_-{\overline{\overline{m}}}&
        \overline{\overline{\overline{\H}}} }
    \end{xy}
  \end{equation*}
  commutes. \qed
\end{lem}

\subsection{The Tensor Map}
\label{sec:tenmap}

Given that by lemma \ref{lem:1cart} we have a 1-Cartesian map
$\pair{d_0^j,d_1^j}\overline{\overline{\overline{\H}}}\to{}P^2(\H)$ we
consider the following diagram in $\Gray\Cat_{\fQ^1}$
\begin{equation}
  \label{eq:tenmapind}
  \begin{xy}
    \xyboxmatrix{
      \overline{\overline{\H}}\times_{\overline{d_0},\overline{d_1}}\overline{\overline{\H}}\ar[drr]^{\pair{h_\ell,h_r}}|-*\dir{/}
      \ar@{..>}[dr]_{t}&&\\
      {}&\overline{\overline{\overline{\H}}}\ar[r]_-{\pair{d_0^j,d_1^j}}&P^2(\H)\\
    }
  \end{xy}
\end{equation}
where $h_\ell$ and $h_r$ are given by (\ref{eq:lhcompdef}) and
(\ref{eq:rhcompdef}) respectively. By (\ref{eq:hcompglob}) we know
that $(h_\ell, h_r)$ is a source for (\ref{eq:par2celdef}) hence we
obtain $\pair{h_\ell,h_r}$.

There is a map
$t_1\from(\overline{\overline{\H}}\times_{\overline{d_0},\overline{d_1}}\overline{\overline{\H}})_1
\to{}(\overline{\overline{\overline{\H}}})_1$
in $\Cat_{\fQ^1}$ given by:
\begin{multline*}
  (g,f)=((g_2;g_0,g_1),(f_2;f_0,f_1))=
  \begin{pmatrix}
    \begin{xy}
      \xyboxmatrix{
        \ar@/^1.5pc/[r]^{f_0}="x"\ar@/_1.5pc/[r]_{f_1}="y"&\ar@/^1.5pc/[r]^{g_0}="x1"\ar@/_1.5pc/[r]_{g_1}="y1"&{}
        \POS\ar@2 "x";"y"**{}?(.2);?(.8)_{f_2} \POS\ar@2
        "x1";"y1"**{}?(.2);?(.8)_{g_2} }
    \end{xy}
  \end{pmatrix}
  \\\mapsto (g_2\ten{}f_2;g_2\lhc{}f_2,g_2\rhc{}f_2;g_0\#_0f_0,g_1\#_0f_1)=
  \begin{pmatrix}
    \begin{xy}
      \xyboxmatrix"A"{
        \ar@/^1.5pc/[r]^{g_0\#_0f_0}="x"\ar@/_1.5pc/[r]_{g_1\#_0f_1}="y"&{}
        \POS\ar@2 "x";"y"**{}?(.2);?(.8)|{g_2\lhc{}f_2} } ,(30,0)
      \xyboxmatrix"B"{
        \ar@/^1.5pc/[r]^{g_0\#_0f_0}="x"\ar@/_1.5pc/[r]_{g_1\#_0f_1}="y"&{}
        \POS\ar@2 "x";"y"**{}?(.2);?(.8)|{g_2\rhc{}f_2} }
      \POS\ar@3"A";"B"^{g_2\ten{}f_2}
    \end{xy}
  \end{pmatrix}
\end{multline*}
and
\begin{multline*}
  \begin{pmatrix}
    (k,h)\from(g,f)\to(g',f')
  \end{pmatrix}
  =
  \begin{pmatrix}
    (k_4;k_2,k_3;h_1,k_1),\\
    (h_4;h_2,h_3;h_0,h_1)
  \end{pmatrix}
  =
  \begin{pmatrix}
    \begin{xy}
      \xyboxmatrix"A"{
        {}\ar@/^1pc/[r]^{f_0}="r"\ar@/_1pc/[r]_{f_1}="s" \ar[d]_{h_0}
        & {}\ar@/^1pc/[r]^{g_0}="r1"\ar@/_1pc/[r]_{g_1}="s1"
        \ar[d]|{h_1} \ar@<+1.5ex>@2[dl] **{}?(.3);?(.7)|{h_3} &
        {}\ar[d]^{k_1}="x1" \ar@<+1.5ex>@2[dl] **{}?(.3);?(.7)|{k_3}
        \\
        \ar@/_1pc/[r]_{f'_1} & \ar@/_1pc/[r]_{g'_1} & {} \POS\ar@2 "r"
        ; "s" **{}?(.3);?(.7)_{f_2} \POS\ar@2 "r1" ; "s1"
        **{}?(.3);?(.7)_{g_2} 
      }
      ,(0,-40)
      \xyboxmatrix"B"{
        {}\ar@/^1pc/[r]^{f_0} \ar[d]_{h_0}
        & {}\ar@/^1pc/[r]^{g_0}
        \ar[d]|{h_1} \ar@<-1.5ex>@2[dl] **{}?(.3);?(.7)|{h_2} &
        {}\ar[d]^{k_1}="x1"\ar@<-1.5ex>@2[dl] **{}?(.3);?(.7)|{k_2}
        \\
        \ar@/_1pc/[r]_{f'_1}="s"\ar@/^1pc/[r]^{f'_0}="r"&
        \ar@/_1pc/[r]_{g'_1}="s1"\ar@/^1pc/[r]^{g'_0}="r1" &
        {}
        \POS\ar@2 "r";"s"**{}?(.3);?(.7)_{f'_2} 
        \POS\ar@2 "r1";"s1"**{}?(.3);?(.7)_{g'_2} 
      }
      \POS\ar@3"A";"B"^{(k_4,h_4)}
    \end{xy}
  \end{pmatrix}
  \\\mapsto 
  \begin{pmatrix}
    \omega_1,\omega_2;(g'_0\#_0h_2)\#_1(k_2\#_0f_0),\\(g'_1\#_0h_3)\#_1(k_3\#_0f_1);h_0,k_1
  \end{pmatrix}
\,,
\end{multline*}
where $\omega_1$ and $\omega_2$ are defined as the vertical composites
in (\ref{eq:tenmapdef1cell}), by definition these constitute the
components of a 1-cell in $\overline{\overline{\overline{\H}}}$.

\begin{sidewaysfigure}
  \begin{equation}
    \label{eq:tenmapdef1cell}
    \begin{xy}
      \save
      (50,0):(0,-1)::
      ,(0,0)="A11"
      ,(1,0)="A21"
      ,(2,0)="A31"
      ,(0,1)="A12"
      ,(1,1)="A22"
      ,(2,1)="A32"
      ,(0,2)="A13"
      ,(1,2)="A23"
      ,(2,2)="A33"
      \restore
      \POS
      ,"A11"
      \xyboxmatrix"A11"@+.5cm{
        \ar@/^1pc/[r]^{f_0}="x1"\ar@/_1pc/[r]_{f_1}="y1"\ar[d]_{h_0}&
        \ar@/^1pc/[r]^{g_0}="x2"\ar@/_1pc/[r]_{g_1}="y2"\ar[d]|{h_0}\ar@2@<+1.5ex>[dl]**{}?(.3);?(.7)|{h_3}&
        \ar@2@<+1.5ex>[dl]**{}?(.3);?(.7)|{k_3}\ar[d]^{k_1}\\
        \ar@/_1pc/[r]_{f'_1}="w1"&
        \ar@/_1pc/[r]_{g'_1}="w2"&
        {}
        \POS
        \ar@2"x1";"y1"**{}?(.3);?(.7)|{f_2}="r"
        \ar@2"x2";"y2"**{}?(.3);?(.7)|{g_2}="s"
        \tria"r";"s"
      }
      ,"A31"      
      \xyboxmatrix"A31"@+.5cm{
        \ar@/^1pc/[r]^{f_0}="x1"\ar@/_1pc/[r]_{f_1}="y1"\ar[d]_{h_0}&
        \ar@/^1pc/[r]^{g_0}="x2"\ar@/_1pc/[r]_{g_1}="y2"\ar[d]|{h_0}\ar@2@<+1.5ex>[dl]**{}?(.3);?(.7)|{h_3}&
        \ar@2@<+1.5ex>[dl]**{}?(.3);?(.7)|{k_3}\ar[d]^{k_1}\\
        \ar@/_1pc/[r]_{f'_1}="w1"&
        \ar@/_1pc/[r]_{g'_1}="w2"&
        {}
        \POS
        \ar@2"x1";"y1"**{}?(.3);?(.7)|{f_2}="r"
        \ar@2"x2";"y2"**{}?(.3);?(.7)|{g_2}="s"
        \tria"s";"r"
      }
      ,"A22"      
      \xyboxmatrix"A22"@+.5cm{
        \ar@/^1pc/[r]^{f_0}="x1"\ar@/_1pc/[r]_{f_1}="y1"\ar[d]_{h_0}&
        \ar@/^1pc/[r]^{g_0}="x2"\ar@/_1pc/[r]_{g_1}="y2"\ar[d]|{h_0}\ar@2@<+1.5ex>[dl]**{}?(.3);?(.7)|{h_3}&
        \ar@2@<+1.5ex>[dl]**{}?(.3);?(.7)|{k_3}="X"\ar[d]^{k_1}\\
        \ar@/_1pc/[r]_{f'_1}="w1"&
        \ar@/_1pc/[r]_{g'_1}="w2"&
        {}
        \POS
        \ar@2"x1";"y1"**{}?(.3);?(.7)|{f_2}="r"
        \ar@2"x2";"y2"**{}?(.3);?(.7)|{g_2}="s"
        \tria"r";"X"
      }
      ,"A13" 
      \xyboxmatrix"A13"@+.5cm{
        \ar@/^1pc/[r]^{f_0}="x1"\ar[d]_{h_0}&
        \ar@/^1pc/[r]^{g_0}="x2"\ar[d]|{h_0}\ar@2@<-1.5ex>[dl]**{}?(.3);?(.7)|{h_2}&
        \ar@2@<-1.5ex>[dl]**{}?(.3);?(.7)|{k_2}\ar[d]^{k_1}\\
        \ar@/^1pc/[r]^{f'_0}="z1"\ar@/_1pc/[r]_{f'_1}="w1"&
        \ar@/^1pc/[r]^{g'_0}="z2"\ar@/_1pc/[r]_{g'_1}="w2"&
        {}
        \POS
        \ar@2"z1";"w1"**{}?(.3);?(.7)|{f'_2}="r1"
        \ar@2"z2";"w2"**{}?(.3);?(.7)|{g'_2}="s1"
        \tria"r1";"s1"
      }
      ,"A33"
      \xyboxmatrix"A33"@+.5cm{
        \ar@/^1pc/[r]^{f_0}="x1"\ar[d]_{h_0}&
        \ar@/^1pc/[r]^{g_0}="x2"\ar[d]|{h_0}\ar@2@<-1.5ex>[dl]**{}?(.3);?(.7)|{h_2}&
        \ar@2@<-1.5ex>[dl]**{}?(.3);?(.7)|{k_2}\ar[d]^{k_1}\\
        \ar@/^1pc/[r]^{f'_0}="z1"\ar@/_1pc/[r]_{f'_1}="w1"&
        \ar@/^1pc/[r]^{g'_0}="z2"\ar@/_1pc/[r]_{g'_1}="w2"&
        {}
        \POS
        \ar@2"z1";"w1"**{}?(.3);?(.7)|{f'_2}="r1"
        \ar@2"z2";"w2"**{}?(.3);?(.7)|{g'_2}="s1"
        \tria"s1";"r1"
      }
      \POS
      \ar@3"A11";"A31"^{(g'_1\#_0h_3)\\\#_1(k_3\#_0f_1)\\\#_1((k_1\#_0g_2)\ten{}f_2)}
      \ar@3"A11";"A22"^{(g'_1\#_0h_3)\\\#_1((k_3\#_1(k_1\#_0g_2))\ten{}f_2)}
      \ar@3"A31";"A22"^{(g'_1\#_0h_3)\\\#_1(k_3\ten{}f_2)\\\#_1((k_1\#_0g_2)\#_0f_0)}
      \ar@3"A22";"A33"^{(g'_1\#_1\underline{h_4})\#_1(\underline{k_4}\#_0f_0)}
      \ar@3"A13";"A33"_{(g'_2\ten{}(f'_2\#_0h_0))\\\#_1(g'_0\#_0h_2)\\\#_1(k_2\#_0f_0)}
      \ar@3"A11";"A13"_{\omega_1}
      \ar@3"A31";"A33"^{\omega_2}
    \end{xy}
  \end{equation}
\end{sidewaysfigure}

such that
\begin{lem}
  \label{lem:tenmapdefcomm}
  $\pair{h_\ell,h_r}_1= \pair{d_0^j,d_1^j}_1t_1$ in $\RGrph$.
\end{lem}
\begin{prf}
  One checks that $(h_\ell)_1=(d_0^jt)_1$ and $(h_r)_1=(d_1^jt)_1$ as
  graph maps using definitions (\ref{eq:lhcompdef}) and (\ref{eq:rhcompdef}). 
\end{prf}

\begin{lem}
  \label{lem:p2gray}
  The 3-globular set 
  \begin{equation*}
    \begin{xy}
      \xyboxmatrix{
        P^2(\H)
        \ar@<-1ex>[r]_-{p_0}\ar@<+1ex>[r]^-{p^1}&
        \overline{\overline{\H}}\ar@<-1ex>[r]_{d_0}\ar@<+1ex>[r]^{d^1}\ar[l]|-\Delta&
        \pathspc{\H}\ar@<-1ex>[r]_{d_0}\ar@<+1ex>[r]^{d^1}\ar[l]|i&
        \H\ar[l]|i&
      }
    \end{xy}
  \end{equation*}
  is an internal $\Gray$-category.
\end{lem}
\begin{prf}
  We already know that its three lower stages constitute a
  sesqui-catgory.  The three top
  parts are trivially a 2-category. The tensor map is given by
  \begin{equation*}
    \begin{xy}
      \xyboxmatrix{
        \overline{\overline{\H}}\times_{\overline{d_0},\overline{d_1}}\overline{\overline{\H}}\ar[r]^-{\pair{h_\ell,h_r}}|-*\dir{/}&
        P^2(\H)\\
      }
    \end{xy}
  \end{equation*}
  which satisfies the tensor axioms by construction.
\end{prf}

We can finally prove our desired theorem:
\begin{thm}
  \label{thm:intgrcat}
  Given a $\Gray$-category $\H$ there is an internal $\Gray$-category
  in $\Gray\Cat_{\fQ^1}$
  \begin{equation}
    \label{eq:mainthmintgrcat}
    \begin{xy}
      \xyboxmatrix{
        \overline{\overline{\overline{\H}}}\ar@<-1ex>[r]_{d_0}\ar@<+1ex>[r]^{d^1}&
        \overline{\overline{\H}}\ar@<-1ex>[r]_{d_0}\ar@<+1ex>[r]^{d^1}\ar[l]|i&
        \pathspc{\H}\ar@<-1ex>[r]_{d_0}\ar@<+1ex>[r]^{d^1}\ar[l]|i&
        \H\ar[l]|i&
      }
    \end{xy}
  \end{equation}
  with composition operations $m, \overline{m},
  \overline{\overline{m}}, w_\ell, w_r, \overline{w_\ell},
  \overline{w_r}, \tilde{w_\ell}, \tilde{w_r}, $ and tensor $t$.
\end{thm}
\begin{prf}
  We have a globular map
  \begin{equation*}
    \begin{xy}
      \xyboxmatrix{
        \overline{\overline{\overline{\H}}}\ar@<-1ex>[r]_{d_0}\ar@<+1ex>[r]^{d^1}\ar[d]_{\pair{d_0^j,d_1^j}}&
        \overline{\overline{\H}}\ar@<-1ex>[r]_{d_0}\ar@<+1ex>[r]^{d^1}\ar[l]|i\ar[d]&
        \pathspc{\H}\ar@<-1ex>[r]_{d_0}\ar@<+1ex>[r]^{d^1}\ar[l]|i\ar[d]&
        \H\ar[l]|i\ar[d]\\
        P^2(\H)
        \ar@<-1ex>[r]_-{p_0}\ar@<+1ex>[r]^-{p^1}&
        \overline{\overline{\H}}\ar@<-1ex>[r]_{d_0}\ar@<+1ex>[r]^{d^1}\ar[l]|-\Delta&
        \pathspc{\H}\ar@<-1ex>[r]_{d_0}\ar@<+1ex>[r]^{d^1}\ar[l]|i&
        \H\ar[l]|i&
      }
    \end{xy}
  \end{equation*}
  This globular map is an internal sesqui-functor in the lower and at
  the upper degrees, and by (\ref{eq:tenmapind}) it preverses the tensor:
  \begin{equation*}
    \begin{xy}
      \xyboxmatrix{
        \overline{\overline{\H}}\times_{\overline{d_0},\overline{d_1}}\overline{\overline{\H}}\ar[r]^-{t}|-*\dir{/}
        \ar[d]&
        \overline{\overline{\overline{\H}}}\ar[d]^{\pair{d_0^j,d_1^j}}\\
         \overline{\overline{\H}}\times_{\overline{d_0},\overline{d_1}}\overline{\overline{\H}}\ar[r]^-{\pair{h_\ell,h_r}}|-*\dir{/}&
         P^2(\H)
      }
    \end{xy}
  \end{equation*}
  Using the results of sections \ref{sec:compaths} and \ref{sec:hcells}
  this proves that (\ref{eq:mainthmintgrcat}) is an internal
  $\Gray$-category.
\end{prf}

\begin{lem}
  \label{lem:strfunopscom}
  The operations $\overline{\overline{m}}$, ${w_\ell}$, ${w_r}$,
  ${\tilde{w_\ell}}$, ${\tilde{w_r}}$,${\overline{w_\ell}}$,
  ${\overline{w_r}}$ and $t$ are natural with respect to strict
  $\Gray$-functors.
\end{lem}
\begin{prf}
  This can be shown using the universality of the respective
  constructions and the fact that $m$ is natural with respect to
  strict $\Gray$-functors, i.~e.\@ lemma \ref{lem:multstrictnat}.
\end{prf}

\section{The Internal Hom Functor}
\label{sec:inthomfct}

We can finally define the internal hom of $\Gray\Cat_{\fQ^1}$
\begin{multline}
  \label{eq:inthomfctdef}
  [\G,\H]\\=
  \begin{pmatrix}
    \begin{xy}
      \xyboxmatrix{
        \Gray\Cat_{\fQ^1}(\G,\overline{\overline{\overline\H}})\ar@<+1.5ex>[r]^{{d_1}_*}\ar@<-1.5ex>[r]_{{d_0}_*}&
        \Gray\Cat_{\fQ^1}(\G,\overline{\overline{\H}})\ar@<+1.5ex>[r]^{{d_1}_*}\ar@<-1.5ex>[r]_{{d_0}_*}\ar[l]|{i_*}&
        \Gray\Cat_{\fQ^1}(\G,\pathspc{{\H}})\ar@<+1.5ex>[r]^{{d_1}_*}\ar@<-1.5ex>[r]_{{d_0}_*}\ar[l]|{i_*}&
        \Gray\Cat_{\fQ^1}(\G,{{\H}})\ar[l]|{i_*} }
    \end{xy}
  \end{pmatrix}
\end{multline}
by applying $\Gray\Cat_{\fQ^1}(\G,-)$ to the diagram
(\ref{eq:mainthmintgrcat}), where the lower star means action by
post-composition in the co-Kleisli sense. This includes the various
induced composition operations $m_*$, $\overline{m}_*$,
$\overline{\overline{m}}_*$, ${w_\ell}_*$, ${w_r}_*$,
${\tilde{w_\ell}}_*$, ${\tilde{w_r}}_*$,${\overline{w_\ell}}_*$,
${\overline{w_r}}_*$ and $t_*$. Because $\Gray\Cat_{\fQ^1}(\G,-)$ by
definition preserves limits in the second variable, it takes internal
$\Gray$-categories in $\Gray\Cat_{\fQ^1}$ to such in $\Set$, that is,
to ordinary $\Gray$-categories. In analogy with our earlier notation
we write the compositions on $[\G,\H]$ as $*_n$ where $n$ is the
dimension of the incident cell, we use $*$ for the tensor of
transformations incident on a functor.

Explicitly, for example, given
\begin{equation*}
  \begin{xy}
    \xyboxmatrix{
      \G\ar@/^2pc/[r]^{G}|*\dir{/}="X"\ar[r]|(.3)*\dir{/}|{H}="Y"\ar@/_2pc/[r]_{K}|*\dir{/}="Z"&\H
      \ar@2"X";"Y"**{}?(.2);?(.8)|\alpha
      \ar@2"Y";"Z"**{}?(.2);?(.8)|\beta
    }
  \end{xy}
\end{equation*}
the composite $\beta*_0\alpha$ is defined as
\begin{equation*}
  \begin{xy}
    \xyboxmatrix{
      \G\ar[r]^-{\pair{\beta,\alpha}}|-*\dir{/}&\pathspc\H\times_{d_0,d_1}\pathspc\H\ar[r]^-{m}|-*\dir{/}&\pathspc\H
    }
  \end{xy}\,.
\end{equation*}
that is, $\beta*_0\alpha=m\fQ^1\pair{\beta,\alpha}d$.

To be slightly more explicit, at the level if 0-, and 1-cells of
$[\G,\H]$, that is, pseudo-functors and transformations the
composition  works as follows:
\begin{equation*}
  \begin{xy}
    \xymatrix@C-1.8cm{
      {}&{}&\Gray\Cat_{\fQ^1}(\G,\pathspc\H)\times_{{d_0}_*,{d_1}_*}\Gray\Cat_{\fQ^1}(\G,\pathspc\H)\ar[d]^{\iso}&{}&{}\\
      {}&{}&\Gray\Cat_{\fQ^1}(\G,\pathspc\H\times_{d_0,d_1}\pathspc\H)\ar[dl]\ar[dr]\ar[d]^{m_*}&{}&{}\\
      {}&\Gray\Cat_{\fQ^1}(\G,\pathspc\H)\ar[dl]_{{d_0}_*}\ar[dr]^(.7){{d_1}_*}&
      \Gray\Cat_{\fQ^1}(\G,\pathspc\H)\ar[dll]^(.7){{d_0}_*}|!{"2,2";"3,3"}\hole\ar[drr]_(.7){{d_1}_*}|!{"2,4";"3,3"}\hole&
      \Gray\Cat_{\fQ^1}(\G,\pathspc\H)\ar[dl]_(.7){{d_0}_*}\ar[dr]^{{d_1}_*}&\\
      \Gray\Cat_{\fQ^1}(\G,\H)&&\Gray\Cat_{\fQ^1}(\G,\H)&&\Gray\Cat_{\fQ^1}(\G,\H)
    }
  \end{xy}
\end{equation*}

\begin{rem}
  The $\Gray$-category $[\G,\H]$ is a $\Gray$-groupoid if $\H$ is one.
\end{rem}

\begin{thm}
  \label{thm:mapspcfunc}
  Given a morphism $F\from\G'\laxto\G$ in $\Gray\Cat_{\fQ^1}$, the map
  $$F^*=[F,\H]\from [\G,\H]\to[\G',\H]$$ acting by pre-composition in
  the co-Kleisli sense is
  a  $\Gray$-functor, that is, a strict morphism.
\end{thm}
\begin{prf}
  Assume a situation 
  \begin{xy}
    \xymatrix{
      \G'\ar[r]|*\dir{/}^F&\G\ar@/^2pc/[r]^{G}|*\dir{/}="X"\ar[r]|(.3)*\dir{/}|{H}="Y"\ar@/_2pc/[r]_{K}|*\dir{/}="Z"&\H
      \ar@2"X";"Y"**{}?(.2);?(.8)|\alpha
      \ar@2"Y";"Z"**{}?(.2);?(.8)|\beta
    }
  \end{xy}
  then we have
  \begin{multline*}
    F^*(\beta*_0\alpha)=(\beta*_0\alpha)F=m\pair{\beta,\alpha}F\\
    =m\pair{\beta{}F,\alpha{}F}=(\beta{}F)*_0(\alpha{}F)=(F^*\beta{})*_0(F^*\alpha{})\,.
  \end{multline*}
  Also, for identity transformations we have:
  \begin{equation*}
    F^*\id_G=iGF=\id_{GF}\,,
  \end{equation*}
  hence $F^*$ is a functor. By the same reasoning the higher
  operations including the tensor, are preserved as well. 
\end{prf}

\begin{rem}
  \label{rem:homstropact}
  This way $[-,\H]\from\Gray\Cat_{\fQ^1}^\op\to\Gray\Cat_{\fQ^1}$ is a
  functor for each $\H$.
\end{rem}

\begin{thm}
  \label{thm:mapspcfuncpost}
  Given a strict morphism $F\from\H\to\H'$ in $\Gray\Cat$, the map
  $$F_*=[\G,F]\from [\G,\H]\to[\G,\H']$$ acting by post-composition is
  a  $\Gray$-functor, that is, a strict morphism.
\end{thm}
\begin{prf}
  Assume a situation 
  \begin{xy}
    \xymatrix{
      \G\ar@/^2pc/[r]^{G}|*\dir{/}="X"\ar[r]|(.3)*\dir{/}|{H}="Y"\ar@/_2pc/[r]_{K}|*\dir{/}="Z"&
      \H\ar[r]|*\dir{/}^F&\H'
      \ar@2"X";"Y"**{}?(.2);?(.8)|\alpha
      \ar@2"Y";"Z"**{}?(.2);?(.8)|\beta }
  \end{xy}
  then we have
  \begin{multline*}
    F*(\beta*_0\alpha)=\pathspc{F}m\fQ^1\pair{\beta,\alpha}d=m\fQ^1(\pathspc{F}\times\pathspc{F})\fQ^1\pair{\beta,\alpha}d\\
    =m\fQ^1(\pair{\pathspc{F}\beta,\pathspc{F}\alpha})d= (F*\beta)*_0(F*\alpha)\,,
   \end{multline*}
   where we use lemma \ref{lem:strfunopscom}. 
   Also, for identity transformations we have:
   \begin{equation*}
    \pathspc{F}*\id_G=\pathspc{F}iG=iFG=\id_{F*G}
  \end{equation*}
  hence $F^*$ is a functor. 

  The other operations are preserved similarly by applying
  lemma \ref{lem:strfunopscom}. 
\end{prf}

We now proceed to constructing the restricted mapping space $\{\G,\H\}$.
We pull back all the parts of (\ref{eq:inthomfctdef}) along $e^*$
given in (\ref{eq:strictinclude}) to obtain
\begin{equation}
  \label{eq:malleable}
    \begin{xy}
      \xyboxmatrix{
        \{\G,\H\}_3\ar@<+1.5ex>[r]^{{d_1}_*}\ar@<-1.5ex>[r]_{{d_0}_*}\ar[d]_{\overline{e^*}}&
        \{\G,\H\}_2\ar@<+1.5ex>[r]^{{d_1}_*}\ar@<-1.5ex>[r]_{{d_0}_*}\ar[l]|{i_*}\ar[d]_{\overline{e^*}}&
        \{\G,\H\}_1\ar@<+1.5ex>[r]^{{d_1}_*}\ar@<-1.5ex>[r]_{{d_0}_*}\ar[l]|{i_*}\ar[d]_{\overline{e^*}}&
        \GC(\G,\H)\ar[l]|{i_*} \ar[d]^{e^*}\\
        \Gray\Cat_{\fQ^1}(\G,\overline{\overline{\overline\H}})\ar@<+1.5ex>[r]^{{d_1}_*}\ar@<-1.5ex>[r]_{{d_0}_*}&
        \Gray\Cat_{\fQ^1}(\G,\overline{\overline{\H}})\ar@<+1.5ex>[r]^{{d_1}_*}\ar@<-1.5ex>[r]_{{d_0}_*}\ar[l]|{i_*}&
        \Gray\Cat_{\fQ^1}(\G,\pathspc{{\H}})\ar@<+1.5ex>[r]^{{d_1}_*}\ar@<-1.5ex>[r]_{{d_0}_*}\ar[l]|{i_*}&
        \Gray\Cat_{\fQ^1}(\G,{{\H}})\ar[l]|{i_*} 
      }
    \end{xy}\,,
\end{equation}
and we set $\{\G,\H\}_0=\GC(\G,\H)$. We call $\{\G,\H\}_1$ the set of
malleable transformations, c.~f.\@ definition \ref{defn:malltranfs}. Obviously the
left and right actions of strict functors described in theorems
\ref{thm:mapspcfuncpost} and \ref{thm:mapspcfunc} restrict to the
restricted mapping space. 

Hence for strict morphisms $F\from\G'\to\G$ and $G\from\H\to\H'$ we
get a commuting square of $\Gray$-functors
\begin{equation*}
  \begin{xy}
    \xyboxmatrix{
      \{\G,\H\}\ar[r]^{F^*}\ar[d]_{G_*}&\{\G',\H\}\ar[d]^{G_*}\\
      \{\G,\H'\}\ar[r]_{F^*}&\{\G',\H'\}\\
    }
  \end{xy}\,.
\end{equation*}

In conclusion, we get the following interesting structure on
$\Gray\Cat$, and leave the question as to further, higher structure
open:
\begin{thm}
  \label{thm:graycatsesqui}
  The category $\Gray\Cat$ of $\Gray$-categories, strict
  $\Gray$-functors and malleable transformations is a sesquicategory. \qed
\end{thm}

\begin{rem}
  \label{rem:mallinclude}
  By section \ref{sec:basfib} $\{\G,\H\}$ is a $\Gray$-category and
  $\overline{e^*}\from\{\G,\H\}\to[\G,\H]$ is a strict
  $\Gray$-functor.
\end{rem}

For $\G$ free up to order 1 the maps $e$ and $k$ discussed in
(\ref{eq:freeorder1split}) give natural transformations
\begin{equation*}
  \begin{xy}
    \xyboxmatrix{
      \GC(\G,\_)\ar[r]^{e^*}\ar@/_1.5pc/[rr]_{\GC(\G,\_)}&\GC_{\fQ^1}(\G,\_)\ar[r]^{k^*}&\GC(\G,\_)
    }
  \end{xy}\,,
\end{equation*}
where the maps act by precomposition in $\Gray\Cat$.

\begin{lem}
  Given a $\Gray$-category $\G$ free up to order 1 there are
  canonical transformations
  \begin{equation*}
    \begin{xy}
      \xyboxmatrix{
        \G\ar@/_1.5pc/[r]_{Fk}="x"\ar@/^1.5pc/[r]^{F}="y"|*\dir{/}&\H
        \POS\ar@2 "y";"x"**{}?(.3);?(.7)^\rho
      }
    \end{xy}
  \end{equation*}
  that is the identity on objects.\footnote{I.~e.\@ basically icons in the sense
  of \citet{lack2007}, except our constraint 2-cell points the other way.}
\end{lem}
\begin{prf}
  We need to give a $\fQ^1$ graph map  $\rho\from\G\laxto\pathspc\H$
  with $d_1\rho=Fke$ and $d_0\rho=F$:
  \begin{enumerate}
  \item 0-cells
    \begin{equation*}
      x\mapsto \xymatrix@1{x\ar[r]^{\id_x}&x}
    \end{equation*}
  \item 1-cells
    \begin{equation*}
      f=[f_1,\ldots,f_n]\mapsto
    \begin{xy}
      \xyboxmatrix@+1.5cm{
        \ar[r]^{\id_x}\ar[d]_{Ff}&\ar[d]^{F[f_1]\#_0\cdots\#_0F[f_n]}
        \ar@2[dl]**{}?(.2);?(.8)|{F^2_{[f_1],\ldots,[f_n]}}\\
        \ar[r]_{\id_y}&
      }
    \end{xy}
  \end{equation*}
  \item 2-cells
    \begin{equation*}
    (\alpha\from{}f\To{}f')\mapsto
    \begin{xy}
      \save
      (70,0):
      (0,0)="A"
      ,(1,0)="B"
      \restore\POS
      ,"A"
      \xyboxmatrix"A"@+1.5cm{
        \ar[r]^{\id_x}\ar[d]|{Ff}="x"\ar@/_4pc/[d]_{Ff'}="y"&
        \ar[d]^{F[f_1]\\\#_0\cdots\\\#_0F[f_n]}
        \ar@2[dl]**{}?(.2);?(.8)|{F^2_{[f_1],\ldots,[f_n]}}\\
        \ar[r]_{\id_y}&{}
        \ar@2"x";"y"**{}?(.3);?(.7)_{F\alpha}
      }
      ,"B"
      \xyboxmatrix"B"@+1.5cm{
        \ar[r]^{\id_x}\ar[d]_{Ff'}&
        \ar[d]|{F[f'_1]\\\#_0\cdots\\\#_0F[f'_{n'}]}="y"\ar@/^4pc/[d]^{F[f_1]\\\#_0\cdots\\\#_0F[f_n]}="x"
        \ar@2[dl]**{}?(.2);?(.8)|(.6){F^2_{[f_1],\ldots,[f_n]}}\\
        \ar[r]_{\id_y}&
        {} \ar@2"x";"y"**{}?(.3);?(.7)_{\omega}="X"
      }
      \POS\ar@{=}"A";"B"
    \end{xy}
  \end{equation*}
    where
  $\omega$ is
    $\overline{F^2_{[f'_1],\ldots,[f'_{n'}]}}\#_1F\alpha\#_1F^2_{[f_1],\ldots,[f_{n}]}$.
  \item 3-cells
    \begin{equation*}
      (\Gamma\from\alpha\Tto\alpha')\mapsto
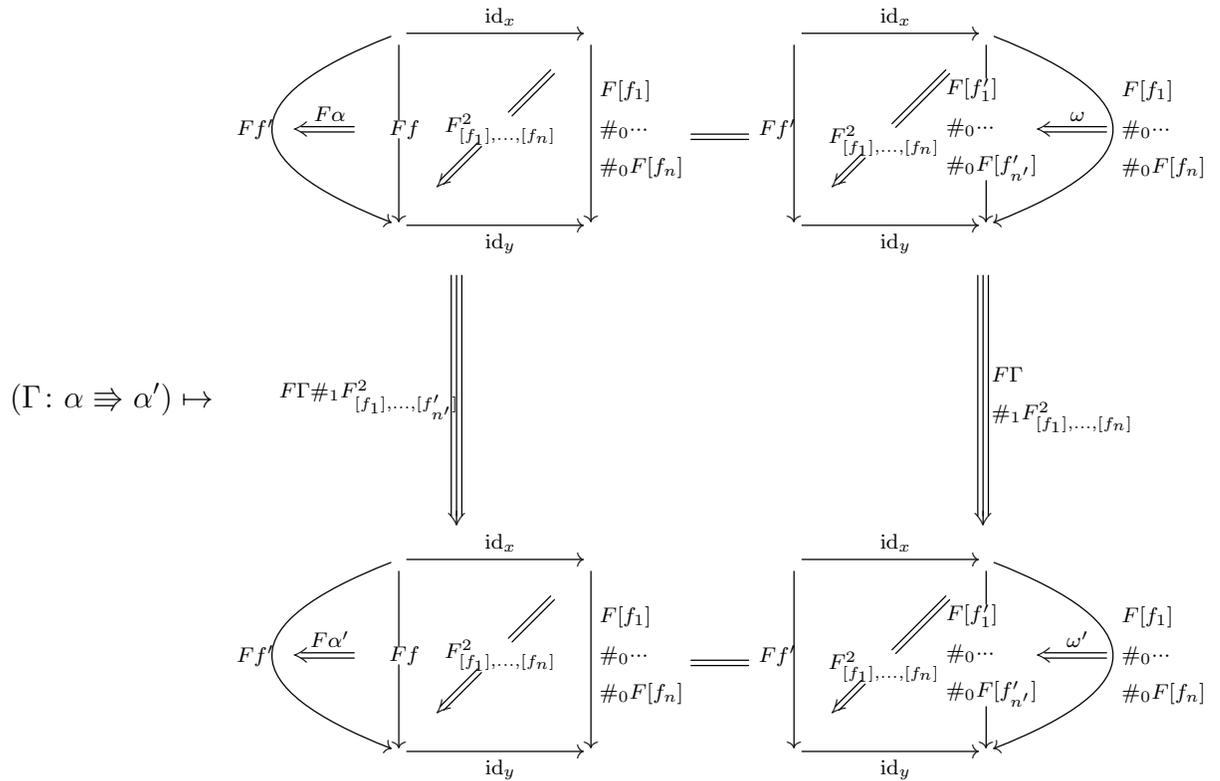
      \begin{matrix}
        \begin{xy}
          \save (70,0):(0,-1):: (0,0)="A" ,(1,0)="B" ,(0,1)="A1"
          ,(1,1)="B1" \restore\POS ,"A" \xyboxmatrix"A"@+1.5cm{
            \ar[r]^{\id_x}\ar[d]|{Ff}="x"\ar@/_4pc/[d]_{Ff'}="y"&
            \ar[d]^{F[f_1]\\\#_0\cdots\\\#_0F[f_n]}
            \ar@2[dl]**{}?(.2);?(.8)|{F^2_{[f_1],\ldots,[f_n]}}\\
            \ar[r]_{\id_y}&{} \POS\ar@2"x";"y"**{}?(.3);?(.7)_{F\alpha}
          } 
          ,"B" 
          \xyboxmatrix"B"@+1.5cm{ \ar[r]^{\id_x}\ar[d]_{Ff'}&
            \ar[d]|{F[f'_1]\\\#_0\cdots\\\#_0F[f'_{n'}]}="y"\ar@/^4pc/[d]^{F[f_1]\\\#_0\cdots\\\#_0F[f_n]}="x"
            \ar@2[dl]**{}?(.2);?(.8)|(.6){F^2_{[f_1],\ldots,[f_n]}}\\
            \ar[r]_{\id_y}& {}\POS
            \ar@2"x";"y"**{}?(.3);?(.7)_{\omega}="X" 
          }
          ,"A1"
          \xyboxmatrix"A1"@+1.5cm{
            \ar[r]^{\id_x}\ar[d]|{Ff}="x"\ar@/_4pc/[d]_{Ff'}="y"&
            \ar[d]^{F[f_1]\\\#_0\cdots\\\#_0F[f_n]}
            \ar@2[dl]**{}?(.2);?(.8)|{F^2_{[f_1],\ldots,[f_n]}}\\
            \ar[r]_{\id_y}&{} \POS\ar@2"x";"y"**{}?(.3);?(.7)_{F\alpha'}
          } 
          ,"B1" 
          \xyboxmatrix"B1"@+1.5cm{ 
            \ar[r]^{\id_x}\ar[d]_{Ff'}&
            \ar[d]|{F[f'_1]\\\#_0\cdots\\\#_0F[f'_{n'}]}="y"
            \ar@/^4pc/[d]^{F[f_1]\\\#_0\cdots\\\#_0F[f_n]}="x"
            \ar@2[dl]**{}?(.2);?(.8)|(.6){F^2_{[f_1],\ldots,[f_n]}}\\
            \ar[r]_{\id_y}& {} \POS
            \ar@2"x";"y"**{}?(.3);?(.7)_{\omega'}="X" 
          }
          \POS\ar@{=}"A";"B"\ar@{=}"A1";"B1"
          \ar@3"A";"A1"_{F\Gamma\#_1F^2_{[f_1],\ldots,[f'_{n'}]}}
          \ar@3"B";"B1"^{F\Gamma\\\#_1F^2_{[f_1],\ldots,[f_{n}]}}
        \end{xy}
      \end{matrix}
    \end{equation*}
  where $\omega=\overline{F^2_{[f'_1],\ldots,[f'_{n'}]}}\#_1F\alpha\#_1F^2_{[f_1],\ldots,[f_{n}]}$
    and
    $\omega'=\overline{F^2_{[f'_1],\ldots,[f'_{n'}]}}\#_1F\alpha'\#_1F^2_{[f_1],\ldots,[f_{n}]}$.
  \item For a composable pair of 1-cells $f', f$ a 2-cocycle element
    \begin{equation*}
      \begin{xy} 
        \xyboxmatrix"A"@+1.5cm{
          \ar[r]^{\id_x}\ar[d]|{Ff}="x"\ar@/_4pc/[dd]_{F(f'\#_0f)}|{}="y"&
          {}\ar[d]^{F[f_1]\\\#_0\cdots\\\#_0F[f_{n}]}\ar@2[dl]**{}?(.2);?(.8)|{F^2_{[f_1],\ldots,[f_{n}]}}
          \\
          {}\ar[r]_{\alpha_y}\ar[d]|{Ff'}&
          {}\ar@2[dl]**{}?(.2);?(.8)|{F^2_{[f'_1],\ldots,[f'_{n'}]}}\ar[d]^{F[f'_1]\\\#_0\cdots\\\#_0F[f'_{n'}]}
          \\
          \ar[r]_{\alpha_z}&{}\ar@2"2,1";"y"**{}?<(.2);?(.8)^{F^2_{f',f}}
        }
        \POS + (80,0)
        \xyboxmatrix"B"@+1.5cm{\ar[r]^{\alpha_x}\ar[dd]_{F(f'\#_0f)}&
          \ar[dd]|{F[f'_1]\#_0\cdots\#_0F[f'_{n'}]\\\#_0F[f_1]\#_0\cdots\#_0F[f_n]}="x"
          \ar@2[ddl]**{}?(.2);?(.8)|(.65){F^2_{[f'_1],\ldots,[f'_{n'}],[f_1],\ldots,[f_{n}]}}
          \ar[dr]^{F[f_1]\#_0\cdots\#_0F[f_n]}&\\ 
          {}&{}&\ar[dl]^{F[f'_1]\#_0\cdots\#_0F[f'_{n'}]}\\
          {}\ar[r]_{\alpha_z}&{}
          \POS\ar@{=}@2"2,3";"x"**{}?<(.2);?>(.8)_(.8){}
        }
        \ar@{=} "A";"B"
      \end{xy}\,.
    \end{equation*}
    The equation holds by \ref{eq:tentriv} and \ref{eq:cocycle}.
  \end{enumerate}
    The verification that this is a $\fQ^1$-graph map is straightforward.
\end{prf}

\section{Putting it all together}
\label{sec:together}

\begin{defn}\label{defn:ltransf} A \defterm{lax transformation}
  $\alpha\from{}F\to{}G$ between pseudo-functors $F,G\from\G\laxto\H$
  of $\Gray$-categories is a pseudo-functor
  $\alpha\from\G\laxto\pathspc\H$ such that $d_0\alpha=F$ and
  $d_1\alpha=G$.
\end{defn}
\begin{defn}
  \label{defn:malltranfs} A \defterm{malleable transformation}
  $\alpha\from{}F\to{}G$ between strict functors $F,G\from\G\to\H$
  of $\Gray$-categories is a pseudo-functor
  $\alpha\from\G\laxto\pathspc\H$ such that $d_0\alpha=F$ and
  $d_1\alpha=G$.
\end{defn}
This was introduced in (\ref{eq:malleable}).

\begin{rem} Using the definition of path spaces in definition \ref{defn:pathspc}
  and the characterization of pseudo-maps in definition \ref{defn:psgrmap} we
  note for reference that a lax transformation $\alpha$ is given by
  the following underlying data:
  \begin{enumerate}
  \item for each 0-cell $x$ of $\G$ a 1-cell
    $\alpha_x\from{}Fx\to{}Gx$,
  \item for each 1-cell $f\from{}x\to{}y$ of $\G$ a 2-cell
    \begin{equation*}
      \begin{xy}\xycompile{
          \xyboxmatrix{Fx\ar[r]^{\alpha_x}\ar[d]_{Ff}&
            Gx\ar[d]^{Gf}\ar@2[dl]**{}?(.2);?(.8)^{\alpha_f}\\
            Fy\ar[r]_{\alpha_y}&Gy }}
      \end{xy}
    \end{equation*}
  \item for each 2-cell $g\from{}f\to{}f'$ of $\G$ a 3-cell of $\H$
    \begin{equation*}
      \begin{matrix}
        \begin{xy} \xycompile{
            \xyboxmatrix"A"{Fx\ar[r]^{\alpha_x}\ar[d]|{Ff}="x"\ar@/_3pc/[d]_{Ff'}="y"&
              Gx\ar[d]^{Gf}\ar@2[dl]**{}?(.2);?(.8)|{\alpha_f}\\Fy
              \ar[r]_{\alpha_y}&Gy\POS\ar@2"x";"y"**{}?<(.3);?>(.7)^{Fg}
            } \POS + (50,0)
            \xyboxmatrix"B"{Fx\ar[r]^{\alpha_x}\ar[d]_{Ff'}&
              Gx\ar[d]|{Gf'}="x"\ar@2[dl]**{}?(.2);?(.8)|{\alpha_{f'}}
              \ar@/^3pc/[d]^{Gf}="y"\\Fy\ar[r]_{\alpha_y}&Gy
              \POS\ar@2"y";"x"**{}?<(.3);?>(.7)^{Gg} } \ar@3 "A";"B"
            ^{\alpha_g} }
        \end{xy}
      \end{matrix}
    \end{equation*}
  \item for each pair of composable 1-cells $f\from{}x\to{}y$,
    $f'\from{}y\to{}z$ an invertible 3-cell
    \begin{equation*}
      \begin{xy} \xycompile{
          \xyboxmatrix"A"{Fx\ar[r]^{\alpha_x}\ar[d]|{Ff}="x"\ar@/_4pc/[dd]_{F(f'\#_0f)}|{}="y"&
            Gx\ar[d]^{Gf}\ar@2[dl]**{}?(.2);?(.8)|{\alpha_f}\\Fy
            \ar[r]_{\alpha_y}\ar[d]|{Ff'}&Gy\ar@2[dl]**{}?(.2);?(.8)|{\alpha_{f'}}\ar[d]^{Gf'}\\
            Fz\ar[r]_{\alpha_z}&Gz\POS\ar@2"2,1";"y"**{}?<(.2);?(.8)^{F^2_{f',f}}
          } \POS + (70,0)
          \xyboxmatrix"B"{Fx\ar[r]^{\alpha_x}\ar[dd]_{F(f'\#_0f)}&
            Gx\ar[dd]|{G(f'\#_0f)}="x"\ar@2[ddl]**{}?(.2);?(.8)|(.65){\alpha_{f'\#_0f}}
            \ar[dr]^{Gf}&\\ 
            {}&{}&Gy\ar[dl]^{Gf'}\\
            Fz\ar[r]_{\alpha_z}&Gz
            \POS\ar@2"2,3";"x"**{}?<(.2);?>(.8)_(.8){G^2_{f',f}}
          }
          \ar@3 "A";"B"^-{\alpha^2_{f',f}} 
        }
      \end{xy}\,.
    \end{equation*}
  \end{enumerate}
  
  Furthermore, these data have to satisfy the following equations:
  \begin{enumerate}
  \item On identities of 0-cells:
    \begin{equation*}
      \alpha_{\id_x}=\id_{\alpha_x}
    \end{equation*}
  \item  for each 3-cell
    $\Gamma\from{}g\to{}g'$ the square of 3-cells in $\H$
    \begin{equation*}
      \begin{xy}\xycompile{ \xyboxmatrix"A"{Fx \ar[r]^{\alpha_x}
            \ar[d]|{Ff}="x"\ar@/_3pc/[d]|{Ff'}="y" & Gx\ar[d]^{Gf} \ar@2 [dl] **{}
            ?(.2); ?(.8) |{\alpha_f}\\Fy \ar[r]_{{\alpha_y}}& Gy \POS \ar@2 "x";
            "y" **{} ?(.2); ?(.8) |{Fg} } \POS + (50,0) \xyboxmatrix"B"{Fx
            \ar[r]^{\alpha_x} \ar[d]_{Ff'}& Gx\ar[d]|{Gf'}="x" \ar@2 [dl] **{}
            ?(.2); ?(.8) |{\alpha_{f'}} \ar@/^3pc/[d]|{Gf}="y"\\Fy
            \ar[r]_{{\alpha_y}}& Gy \POS \ar@2 "y"; "x" **{} ?(.2); ?(.8) |{Gg}}
          \POS + (-50,-50) \xyboxmatrix"C"{Fx \ar[r]^{\alpha_x}
            \ar[d]|{Ff}="x"\ar@/_3pc/[d]|{Ff'}="y" & Gx\ar[d]^{Gf} \ar@2 [dl] **{}
            ?(.2); ?(.8) |{\alpha_f}\\Fy \ar[r]_{{\alpha_y}}& Gy \POS \ar@2 "x";
            "y" **{} ?(.2); ?(.8) |{Fg'} } \POS + (50,0) \xyboxmatrix"D"{Fx
            \ar[r]^{\alpha_x} \ar[d]_{Ff'}& Gy\ar[d]|{Gf'}="x" \ar@2 [dl] **{}
            ?(.2); ?(.8) |{\alpha_{f'}} \ar@/^3pc/[d]|{Gf}="y"\\Fy
            \ar[r]_{{\alpha_y}}& Gy \POS \ar@2 "y"; "x" **{} ?(.2); ?(.8) |{Gg'}}
          \ar@3 "A";"B" ^{\alpha_g} \ar@3 "C";"D" _{\alpha_{g'}} \ar@3 "A";"C"
          _{({\alpha_y}\#_0F\Gamma)\#_1\alpha_f} \ar@3 "B";"D"
          ^{\alpha_{f'}\#_1(G\Gamma\#_0\alpha_x)} }\end{xy}
    \end{equation*}commutes. This condition obviously comes from the definition
    of 3-cells in the path space.
  \item  For every pair $g\from{}f\To{}f', g'\from{}f'\To{}f''$:
    \begin{equation*}
      \begin{xy}\xycompile{
          \xyboxmatrix"A"{
            Fx \ar[r]^{\alpha_x}\ar[d]|{Ff}="x"
            \ar@/_2.5pc/[d]|{Ff'}="y"
            \ar@/_5pc/[d]|{Ff''}="z" &
            Gx\ar[d]^{Gf} \ar@2 [dl] **{} ?(.2); ?(.8) |{\alpha_f}\\
            Fy \ar[r]_{\alpha_y}& Gy
            \POS \ar@2 "x"; "y" **{} ?(.2); ?(.8)|{Fg} 
            \POS \ar@2 "y"; "z" **{} ?(.2); ?(.8)|{Fg'} 
          } 
          \POS + (55, 0)
          \xyboxmatrix"B"{
            Fx\ar[r]^{\alpha_x} \ar[d]|{Ff'}="x" \ar@/_2.5pc/[d]|{Ff''}="y" &
            Gx\ar[d]|{Gf'}="w" \ar@/^2.5pc/[d]|{Gf}="z" 
            \ar@2 [dl] **{} ?(.2); ?(.8) |{\alpha_{f'}}\\
            Fx\ar[r]_{\alpha_{y}}& Gy 
            \POS \ar@2 "x"; "y" **{} ?(.2); ?(.8) |{Fg'} 
            \POS \ar@2 "z"; "w" **{} ?(.2); ?(.8) |{Gg} 
          } 
          \POS + (55, 0)
          \xyboxmatrix"C"{
            Fx \ar[r]^{\alpha_x} \ar[d]_{Ff''}  &
            Gx\ar[d]|{Gf''}="z" \ar@2 [dl] **{} ?(.2); ?(.8) |{\alpha_{f''}}
            \ar@/^2.5pc/[d]|{Gf'}="y" \ar@/^5pc/[d]|{Gf}="x" \\
            Fy\ar[r]_{\alpha_y}& Gy 
            \POS \ar@2 "x"; "y" **{} ?(.2); ?(.8)|{Gg} 
            \POS \ar@2 "y"; "z" **{} ?(.2); ?(.8)|{Gg'} 
          }
          \ar@3 "A";"B" |{}="x"\save\POS"x"+(0,13)*!C\labelbox{(\alpha_y\#_0Fg')\#_1\alpha_g}\ar@{.}"x";c\restore
          \ar@3 "B";"C" |{}="x"\save\POS"x"+(0,13)*!C\labelbox{\alpha_{g'}\#_1(Gg\#_0\alpha_x)}\ar@{.}"x";c\restore
          \ar@3@/_5pc/ "A";"C" _{\alpha_{g'\#_1g}}
        }\end{xy}\,,
    \end{equation*}
    and for identity 2-cells $\id_f\from{}f\To{}f$ we have an identity
    3-cell
    \begin{equation*}
      \alpha_{\id_f}=\id_{\alpha_f}\,.
    \end{equation*}
  \item The family of 3-cells has to satisfy a kind of cocycle
    condition: For a composable triple $f,f',f''$ of 1-cells
    $\alpha^2$ has to satisfy equation (\ref{eq:transex2cocy}).
    \begin{sidewaysfigure}
      \begin{equation}
        \label{eq:transex2cocy}
        \begin{xy}
	  \save (85,0):(0,-.7):: 
	  ,(0,0)="A11" 
	  ,(1,0)="A21" 
	  ,(2,0)="A31"
	  ,(3,0)="A41" 
	  ,(0,1)="A12" 
	  ,(1,1)="A22" 
	  ,(2,1)="A32" 
	  ,(3,1)="A42"
	  ,"A11";"A21"**{}?<>(.5)+(0,-.7)="K1"
	  ,"A12";"A22"**{}?<>(.5)+(0,+.7)="K2"
	  \restore 
	  \POS
	  ,"K1"
	  \xyboxmatrix"K1"{
	    Fx\ar[r]^{\alpha_x}\ar[d]|{Ff}="x"\ar@/_4pc/[dd]|{F(f'\#_0f)}="y"\ar@/_6pc/[ddd]|{F(f''\#_0f'\#_0f)}="z"&
	    Gx\ar[d]^{Gf}\ar@2[dl]**{}?(.2);?(.8)|{\alpha_f}\\
	    Fy\ar[r]_{\alpha_y}\ar[d]|{Ff'}&Gy\ar@2[dl]**{}?(.2);?(.8)|{\alpha_{f'}}\ar[d]^{Gf'}\\
	    Fz\ar[r]_{\alpha_z}\ar[d]|{Ff''}&Gz\ar@2[dl]**{}?(.2);?(.8)|{\alpha_{f''}}="Y"
	    \ar[d]^{Gf''}\\
            Fw\ar[r]_{\alpha_w}&Gw
            \ar@2"2,1";"y"**{}?<(.2);?>(.8)^{F^2_{f',f}}="X"
            \ar@2"3,1";"z"**{}?<(.2);?>(.8)^{F^2_{f'',f'\#_0f}} 
            \tria"X";"Y"
          }
          ,"A11" 
          \xyboxmatrix"A11"{
            Fx\ar[r]^{\alpha_x}\ar[d]|{Ff}="x"\ar@/_4pc/[dd]|{F(f'\#_0f)}="y"\ar@/_6pc/[ddd]|{F(f''\#_0f'\#_0f)}="z"&
            Gx\ar[d]^{Gf}\ar@2[dl]**{}?(.2);?(.8)|{\alpha_f}\\
            Fy\ar[r]_{\alpha_y}\ar[d]|{Ff'}&Gy\ar@2[dl]**{}?(.2);?(.8)|{\alpha_{f'}}\ar[d]^{Gf'}\\
            Fz\ar[r]_{\alpha_z}\ar[d]|{Ff''}&Gz\ar@2[dl]**{}?(.2);?(.8)|{\alpha_{f''}}="Y"
            \ar[d]^{Gf''}\\
            Fw\ar[r]_{\alpha_w}&Gw
            \ar@2"2,1";"y"**{}?<(.2);?>(.8)^{F^2_{f',f}}="X"
            \ar@2"3,1";"z"**{}?<(.2);?>(.8)^{F^2_{f'',f'\#_0f}} 
            \tria"Y";"X"
          },"A21"
          \xyboxmatrix"A21"{
            Fx\ar[r]^{\alpha_x}\ar[dd]_{F(f'\#_0f)}="x"\ar@/_5pc/[ddd]|(.4){F(f''\#_0f'\#_0f)}="z"&
            Gx\ar[dd]|(.6){G(f'\#_0f)}="y"\ar@2[ddl]**{}?(.2);?(.8)|(.4){\alpha_{f'\#_0f}}
            \ar[dr]^{Gf}&\\
            {}&&Gy\ar[dl]^{Gf'}\\
            Fz\ar[r]_{\alpha_z}\ar[d]|{Ff''}&Gz\ar@2[dl]**{}?(.2);?(.8)|{\alpha_{f''}}
            \ar[d]^{Gf''}&\\
            Fw\ar[r]_{\alpha_w}&Gw&
            \ar@2"2,3";"y"**{}?<(.2);?>(.8)_(.8){G^2_{f',f}}
            \ar@2"3,1";"z"**{}?<(.2);?>(.8)^{F^2_{f'',f'\#_0f}} 
          },"A31"
          \xyboxmatrix"A31"{
            Fx\ar[r]^{\alpha_x}\ar[ddd]_{F(f''\#_0f'\#_0f)}="x"&
            Gx\ar[ddd]|(.7){G(f''\#_0f'\#_0f)}="z"\ar@2[dddl]**{}?(.2);?(.8)|(.5){\alpha_{f''\#_0f'\#_0f}}
            \ar[dddr]|{G(f'\#_0f)}="y"\ar[r]^{Gf}&Gy\ar[ddd]^{Gf'}\\
            {}&&\\
            &&\\
            Fw\ar[r]_{\alpha_w}&Gw&Gz\ar[l]^{Gf''}
            \ar@2"1,3";"y"**{}?<(.2);?>(.8)_(.5){G^2_{f',f}}
            \ar@2"4,3";"z"**{}?<(.2);?>(.8)|{}="k"\save"k"+(-5,-15)*!C\labelbox{G^2_{f'',f'\#_0f}};"k"**\dir{.}\restore 
          }
          ,"A12" 
          \xyboxmatrix"A12"{
            Fx\ar[r]^{\alpha_x}\ar[d]|{Ff}="x"\ar@/_6pc/[ddd]|{F(f''\#_0f'\#_0f)}="z"&
            Gx\ar[d]^{Gf}\ar@2[dl]**{}?(.2);?(.8)|{\alpha_f}\\
            Fy\ar[r]_{\alpha_y}\ar[d]|{Ff'}\ar@/_4pc/[dd]|{F(f''\#_0f')}="y"&
            Gy\ar@2[dl]**{}?(.2);?(.8)|{\alpha_{f'}}\ar[d]^{Gf'}\\
            Fz\ar[r]_{\alpha_z}\ar[d]|{Ff''}&Gz\ar@2[dl]**{}?(.2);?(.8)|{\alpha_{f''}}
            \ar[d]^{Gf''}\\
            Fw\ar[r]_{\alpha_w}&Gw
            \ar@2"3,1";"y"**{}?<(.2);?>(.8)^{F^2_{f'',f'}}
            \ar@2"2,1";"z"**{}?<(.2);?>(.8)_{F^2_{f''\#_0f',f}} }
          ,"K2"
          \xyboxmatrix"K2"{
            Fx\ar[r]^{\alpha_x}\ar[d]_{Ff}="x"\ar@/_5pc/[ddd]|{F(f''\#_0f'\#_0f)}="z"&
            Gx\ar[d]^{Gf}\ar@2[dl]**{}?(.2);?(.8)|(.4){\alpha_{f}}="X"
            &\\
            Fz\ar[r]_{\alpha_z}\ar[dd]_{F(f''\#_0f')}&Gz\ar@2[ddl]**{}?(.2);?(.8)|(.4){\alpha_{f''\#_0f'}}
            \ar[dd]|(.6){G(f''\#_0f')}="y"\ar[dr]^{Gf'}&\\
            {}&&Gy\ar[dl]^{Gf''}\\
            Fw\ar[r]_{\alpha_w}&Gw&
            \ar@2"3,3";"y"**{}?<(.2);?>(.8)_(.8){G^2_{f',f}}="Y"
            \ar@2"2,1";"z"**{}?<(.2);?>(.8)_(.3){F^2_{f''\#_0f',f}} 
            \tria"Y";"X"
          },"A22"
          \xyboxmatrix"A22"{
            Fx\ar[r]^{\alpha_x}\ar[d]_{Ff}="x"\ar@/_5pc/[ddd]|{F(f''\#_0f'\#_0f)}="z"&
            Gx\ar[d]^{Gf}\ar@2[dl]**{}?(.2);?(.8)|(.4){\alpha_{f}}="X"
            &\\
            Fz\ar[r]_{\alpha_z}\ar[dd]_{F(f''\#_0f')}&Gz\ar@2[ddl]**{}?(.2);?(.8)|(.4){\alpha_{f''\#_0f'}}
            \ar[dd]|(.6){G(f''\#_0f')}="y"\ar[dr]^{Gf'}&\\
            {}&&Gy\ar[dl]^{Gf''}\\
            Fw\ar[r]_{\alpha_w}&Gw&
            \ar@2"3,3";"y"**{}?<(.2);?>(.8)_(.8){G^2_{f',f}}="Y"
            \ar@2"2,1";"z"**{}?<(.2);?>(.8)_(.3){F^2_{f''\#_0f',f}} 
            \tria"X";"Y"
          },"A32"
          \xyboxmatrix"A32"{
            Fx\ar[r]^{\alpha_x}\ar[ddd]_{F(f''\#_0f'\#_0f)}="x"&
            Gx\ar[ddd]|(.3){}="z"\save"z"+(-15,15)*!C\labelbox{G(f''\#_0f'\#_0f)};"z"**\dir{.}\restore
            \ar@2[dddl]**{}?(.2);?(.8)|(.5){\alpha_{f''\#_0f'\#_0f}}
            \ar[r]^{Gf}&Gy\ar[ddd]^{Gf'}\ar[dddl]|{G(f''\#_0f')}="y"\\
            {}&&\\
            &&\\
            Fw\ar[r]_{\alpha_w}&Gw&Gz\ar[l]^{Gf''}
            \ar@2"4,3";"y"**{}?<(.2);?>(.8)^(.3){G^2_{f'',f'}}
            \ar@2"1,3";"z"**{}?<(.2);?>(.8)|{}="k"\save"k"+(-5,15)*!C\labelbox{G^2_{f''\#_0f',f}};"k"**\dir{.}\restore
          }
          \ar@3"A11"+U;"K1"+L^{(\alpha_w\#_0F^2_{f'',f'\#_0f})\\
            \#_1(\overline{\alpha_{f''}\ten{F^2_{f',f}}})\\
            \#_1(Gf''\#_0\alpha_{f'}\#_0Ff)\\
            \#_1(Gf''\#_0Gf'\#_0\alpha_{f})}
          \ar@3"K1"+R;"A21"+U^{(\alpha_w\#_0F^2_{f'',f'\#_0f})\\
            \#_1(\alpha_{f''}\#_0{F(f'\#_0f)})\\
            \#_1(Gf''\#_0\underline{\alpha^2_{f',f}})}
          \ar@3"A12"+D;"K2"+L_{(\alpha_w\#_0F^2_{f''\#_0f',f})
            \\\#_1(\underline{\alpha^2_{f'',f}}\#_0{Ff})
            \\\#_1(Gf''\#_0Gf'\#_0\alpha_f)}
          \ar@3"K2"+R;"A22"+D_{(\alpha_w\#_0F^2_{f''\#_0f',f})\\
            \#_1(\alpha_{f''\#_0f'}\#_0Ff)\\
            \#_1(\overline{G^2_{f''\#_0f'}\ten\alpha_f})}
          \ar@3"A21";"A31"|-{}="here"\ar@{.}"here";p+(0,20)*+\labelbox{
            (\underline{\alpha^2_{f'',f'\#_0f}})\\
            \#_1(Gf''\#_0G^2_{f',f}\#_0\alpha_x)}
          \ar@3"A22";"A32"|-{}="here"\ar@{.}"here";p+(0,-20)*+\labelbox{
            (\underline{\alpha^2_{f''\#_0f',f}})\\\#_1(G^2_{f''\#_0f'}\#_0Gf\#_0\alpha_x)
          }
          \ar@{=}"A11";"A12"
          \ar@{=}"A31";"A32"
        \end{xy}
      \end{equation}
    \end{sidewaysfigure}
    furthermore, $\alpha^2$ has to satisfy the normalization
    condition:
    \begin{equation*}
      \alpha^2_{f',f}=
      \begin{cases}
        \id_{\alpha_{f'}}&\text{if }f'=\id_y\\
        \id_{\alpha_{f}}&\text{if }f=\id_x\\
      \end{cases}
    \end{equation*}
  \item The family of 3-cells $\alpha^2$ has to be compatible with left and
    right whiskering according to (\ref{eq:transex12whiskleft}) and
    (\ref{eq:transex12whiskright}).
    \begin{sidewaysfigure}
      \begin{equation} 
        \label{eq:transex12whiskleft}
        \begin{xy}
          (60,0) : (0,0)
          \xyboxmatrix"A"{
            Fx \ar[r]^{\alpha_x} \ar[d]_{Ff} \ar@/_5pc/[dd]|{F(g'\#_0f)}="w"& Gx\ar[d]^{Gf}
            \ar@2[dl]**{}?(.2);?(.8)|{\alpha_f}\\
            Fy\ar[r]|{\alpha_y} \ar[d]|{Fg}="x"
            \ar@/_2.5pc/[d]|{Fg'}="y"&
            Gy\ar[d]^{Gg}\ar@2[dl]**{}?(.2);?(.8)|{\alpha_g}\\
            Fz\ar[r]_{\alpha_z} & Gz
            \POS\ar@2"x";"y" **{} ?(.2); ?(.8)|{F\gamma}
            \ar@2"2,1";"w" **{} ?(.2); ?>(.8)_{F^2_{g',f}}
          }
          +(1,0)
          \xyboxmatrix"B"{
            Fx \ar[r]^{\alpha_x} \ar[d]_{Ff}
            \ar@/_4pc/[dd]|{F(g'\#_0f)}="w"&
            Gx\ar[d]^{Gf}
            \ar@2[dl]**{}?(.2);?(.8)|{\alpha_f}="X"\\
            Fy\ar[r]|{\alpha_y} \ar[d]_{Fg'} & 
            Gy\ar[d]|{Gg'}="y" \ar@2[dl]**{}?(.2);?(.8)|{\alpha_{g'}} 
            \ar@/^2.5pc/[d]|{Gg}="x"\\
            Fz\ar[r]_{\alpha_z} & Gz
            \POS\ar@2"x";"y" **{} ?(.2); ?(.8)|{G\gamma}="Y"
            \tria"X";"Y"
            \ar@2"2,1";"w" **{} ?(.2); ?>(.8)_{F^2_{g',f}}
          }
          +(1,0)
          \xyboxmatrix"C"{
            Fx \ar[r]^{\alpha_x} \ar[d]_{Ff}
            \ar@/_4pc/[dd]|{F(g'\#_0f)}="w"&
            Gx\ar[d]^{Gf}
            \ar@2[dl]**{}?(.2);?(.8)|{\alpha_f}="X"\\
            Fy\ar[r]|{\alpha_y} \ar[d]_{Fg'} & 
            Gy\ar[d]|{Gg'}="y" \ar@2[dl]**{}?(.2);?(.8)|{\alpha_{g'}} 
            \ar@/^2.5pc/[d]|{Gg}="x"\\
            Fz\ar[r]_{\alpha_z} & Gz
            \POS\ar@2"x";"y" **{} ?(.2); ?(.8)|{G\gamma}="Y"
            \tria"Y";"X"
            \ar@2"2,1";"w" **{} ?(.2); ?>(.8)_{F^2_{g',f}}
          }
          +(1,0)
          \xyboxmatrix"D"{
            Fx \ar[r]^{\alpha_x} \ar[dd]_{F(g'\#_0f)}&
            Gx\ar[dd]|(.6){G(g'\#_0f)}="w"\ar[dr]^{Gf}
            \ar@2[ddl]**{}?(.2);?(.8)|(.4){\alpha_{g'\#_0f}}&\\
            &&
            Gy\ar[dl]|{Gg'}="y" 
            \ar@/^2.5pc/[dl]|{Gg}="x"\\
            Fz\ar[r]_{\alpha_z} & Gz
            \POS\ar@2"x";"y" **{} ?(.2); ?(.8)|{G\gamma}="Y"
            \ar@2"2,3";"w" **{} ?(.2); ?>(.8)_(.7){G^2_{g',f}}
          }
          ,(0,-1)
          \xyboxmatrix"E"{
            Fx \ar[r]^{\alpha_x} \ar[d]_{Ff}
            \ar@/_8pc/[dd]|{F(g'\#_0f)}="w"\ar@/_4pc/[dd]|{F(g\#_0f)}="y"&
            Gx\ar[d]^{Gf}
            \ar@2[dl]**{}?(.2);?(.8)|{\alpha_f}\\
            Fy\ar[r]|{\alpha_y} \ar[d]_{Fg}="x"
            &
            Gy\ar[d]^{Gg}\ar@2[dl]**{}?(.2);?(.8)|{\alpha_g}\\
            Fz\ar[r]_{\alpha_z} & Gz
            \POS\ar@2"2,1";"y" **{} ?<; ?>_{F^2_{g,f}}
            \ar@2"y";"w" **{} ?<; ?>_{F(\gamma\#_0f)}
          }
          ;p+(3,0)="G"**{}?(.5)
          \xyboxmatrix"F"{
            Fx \ar[r]^{\alpha_x} \ar[dd]|{F(g\#_0f)}="x"\ar@/_4pc/[dd]|{F(g'\#_0f)}="y"&
            Gx\ar[dd]|(.6){G(g\#_0f)}="w"\ar[dr]^{Gf}
            \ar@2[ddl]**{}?(.2);?(.8)|(.35){\alpha_{g\#_0f}}&\\
            &&Gy\ar[dl]^{Gg}\\
            Fz\ar[r]_{\alpha_z} & Gz
            \POS\ar@2"x";"y" **{} ?<;?>_{F(\gamma\#_0f)}="Y"
            \ar@2"2,3";"w" **{} ?(.2); ?>(.8)_(.7){G^2_{g,f}}
          }
          ,"G"
          \xyboxmatrix"G"{
            Fx \ar[r]^{\alpha_x} \ar[dd]|{F(g'\#_0f)}="x"&
            Gx\ar[dd]|(.65){G(g'\#_0f)}="w"\ar[drr]^{Gf}
            \ar@2[ddl]**{}?(.2);?(.8)|(.35){\alpha_{g'\#_0f}}
            \ar@/^3pc/[dd]|(.5){G(g\#_0f)}="y"&&\\
            &&&Gy\ar[dll]^{Gg}\\
            Fz\ar[r]_{\alpha_z} & Gz
            \POS\ar@2"y";"w" **{} ?<;?>|{}="here"\ar@{.}"here";p+(.2,-.2)*\labelbox{G(\gamma\#_0f)}="Y"
            \ar@2"2,4";"y" **{} ?(.2); ?>(.8)|{}="here"\ar@{.}"here";p+(.2,-.2)*\labelbox{G^2_{g,f}}&
          }
          \ar@3"A";"B"^{}="here"\ar@{.}"here";p+(0,.4)*\labelbox{
            (\alpha_z\#_0F^2_{g',f})\\
            \#_1(\underline{\alpha_\gamma}\#_0Ff)\\
            \#_1(Gg\#_0\alpha_{f})
          }
          \ar@3"B";"C"^{}="here"\ar@{.}"here";p+(0,.4)*\labelbox{
            (\alpha_z\#_0F^2_{g',f})\\\#_1(\alpha_{g'}\#_0Ff)\\
            \#_1(G\gamma\ten\alpha_f)
          }
          \ar@3"C";"D"^{}="here"\ar@{.}"here";p+(0,.4)*\labelbox{
            \underline{\alpha^2_{g',f}}\\\#_1(G\gamma\#_0Gf\#_0\alpha_x)
          }
          \ar@3"E";"F"^{}="here"\ar@{.}"here";p+(0,-.4)*\labelbox{
            (\alpha_z\#_0F(\gamma\#_0f))\\\#_1(\underline{\alpha^2_{g\#_0f}})
          }
          \ar@3"F";"G"^{}="here"\ar@{.}"here";p+(0,-.4)*\labelbox{
            (\underline{\alpha_{\gamma\#_0f}})\\\#_1(G^2_{g,f}\#_0\alpha_x)
          }
          \ar@{=}"A";"E"
          \ar@{=}"D";"G"
        \end{xy}
      \end{equation}
      \begin{center}
        Compatibility of the cocycle $\alpha^2$ with left whiskers $\gamma\#_0f$.
      \end{center}
    \end{sidewaysfigure}

    \begin{sidewaysfigure}
      \begin{equation} 
        \label{eq:transex12whiskright}
        \begin{xy}
          (60,0) : (0,0)
          \xyboxmatrix"A"{
            Fx \ar[r]^{\alpha_x} \ar[d]_{Ff}="x"
            \ar@/_5pc/[dd]|{F(g\#_0f')}="w"
            \ar@/_2.5pc/[d]|{Ff'}="y"& 
            Gx\ar[d]^{Gf}
            \ar@2[dl]**{}?(.2);?(.8)|{\alpha_f}\\
            Fy\ar[r]|{\alpha_y} \ar[d]|{Fg}
            &
            Gy\ar[d]^{Gg}\ar@2[dl]**{}?(.2);?(.8)|{\alpha_g}="Y"\\
            Fz\ar[r]_{\alpha_z} & Gz
            \POS\ar@2"x";"y" **{} ?(.2); ?(.8)|{F\delta}="X"
            \ar@2"2,1";"w" **{} ?(.2); ?>(.8)^{F^2_{g,f'}}
            \tria"Y";"X"
          }
          +(1,0)
          \xyboxmatrix"B"{
            Fx \ar[r]^{\alpha_x} \ar[d]_{Ff}="x"
            \ar@/_5pc/[dd]|{F(g\#_0f')}="w"
            \ar@/_2.5pc/[d]|{Ff'}="y"& 
            Gx\ar[d]^{Gf}
            \ar@2[dl]**{}?(.2);?(.8)|{\alpha_f}\\
            Fy\ar[r]|{\alpha_y} \ar[d]|{Fg}
            &
            Gy\ar[d]^{Gg}\ar@2[dl]**{}?(.2);?(.8)|{\alpha_g}="Y"\\
            Fz\ar[r]_{\alpha_z} & Gz
            \POS\ar@2"x";"y" **{} ?(.2); ?(.8)|{F\delta}="X"
            \ar@2"2,1";"w" **{} ?(.2); ?>(.8)^{F^2_{g,f'}}
            \tria"X";"Y"
          }
          +(1,0)
          \xyboxmatrix"C"{
            Fx \ar[r]^{\alpha_x} \ar[d]_{Ff'}
            \ar@/_4pc/[dd]|{F(g\#_0f')}="w"&
            Gx\ar[d]^{Gf'}="y"
            \ar@2[dl]**{}?(.2);?(.8)|{\alpha_{f'}}="X"\ar@/^2.5pc/[d]|{Gf}="x"\\
            Fy\ar[r]|{\alpha_y} \ar[d]_{Fg} & 
            Gy\ar[d]^{Gg} \ar@2[dl]**{}?(.2);?(.8)|{\alpha_{g}} 
            \\
            Fz\ar[r]_{\alpha_z} & Gz
            \POS\ar@2"x";"y" **{} ?(.2); ?(.8)|{G\delta}="Y"
            \ar@2"2,1";"w" **{} ?(.2); ?>(.8)_{F^2_{g,f'}}
          }
          +(1,0)
          \xyboxmatrix"D"{
            Fx \ar[r]^{\alpha_x} \ar[dd]_{F(g\#_0f')}&
            Gx\ar[dd]|(.6){G(g\#_0f')}="w"\ar[dr]|{Gf'}="y" 
            \ar@2[ddl]**{}?(.2);?(.8)|(.4){\alpha_{g\#_0f'}}\ar@/^2.5pc/[dr]|{Gf}="x"&\\
            &&
            Gy\ar[dl]^{Gg}
            \\
            Fz\ar[r]_{\alpha_z} & Gz
            \POS\ar@2"x";"y" **{} ?(.2); ?(.8)|{G\delta}="Y"
            \ar@2"2,3";"w" **{} ?(.2); ?>(.8)_(.7){G^2_{g,f'}}
          }
          ,(0,-1)
          \xyboxmatrix"E"{
            Fx \ar[r]^{\alpha_x} \ar[d]_{Ff}
            \ar@/_8pc/[dd]|{F(g\#_0f')}="w"\ar@/_4pc/[dd]|{F(g\#_0f)}="y"&
            Gx\ar[d]^{Gf}
            \ar@2[dl]**{}?(.2);?(.8)|{\alpha_f}\\
            Fy\ar[r]|{\alpha_y} \ar[d]_{Fg}="x"
            &
            Gy\ar[d]^{Gg}\ar@2[dl]**{}?(.2);?(.8)|{\alpha_g}\\
            Fz\ar[r]_{\alpha_z} & Gz
            \POS\ar@2"2,1";"y" **{} ?<; ?>_{F^2_{g,f}}
            \ar@2"y";"w" **{} ?<; ?>_{F(g\#_0\delta)}
          }
          ;p+(3,0)="G"**{}?(.5)
          \xyboxmatrix"F"{
            Fx \ar[r]^{\alpha_x} \ar[dd]|{F(g\#_0f)}="x"\ar@/_4pc/[dd]|{F(g\#_0f')}="y"&
            Gx\ar[dd]|(.6){G(g\#_0f)}="w"\ar[dr]^{Gf}
            \ar@2[ddl]**{}?(.2);?(.8)|(.35){\alpha_{g\#_0f}}&\\
            &&Gy\ar[dl]^{Gg}\\
            Fz\ar[r]_{\alpha_z} & Gz
            \POS\ar@2"x";"y" **{} ?<;?>_{F(g\#_0\delta)}="Y"
            \ar@2"2,3";"w" **{} ?(.2); ?>(.8)_(.7){G^2_{g,f}}
          }
          ,"G"
          \xyboxmatrix"G"{
            Fx \ar[r]^{\alpha_x} \ar[dd]_{F(g\#_0f')}="x"&
            Gx\ar[dd]|(.65){G(g\#_0f')}="w"\ar[drr]^{Gf}
            \ar@2[ddl]**{}?(.2);?(.8)|(.35){\alpha_{g\#_0f'}}
            \ar@/^3pc/[dd]|(.5){G(g\#_0f)}="y"&&\\
            &&&Gy\ar[dll]^{Gg}\\
            Fz\ar[r]_{\alpha_z} & Gz
            \POS\ar@2"y";"w" **{} ?<;?>|{}="here"\ar@{.}"here";p+(.2,-.2)*\labelbox{G(g\#_0\delta)}="Y"
            \ar@2"2,4";"y" **{} ?(.2); ?>(.8)|{}="here"\ar@{.}"here";p+(.2,-.2)*\labelbox{G^2_{g,f}}&
          }
          \ar@3"A";"B"^{}="here"\ar@{.}"here";p+(0,.4)*\labelbox{
            (\alpha_z\#_0F^2_{g,f'})\\
            \#_1(\overline{\alpha_g\ten{}F\delta})\\
            \#_1(Gg\#_0\alpha_f)
          }
          \ar@3"B";"C"^{}="here"\ar@{.}"here";p+(0,.4)*\labelbox{
            (\alpha_z\#_0F^2_{g,f'})\\
            \#_1(\alpha_g\#_0Ff')\\
            \#_1(Gg\#_0\underline{\alpha_\delta})
          }
          \ar@3"C";"D"^{}="here"\ar@{.}"here";p+(0,.4)*\labelbox{
            \underline{\alpha^2_{g,f'}}\\
            \#_1(Gg\#_0G\delta\#_0\alpha_x)
          }
          \ar@3"E";"F"^{}="here"\ar@{.}"here";p+(0,-.4)*\labelbox{
            (\alpha_z\#_0F(g\#_0\delta)\\
            \#_1(\underline{\alpha^2_{g\#_0f}})
          }
          \ar@3"F";"G"^{}="here"\ar@{.}"here";p+(0,-.4)*\labelbox{
            (\underline{\alpha_{g_0\delta}})\\
            \#_1(G^2_{g,f}\#_0\alpha_x)
          }
          \ar@{=}"A";"E"
          \ar@{=}"D";"G"
        \end{xy}
      \end{equation}
      \begin{center}
        Compatibility of the cocycle $\alpha^2$ with right whiskers $g\#_0\delta$.
      \end{center}
    \end{sidewaysfigure}

  \end{enumerate}
  These conditions are derived from the ones in the definition of
  pseudo-$\Gray$-functors \ref{defn:psgrmap}. Note how conditions 4,
  5, 6 of definition \ref{defn:psgrmap} are trivially satisfied for
  transformations. 
\end{rem}

\begin{defn}
  A transformation $\alpha\from{}F\to{}G$ where the cocycle $\alpha^2$
  has only trivial components we call a \defterm{stiff transformation}.
\end{defn}

\begin{lem}
  A stiff transformation $\alpha\from{}F\to{}G$ with $F$ and $G$
  strict $\Gray$-functors is a 1-transfor in the sense of
  \citep{crans}. \qed
\end{lem}

\begin{rem} Given two lax-transformations
  $\xymatrix@1{F\ar[r]^\alpha&G\ar[r]^\beta&H}$ their composite
  $\beta*\alpha$ given by $m\langle\beta,\alpha\rangle$ and has the
  following components:
  \begin{enumerate}
  \item for each 0-cell $x$ of $\G$ the 1-cell
    \begin{equation*}
      \begin{matrix}
        \begin{xy} \xymatrix{Fx\ar[r]^{(\beta*\alpha)_x}&Hx}
        \end{xy}
      \end{matrix}=
      \begin{matrix}
        \begin{xy}
          \xymatrix{Fx\ar[r]^{\alpha_x}&Gx\ar[r]^{\beta_x}&Hx}
        \end{xy}
      \end{matrix},
    \end{equation*}
  \item for each 1-cell $f\from{}x\to{}y$ of $\G$ the 2-cell
    \begin{equation*}
      \begin{matrix}
        \begin{xy}\xycompile{
            \xyboxmatrix{Fx\ar[r]^{(\beta*\alpha)_x}\ar[d]_{Ff}&
              Hx\ar[d]^{Hf}\ar@2[dl]**{}?(.2);?(.8)|{(\beta*\alpha)_f}\\
              Fy\ar[r]_{(\beta*\alpha)_y}&Hy }}
        \end{xy}
      \end{matrix} =\begin{matrix}
        \begin{xy}\xycompile{
            \xyboxmatrix{Fx\ar[r]^{\alpha_x}\ar[d]_{Ff}&
              Gx\ar[d]|{Gf}\ar@2[dl]**{}?(.2);?(.8)|{\alpha_f}\ar[r]^{\beta_x}&
              Hx\ar[d]^{Hf}\ar@2[dl]**{}?(.2);?(.8)|{\beta_f}\\
              Fy\ar[r]_{\alpha_y}&Gy\ar[r]_{\beta_y}&Hy }}
        \end{xy}
      \end{matrix}
    \end{equation*}
  \item for each 2-cell $g\from{}f\to{}f'$ of $\G$ the 3-cell of $\H$ shown in \eqref{eq:exptranscomp2} 
    \begin{sidewaysfigure}
      \begin{multline}\label{eq:exptranscomp2}
        \begin{matrix}
          \begin{xy} \xycompile{
              \xyboxmatrix"A"{Fx\ar[r]^{\alpha_x}\ar[d]|{Ff}="x"\ar@/_3pc/[d]_{Ff'}="y"&
                Hx\ar[d]^{Hf}\ar@2[dl]**{}?(.2);?(.8)|{(\beta*\alpha)_f}\\Fy
                \ar[r]_{(\beta*\alpha)_y}&Hy
                \POS\ar@2"x";"y"**{}?<(.3);?>(.7)^{Fg} } \POS + (55,0)
              \xyboxmatrix"B"{Fx\ar[r]^{(\beta*\alpha)_x}\ar[d]_{Ff'}&
                Hx\ar[d]|{Hf'}="x"\ar@2[dl]**{}?(.2);?(.8)|{(\beta*\alpha)_{f'}}
                \ar@/^3pc/[d]^{Hf}="y"\\Fy\ar[r]_{(\beta*\alpha)_y}&Hy
                \POS\ar@2"y";"x"**{}?<(.3);?>(.7)^{Hg} } \ar@3 "A";"B"
              ^{(\beta*\alpha)_g} }
          \end{xy}
        \end{matrix}\\=
        \begin{matrix}
          \begin{xy} \xycompile{ \xyboxmatrix"A"{Fx \ar[r]^{\alpha_x}
                \ar[d]|{Ff}="x"\ar@/_3pc/[d]_{Ff'}="y" &
                Gx\ar[d]|{Gf}\ar@2[dl]**{}?(.2);?(.8)|{\alpha_f}{}\ar[r]^{\beta_x}&
                Hx\ar[d]^{Hf}\ar@2[dl]**{}?(.2);?(.8)|{\beta_f}\\
                Fy\ar[r]_{\alpha_y}&Gy\ar[r]_{\beta_y}&Hy
                \POS\ar@2"x";"y"**{}?<(.3);?>(.7)^{Fg} } \POS+(85,0)
              \xyboxmatrix"B"@+.7cm{ Fx\ar[r]^{\alpha_x}\ar[d]_{Ff'}&
                Gx\ar@/_1.5pc/[d]|{Gf'}="x"
                \ar@2[dl]**{}?(.4);?(.6)|{\alpha_{f'}}
                \ar@/^1.5pc/[d]|{Gf}="y"\ar[r]^{\beta_x}&Hx\ar[d]^{Hf}
                \ar@2[dl]**{}?(.4);?(.6)|{\beta_f} \\
                Fy\ar[r]_{\alpha_y}&Gy\ar[r]_{\beta_y}&Hy
                \POS\ar@2"y";"x"**{}?<(.3);?>(.7)^{Gg} }
              \POS\ar@3 "B";"A"**{}?<(.4);?>(.6)^{(\beta_y\#_0\alpha_g)\\
                \;\#_1(\beta_f\#_0\alpha_x)}
              \POS+(85,0) 
              \xyboxmatrix"C"{
                Fx\ar[r]^{\alpha_x}\ar[d]_{Ff'}&
                Gx\ar[d]|{Gf'}\ar@2[dl]**{}?(.2);?(.8)|{\alpha_{f'}}\ar[r]^{\beta_x}&
                Hx\ar[d]|{Hf'}="x"\ar@2[dl]**{}?(.2);?(.8)|{\beta_{f'}}\ar@/^3pc/[d]^{Hf}="y"\\
                Fy\ar[r]_{\alpha_y}&Gy\ar[r]_{\beta_y}&Hy
                \POS\ar@2"y";"x"**{}?<(.3);?>(.7)^{Hg} }
              \POS\ar@3 "B";"C"**{}?<;?>^{(\beta_y\#_0\alpha_{f'})\\
                \;\#_1(\beta_g\#_0\alpha_x)} }
          \end{xy}
        \end{matrix}
      \end{multline}
    \end{sidewaysfigure}
  \item for each pair of composable 1-cells $f\from{}x\to{}y$,
    $f'\from{}y\to{}z$ a 3-cell shown in \eqref{eq:exptranscomp2cocycle}
    \begin{sidewaysfigure}
      \begin{multline}\label{eq:exptranscomp2cocycle}
        \begin{matrix}
          \begin{xy} \xycompile{ 
              \xyboxmatrix"A"{
                Fx\ar[r]^{(\beta*\alpha)_x}\ar[d]|{Ff}="x"\ar@/_4pc/[dd]_{F(f'\#_0f)}|{}="y"&
                Hx\ar[d]^{Hf}\ar@2[dl]**{}?(.2);?(.8)|{(\beta*\alpha)_f}\\Fy
                \ar[r]|{(\beta*\alpha)_y}\ar[d]|{Ff'}&Hy\ar@2[dl]**{}?(.2);?(.8)|{(\beta*\alpha)_{f'}}\ar[d]^{Hf'}\\
                Fz\ar[r]_{(\beta*\alpha)_z}&=Hz\POS\ar@2"2,1";"y"**{}?<(.2);?(.8)^{F^2_{f',f}}
              } \POS + (75,0) 
              \xyboxmatrix"B"{
                Fx\ar[r]^{(\beta_\alpha)_x}\ar[dd]_(.4){F(f'\#_0f)}&
                Hx\ar[dd]|{H(f'\#_0f)}="x"\ar@2[ddl]**{}?(.2);?(.8)|(.65){(\beta*\alpha)_{f'\#_0f}}
                \ar@/^4pc/[dd]^{Hf'\#_0Hf}="y"\\
                {}&{}\\Fy\ar[r]_{(\beta*\alpha)_z}&Hy
                \POS\ar@2"y";"2,2"**{}?<(.3);?>(.7)^{H^2_{f',f}} }
              \ar@3"A";"B"^-{(\beta*\alpha)^2_{f',f}} }
          \end{xy}
        \end{matrix}
        \\=
        \begin{matrix}
          \begin{xy} 
              \xyboxmatrix"Z"{
                Fx\ar[r]^{\alpha_x}\ar[d]|{Ff}="x"\ar@/_4pc/[dd]_{F(f'\#_0f)}|{}="y"&
                Gx\ar[r]^{\beta_x}\ar[d]|{Gf}\ar@2[dl]**{}?(.2);?(.8)|{\alpha_f}="k"&
                Hx\ar[d]^{Hf}\ar@2[dl]**{}?(.2);?(.8)|{\beta_f}\\
                Fy\ar[r]|{\alpha_y}\ar[d]|{Ff'}&
                *i{Gy}\ar[r]|{\beta_y}\ar[d]|{Gf'}\ar@2[dl]**{}?(.2);?(.8)|{\alpha_{f'}}&
                Hy\ar@2[dl]**{}?(.2);?(.8)|{\beta_{f'}}="l"\ar[d]^{Hf'}\\
                Fz\ar[r]_{\alpha_z}&Gz\ar[r]_{\beta_z}&Hz
                \POS\ar@2"2,1";"y"**{}?<(.2);?(.8)^{F^2_{f',f}}
                \tria "k";"l"
              }
              \POS (125,0)
              \xyboxmatrix"A"@+1.5cm{
                Fx \ar[r]^{\alpha_x}
                \ar[d]|{Ff'\#_0Ff}="x"\ar@/_4pc/[d]_{F(f'\#_0f)}="y" &
                Gx\ar[d]|{Gf'\#_0Gf}\ar@2[dl]**{}?(.2);?(.8)|{(\alpha_{f'}\#_0Ff)\\
                  \;\#_1(Gf'\#_0\alpha_f)}{}\ar[r]^{\beta_x}&
                Hx\ar[d]|{Hf'\#_0Hf}\ar@2[dl]**{}?(.2);?(.8)|{(\beta_{f'}\#_0Gf)\\
                  \;\#_1(Hf'\#_0\beta_f)}\\
                Fz\ar[r]_{\alpha_z}&Gz\ar[r]_{\beta_z}&Hz
                \POS\ar@2"x";"y"**{}?<(.2);?>(.8)^{F^2_{f',f}} 
              }
              \POS (0,-60) 
              \xyboxmatrix"B"@+1.5cm{ 
                Fx\ar[r]^{\alpha_x}\ar[d]|{F(f'\#_0f)}&
                Gx\ar@/_1.7pc/[d]|(.7){G(f'\#_0f)}="x" 
                \ar@2@<-0em>[dl]**{}?(.2);?(.8)|{\alpha_{f'\#_0f}} 
                \ar@/^1.7pc/[d]|(.3){G(f'\#_0f)}="y"\ar[r]^{\beta_x}&Hx\ar[d]^{Hf'\#_0Hf}
                \ar@2@<+1em>[dl]**{}?(.2);?(.8)|{(\beta_{f'}\#_0Gf)\\
                  \;\#_1(Hf'\#_0\beta_f)} \\
                Fz\ar[r]_{\alpha_z}&Gz\ar[r]_{\beta_z}&Hz
                \POS\ar@2@<-.5em>"y";"x"**{}?<(.1);?>(.9)^{G^2_{f',f}} 
              }
              \POS (125,-60) 
              \xyboxmatrix"C"@+1.5cm{
                Fx\ar[r]^{\alpha_x}\ar[d]_{F(f'\#_0f)}&
                Gx\ar[d]|{G(f'\#_0f)}\ar@2[dl]**{}?(.2);?(.8)|{\alpha_{f'\#_0f}}\ar[r]^{\beta_x}&
                Hx\ar[d]|{H(f'\#_0f)}="x"\ar@2[dl]**{}?(.2);?(.8)|{\beta_{f'\#_0f}}\ar@/^4pc/[d]^{Hf'\#_0Hf}="y"\\
                Fz\ar[r]_{\alpha_z}&Gz\ar[r]_{\beta_z}&Hz
                \POS\ar@2@<-.7ex>"y";"x"**{}?<(.2);?>(.8)^{H^2_{f',f}} 
              } 
              \POS\ar@3 "Z";"A"**{}?<;?>^{(\beta_z\#_0\alpha_z\#_0F^2_{f',f})\\
                \;\#_1(\beta_z\#_0\alpha_{f'}\#_0Ff)\\
                \;\#_1(\beta_{f'}\ten\alpha_f)\\
                \;\#_1(Hf'\#_0\beta_f\#_0\alpha_x)}
              \POS\ar@3 "A";"B"**{}?<;?>|{(\beta_z\#_0\alpha^2_{f',f})\\
                \;\#_1(((\beta_{f'}\#_0Gf)\\
                \;\#_1(Hf'\#_0\beta_f))\#_0\alpha_x)}
              \POS\ar@3 "B";"C"**{}?<;?>_{ (\beta_z\#_0\alpha_{f'\#_0f})\\
                \;\#_1(\beta^2_{f',f}\#_0\alpha_x)}
          \end{xy}
        \end{matrix}
      \end{multline}
    \end{sidewaysfigure}
  \end{enumerate}
\end{rem}

\begin{defn}\label{defn:modif} Assuming $\alpha$ and $\beta$ are as in
  definition \ref{defn:ltransf} and $F$ and $G$ are pseudo-functors
  $\G\laxto\H$, a \defterm{modification}
  $A\from\alpha\to\beta\from{}F\to{}G$ is a pseudo-functor
  $A\from\G\laxto\dblbarspc\H$, such that $d_0A=\alpha$ and
  $d_1A=\beta$.
\end{defn}

\begin{rem}
  A modification $A\from\alpha\to\beta$ according to definitions
  \ref{defn:modif} and \ref{defn:psgrmap} is given by the following
  data:
  \begin{enumerate}
  \item For every 0-cell $x$ in $\G$ a 2-cell 
    \begin{equation*}
      \begin{xy}
        \xyboxmatrix{
          Fx\ar@/^1.5pc/[r]^{\alpha_x}="a"\ar@/_1.5pc/[r]_{\beta_x}="b"&Gx
          \ar@2"a";"b"**{}?(.2);?(.8)|{A_x}
        }
      \end{xy}
    \end{equation*}
  \item For every 1-cell $f\from{}x\to{}y$ a 3-cell in $\H$
    \begin{equation*}
      \begin{xy}
        \xyboxmatrix"A"{
          Fx\ar@/^1.5pc/[r]^{\alpha_x}="a"\ar@/_1.5pc/[r]_{\beta_x}="b"\ar[d]_{Ff}&
          Gx\ar[d]^{Gf}\ar@<+1em>@2[dl]**{}?(.5);?(.9)|{\beta_f}\\
          Fy\ar@/_1.5pc/[r]_{\beta_y}="b1"&
          Gy
          \ar@2"a";"b"**{}?(.2);?(.8)|{A_x}
        }
        +(40,0)
        \xyboxmatrix"B"{
          Fx\ar@/^1.5pc/[r]^{\alpha_x}="a"\ar[d]_{Ff}&
          Gx\ar[d]^{Gf}\ar@<-1em>@2[dl]**{}?(.1);?(.5)|{\alpha_f}\\
          Fy\ar@/^1.5pc/[r]^{\alpha_y}="a1"\ar@/_1.5pc/[r]_{\beta_y}="b1"&
          Gy
          \ar@2"a1";"b1"**{}?(.2);?(.8)|{A_y}
        }
        \ar@3"A";"B"**{}?<;?>^{A_f}
      \end{xy}
    \end{equation*}
  \end{enumerate}
  This data has to satisfy the following conditions:
  \begin{enumerate}
  \item Units are preserved:
    \begin{equation*}
      A_{\id_{x}}=\id_{A_x}
    \end{equation*}
  \item Compatibility with the cocycles of $F,G,\alpha,\beta$
    according to (\ref{eq:modifexcoccomp})
    \begin{sidewaysfigure}
      \begin{equation}
        \label{eq:modifexcoccomp}
        \begin{xy}
          \save (55,0):(0,-1):: 
          ,(0,0)="A11" 
          ,(1,0)="A21" 
          ,(2,0)="A31" 
          ,(3,0)="A41" 
          ,(0,1)="A12"
          ,(3,1)="A32" 
          ,"A12";"A32"**{};?(.5)="A22" 
          ,(0,2)="A13" 
          ,(1,2)="A23" 
          ,(2,2)="A33"
          ,(3,2)="A43"
          ,(0,3)="A14" 
          ,(1,3)="A24" 
          ,(2,3)="A34" 
          ,(3,3)="A44" 
          \restore 
          \POS,"A11"
          \xyboxmatrix"A11"{
            Fx\ar@/^1.5pc/[r]^{\alpha_x}="a"\ar@/_1.5pc/[r]|{\beta_x}="b"\ar[d]_{Ff}
            \ar@/_3pc/[dd]_{F(f'\#_0f)}="z"&
            Gx\ar[d]^{Gf}\ar@<+1em>@2[dl]**{}?(.5);?(.9)|{\beta_f}\\
            Fy\ar@/_1.5pc/[r]|{\beta_y}="b1"\ar[d]_{Ff'}&
            Gy\ar[d]^{Gf'}\ar@<+1em>@2[dl]**{}?(.5);?(.9)|{\beta_{f'}}\\
            Fz\ar@/_1.5pc/[r]_{\beta_z}="b2"&Gz
            \ar@2"a";"b"**{}?(.2);?(.8)|{A_x}
            \ar@2"2,1";"z"**{}?<(.2);?>(.8)_{F^2_{f',f}}
          },"A12"
          \xyboxmatrix"A12"{
            Fx\ar@/^1.5pc/[r]^{\alpha_x}="a"\ar@/_1.5pc/[r]|{\beta_x}="b"\ar[dd]_{F(f'\#_0f)}&
            Gx\ar[dd]|{G(f'\#_0f)}="z"\ar[dr]^{Gf}\ar@<+1em>@2[ddl]**{}?(.5);?(.9)|{\beta_{f'\#_0f}}\\
            {}&{}&{}\ar[dl]^{Gf'}\\
            Fz\ar@/_1.5pc/[r]_{\beta_z}="b2"&Gz
            \ar@2"a";"b"**{}?(.2);?(.8)|{A_x}="K"
            \ar@2"2,3";"z"**{}?<(.2);?>(.8)|{}="here"\ar@{.}"here";p+(-5,20)*!C\labelbox{G^2_{f',f}}
            \tria"K";"here"
          }
          ,"A22"
          \xyboxmatrix"A22"{
            Fx\ar@/^1.5pc/[r]^{\alpha_x}="a"\ar@/_1.5pc/[r]|{\beta_x}="b"\ar[dd]_{F(f'\#_0f)}&
            Gx\ar[dd]|{G(f'\#_0f)}="z"\ar[dr]^{Gf}\ar@<+1em>@2[ddl]**{}?(.5);?(.9)|{\beta_{f'\#_0f}}\\
            {}&{}&{}\ar[dl]^{Gf'}\\
            Fz\ar@/_1.5pc/[r]_{\beta_z}="b2"&Gz
            \ar@2"a";"b"**{}?(.2);?(.8)|{A_x}="K"
            \ar@2"2,3";"z"**{}?<(.2);?>(.8)|{}="here"\ar@{.}"here";p+(-5,20)*!C\labelbox{G^2_{f',f}}
            \tria"here";"K"
          }
          ,"A21"
          \xyboxmatrix"A21"{
            Fx\ar@/^1.5pc/[r]^{\alpha_x}="a"\ar[d]_{Ff}
            \ar@/_3pc/[dd]_{F(f'\#_0f)}="z"&
            Gx\ar[d]^{Gf}\ar@<-1em>@2[dl]**{}?(.1);?(.5)|{\alpha_f}\\
            Fy\ar@/_1.5pc/[r]|{\beta_y}="b1"\ar@/^1.5pc/[r]|{\alpha_y}="a1"\ar[d]_{Ff'}&
            Gy\ar[d]^{Gf'}\ar@<+1em>@2[dl]**{}?(.5);?(.9)|{\beta_{f'}}\\
            Fz\ar@/_1.5pc/[r]_{\beta_z}&Gz
            \ar@2"a1";"b1"**{}?(.2);?(.8)|{A_y}
            \ar@2"2,1";"z"**{}?<(.2);?>(.8)_{F^2_{f',f}}
          }
          ,"A31"
          \xyboxmatrix"A31"{
            Fx\ar@/^1.5pc/[r]^{\alpha_x}="a"\ar[d]_{Ff}
            \ar@/_3pc/[dd]_{F(f'\#_0f)}="z"&
            Gx\ar[d]^{Gf}\ar@<-1em>@2[dl]**{}?(.1);?(.5)|{\alpha_f}\\
            Fy\ar@/^1.5pc/[r]|{\alpha_y}="a1"\ar[d]_{Ff'}&
            Gy\ar[d]^{Gf'}\ar@<-1em>@2[dl]**{}?(.1);?(.5)|{\alpha_{f'}}\\
            Fz\ar@/^1.5pc/[r]|{\alpha_z}="a2"\ar@/_1.5pc/[r]_{\beta_z}="b2"&Gz
            \ar@2"a2";"b2"**{}?(.2);?(.8)|{A_z}="Y"
            \ar@2"2,1";"z"**{}?<(.2);?>(.8)_{F^2_{f',f}}="X"
            \tria"Y";"X"
          }
          ,"A41"
          \xyboxmatrix"A41"{
            Fx\ar@/^1.5pc/[r]^{\alpha_x}="a"\ar[d]_{Ff}
            \ar@/_3pc/[dd]_{F(f'\#_0f)}="z"&
            Gx\ar[d]^{Gf}\ar@<-1em>@2[dl]**{}?(.1);?(.5)|{\alpha_f}\\
            Fy\ar@/^1.5pc/[r]|{\alpha_y}="a1"\ar[d]_{Ff'}&
            Gy\ar[d]^{Gf'}\ar@<-1em>@2[dl]**{}?(.1);?(.5)|{\alpha_{f'}}\\
            Fz\ar@/^1.5pc/[r]|{\alpha_z}="a2"\ar@/_1.5pc/[r]_{\beta_z}="b2"&Gz
            \ar@2"a2";"b2"**{}?(.2);?(.8)|{A_z}="Y"
            \ar@2"2,1";"z"**{}?<(.2);?>(.8)_{F^2_{f',f}}="X"
            \tria"X";"Y"
          }
          ,"A32"
          \xyboxmatrix"A32"{
            Fx\ar@/^1.5pc/[r]^{\alpha_x}="a"\ar[dd]_{F(f'\#_0f)}&
            Gx\ar[dd]|{G(f'\#_0f)}="z"\ar@<-1em>@2[ddl]**{}?(.1);?(.5)|{\alpha_{f'\#_0f}}\ar[dr]^{Gf}&\\
            {}&{}&\ar[dl]^{Gf'}\\
            Fz\ar@/^1.5pc/[r]|{\alpha_z}="a2"\ar@/_1.5pc/[r]_{\beta_z}="b2"&Gz&
            \ar@2"a2";"b2"**{}?(.2);?(.8)|{A_z}
            \ar@2"2,3";"z"**{}?<(.2);?>(.8)|{}="here"\ar@{.}"here";p+(-5,20)*!C\labelbox{G^2_{f',f}}
          }
          \ar@3"A11";"A21"|{}="here"\ar@{.}"here";p+(0,30)*!C\labelbox{
            (\alpha_z\#_0F^2_{f',f})\\
            \#_1(\beta_{f'}\#_0Ff)\\
            \#_1(Gf'\#_0\underline{A_f})}
          \ar@3"A21";"A31"|{}="here"\ar@{.}"here";p+(0,30)*!C\labelbox{
            (\beta_z\#_0F^2_{f',f})\\
            \#_1(\underline{A_{f'}}\#_0Ff)\\
            \#_1(Gf'\#_0\alpha_f)}
          \ar@3"A31";"A41"|{}="here"\ar@{.}"here";p+(0,30)*!C\labelbox{
          \overline{A_z\ten{}F^2_{f',f}}\\\#_1(\alpha_{f'}\#_0Ff)\\\#_1(Gf'\#_0\alpha_f)}
          \ar@3"A41";"A32"|{}="here"\ar@{.}"here";p+(30,0)*!C\labelbox{
            (A_z\#_0F(f'\#_0f))\\
            \#_1(\underline{\alpha^2_{f',f}})}
          \ar@3"A11";"A12"|{}="here"\ar@{.}"here";p+(-20,0)*!C\labelbox{
            (\underline{\beta^2_{f',f}})\\\#_1(Gf'\#_0Gf\#_0A_x)}
          \ar@3"A12";"A22"|{}="here"\ar@{.}"here";p+(0,-20)*!C\labelbox{
            \beta_{f'\#_0f}\\\#_1(G^2_{f',f}\ten{}A_x)}
          \ar@3"A22";"A32"|{}="here"\ar@{.}"here";p+(0,-20)*!C\labelbox{
            \underline{A_{f'\#_0f}}\\\#_1(G^2_{f',f}\#_0\alpha_x)}
        \end{xy}
      \end{equation}
      \begin{center}
        Compatibility of the modification $A$ with the cocycles of
        $F,G,\alpha,\beta$
      \end{center}

    \end{sidewaysfigure}
  \item For 2-cells $g\from{}f\To{}f'$ in $\G$ the images under $F$
    and $G$ as well the data of $A$, $\alpha$ and $\beta$ are
    compatible as shown in (\ref{eq:modifex2cell})
    \begin{sidewaysfigure}
      \begin{equation}
        \label{eq:modifex2cell}
        \begin{xy}
          \save (55,0):(0,-1):: 
          ,(0,0)="A11" 
          ,(1,0)="A21" 
          ,(2,0)="A31" 
          ,(3,0)="A41" 
          ,(0,1)="A12"
          ,(1,1)="A22"
          ,(2,1)="A32" 
          ,(3,1)="A42" 
          ,(0,2)="A13" 
          ,(1,2)="A23" 
          ,(2,2)="A33"
          ,(3,2)="A43"
          ,(0,3)="A14" 
          ,(1,3)="A24" 
          ,(2,3)="A34" 
          ,(3,3)="A44" 
          \restore 
          \POS,"A11"
          \xyboxmatrix"A11"{
            Fx\ar@/^1.5pc/[r]^{\alpha_x}="a"\ar@/_1.5pc/[r]_{\beta_x}="b"
            \ar[d]|{Ff}="a2"\ar@/_3pc/[d]_{Ff'}="b2"&
            Gx\ar[d]^{Gf}\ar@<+1em>@2[dl]**{}?(.5);?(.9)|{\beta_f}\\
            Fy\ar@/_1.5pc/[r]_{\beta_y}="b1"&
            Gy
            \ar@2"a";"b"**{}?(.2);?(.8)|{A_x}="Y"
            \ar@2"a2";"b2"**{}?(.2);?(.8)|{Fg}="X"
          }
          ,"A21"
          \xyboxmatrix"A21"{
            Fx\ar@/^1.5pc/[r]^{\alpha_x}="a"
            \ar[d]|{Ff}="a2"\ar@/_3pc/[d]_{Ff'}="b2"&
            Gx\ar[d]^{Gf}\ar@<-1em>@2[dl]**{}?(.1);?(.5)|{\alpha_f}\\
            Fy\ar@/^1.5pc/[r]^{\alpha_y}="a1"\ar@/_1.5pc/[r]_{\beta_y}="b1"&
            Gy
            \ar@2"a1";"b1"**{}?(.2);?(.8)|{A_y}="Y"
            \ar@2"a2";"b2"**{}?(.2);?(.8)|{Fg}="X"
            \tria"Y";"X"
          }
          ,"A31"
          \xyboxmatrix"A31"{
            Fx\ar@/^1.5pc/[r]^{\alpha_x}="a"
            \ar[d]|{Ff}="a2"\ar@/_3pc/[d]_{Ff'}="b2"&
            Gx\ar[d]^{Gf}\ar@<-1em>@2[dl]**{}?(.1);?(.5)|{\alpha_f}\\
            Fy\ar@/^1.5pc/[r]^{\alpha_y}="a1"\ar@/_1.5pc/[r]_{\beta_y}="b1"&
            Gy
            \ar@2"a1";"b1"**{}?(.2);?(.8)|{A_y}="Y"
            \ar@2"a2";"b2"**{}?(.2);?(.8)|{Fg}="X"
            \tria"X";"Y"
          }
          ,"A41"
          \xyboxmatrix"A41"{
            Fx\ar@/^1.5pc/[r]^{\alpha_x}="a"\ar[d]_{Ff'}&
            Gx\ar[d]|{Gf'}="a2"\ar@<-1em>@2[dl]**{}?(.1);?(.5)|{\alpha_{f'}}
            \ar@/^3pc/[d]^{Gf}="b2"\\
            Fy\ar@/^1.5pc/[r]^{\alpha_y}="a1"\ar@/_1.5pc/[r]_{\beta_y}="b1"&
            Gy
            \ar@2"a1";"b1"**{}?(.2);?(.8)|{A_y}="Y"
            \ar@2"b2";"a2"**{}?(.2);?(.8)|{Gg}="X"
          },"A12"
          \xyboxmatrix"A12"{
            Fx\ar@/^1.5pc/[r]^{\alpha_x}="a"\ar@/_1.5pc/[r]_{\beta_x}="b"
            \ar[d]|{Ff}="a2"\ar@/_3pc/[d]_{Ff'}="b2"&
            Gx\ar[d]^{Gf}\ar@<+1em>@2[dl]**{}?(.5);?(.9)|{\beta_{f'}}\\
            Fy\ar@/_1.5pc/[r]_{\beta_y}="b1"&
            Gy
            \ar@2"a";"b"**{}?(.2);?(.8)|{A_x}="Y"
            \ar@2"a2";"b2"**{}?(.2);?(.8)|{Fg}="X"
          }
          ,"A22"
          \xyboxmatrix"A22"{
            Fx\ar@/^1.5pc/[r]^{\alpha_x}="a"\ar@/_1.5pc/[r]_{\beta_y}="b"
            \ar[d]_{Ff'}&
            Gx\ar[d]|{Gf'}="a2"\ar@/^3pc/[d]^{Gf}="b2"
            \ar@<+1em>@2[dl]**{}?(.6);?(.9)|{\beta_f}\\
            Fy\ar@/_1.5pc/[r]_{\beta_y}="a1"&
            Gy
            \ar@2"a";"b"**{}?(.2);?(.8)|{A_x}="Y"
            \ar@2"b2";"a2"**{}?(.2);?(.8)|{Gg}="X"
            \tria"Y";"X"
          }
          ,"A32"
          \xyboxmatrix"A32"{
            Fx\ar@/^1.5pc/[r]^{\alpha_x}="a"\ar@/_1.5pc/[r]_{\beta_y}="b"
            \ar[d]_{Ff'}&
            Gx\ar[d]|{Gf'}="a2"\ar@/^3pc/[d]^{Gf}="b2"
            \ar@<+1em>@2[dl]**{}?(.6);?(.9)|{\beta_{f'}}\\
            Fy\ar@/_1.5pc/[r]_{\beta_y}="a1"&
            Gy
            \ar@2"a";"b"**{}?(.2);?(.8)|{A_x}="Y"
            \ar@2"b2";"a2"**{}?(.2);?(.8)|{Gg}="X"
            \tria"X";"Y"
          }
          ,"A42"
          \xyboxmatrix"A42"{
            Fx\ar@/^1.5pc/[r]^{\alpha_x}="a"\ar[d]_{Ff'}&
            Gx\ar[d]|{Gf'}="a2"\ar@<-1em>@2[dl]**{}?(.1);?(.5)|{\alpha_{f'}}
            \ar@/^3pc/[d]^{Gf}="b2"\\
            Fy\ar@/^1.5pc/[r]^{\alpha_y}="a1"\ar@/_1.5pc/[r]_{\beta_y}="b1"&
            Gy
            \ar@2"a1";"b1"**{}?(.2);?(.8)|{A_y}="Y"
            \ar@2"b2";"a2"**{}?(.2);?(.8)|{Gg}="X"
          }
          \ar@{=}"A11";"A12"
          \ar@{=}"A41";"A42"
          \ar@3"A11";"A21"**{}?<;?>^{(\beta_y\#_0Fg)\\\#_1\underline{A_f}}
          \ar@3"A21";"A31"**{}?<;?>^{\overline{A_y\ten{}Fg}\\\#_1\alpha_f}
          \ar@3"A31";"A41"**{}?<;?>^{(A_y\#_0Ff')\\\#_1\underline{\alpha_{g}}}
          \ar@3"A12";"A22"**{}?<;?>_{\underline{\beta_g}\\\#_1(Gf\#_0A_x)}
          \ar@3"A22";"A32"**{}?<;?>_{\beta_{f'}\\\#_1(Gg\ten{}A_x)}
          \ar@3"A32";"A42"**{}?<;?>_{\underline{A_{f'}}\\\#_1(Gg\#_0\alpha_x)}
        \end{xy}
      \end{equation}
      \begin{center}
        Compatibility of 2-cells with $A$, $\alpha$ and $\beta$
      \end{center}
    \end{sidewaysfigure}
  \end{enumerate}
\end{rem}

\begin{lem}
  \label{lem:modif2trsf}
  A transformation $A\from\alpha\to\beta$ where
  $\alpha,\beta\from{}F\to{}G$ are stiff and $F, G$ are strict is a
  2-transfor in the sense of \citep{crans}.\qed
\end{lem}

\begin{defn}
  \label{defn:pertdef} Given modifications $A,B\from\alpha\to\beta$ a
  \defterm{perturbation} is a pseudo-$\Gray$-functor
  $\sigma\from\G\laxto\overline{\overline{\overline{\H}}}$ such that
  $d_0\sigma=A$ and $d_1\sigma=B$.
\end{defn}

\begin{rem}
  According to definition \ref{defn:pertdef} a perturbation is given by
  a 3-cell in $\H$ 
  \begin{equation*}
    \begin{xy}
      \save (40,0):(0,-1):: 
      ,(0,0)="A11" 
      ,(1,0)="A21" 
      \restore
      \POS,"A11"
      \xyboxmatrix"A11"{
        Fx\ar@/^1.5pc/[r]^{\alpha_x}="a"
        \ar@/_1.5pc/[r]_{\beta_x}="b"& Gx
        \ar@2"a";"b"**{}?(.2);?(.8)|{A_x}
      }
      ,"A21"
      \xyboxmatrix"A21"{
        Fx\ar@/^1.5pc/[r]^{\alpha_x}="a"
        \ar@/_1.5pc/[r]_{\beta_x}="b"& Gx
        \ar@2"a";"b"**{}?(.2);?(.8)|{B_x}
      }
      \ar@3"A11";"A21"**{}?<;?>^{\sigma_x}
    \end{xy}
  \end{equation*}
  for each 0-cell $x$ in $\G$ such that
  \begin{equation*}
    \begin{xy}
          \save (45,0):(0,-1):: 
          ,(0,0)="A11" 
          ,(1,0)="A21" 
          ,(0,1)="A12"
          ,(1,1)="A22"
          \restore 
          \POS,"A11"
          \xyboxmatrix"A11"{
            Fx\ar@/^1.5pc/[r]^{\alpha_x}="a"\ar@/_1.5pc/[r]_{\beta_x}="b"
            \ar[d]_{Ff}="a2"&
            Gx\ar[d]^{Gf}\ar@<+1em>@2[dl]**{}?(.5);?(.9)|{\beta_f}\\
            Fy\ar@/_1.5pc/[r]_{\beta_y}="b1"&
            Gy
            \ar@2"a";"b"**{}?(.2);?(.8)|{A_x}="Y"
          }
          ,"A12"
          \xyboxmatrix"A12"{
            Fx\ar@/^1.5pc/[r]^{\alpha_x}="a"
            \ar[d]_{Ff}&
            Gx\ar[d]^{Gf}\ar@<-1em>@2[dl]**{}?(.1);?(.5)|{\alpha_f}\\
            Fy\ar@/^1.5pc/[r]^{\alpha_y}="a1"\ar@/_1.5pc/[r]_{\beta_y}="b1"&
            Gy
            \ar@2"a1";"b1"**{}?(.2);?(.8)|{A_y}="Y"
          }   
          ,"A21"
          \xyboxmatrix"A21"{
            Fx\ar@/^1.5pc/[r]^{\alpha_x}="a"\ar@/_1.5pc/[r]_{\beta_x}="b"
            \ar[d]_{Ff}="a2"&
            Gx\ar[d]^{Gf}\ar@<+1em>@2[dl]**{}?(.5);?(.9)|{\beta_f}\\
            Fy\ar@/_1.5pc/[r]_{\beta_y}="b1"&
            Gy
            \ar@2"a";"b"**{}?(.2);?(.8)|{B_x}="Y"
          }
          ,"A22"
          \xyboxmatrix"A22"{
            Fx\ar@/^1.5pc/[r]^{\alpha_x}="a"
            \ar[d]_{Ff}&
            Gx\ar[d]^{Gf}\ar@<-1em>@2[dl]**{}?(.1);?(.5)|{\alpha_f}\\
            Fy\ar@/^1.5pc/[r]^{\alpha_y}="a1"\ar@/_1.5pc/[r]_{\beta_y}="b1"&
            Gy
            \ar@2"a1";"b1"**{}?(.2);?(.8)|{B_y}="Y"
          }
          \ar@3"A11";"A21"**{}?<;?>^{\beta_f\\\#_1(Gf\#_0\sigma_x)}
          \ar@3"A21";"A22"**{}?<;?>^{B_f}
          \ar@3"A11";"A12"**{}?<;?>_{A_f}
          \ar@3"A12";"A22"**{}?<;?>_{(\sigma_y\#_0Ff)\\\#_1\alpha_f}      
    \end{xy}
  \end{equation*}
  commutes.
\end{rem}

\begin{lem}
  \label{lem:pert3trsf}
  A perturbation $\sigma\from{}A\to{}B$ fulfilling the conditions of lemma
  \ref{lem:modif2trsf} is a 3-transfor in the sense of
  \citep{crans}.\qed 
\end{lem}

\appendix
\section{Adjunctions}
\label{sec:adjunctions}

We can embed the ideas developed in section \ref{sec:cofrpl} in a more global
picture. The functor $\fQ^1\from\Gray\Cat\to\Gray\Cat$ is part of the
following adjunction of fibered categories:
\begin{equation*}
  \xymatrix{
    F^*(\Gray\Cat)\ar@/^/[r]^-{{(\_)_1}^*(F)}="x1"\ar[d]_{F^*((\_)_1)}&
    \Gray\Cat\ar[d]^{(\_)_1}\ar@/^/[l]^-{\underline{U}}="y1"\\
    \RGrph\ar@/^/[r]^-{F}="x"&\Cat\ar@/^/[l]^-{U}="y"
    \ar@{|- }"y1";"x1"**{}?(.4);?(.6)
    \ar@{|- }"y";"x"**{}?(.4);?(.6)
  }
\end{equation*}
where $F$ means ``free category over a reflexive graph'' and $U$ means
``underlying reflexive graph of a category'', $(\_)_1$ means
``underlying category of a $\Gray$-category. According to
\cite[4.1]{hermida} the adjunction $F\ladj U$ lifts canonically to an
adjunction $({(\_)_1}^*(F),F)\ladj(\underline{U},U)$ of fibered
categories. Which means in particular that
${(\_)_1}^*(F)\ladj\underline{U}$ is an adjunction and our $\fQ^1$ can
be defined as ${(\_)_1}^*(F)\underline{U}$.

The objects of $\Graph\times\Gray\Cat$ might be called 1-free
$\Gray$-categories. 

\begin{rem}
  Let $P\from\mathcal{E}\to\mathcal{B}$ be a 2-fibration in the sense
  of \cite{hermida}. Given $u\from I\to PX$ and $u'\from I'\to PX$ for
  $X$ an object in $\mathcal E$; and an equivalence $h\from I\to I'$
  such that $u'h=u$. Then the unique filler $\hat{h}$ over $h$ is an
  equivalence as well.

  In particular, given the comparison functor $K\from\X_{FU}\to\A$ for
  the comonad induced by $F\ladj U\from\A\to\X$ lifts to a comparison
  functor $\hat{K}$.
\end{rem}

\begin{lem}
  If $F$ is comonadic, then so is $({(\_)_1}^*(F),F)$.
\end{lem}

\label{sec:bibl} \bibliographystyle{abbrvnat-B} \bibliography{mapspc}

\end{document}
